\documentclass[12pt]{report}
\usepackage{amsmath,amsfonts,amssymb,amsthm,mathrsfs}
\usepackage{wasysym}
\usepackage{minitoc}
\usepackage{endnotes}
\usepackage[dvipsnames]{xcolor}
\usepackage[a4paper,vmargin={3.5cm,3.5cm},hmargin={2.5cm,2.5cm}]{geometry}

\usepackage{graphicx,graphics}
\usepackage{epsfig}
\usepackage{latexsym}
\usepackage[applemac]{inputenc}
\usepackage{ae,aecompl}

\usepackage[english]{babel}
 \usepackage[colorlinks=true]{hyperref}
\usepackage{pstricks}
\usepackage{enumerate}
\newcommand{\subgraph}{\sqsubset}
\renewcommand{\leq}{\leqslant}
\renewcommand{\geq}{\geqslant}
\usepackage{baskervald}
\usepackage[baskervaldx]{newtxmath}
\usepackage[font=sf, labelfont={sf,bf}, margin=1cm]{caption}

\usepackage{titlesec}
\usepackage{titletoc}
\usepackage{color}

\titleformat{\part}[display]{\beaupetit}{\Huge\textsf{Part \textsc{\Roman{part}}:}    }{0pt}{}[]

\DeclareFixedFont{\chapterNumberFont}{OT1}{fsk
}{b}{n}{5cm} 
\DeclareFixedFont{\classik}{OT1}{ptm}{b}{sc}{1cm} 
\definecolor{gris75}{gray}{0.75} 
\definecolor{gris25}{gray}{0.15}

\titleformat{\section}[block]{\Large\sffamily}{\noindent \bfseries \thesection}{1em}{} 
\titleformat{\subsection}[block]{\large\sffamily}{\thesubsection}{1em}{} 
\titleformat{\subsubsection}[block]{\sffamily}{}{}{} 
\titleformat{\paragraph}[runin]{\sffamily}{}{}{} 
\titlespacing{\section} {0pc}{3.5ex plus .1ex minus .2ex}{1.5ex minus .1ex}

\titleformat{\chapter}[hang]{\bfseries\Large
}{\textsf{
\textsc{\Roman{chapter}}:}  }{0pt}{}[\titlerule 
]

\newcommand{\W}{\mathcal{W}}

\newcommand{\noeud}[1]{\raisebox{.5pt}{\textcircled{\raisebox{-.9pt} {$#1$}}}}
  \DeclareFixedFont{\beaupetit}{T1}{ftp}{b}{n}{2cm} 

\usepackage{thmtools}

\declaretheorem[thmbox=M,name=Theorem,parent=chapter]{theorem}
\declaretheorem[name=Proposition,sibling=theorem]{proposition}
\declaretheorem[name=Lemma,sibling=theorem]{lemma}
\declaretheorem[name=Corollary,sibling=theorem]{corollary}

\declaretheorem[shaded={bgcolor= {rgb}{0.8,0.85,1}},parent=chapter]{definition}

\declaretheorem[parent=chapter,style=remark,name=Remark,parent=chapter]{remark}
\declaretheorem[parent=chapter,style=remark,name=Example,parent=chapter]{example}
\declaretheorem[parent=chapter,style=remark,name=Exercise,parent=chapter]{exo}

\def\llbrack{\{\hspace{-.25em} \{ }
\def\rrbrack{ \} \hspace{-.25em}\}}

\def\build#1_#2^#3{\mathrel{
\mathop{\kern 0pt#1}\limits_{#2}^{#3}}}

             \title{\beaupetit A random Walk among \\
             random Graphs}
             \author{\textsc{Nicolas CURIEN}}
             
\begin{document}
             \maketitle

                          \clearpage
         \section*{A random walk among random graphs}   
         

The theory of random graphs is now ubiquitous in probability theory, and there are already many comprehensive textbooks (to name just a few \cite{RemcoRGI,RemcoRGII,bollobas2001random,janson2011random,drmota2009random,durrett2010random}) dealing with the numerous models of random graphs invented over the last decades. The goal of these lecture notes is to give a glimpse of a few models of random graphs together with some of the probabilistic tools used to study them. It is intended for master or PhD  students in probability theory. I chose the models of random graphs mainly by taste and by the will to cover different types of probabilistic arguments. This document should not be seen as an authoritative reference but rather as a recreational (random) walk in the wonderland of random graph theory. Several exercises of varying difficulty (most of them being non trivial) are scattered along the text and each chapter is ended with bibliographical pointers. Here are the main topics covered in the lecture notes together with the \textit{mathematical tools they introduce}:\medskip 

\begin{itemize} 
\item \textsc{Chapter I:} Basic of (bond) percolation. Phase transition. The Rado graph. \\ \indent \textit{Graph theory, First and second moment, duality}.
\item \textsc{Chapter II:} One-dimensional random walk, Recurrence/transience, Oscillation/drift. \\ \indent \textit{Law of large numbers and its reciproque, Fourier transform}.
\item \textsc{Chapter III:} Skip-free random walk, duality and cycle lemma. Applications: Kemperman formula, Ballot theorem, parking on the line.\\ \noindent \textit{Feller combinatorial cyclic lemma}.
\item \textsc{Chapter IV:} Bienaymé-Galton-Watson trees, {\L}ukasiewicz encoding, Enumeration. \\ \noindent \textit{Formal series, Neveu's plane tree formalism}.
\item \textsc{Chapter V:} Sharp threshold for graph properties on the Erd{\H{o}}s--R\'enyi: connectedness, clique number, diameter, cycle. Convergence of the spectrum. \\ \noindent \textit{First and second moment method, method of moments, Poisson paradigm}.
\item \textsc{Chapter VI:} Phase transition for the giant component I. \\ \noindent \textit{$ \varepsilon$-cut, first moment method, sprinkling, multiplicative coalescent}.
\item \textsc{Chapter VII:} Phase transition for the giant component II. \\ \noindent \textit{Markovian exploration, differential equation method}.
\item \textsc{Chapter VIII:} Phase transition for the giant component III.  \\ \noindent \textit{Poissonization, Bin counting processes, Brownian asymptotics}.
\item \textsc{Chapter IX:} (Uniform) random permutations. Poisson-Dirichlet distribution and Dyckman function for large cycles. Poisson limit for small cycle counts. \\ \noindent \textit{Feller's coupling, Randomization, Stick breaking construction}.
\item \textsc{Chapter X:} Random recursive tree (and random permutations). \\ \noindent \textit{Chinese restaurant process, Recursive distributional equation, Polya urn scheme}.
\item \textsc{Chapter XI:} Continuous-time embedding and applications. \\ \noindent \textit{Athreya-Karling embedding of Markov chains, convergence of Yule processes and links between exponential and Poisson processes}.
\item \textsc{Chapter XII:} Spine decomposition and applications. \\ \noindent \textit{Martingale transform, spine decomposition, many-to-one formulas}.
\item \textsc{Chapter XIII:} Barab\'asi--Albert random tree. \\ \noindent \textit{Preferential attachment mechanism, scale-free random networks}.

 \end{itemize} 

Many thanks go to the students that attended the ``random graph'' master course I gave in 2019-2025 at Orsay. They contributed to the development of the material and spotted many typos. I am particularly grateful to Alice Contat, Baojun Wu (promotion 2019),  Guillaume Blanc, Maude Bellugeon, Elie Khalfallah (promotion 2020), Tanguy Lions, Francisco Calvillo (promotion 2021), Corentin Correia, Lo\"ic Gassmann (promotion 2022), Nathan de Montgolfier, Laureline Legros, Emile Averous (promotion 2023), Remi Bernard, 
Simone Maria Giancola (promotion 2024). Special thanks go to Damian Cid  for spotting (so)many typoes and inaccuracies and for his participation to Chapter \ref{chap:poissonER}. I am also grateful to Serte Donderwinkel for many useful comments. 

\clearpage 
\section*{Notations}
We list here the (perhaps non-standard) notation we use through the lecture notes:\\

\noindent \begin{tabular}{cl}
e.g. & for example (\textit{exempli gratia})\\
i.e. & namely (\textit{id est})\\
a.s. & almost surely\\
i.o. & infinitely often\\
$ \mathbb{Z}_{>0}$ & $=\{1,2,3,\cdots\}$\\
$ \mathbb{Z}_{\geq0}$ & $=\{0,1,2,3,\cdots\}$\\
$ \mathbb{Z}_{<0}$ & $=\{\cdots,-3,-2,-1\}$\\
$ \mathbb{Z}_{\leq0}$ & $=\{\cdots,-3,-2,-1,0\}$\\
 $ \equiv$ & \mbox{gives a shorter and temporary notation for an object}\\
$\#E$ & cardinality of the set $E$\\
$ [z^{n}]f(z)$ & $=f_{n}$ when $f(z) = \sum_{i \geq0} f_{i}z^{i} \in \mathbb{C}[[X]]$
\end{tabular}
\medskip 

\noindent For an asymptotically positive function $f(n)$  and random variables $X_{n} : n \geq 0$ we write

\noindent \begin{tabular}{cl}
$ X_{n} \sim_{ \mathbb{P}} f(n)$ & if  $\frac{X_{n}}{f(n)} \xrightarrow[n\to\infty]{( \mathbb{P})} 1$\\
$ X_{n}  = o_{ \mathbb{P}}( f(n))$ & if  $\frac{X_{n}}{f(n)} \xrightarrow[n\to\infty]{( \mathbb{P})} 0$\\
$ X_{n} = O_{ \mathbb{P}}( f(n))$ & if $(X_{n}/f(n) : n \geq 1)$ is tight\\
\end{tabular}

\medskip 

\noindent If furthermore the variables $X_{n}$ are coupled and form a sequence $(X_{n} : n \geq 0)$ then we write

\noindent \begin{tabular}{cl}
$ X_{n} \sim_{ a.s.} f(n)$ & if  $\frac{X_{n}}{f(n)} \xrightarrow[n\to\infty]{a.s.} 1$\\
$ X_{n}  = o_{ a.s.} (f(n))$ & if  $\frac{X_{n}}{f(n)} \xrightarrow[n\to\infty]{a.s.} 0$\\
$ X_{n} = O_{ a.s.}( f(n))$ & if $(X_{n}/f(n) : n \geq 1)$ is bounded above
\end{tabular}

\medskip 

\noindent We use standard notation for several (laws of) random variables:
\medskip 

\noindent \begin{tabular}{ll}
$  \mathcal{N}(m, \sigma^2)$ & real Gaussian law with mean $m$ and variance $\sigma^2$\\
$ \mathcal{E}(\alpha)$ & exponential variable with mean $1/\alpha$\\
$ (B_t : t \geq 0)$& standard linear Brownian motion issued from $0$\\
$ ( \mathfrak{P}(t) : t \geq 0)$& unit rate Poisson counting process, \\
& in particular $ \mathfrak{P}(t)$ is a Poisson random variable with mean $t$\\
$G(n, p)$ & Erd{\H{o}}s--R\'enyi random graph with $n$ vertices and edge parameter $p$\\
$(S_n : n \geq 0)$ & random walk with i.i.d.~increments (see context for the law of increments)\\
\end{tabular}

\medskip 

\noindent Graph notation: 

\noindent \begin{tabular}{ll}
$  \mathrm{V}(\mathfrak{g}),\mathrm{E}(\mathfrak{g})$& vertex and edge sets of a graph $ \mathfrak{g}$\\
$x \sim y$ & vertices $x,y$ are neighbors in the underlying graph $ \mathfrak{g}$\\
$x \leftrightarrow y$ & vertices $x,y$ are in the same connected component in the underlying graph $ \mathfrak{g}$\\
$ \mathrm{deg}_{ \mathfrak{g}}(x)$ or $ \mathrm{deg}(x)$ & degree of the vertex $x \in \mathrm{V}( \mathfrak{g})$\\
$\mathfrak{g}' \sqsubset \mathfrak{g}$ & $ \mathfrak{g'}$ is a subgraph of $ \mathfrak{g}$\\
$ \mathfrak{g}[V]$ & graph induced by $ \mathfrak{g}$ on the vertices $V$\\
$ \mathfrak{g}\simeq \mathfrak{g}'$ & two isomorphic graphs\\
$ \mathrm{d}_{ \mathrm{gr}}^\mathfrak{g}$ or $ \mathrm{d_{gr}}$ & graph distance on $  \mathrm{V}(\mathfrak{g})$\\
$ \mathbb{G}_n$ & set of all simple graphs on the vertex set $\{1,2, \dots , n \}$
\end{tabular}

\medskip 

\noindent Tree notation:
\medskip 

\noindent \begin{tabular}{ll}
$ (T_n : n \geq 0)$ & is the random recursive tree or uniform attachment chain \\
$ ( \mathsf{T}_n : n \geq 1)$ & is the Barab\'asi--Albert or linear preferential attachment chain \\
$( [\mathbb{T}]_t : t \geq 0)$ & is a standard Yule tree (rate $1$ and usually order $2$) process\\
$ \mathcal{T}$ & is a Bienaym\'e--Galton--Watson tree \\ & whose offspring distribution should be clear from the context\\
 \end{tabular}

             \tableofcontents
             


\bigskip 



\chapter{Basics of percolation}
\hfill An appetizer.

\bigskip 

In this introductory chapter we present the model of Bernoulli bond \textbf{percolation}. This is a way to generate a random graph from a deterministic graph by keeping some of its edges at random. The random graphs studied in part I (Bienaym\'e--Galton--Watson trees) and in part II (Erd{\H{o}}s--R\'enyi random graphs) can be seen as percolation models on some special graphs. Our goal here is only to present the main features of the Bernoulli percolation model focusing on the  \textbf{phase transition} for the existence of an infinite cluster.

\begin{figure}[!h]
 \begin{center}
 \includegraphics[width=3.7cm]{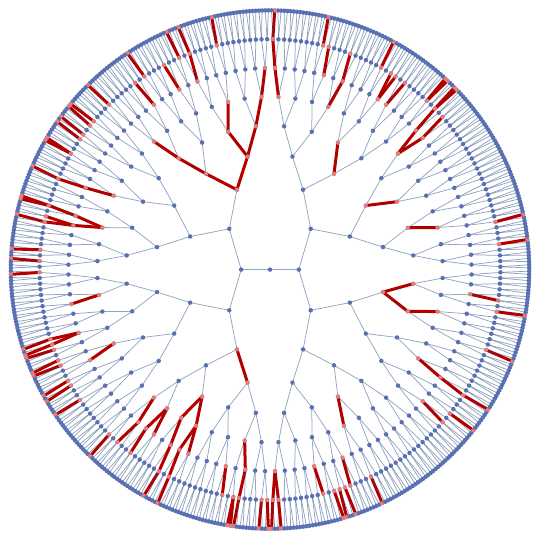}
  \includegraphics[width=3.7cm]{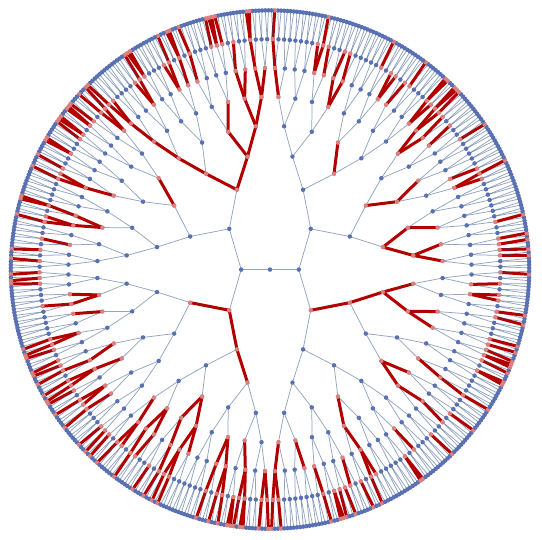}
   \includegraphics[width=3.7cm]{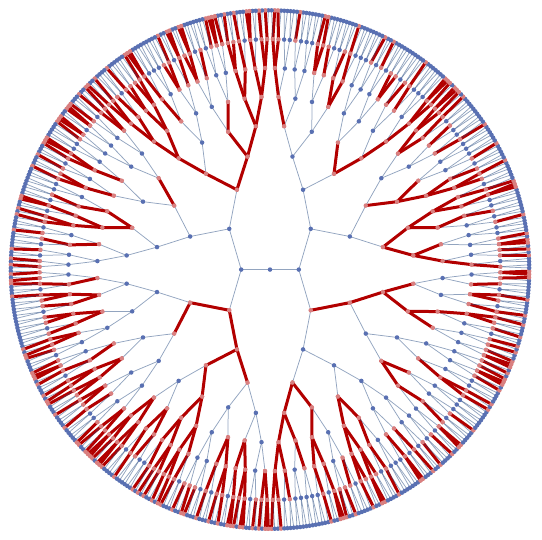}
    \includegraphics[width=3.7cm]{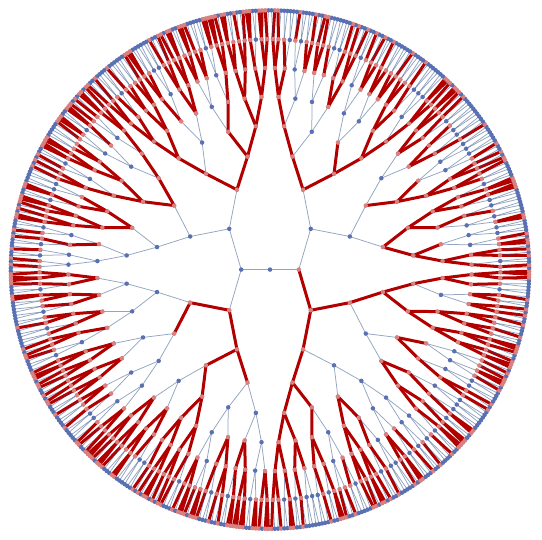}
 \caption{ Increasing Bernoulli percolation on a complete binary tree up to level $10$ with parameters $p=0.2$, $0.4$, $0.5$ and $p=0.6$ from left to right.}
 \end{center}
 \end{figure}
 
 \begin{figure}[!h]
 \begin{center}
 \includegraphics[width=3.7cm]{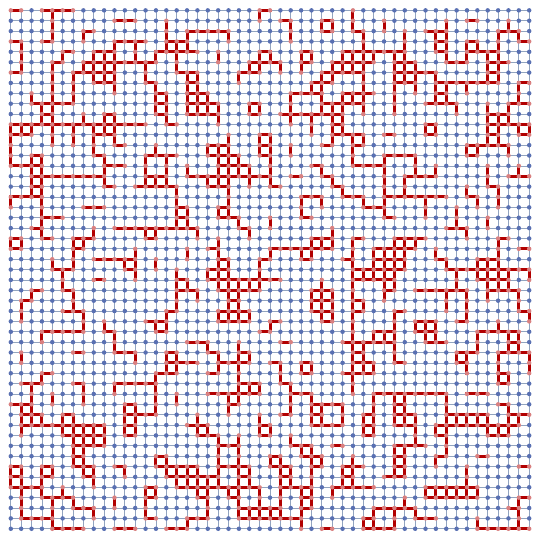}
  \includegraphics[width=3.7cm]{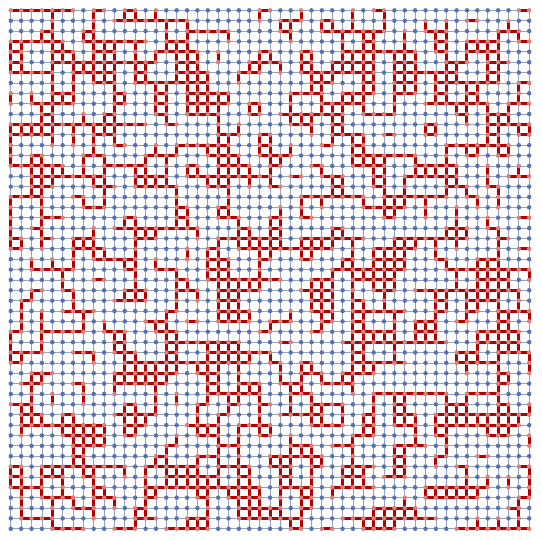}
   \includegraphics[width=3.7cm]{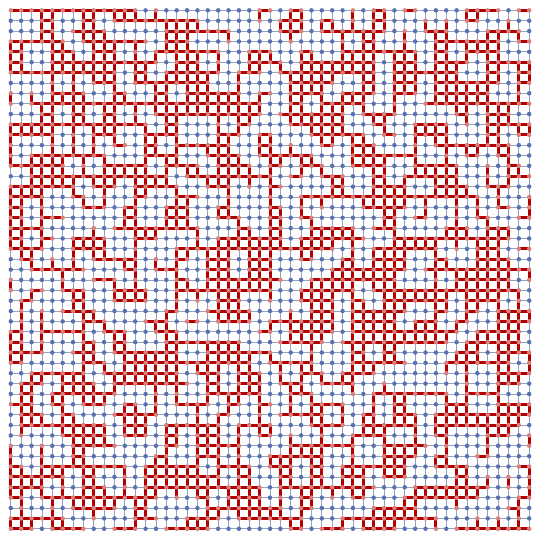}
    \includegraphics[width=3.7cm]{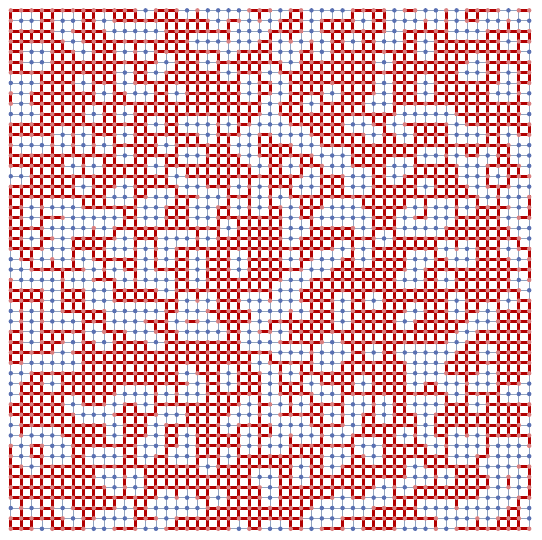}
 \caption{Increasing Bernoulli percolation on a $50\times 50$ grid with parameters $p=0.35$, $p=0.45$, $p=0.55$ and $p=0.65$. Notice the appearance of an ubiquitous cluster between the second and the third picture.}
 \end{center}
 \end{figure}

\section{Basics on graphs}
 A \textbf{graph}\footnote{more formally, a non-oriented multi-graph} $ \mathfrak{g}$ is a  pair $ \mathfrak{g}=(\mathrm{V}(\mathfrak{g}),  \mathrm{E}(\mathfrak{g}))$, where $ V= \mathrm{V}(\mathfrak{g})$ is the set of \textbf{vertices} of $\mathfrak{g}$ and $E=\mathrm{E}(\mathfrak{g})$ is the set of \textbf{edges} of $\mathfrak{g}$ which is a multiset ( i.e.~where repetitions are allowed) over the set $\{V^{2}\}$ of all unordered pairs of elements of $V$.  The graph is $\textbf{simple}$ if they are no multiple edges nor loops (an edge with confounded end vertices).
\begin{figure}[!h]
 \begin{center}
\includegraphics[height=4cm]{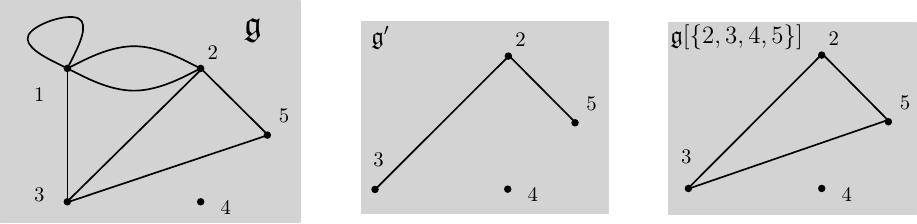}
\caption{(Left) An example of a graph $\mathfrak{g}=(V,E)$ with vertex set  
$  \displaystyle V= \{1,2,3,4,5\}$ and edge set  $ \displaystyle E = \{\hspace{-1mm}\{\{1,1 \},\{1,2 \},\{ 1,2\},\{1,3 \},\{3,2 \},\{2,5 \},\{3,5 \}\}\hspace{-1mm}\} $. The vertex degrees of $1,2,3,4$ in $ \mathfrak{g}$ are respectively $5,3,2,0$. (Center) An example of a subgraph $ \mathfrak{g}'\subgraph \mathfrak{g}$ and (Right) the subgraph induced on the vertices $2,3,4,5$.}
 \end{center}
 \end{figure}
If $x,y \in V$ and $\{x,y\} \in E$ we say that $x$ and $y$ are \textbf{neighbors} and we write $ x \sim y$. We say that an edge is \textbf{adjacent} to a vertex if it is one of its endpoints, and two edges are adjacent if they are adjacent to a common vertex. The \textbf{degree} of a vertex $x \in V$ denoted by $ \mathrm{deg}_{ \mathfrak{g}}(x)$ (or $ \mathrm{deg}(x)$ if this no ambiguity) is the number of half-edges adjacent to $x$, otherwise said it is the number of edges adjacent to $x$ where loops are counted twice. A \textbf{subgraph} of $ \mathfrak{g}$ is a graph $ \mathfrak{g'}$ such that $ \mathrm{V}( \mathfrak{g}') \subset \mathrm{V}( \mathfrak{g})$ and where $ \mathrm{E} ( \mathfrak{g}') \subset \mathrm{E}( \mathfrak{g})$. We shall write $ \mathfrak{g}' \sqsubset \mathfrak{g}$ in this case.  If $ V' \subset \mathrm{V}( \mathfrak{g})$ the subgraph \textbf{graph induced} by $ \mathfrak{g}$ on $V'$ is the graph with vertex set $V'$ obtained by keeping only the edges of $ \mathrm{E}( \mathfrak{g})$ whose endpoints are in $V'$. It is denoted by $ \mathfrak{g}[V']$, note that $ \mathfrak{g}[V'] \subgraph \mathfrak{g}$. 

\paragraph{Graph equivalence.}
If $\mathfrak{g}$ and $\mathfrak{g}'$ are two graphs we say that $\mathfrak{g}$ and $\mathfrak{g}'$ are \textbf{equivalent} if they represent the same graph up to renaming the vertex set. Formally this means that there exists a bijection $\phi :  \mathrm{V}(\mathfrak{g}) \to  \mathrm{V}(\mathfrak{g}')$ which maps the multi-set $ \mathrm{E}(\mathfrak{g})$ to $ \mathrm{E}(\mathfrak{g}')$: such a function is called a homomorphism of graph (automorphism if $ \mathfrak{g}= \mathfrak{g}'$) and we write $ \mathfrak{g} \simeq \mathfrak{g}'$. In this course we shall often implicitly identify two equivalent\footnote{although the space of equivalence classes of all finite connected countable graphs is  monstrous, see \cite{vershik1998universal}} graphs. 
\clearpage

\begin{center} \hrulefill \textit{Convention} \hrulefill  \end{center}
Unless explicitly specified, we shall always suppose that $ \mathrm{E}( \mathfrak{g})$ is finite or countable and that $ \mathfrak{g}$ is \textbf{locally finite} i.e.~that the vertex degrees are all finite (no vertices of infinite degree). 
\begin{center} \hrulefill  \end{center}

\paragraph{Connected graphs.} A \textbf{path} $\gamma = (e_1, e_2, \dots)$ is a sequence of adjacent edges in the graph, its \textbf{length} is the number of edges it contains. If the starting and endpoint points of $\gamma$ are the same it is called a \textbf{cycle}. The path $\gamma$ is \textbf{self-avoiding} if $e_i$ and $e_j$ are not adjacent when $|i-j|>1$.   The \textbf{graph distance} on $\mathfrak{g}$ is denoted by $ \mathrm{d}_{ \mathrm{gr}}^\mathfrak{g}$ or $ \mathrm{d_{gr}}$ when there is no ambiguity, and is defined for $x,y \in \mathrm{V}( \mathfrak{g})$ by 
$$ \mathrm{d_{gr}}(x,y) = \mbox{minimal length of a path $\gamma$ going from  }x \mbox{ to }y.$$
By convention we put $ \mathrm{d_{gr}}(x,y)=\infty$ if there is no path linking $x$ to $y$ in $\mathfrak{g}$. The equivalence classes from the relation $ x \leftrightarrow y \iff \mathrm{d_{gr}}(x,y) < \infty$ are the connected components of $\mathfrak{g}$.  If the connected component of $v_0 \in \mathrm{V}( \mathfrak{g})$ is infinite we write $ v_0\leftrightarrow \infty$. We say that $\mathfrak{g}$ is \textbf{connected} if it has only one connected component.  The connected graphs with a minimal number of edges are famously called \textbf{trees}:
\begin{proposition}[Tree] \label{prop:tree} Let $\mathfrak{g}=(V,E)$ be a connected graph on $n$ vertices. Then we must have $ \# E \geq n-1$. If $\# E=n-1$ then $\mathfrak{g}$ is a \textbf{tree}, meaning that is has no non trivial cycle.
\end{proposition}
\noindent \textbf{Proof.}  We can suppose that the vertex set of $\mathfrak{g}$ is $\{1,2,3,\dots ,n\}$. We start with the vertex $1$. Since $\mathfrak{g}$ is connected there exists an edge adjacent to $1$ of the form $\{1, i_{1}\}$. If $i_{1} =1$ then this edge is a loop and  otherwise $i_{1} \ne 1$. We then throw this edge away and pick a new edge adjacent to either $1$ or $i_{1}$. Iteratively, after having explored  $k$ edges, we have discovered a part of the connected component of $1$ which has at most $k+1$ vertices. Since $\mathfrak{g}$ is connected it follows that $ \# E \geq n-1$. In case of equality this means that during the exploration process we have never found an edge linking two vertices already explored, in other words, no non trivial cycle has been created and $\mathfrak{g}$ is thus a tree. \qed \medskip 

We record here a useful property (whose proof is left as an exercise) known as K\"onig's lemma  which characterizes infinite connected components via existence of infinite self-avoiding paths:
\begin{lemma}[K\"onig's lemma] \label{lem:konig} Let $ \mathfrak{g}$ be a locally finite graph and let $v_0 \in \mathrm{V}( \mathfrak{g})$. Then the following propositions are equivalent 
\begin{enumerate}[(i)]
\item The connected component of $v_0$ is infinite, i.e.~$v_0 \leftrightarrow \infty$,
\item There is a self-avoiding infinite path starting from $v_0$,
\item For every $ n \geq 1$, there is a self-avoiding path starting from $v_0$ and of length $n$.
\end{enumerate}
\end{lemma}



\section{Percolation} \label{sec:percoabstrait}

\begin{definition}[Bernoulli bond percolation] Fix a countable graph $ \mathfrak{g}$ and a parameter $p \in [0,1]$. The Bernoulli bond percolation on $ \mathfrak{g}$ with parameter $p$ is the  random graph $$ \mathrm{Perc}( \mathfrak{g},p)$$ whose vertex set is $ \mathrm{V}( \mathfrak{g})$ and where each edge $e \in \mathrm{E}( \mathfrak{g})$ is kept independently of each other with probability $p$. The edges kept are called ``open" and those discarded are called ``closed".
\end{definition} 

Obviously, for each $p \in [0,1]$, the random graph $\mathrm{Perc}( \mathfrak{g},p)$ is a subgraph of $ \mathfrak{g}$ which is bigger and bigger as $p$ increases. To make this statement formal, it is useful to \textbf{couple}~i.e.~to realize on the same probability space, all graphs $\mathrm{Perc}( \mathfrak{g},p)$ for $p \in [0,1]$. A natural way to do this is to consider a probability space $(\Omega, \mathcal{F}, \mathbb{P})$ which supports i.i.d.~random variables $ (U_e : e \in \mathrm{E}( \mathfrak{g}))$ which are uniformly distributed on $[0,1]$ --this is possible since we supposed that $ \mathrm{E}( \mathfrak{g})$ is at most countable--. It is now clear that if we set 
$$ \mathrm{Perc}( \mathfrak{g},p) = \Big( \mathrm{V}( \mathfrak{g}) ;  \left\{ e \in \mathrm{E}( \mathfrak{g}) : U_e \leq p \right\} \Big),$$
then for each $p$, the random graph $\mathrm{Perc}( \mathfrak{g},p)$ is indeed distributed as a percolation on $ \mathfrak{g}$ with parameter $p$ and furthermore $p \mapsto \mathrm{Perc}( \mathfrak{g},p)$ is increasing (for the inclusion of edges).  The connected components  of $\mathrm{Perc}( \mathfrak{g},p)$ are called \textbf{clusters}.


\begin{remark}[History, see \cite{Mendes}] Percolation was designed to model the porosity of coal (used for gas masks during the second world war). In 1942, Rosalind Franklin (famous later for participating to  the discovery of DNA structure) working for the  British Coal Utilisation Research Association remarked that the porosity of coal depends on the size of the molecules of the gas and on the temperature at which the coal was formed. Later on, in the 50's, Simon Broadbent also working at BCURA as a statistician, together with the mathematician  John Hammersley, introduced  Bernoulli bond percolation on the grid to model these phenomena. \end{remark}

\section{Phase transition}
In the rest of the chapter, we focus on graphs $ \mathfrak{g}$ which are infinite and connected.
Much of the theory of percolation is focused on the existence of large clusters in $\mathrm{Perc}( \mathfrak{g},p)$. More precisely, if $ \mathfrak{g}$ is an infinite connected graph, one can ask whether for some parameter $p$, the random graph $\mathrm{Perc}( \mathfrak{g},p)$ has an infinite cluster\footnote{using Proposition \ref{lem:konig} one can prove that this event is indeed measurable for the the $\sigma$-field generated by the variables $ \mathbf{1}_{e \mathrm{ \ is \ open }}$ for $e \in \mathrm{E}( \mathfrak{g})$}. More precisely, the function 
$$ p \mapsto \mathbb{P}( \mathrm{Perc}( \mathfrak{g},p) \mbox{ contains an infinite cluster} ),$$ is easily seen to be increasing using the coupling of Section \ref{sec:percoabstrait}. Since the existence of an infinite cluster in  $\mathrm{Perc}( \mathfrak{g},p)$ is an event which is independent of the status of any finite number of edges, it has probability $0$ or $1$ by Kolmogorov $0-1$ law. We say that there is a phase transition, if this probability does not depend trivially on $p$:
\begin{definition}[$p_c$ and phase transition] We define the \textbf{critical parameter} $p_c( \mathfrak{g})$ as 
$$ p_c( \mathfrak{g})= \inf \{ p \in [0,1] :  \mathbb{P}( \mathrm{Perc}( \mathfrak{g},p) \mbox{ has an infinite cluster})=1 \}.$$
If $p_c( \mathfrak{g}) \in (0,1)$ we say that there is a non trivial phase transition for percolation on $ \mathfrak{g}$.
\end{definition}
For example, the line graph $ \mathfrak{z}$ whose vertex set is $\mathbb{Z}$ with the edges $\{\{i,i+1\} : i \in \mathbb{Z}\}$ has no phase transition since $p_c(  \mathfrak{z}) = 1$. Similarly, the (non-locally finite) graph made of a star with infinite degree has  $p_c( \mathrm{star})=0$. We will see in Proposition \ref{prop:lowerperco} that having vertices with large degrees is the only way to achieve  $p_c=0$. \\ 
Knowing whether or not there is an infinite cluster at the critical threshold $p_{c}$ is one of the main open question in the area: it is widely believed that for ``homogeneous'' graphs there is no infinite cluster at the critical point. 
\begin{remark} The terminology ``phase transition" comes from the fact that around the critical parameter $p_c$, a slight variation of the parameter $p$ induces dramatic changes in the large scale geometry of the random graph $ \mathrm{Perc}( \mathfrak{g},p)$. This can be  used  to model physical phase transitions (such as the transformation of water into ice when the temperature drops below $0^{\circ}$C).
\end{remark}

\section{Two examples}
  In the rest of this section we shall prove the existence of a non-trivial phase transition for percolation on two infinite graphs: the infinite binary tree and the cubic planar lattice. We shall use the so-called first and second moment method which is (a sometimes subtle) application of Markov\footnote{\raisebox{-5mm}{\includegraphics[width=1cm]{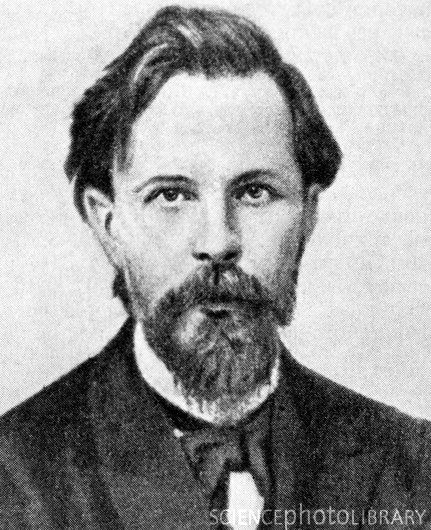}} Andreï Andreïevitch Markov  (1856--1922), Russian}  and Cauchy--Schwarz\footnote{\raisebox{-5mm}{\includegraphics[width=1cm]{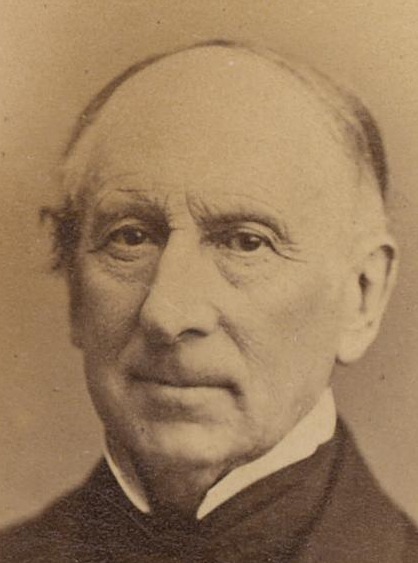}} Augustin Louis Cauchy (1789--1857), French  \raisebox{-5mm}{\includegraphics[width=1cm]{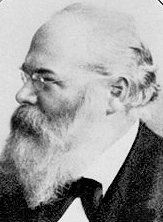}} Hermann Amandus Schwarz (1843--1921), German} inequalities and which will accompany us all along this course. Our first proposition shows that the critical parameter must be positive as long as the underlying graph $ \mathfrak{g}$ has bounded degree.
\begin{proposition}  \label{prop:lowerperco} Let $ \mathfrak{g}$ be an (infinite connected countable) graph such that $$ \max_{v\in \mathrm{V}( \mathfrak{g})} \mathrm{deg}(v) \leq M.$$ Then we have 
$ \mathbb{P}( \exists \mbox{ infinite cluster in }\mathrm{Perc}( \mathfrak{g},p)) =0$ as long as $p (M-1)< 1$.
\end{proposition}

The proof of this proposition is our first application of the  \textbf{first moment method} which we single out as a lemma: 

\begin{lemma}[First moment method]\label{def:first-moment} Let $X \in \{0,1,2,\dots\}$ be a non-negative integer valued random variable. Then we have
 $$ \mathbb{P}(X \geq 1) \leq \mathbb{E}[X].$$
 \end{lemma}
  \noindent \textbf{One-line proof:}   Since $X \in  \mathbb{Z}_{\geq0}$ we have $ \mathbb{P}(X\geq 1) =\mathbb{E}[ \mathbf{1}_{X>0}]\leq \mathbb{E}[X \mathbf{1}_{X>0}] = \mathbb{E}[X]$. \qed \medskip

\noindent \textbf{Proof of Proposition \ref{prop:lowerperco}.}  Let us consider a reference vertex $v_0$ in $ \mathfrak{g}$ and let $X(p)=  \mathbf{1}_{v_0 \leftrightarrow \infty}$. Our goal is to show that $X(p)=0$  for $p$ small. For this, we shall use the proxy random variables $X_n(p)$ counting the number of self-avoiding paths starting from $v_0$ of length $n$ and made of open edges in $\mathrm{Perc}( \mathfrak{g},p)$. Clearly, since the degree of each vertex in $ \mathfrak{g}$ is bounded above by $M$,  there is at most $ M\cdot (M-1)^{n-1}$ non-backtracking paths of length $n$ starting from $v_0$ in  $ \mathfrak{g}$. Since there are more non-backtracking paths than self-avoiding paths, by independence of the status of the edges we have 
$$ 	\mathbb{E}[ X_n(p)] \leq M \cdot (M-1)^{n-1} \cdot p^n.$$
Lemma \ref{lem:konig} shows that $v_{0}$ is in an infinite cluster if and only if there is a self-avoiding path of arbitrary length starting from $v_{0}$. We deduce that 
 \begin{eqnarray*}  \mathbb{P}( v_0 \leftrightarrow \infty  \mbox{ in }\mathrm{Perc}( \mathfrak{g},p)) & \underset{ \mathrm{Lem.\ } \ref{lem:konig} }{=}& \mathbb{P}( X_n(p) \geq 1, \forall n \geq 1)\\
 & \leq & \inf_{n\geq 1} \mathbb{P}(X_n(p) \geq 1)\\
 & \underset{ \mathrm{First\ Moment }}{\leq} & \inf_{n \geq 1} \mathbb{E}[ X_n(p)] \\ & \leq & \inf_{n \geq 1 }M \cdot (M-1)^{n-1} \cdot p^n. \end{eqnarray*} Hence if $p  (M-1)<1$ the above probability is $0$. By countable union over all $v_0 \in \mathfrak{g}$, the probability that there exists an infinite cluster (at all) is also zero in this regime. \qed

\subsection{Regular $d$-ary tree}
Fix $d \geq 3$. Let us suppose in this section that $ \mathfrak{g}$ is the infinite $(d-1)$-ary tree $ \mathfrak{t}_d$ where all vertices have degree $d$ except for the origin vertex $v_0$ which has degree $d-1$ (so that there are exactly $(d-1)^n$ vertices at distance $n$ from $v_0$). By Proposition \ref{prop:lowerperco} we have $p_c(  \mathfrak{t}_d) \geq 1/ (d-1)$ and in fact this lower bound is sharp:
\begin{proposition} \label{prop:binarythreshold}
We have $p_c(  \mathfrak{t}_d)= \frac{1}{d-1}$.
\end{proposition}

To prove the proposition we shall now use the \textbf{second moment method}:

\begin{lemma}[Second moment method]\label{def:second-moment} Let $X \in \{0,1,2,\dots\}$ be a non-negative integer valued random variable which is not constant equal to $0$. Then we have
 $$ \mathbb{P}(X \geq 1) \geq \frac{\mathbb{E}[X]^{2}}{ \mathbb{E}[X^{2}]}.$$
 \end{lemma}
 \noindent \textbf{One-line proof:} Use Cauchy-Schwarz  $\mathbb{E}[X]^{2}= \mathbb{E}[X \mathbf{1}_{X>0}]^{2} \leq \mathbb{E}[X^{2}] \mathbb{P}(X>0)$. \qed \medskip

\noindent \textbf{Proof of Proposition \ref{prop:binarythreshold}.} Let us focus on the case $d=3$ to ease notation. 
 Recall from the proof of Proposition \ref{prop:lowerperco}  in the case when $v_0$ is the origin of $ \mathfrak{t}_3$ that $X_n(p)$ is the number of open paths in $\mathrm{Perc}( \mathfrak{t}_3,p)$ starting at $v_0$ and reaching level $n$. When $p > 1/2$ we know that $ \mathbb{E}[X_n(p)] = (2p)^n$ tends to infinity, but \textbf{that does not imply} that $X_n(p) \geq 1$ with large probability. To ensure this, we shall compute the second moment of $X_n(p)$:
 \begin{eqnarray*}\mathbb{E}[\left(X_n(p)\right)^2]&=& \mathbb{E}\left[\left(\sum_{x : \mathrm{d_{gr}}(x,v_0)=n}  \mathbf{1}_{v_0 \leftrightarrow x  \mbox{ in } \mathrm{Perc}(  \mathfrak{t}_3,p)} \right)^2\right]\\
 &=& \sum_{x,y : \mathrm{d_{gr}}(x,v_0)=\mathrm{d_{gr}}(y,v_0)=n } \mathbb{P}(v_0 \leftrightarrow x \mbox{ and } v_0 \leftrightarrow y \mbox{ in } \mathrm{Perc}( \mathfrak{t}_3,p))\\
 &=& (2p)^n \left( 1 + p+ 2p^2 + 4 p^3+\dots +2^{n-1}p^n \right) \sim \frac{p}{2p-1} (2p)^{2n},   \end{eqnarray*} as $n \to \infty$ for $p >1/2$.
We thus find that the second moment of $X_n(p)$ is of the same order as the first moment squared. Applying Lemma \ref{def:second-moment} we deduce that $ \mathbb{P}(X_n(p)>0) \geq \mathbb{E}[X_n(p)]^2/\mathbb{E}[X_n(p)^2] \geq \frac{2p-1}{p}$ asymptotically. We deduce as in the proof of Proposition \ref{prop:lowerperco} that 
$$ \mathbb{P}( v_0 \leftrightarrow \infty \mbox{ in } \mathrm{Perc}( \mathfrak{t}_3,p)) = \inf_{n \geq 1} \mathbb{P}(X_n(p)\geq 1) \geq \frac{2p-1}{p}.$$ By the $0-1$-law  there is an infinite cluster in $\mathrm{Perc}( \mathfrak{t}_3,p)$ with probability $1$ when $p >1/2$. 
\qed \medskip

 \begin{exo} Show that there is no infinite cluster in $ \mathrm{Perc}( \mathfrak{t}_{d},p)$ at $p= \frac{1}{d-1}$. \label{exo:11}
 \end{exo}
 
Of course, the knowledgeable reader may have noticed that the open subtree of the origin in $ \mathrm{Perc}( \mathfrak{t}_d,p)$ is a Bienaym\'e--Galton--Watson tree with offspring distribution $ \mathrm{Bin}(d-1,p)$. The phase transition for the existence of an infinite cluster happens when $p(d-1)$, the mean number of children in the Bienaym\'e--Galton--Watson tree, is larger than $1$. We shall study in more details Bienaym\'e--Galton--Watson trees in Part \ref{part:trees} and in particular get a new proof of the above proposition. 

\subsection{Cubic lattice $ \mathbb{Z}^2$}
Let us now focus  on the case when $ \mathfrak{g}$ is the standard Manhattan lattice i.e.~the cubic lattice in dimension $2$. This is the usual grid graph, whose vertex set is $ \mathbb{Z}^2$ and where an edge joins the point $x$ to the point $x+ e$ for $e \in \{(0,\pm1), (\pm1,0)\}$. Let us denote this graph by $ \mathfrak{z}_2$. We know from Proposition \ref{prop:lowerperco} that $p_c \geq 1/3$, but the second moment method does not work well in this setting since two paths of length $n$ may have a very complicated structure. To show that $p_c <1$, we shall rely on another argument specific to planar lattices.

\begin{figure}[!h]
 \begin{center}
 \includegraphics[width=7cm]{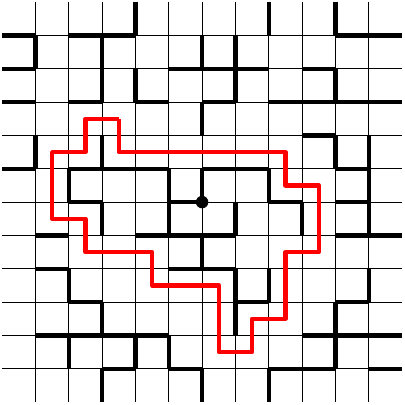}
 \caption{ \label{fig:blocking} If the cluster of the origin is finite, then it is surrounded by a blocking dual self-avoiding cycle of length at least $4$.}
 \end{center}
 \end{figure}

\begin{proposition} \label{prop:Z2} We have 
$ 0 < p_c( \mathfrak{z}_2) <1.$
\end{proposition}
\noindent\textbf{Proof.} The idea is to use plane duality. More precisely, if the cluster of the origin vertex $(0,0)$ is finite in $ \mathrm{Perc}( \mathfrak{z}_2,p)$ this forces the existence of a \textbf{blocking self-avoiding cycle} in the dual graph, see Figure \ref{fig:blocking}.

 Since the dual graph of $ \mathfrak{z}_2$ is $ \mathfrak{z}_2$ itself, there are at most $4\cdot3^{n-1}$ dual cycles of length $n$ starting from the origin, and at most $n \cdot 4 \cdot 3^{n-1}$ such cycles blocking the origin (re-root at its first intersection with the positive origin axis which must be at distance less than $n$). We can now use the first-moment method on these blocking cycles: the probability that there exists such a cycle in the dual graph is upper bounded by its expectation and so by

$$ \mathbb{P}( \exists \mbox{ blocking cycle}) \leq \sum_{n \geq 4} n 4^{n} (1-p)^n.$$
The above sum can be made smaller than $1$ if $p$ is close enough to $1$. In this case we get $ \mathbb{P}( (0,0) \leftrightarrow \infty \mbox{ in }  \mathrm{Perc}( \mathfrak{z}_2,p))>0$ and so $p_c( \mathfrak{z}_2) <1$.
\qed

\begin{remark} In essence, the duality argument shows that the percolation on $ \mathfrak{z}_2$ is self-dual at $p=1/2$ and this is one of the key ingredients to prove that $p_c( \mathfrak{z}_2)= 1/2$ (a result due to Kesten). 
\end{remark}
 
 Since the two-dimensional cubic lattice is included in its higher dimensional analog in a trivial fashion, we deduce that there is a non-trivial phase transition in $ \mathfrak{z}_{d}$ for any $d \geq 2$. The nature of the phase transition in low dimensions  $3,4,5,6\dots$ (and also in dimension $2$ to some extent) is still elusive.
\section{Mean-field regime}
In Part \ref{part:ER}, we will study the case when the underlying graph $ \mathfrak{g}$ is the complete graph $ \mathbb{K}_n$ on $n$ vertices. This is the graph made of the vertices $1,2, \dots , n$ and where there is an edge between two distinct vertices (no loops).  One of the main objects in this course is obtained by studying $ \mathrm{Perc}( \mathbb{K}_{n},p)$ when $p$ may vary with $n$. This random graph, usually referred to as the \textbf{Erd{\H{o}}s--R\'enyi random graph} will be denoted by $G(n,p)$.   This  model was introduced\footnote{Actually this definition of random graph is not really due to Erd{\H{o}}s and R\'enyi  who considered a random graph on $n$ vertices with a fixed number $m \leq  {n \choose 2}$ of edges. However, once conditioned on the number of edges the two models are equivalent and we shall use the name Erd{\H{o}}s--R\'enyi instead of Edgar Gilbert who introduced this variant.} by Erd{\H{o}}s and R\'enyi\footnote{\raisebox{-5mm}{\includegraphics[width=1cm]{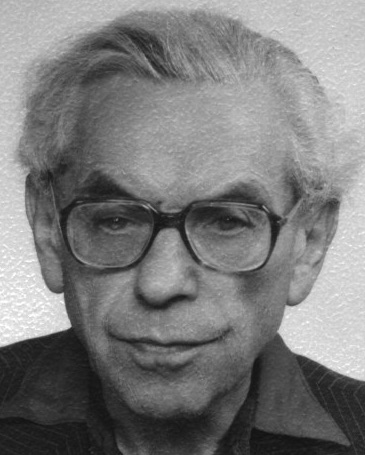}} Paul Erd{\H{o}}s (1913--1996), Hungarian and \raisebox{-5mm}{\includegraphics[width=1cm]{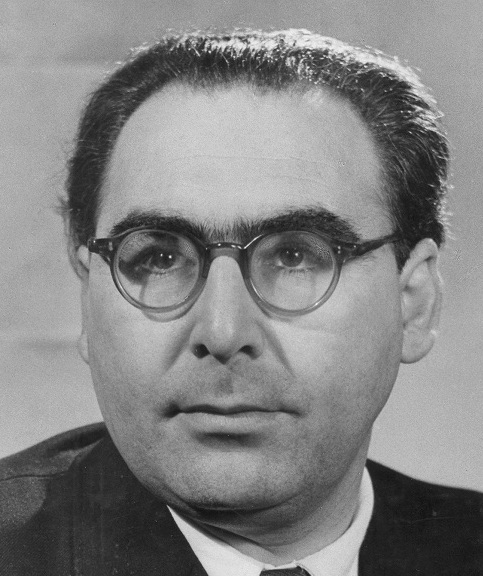}} Alfr\'ed R\'enyi (1921--1970), Hungarian} in 1959 who wanted to probe randomly  a graph with $n$ (labeled) vertices. This random graph model has become ubiquitous in probability and commonly referred to as the ``\textbf{mean field model}''. This means that the initial geometry of the model is trivial: one could permute all the vertices and get the same model. \medskip 

There is a convenient way to couple all these realizations including the case $n= \infty$: consider the \textbf{complete} graph $ \mathbb{K}_{\infty}$ whose vertex set is $ \mathbb{Z}_{>0} = \{1,2,3, \dots \}$ and whose edge set is $ \{ \{i,j\} : i \ne j \in \mathbb{Z}_{>0}\}$ (hence, an edge between any possible pair of distinct vertices). This graph is connected and countable although it is not locally finite. We can then consider for each edge $e= \{i,j\}$ and independent uniform random variable $U_{ e} \in [0,1]$ and set for each $n \in \{1,2, \dots\}$ and $p \in [0,1]$
  \begin{eqnarray} \label{def:erdosrenyicoupled} G(n, p) = \Big( \underbrace{\big\{ i \in \mathbb{Z}_{>0} : i \leq n\big\}}_{ \mathrm{vertex\ set}},  \underbrace{\big\{ \{i,j\}  : i,j \leq n \mbox{ and } U_{\{i,j\}} \leq p\big\}}_{ \mathrm{edge\ set}} \Big).  \end{eqnarray}
Once again, studying the properties of $G(n, p_{n})$ for finite $n'$s and varying parameter $p \equiv p_{n}$ will be the subject of the whole Part \ref{part:ER}. To conclude this chapter, let us  focus on the case when $n = \infty$ and $p>0$. It should be clear to the reader that $G(\infty, p)$ is almost surely connected, but the following result might come as a surprise:

\begin{theorem}[Erd{\H{o}}s--R\'enyi (1963)] \label{thm:rado} For any $p,p' \in (0,1)$ almost surely $G( \infty,p)$ and $ G(\infty,p')$ are equivalent. In particular, the equivalence class of non trivial Bernoulli bond-percolation on $ \mathbb{K}_\infty$ is almost surely constant: this is the \textbf{Rado graph}.
\end{theorem}
 Recall that two graph $ \mathfrak{g}, \mathfrak{g}'$ are equivalent if there is a bijection of their vertex sets which preserve the adjacency properties.  The proof of the result is easy once we know the following characteristic property of the  Rado\footnote{\raisebox{-5mm}{\includegraphics[width=1cm]{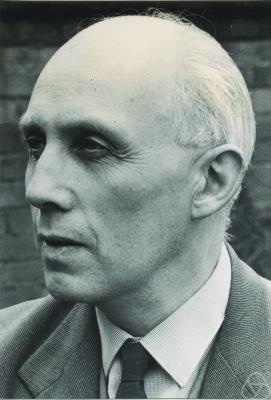}} Richard Rado (1906--1989), German } graph (over the vertex set $ \mathbb{ Z}_{>0}$): it is the only (equivalence class of) countable simple graph such that for any finite disjoint subsets $U,V \subset \mathbb{Z}_{>0}$, there exists $v$ outside of $U$ and $V$ such that $v$ is neighbor to all vertices in $U$ and none of $V$. The previous property is called the \textbf{extension property} and can be used to prove by induction that any finite or countable graph $ \mathfrak{g}$ can be embedded inside the Rado graph. See the excellent wikipedia article on the Rado graph for more details. \medskip 
 
 \noindent \textbf{Proof of Theorem \ref{thm:rado}.} Let us check that $G( \infty,p)$ with $p \in (0,1)$ almost surely satisfies the extension property.  Fix $U,V \subset  \mathbb{Z}_{>0}$. For $v  \in   \mathbb{Z}_{>0} \backslash (U \cup V)$ the probability that $v$ is connected to all vertices of $U$ and none of $V$ is $ p^{\# U} (1-p)^{\# V}>0$. By independence and the Borel--Cantelli lemma, there exists $v$ connected to all the vertices of $U$ and none of $V$ in $G( \infty,p)$ with probability one. The property holds true for all finite subsets $U,V \subset \mathbb{Z}_{>0}$ simultaneously by countable union. \qed  \bigskip

\noindent\textbf{Bibliographical notes.} Percolation theory is a very broad and vivid area in nowadays probability theory,  \cite{grimmett1997percolation,werner2007lectures,duminil2017sixty}. When the underlying graph has strong geometric constraints (e.g.~the cubic lattices in $ \mathbb{Z}^{d}$ for $d \geq 2$) then the study of the phase transition and in particular of the critical behavior is still a challenge for mathematicians. Theorem \ref{thm:rado} is proved in \cite{erdos1963asymmetric}. For more about the links between the geometry of the graph and the behavior of Bernoulli percolation we advise the reading of the influential paper \cite{benjamini1996percolation}. \medskip 

\noindent{\textbf{Hints for Exercises.}}\ \\
\noindent Exercise \ref{exo:11}:  With the notation above the exercise, show that $ \mathbb{E}[X_{n}(p) \mid X_{n}(p) \geq 1] \to \infty$ as $n \to \infty$ when $p = \frac{1}{d-1}$.

\part[Bienaym\'e--Galton--Watson trees]{Bienaym\'e-Galton-Watson trees 
                             \\ \\ 
\label{part:trees}
  \begin{center}
                     \begin{minipage}[l]{15cm}
       \normalsize 
       In this part we study the model of Bienaym\'e--Galton--Watson (BGW) tree, or discrete branching process, modeling the genealogy of an asexual population where individuals reproduce independently of each other according to the same offspring distribution. The main tool to study such objects is their encodings by one-dimensional random walks. 
                     \end{minipage}
                  \end{center}
                  \vspace{1cm}
                   \begin{figure}[!h]
  \begin{center}
 \includegraphics[height=8cm]{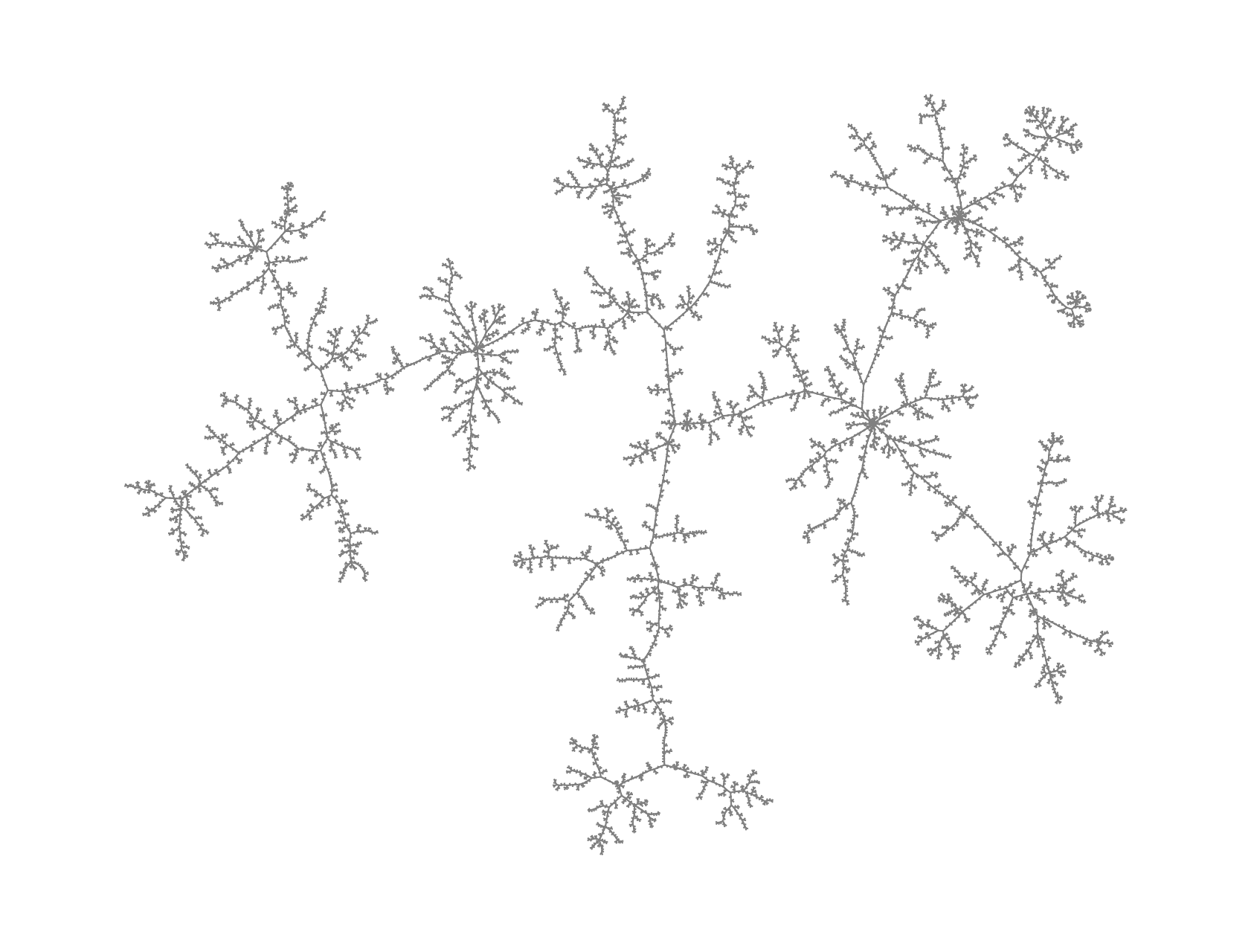}
         \caption{A large critical Bienaym\'e--Galton--Watson tree with finite variance \label{fig:CRT}}
  \end{center}
  \end{figure}
                }

\chapter{One-dimensional random walks}
\hfill Back to basics.

\label{chap:generalRW}
In this chapter we consider the following object:
                       \begin{definition}[One-dimensional random walk] Let $\mu = ( \mu_k : k \in \mathbb{Z})$ be a probability distribution on $ \mathbb{Z}$ with $\mu(  \mathbb{Z}_{>0})>0$ as well as $\mu( \mathbb{Z}_{<0}) >0$. Consider $X_{1}, X_{2}, \ldots$ i.i.d.\,copies of law $\mu$ which we see as the increments of the process $(S) \equiv (S_{n}: n \geq 0)$ on $ \mathbb{Z}$ defined as follows : $S_{0}=0$ and for $n \geq 1$
$$ S_{n} = X_{1}+ \cdots + X_{n}.$$
We say that $(S)$ is a one-dimensional random walk with step distribution $\mu$ (or $\mu$-random walk for short).
\end{definition}

                   \begin{figure}[!h]
  \begin{center}
 \includegraphics[height=4cm]{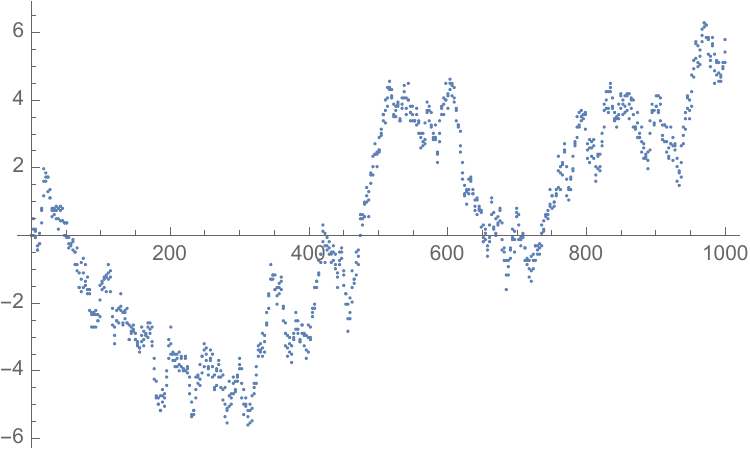}  \hspace{1cm}   \includegraphics[height=4cm]{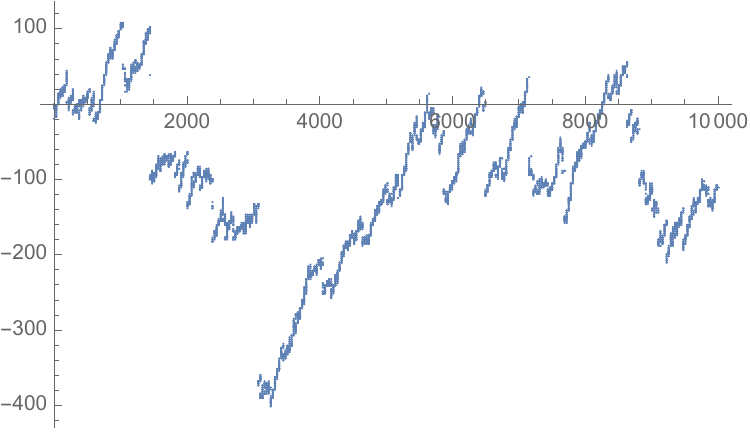}
  \caption{Two samples of one-dimensional random walks with different step distributions. The first one seems continuous at large scales whereas the second one displays macroscopic jumps.}
  \end{center}
  \end{figure}
Notice that we restrict (for simplicity)  to the \textbf{lattice case} by demanding that the support of $\mu$ be included in $ \mathbb{Z}$ and that we excluded the monotone situation since the support of $\mu$ contains both positive and negative integers. Of course, the behavior of a one-dimensional random walk depends on the step distribution $ \mu$ in a non-trivial way as we will see. We first recall the general background on such objects before moving to  \textbf{skip-free random walks} which can only make negative jumps of size  $-1$ and which will be used in the next chapters to study random trees and graphs. 
\section{General theory}
In this section we gather a few general results on one-dimensional random walks and start with the applications of discrete Markov chain theory since   a one-dimensional random walk is clearly a very particular case of Markov chain in discrete time with a discrete state space. 
\subsection{Reminder on Markov chains}
We start with the parity consideration:
\begin{proposition} The Markov chain $(S_n: n \geq 0)$ is 
\begin{itemize}
\item \textbf{irreducible} if $ \mathrm{Supp}(\mu)$ is not included in $ a \mathbb{Z}$ for some $a>1$,
\item It is furthermore \textbf{aperiodic} if $ \mathrm{Supp}(\mu)$ is not included in $ b + a \mathbb{Z}$ for some $a>1$ and $ b \in \mathbb{Z}$.
\end{itemize}
\end{proposition}
\noindent \textbf{Proof.} Using the fact that the walk is not monotone, it is an exercise to check that the chain can come back to $0$ with positive probability and so the set of integers accessible by the chain starting from $0$ is a subgroup of $ \mathbb{Z}$. Writing Bezout relation we can find non-negative integers $\alpha_{1}, \dots , \alpha_{k}, \alpha'_{1}, \dots , \alpha'_{k}$ and $\ell_{1}, \dots , \ell_{k},\ell'_{1}, \dots , \ell'_{k'} \in \mathrm{Supp}( \mu)$ so that
$$ \alpha_{1}\ell_{1} + \dots + \alpha_{k} \ell_{k }= \alpha'_{1}\ell'_{1} + \dots + \alpha'_{k'} \ell'_{k' } +\mathrm{gcd}( \mathrm{Supp}(\mu)).$$
Hence, by the above consideration $\mathrm{gcd}( \mathrm{Supp}(\mu))$ is an accessible value for the walk and the first point is proved.
For the second point, notice that if $\mathrm{Supp}(\mu) \subset b + a \mathbb{Z}$ then $$S_n \equiv bn \quad [a]$$ and so the chain cannot be aperiodic if $a \notin\{-1,+1\}$. If $\mathrm{Supp}(\mu)$ is not included in $b + a \mathbb{Z}$ for $|a|>1$, then we have $ \mathrm{gcd} (\mathrm{Supp}( \mu) - k_0) =1$ where $k_0$ is any integer in the support of $ \mu$. We pick $k_0 \in \mathrm{Supp}( \mu)$ such that the measure $\tilde{\mu}(\cdot) = \mu(k_0 + \cdot)$ does not put all its mass on $ \mathbb{Z}_{\geq0}$ nor on $ \mathbb{Z}_{\leq 0}$. It is possible since otherwise $ \mathrm{Supp}( \mu) = \{ \alpha,\beta\}$ with $\alpha < 0 <\beta$ and so $ \mathrm{Supp}(\mu) \subset	\alpha +(\beta-\alpha) \mathbb{Z}$ which is excluded by hypothesis. Then, by the first point of the proposition, we can find $n_0\geq 0$ large enough so that a $\tilde{\mu}$-random walk satisfies $ \mathbb{P}( \widetilde{S}_{n_0} = k_0) >0$ which means that 
$$ \mathbb{P}(S_{n_0} = (n_{0}+1) k_0) >0.$$
Combining this with the trivial point $ \mathbb{P}(S_{n_0+1} = (n_0+1) k_0)  \geq (\mu_{k_{0}})^{{n_{0}+1}} >0$ we deduce that the integer $(n_0+1) k_0$ is accessible both at time $n_0$ and time $n_0+1$ for the chain. By standard results on Markov chains this implies aperiodicity.  \qed \medskip 

\begin{example} Simple random walk on $ \mathbb{Z}$ with $ \mathbb{P}(S_1= \pm 1)= \frac{1}{2}$ is irreducible but not aperiodic. \end{example}

The counting measure on $ \mathbb{Z}$ is clearly an invariant measure for any $\mu$-random walk (beware, it is not usually reversible, and it might not be the only invariant measure up to multiplicative constant in the transient case). Due to homogeneity of the process, the Markov property takes a nice form in our setup: as usual $ \mathcal{F}_{n} = \sigma ( X_{1}, \dots , X_{n})$ is the natural filtration generated by the walk $(S)$ up to time $n$ and a  \textbf{stopping time} is a random variable $\tau \in \{0,1,2, \dots \} \cup\{\infty\}$ such that for each $n \geq 0$ the event
$\{ \tau =n \}$ is measurable with respect to $ \mathcal{F}_{n}$. 
\begin{proposition}[Strong Markov property]  \label{prop:strongmarkov}If $\tau$ is a stopping time then conditionally on $\{\tau < \infty\}$ (implicitly of positive probability) the process $(S^{(\tau)}_{n})_{n \geq 0}=(S_{\tau+n}- S_{\tau})_{n \geq 0}$ is independent of $( S_{n})_{ 0 \leq n\leq \tau}$ and is distributed as the initial walk $(S_{n})_{n \geq 0}$.
\end{proposition}
\noindent \textbf{Proof.} Let $f,g$ be two positive measurable functions and let us compute 
 \begin{eqnarray*} \mathbb{E}\left[f\left((S_{n})_{ 0 \leq n \leq \tau}\right)g\left((S^{(\tau)}_{n})_{n \geq 0}\right) \mathbf{1}_{\tau < \infty}\right] &\underset{\tau < \infty}{=}& \sum_{ k=0}^{\infty} \mathbb{E}\left[  \mathbf{1}_{\tau=k} f\left((S_{n})_{ 0 \leq n \leq k}\right)g\left((S^{(k)}_{n})_{n \geq 0}\right) \right] \\
 & \underset{\mathrm{indep}.}{=} & \sum_{ k=0}^{\infty} \mathbb{E}\left[  \mathbf{1}_{\tau=k} f\left((S_{n})_{ 0 \leq  n \leq k}\right)\right] \mathbb{E}\left[g\left((S^{(k)}_{n})_{n \geq 0}\right) \right]\\
  & \underset{\mathrm{stat}.}{=} & \sum_{ k=0}^{\infty} \mathbb{E}\left[  \mathbf{1}_{\tau=k} f\left((S_{n})_{ 0 \leq n \leq k}\right)\right] \mathbb{E}\left[g\left((S_{n})_{n \geq 0}\right) \right]\\
  &=&  \mathbb{E}\left[  f\left((S_{n})_{ 0 \leq n \leq \tau}\right) \mathbf{1}_{\tau < \infty}\right] \mathbb{E}\left[g\left((S_{n})_{n \geq 0}\right) \right]. \end{eqnarray*}
This proves the proposition. \qed 
  
  \subsection{$0-1$ laws}
  In the study of random walks, one often uses $0-1$ laws when dealing with asymptotic events such as $\{S_{n} \to \infty\}$. The most well-known of such laws is Kolmogorov's\footnote{\raisebox{-5mm}{\includegraphics[width=1cm]{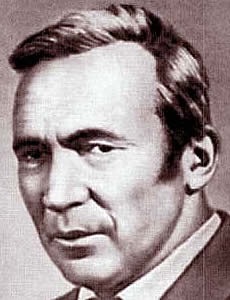}}Andre\"i Nikola\"ievitch Kolmogorov (1903--1987), Russian} $0-1$ law which states that if $(X_{i})_{i \geq 0}$ are independent random variables (not necessarily identically distributed), then any event $ \mathcal{A}$ measurable with respect to $\sigma( X_{i}: i \geq 0)$ and which is independent of $(X_{1}, \dots , X_{n_{0}})$ for any $n_{0}$ has measure $ \mathbb{P}( \mathcal{A}) \in \{0,1\}$. Let us present a stronger version of Kolmogorov $0-1$ law in the case of \textbf{i.i.d.}\ increments. This $0-1$-law, due to Hewitt \& Savage\footnote{\raisebox{-5mm}{\includegraphics[width=1cm]{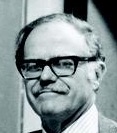}} Edwin Hewitt (1920--1999), \raisebox{-5mm}{\includegraphics[width=1cm]{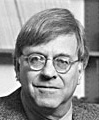}} Leonard Savage (1917--1971), American}, has many applications in the random walk setting: 
  
  \clearpage 
  
  \begin{theorem}[Hewitt--Savage exchangeable $0-1$ law]\label{thm:hewitt-savage}Let $(X_{i})_{ i \geq 1}$ be a sequence of independent and \underline{identically distributed} random variables with values in a Polish space $(E,d)$. Suppose that $ \mathcal{A}$ is a measurable event with respect to $ \sigma (X_{i}: i \geq 1)$ which is invariant (up to negligible events) by any  permutation of the $(X_{i} : i \geq 1)$ with finite support. Then $ \mathbb{P}( \mathcal{A}) \in \{0,1 \}$.
  \end{theorem}
  
  \noindent \textbf{Proof.} Let $ \mathcal{A} \in \sigma (X_{i} : i \geq 1)$ be invariant by any permutations of the $X_{i}$ with finite support (i.e.~only finitely many terms are permuted). By a standard measure-theory argument (see \cite[Lemma 3.16]{Kal07})  one can approximate $ \mathcal{A}$ by a sequence of  events $ \mathcal{A}_{{n}} \in \sigma( X_{1}, \dots , X_{n})$ in the sense that 
  $$ \mathbb{P}( \mathcal{A} \Delta \mathcal{A}_{n})  \xrightarrow[n\to\infty]{} 0.$$
By definition, any event $ \mathcal{E} \in \sigma (X_{i} : i \geq 1)$ can be written $ \mathcal{E} = \mathbf{1}_{ (X_{i} : i \geq 1) \in \tilde{\mathcal{E}}}$ where $ \tilde{\mathcal{E}}$ is an event of the Borel cylindric $\sigma$-field on $ E^{ \mathbb{Z}_{> 0}}$. We can thus consider the function $\psi_{n}$ acting on events  $  \mathcal{E} \in \sigma( X_{i}: i \geq 1)$ by swapping $X_{1}, \dots , X_{n}$ with $ X_{n+1}, \dots , X_{2n}$~i.e.~
  $$ \psi_{n}( \mathcal{E}) = \mathbf{1}_{X_{n+1}, \dots ,X_{2n}, X_{1}, \dots , X_{n}, X_{2n+1}, \dots   \in \tilde{ \mathcal{E}}} \in \sigma( X_{i}: i \geq 1).$$  Since the $X_{i}$ are i.i.d.\, we have $ \mathbb{P}( \psi_{n} (\mathcal{E})) = \mathbb{P} ( \mathcal{E})$ for any event $ \mathcal{ E}$ and also $ \psi_{n}( \mathcal{A}_{n})$ is independent of $ \mathcal{A}_{n}$. Using this we have 
$$0 \xleftarrow[n \to \infty]{} \mathbb{P}(\mathcal{A} \Delta  \mathcal{A}_{n})=\mathbb{P}(\psi_{n}( \mathcal{A} \Delta  \mathcal{A}_{n})) =  \mathbb{P}( \psi_{n}(\mathcal{A}) \Delta \psi_{n}( \mathcal{A}_{n})) = \mathbb{P}(  \mathcal{A} \Delta \psi_{n} (\mathcal{A}_{n})).$$
  We deduce that $ \mathcal{A}$ is both very well approximated by $ \mathcal{A}_{n}$ but also  by $ \psi_{n}( \mathcal{A}_{n})$. Since the last two events are independent we deduce that $ \mathbb{P}( \mathcal{A}) \in \{0,1\}$ because  $$ \mathbb{P}( \mathcal{A}) = \lim_{n \to \infty} \mathbb{P}( \mathcal{A}_{n} \cap \psi_{n} ( \mathcal{A}_{n})) \underset{ \mathrm{indept.}}{=} \lim_{n \to \infty} \mathbb{P}( \mathcal{A}_{n}) \mathbb{P}( \psi_{n}( \mathcal{A}_{n})) \underset{ \mathrm{i.d.}}{=} \lim_{n \to \infty}\mathbb{P}( \mathcal{A}_{n})^{2}  = \mathbb{P}( \mathcal{A})^{2}. $$ \qed \medskip 
  
  \begin{example} \label{ex:O1} If $A \in \mathbb{Z}$ is a measurable subset and $(S)$ a one-dimensional random walk with i.i.d.~increments, we write 
 \begin{eqnarray} \label{eq:IA} \mathcal{I}_{A}:=  \sum_{n = 0}^{\infty} \mathbf{1}_{S_{n} \in A}.  \end{eqnarray} Then the commutativity of $ \mathbb{Z}$ (sic!) shows that the event $ \{  \mathcal{I}_{A} = \infty\}$ is invariant under finite permutations of the $X_{i}$'s (indeed any finite permutation leaves $S_{n}$ invariant for large $n$);  hence it has probability $0$ or $1$. Notice that this cannot be deduced directly from Kolmogorov's $0-1$ law.
 \end{example}
   
 \subsection{Asymptotic behavior}
 Let us denote $$ \overline{S} := \limsup_{n \to \infty} S_n\quad \mbox{ and } \quad  \underline{S} := \liminf_{n \to \infty} S_n.$$ For any $k \in \mathbb{Z}$, the probability that $ \overline{S}$ or $\underline{S}$ is equal to $k$ is null, since otherwise the walk would take the value $k$ an infinite number of times with positive probability: by the classification of states, the walk would be recurrent and so would visit the whole subgroup $  \mathrm{gcd}( \mathrm{Supp}(\mu)) \cdot \mathbb{Z}$ almost surely which is incompatible with a finite $\limsup$ or $\liminf$ (recall that $\mu( \mathbb{Z}_{>0})$ and $\mu( \mathbb{Z}_{<0})$ are positive). This motivates the following definition:
\begin{definition} \label{def:prop:oscillates} A (non-trivial) one-dimensional random walk $(S)$ with i.i.d.\ increments falls into exactly one of the three following categories:
\begin{enumerate}[(i)]
\item Either $\overline{S} = \underline{S} = \infty$, that is $ S_{n} \xrightarrow[n\to\infty]{} \infty$ in which case $(S)$ is said to \textbf{drift towards $\infty$},
\item Or $\overline{S} = \underline{S} = -\infty$, that is $ S_{n} \xrightarrow[n\to\infty]{} - \infty$  in which case $(S)$ is said to \textbf{drift towards $-\infty$},
\item Or $(S)$ \textbf{ oscillates} i.e.~$\limsup_{n \to \infty} S_{n}= +\infty$ and $\liminf_{n \to \infty} S_{n} = -\infty$ almost surely.
\end{enumerate}
\end{definition} 

When a random walk drifts, it is obviously transient, but in the oscillating case, it may be transient or recurrent, see Theorem \ref{prop:discretetransient} for examples.

\begin{remark}  \label{rek:>0pos>0}If the random walk $(S)$ drifts towards $+\infty$, then $$ \mathbb{P}( S_i >0 : \forall i \geq 1)>0.$$ Indeed, if we had $\mathbb{P}( S_i >0 : \forall i \geq 1) =0$ then the stopping time $ \theta = \inf\{ i \geq 1 : S_i \leq 0\}$ would be almost surely finite. Using (iterations of) the Markov property this would imply that $(S)$ visits $ \mathbb{Z}_{\leq 0}$ infinitely often a.s.\  which contradicts the fact that $(S)$ drifts to $+ \infty$.
\end{remark}

\section{Walks with finite mean and the law of large numbers}
In this section we examine the particular case when $\mu$ has finite mean and show  that the walk is recurrent whenever it is centered, otherwise it is transient and drifts. It will be a good opportunity to wiggle around the strong and weak laws of large numbers. We will see in the next chapter a quick proof (Lemma \ref{lem:LLN}) of the strong law of large numbers based on a path transformation called \textbf{duality}.

\subsection{Recurrence/transience}

 Recall that a random walk $(S)$ is recurrent iff one of the following equivalent conditions is satisfied
 \begin{eqnarray} \label{eq:cnsrecurrence}\mathbb{P}(\exists n >0: S_{n}=0)=1 \iff \mathbb{E}\left[ \sum_{n = 0}^{\infty} \mathbf{1}_{S_{n}=0}\right] =\infty \iff  \liminf_{ n \to \infty} |S_{n}| < \infty \quad a.s.  \end{eqnarray}

\begin{theorem}[Dichotomy for walks with finite mean] \label{thm:rec0}Suppose $ \mathbb{E}[|X_{1}|] < \infty$ then 
\begin{enumerate}[(i)]
\item If  $\mathbb{E}[X_{1}] \ne 0$ then $(S)$ is transient and drifts,
\item otherwise if $ \mathbb{E}[X_{1}]=0$ then $(S)$ is recurrent.
\end{enumerate}
\end{theorem}
\noindent \textbf{Proof.} The first point $(i)$ is easy since by the strong law of large numbers we have $ n^{{-1}} S_{n} \to \mathbb{E}[X_{1}]$ almost surely: when $ \mathbb{E}[X_{1}] \ne 0$ this automatically implies that $(S)$ drifts towards $\pm \infty$ depending on the sign of $ \mathbb{E}[X]$.

In the second case we still use the law of large numbers to deduce that $S_{n}/n \to 0$ almost surely as $n \to \infty$. This implies that for any $\varepsilon >0$ we have $ \mathbb{P}(|S_{n}| \leq \varepsilon n) \to 1$ as $n$ tends to infinity. In particular, we have 
 \begin{eqnarray}  \mathbb{E}\left[\sum_{i=0}^{\infty}  \mathbf{1}_{|S_{i}| \leq \varepsilon n} \right] \underset{ i \leq n}{\geq} \sum_{i=0}^{n}  \mathbb{P} \left(|S_{i}| \leq \varepsilon i\right )  \underset{ \mathrm{Ces\`aro}}{ \geq} \frac{n}{2},  
 \label{eq:contradict}   \end{eqnarray}
 eventually. We claim that this inequality is not compatible with transience. Indeed, according to  \eqref{eq:cnsrecurrence}, if the walk $(S)$ is transient then for some constant $C>0$ we have 
$$ \mathbb{E}\left[ \sum_{i=0}^{\infty} \mathbf{1}_{S_{i}=0} \right]	\leq C.$$
If $ k \in \mathbb{Z}$, applying the strong Markov property at the stopping time $\tau_k = \inf\{i \geq 0 : S_{i} = k\}$ we deduce that 
$$ \mathbb{E}\left[ \sum_{i=0}^{\infty} \mathbf{1}_{S_{i} = k} \right] = \mathbb{P}( \tau_k < \infty) \mathbb{E}\left[ \sum_{i=0}^{\infty} \mathbf{1}_{S_{i}=0} \right]	\leq C.$$
Hence, if the walk were transient we would have 
$$\mathbb{E}\left[\sum_{i=0}^{\infty}  \mathbf{1}_{|S_{i}| \leq \varepsilon n}\right ] \leq \sum_{ - \varepsilon n \leq k \leq \varepsilon n} \mathbb{E}\left[ \sum_{i =0}^{\infty} \mathbf{1}_{S_{i} = k}\right] \leq 2 \varepsilon n C,$$ which contradicts  \eqref{eq:contradict} for  $\varepsilon>0$ small enough. Hence the walk cannot be transient. \qed \medskip

Notice that we only use the weak law of large numbers to deduce recurrence: any one-dimensional random walk $(S)$ for which $S_{n}/n \to 0$ in probability is recurrent. There are examples where the step distribution is not integrable, see Exercise \ref{exo:nolfgn}. The theorem above can be seen as a particular example of the Kesten--Spitzer--Whitman theorem (see \cite[Chapter I]{Spi76}) saying that a random walk with independent increments on a group is transient if and only if its range (i.e.~the number of visited vertices) grows linearly with time.

\subsection{Wald's equality}
\begin{theorem}[Wald equality] \label{thm:wald}\noindent Suppose $ \mathbb{E}[|X_{1}|] < \infty$. Let $\tau$ be  a stopping time with finite expectation. Then we have 
$$ \mathbb{E}[\tau] \cdot \mathbb{E}[X_{1}] = \mathbb{E}[S_{\tau}].$$
\end{theorem}

\noindent \textbf{Proof with martingales.} We present a first proof based on martingale techniques. If we denote by $m$ the mean of $X_{1}$ then clearly the process $( S_{n}- nm)_{n \geq 0}$ is a martingale for the canonical filtration $ \mathcal{F}_{n} = \sigma ( X_{1}, \dots, X_{n})$. By the optional sampling theorem we deduce that 
 \begin{eqnarray} \label{eq:butbut} \mathbb{E}[S_{n \wedge \tau}] = m \mathbb{E}[n \wedge \tau].  \end{eqnarray}
Since $\tau$ is almost surely finite, we can let $n \to \infty$ and get by monotone convergence that the right hand side tends to $ m\mathbb{E}[\tau]$. However, to deduce that the left hand side also converges towards $ \mathbb{E}[S_{\tau}]$ one would need a domination... To get this, the trick is to reproduce the argument with the process  $$ Y_{n}=\sum_{i =1}^{n}|X_{i}| - \tilde{m} n,$$ where $\tilde{m} = \mathbb{E}[|X_{1}|]$. Then $(Y_{n})_{n \geq 0}$ is again a martingale for the filtration $( \mathcal{F}_{n})$. Notice that $Y_{n}$ is also a martingale for its own filtration but the previous statement is stronger. We can then apply the optional sampling theorem again for $n \wedge \tau$ and use \textit{monotone} convergence on both sides to get that $ \mathbb{E}[\sum_{i=1}^{\tau}|X_{i}|] = \tilde{m} \mathbb{E}[\tau]$. Clearly the variable $\sum_{i=1}^{\tau}|X_{i}|$ dominates all variables $S_{n \wedge \tau}$ for $n \geq 0$.
One can then use this domination to prove convergence of the left-hand side in \eqref{eq:butbut}. \qed \medskip 

We now give a second proof of Wald's identity based on the less well-known \emph{converse} to the strong law of large numbers (Lemma \ref{lem:recilfgn}):

\noindent \textbf{Proof of Wald's identity with the law of large numbers.}  The idea is to iterate the stopping rule. Let $0 = \tau_{0} \leq \tau=\tau_{1} \leq \tau_{2} \leq \tau_{3} \leq \cdots$ be the successive stopping times obtained formally as 
$$ \tau_{i+1}=\tau_{i+1}\left(S_{n} : n \geq 0\right) = \tau_{i}(S_{n} : n \geq 0) + \tau(S_{n+\tau_{i}}-S_{\tau_i} : {n \geq 0}),$$ for $i \geq 1$ where we see here $\tau$ as a measurable function of the underlying walk\footnote{In particular, if $\tau$ were a stopping time of a larger filtration $ (\mathcal{G}_n : n \geq 0)$ than the filtration $ (\mathcal{F}_n : n \geq 0)$ generated by the walk, then we could not write the previous display in full generality.}.
In particular since $\tau < \infty$ a.s., we deduce by successive applications of the Markov property (Proposition \ref{prop:strongmarkov}) that $\tau_i < \infty$ for all $i \geq 0$ a.s. and that $$(\tau_{i+1}-\tau_{i} ; S_{\tau_{i+1}}-S_{\tau_{i}})_{i \geq 0} \quad \mbox{ are i.i.d.~of law}  \quad (\tau, S_{\tau}).$$ Since $ \tau$ has finite expectation by assumption, the law of large numbers gives 
$$ \frac{\tau_{i}}{i}  \xrightarrow[i\to\infty]{a.s.} \mathbb{E}[\tau].$$
In particular $\tau_i \to \infty$ almost surely and by the law of large numbers applied on the walk $(S)$ (recall that $ \mathbb{E}[|X|] < \infty$) we deduce that 
$$ \frac{S_{\tau_{i}}}{i} = \frac{S_{\tau_{i}}}{\tau_{i}}  \cdot \frac{\tau_{i}}{i}  \xrightarrow[i\to\infty]{a.s.} \mathbb{E}[X_{1}] \cdot \mathbb{E}[\tau].$$
We then use the \textit{converse} to the law of large numbers (Lemma \ref{lem:recilfgn}) to deduce that $ S_{\tau}$ has finite expectation and equal to $\mathbb{E}[\tau] \cdot \mathbb{E}[X_{1}]$ as claimed by Wald\footnote{\raisebox{-5mm}{\includegraphics[width=1cm]{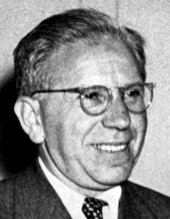}} Abraham Wald (1902--1950), American}.\qed

\begin{lemma}  \label{lem:recilfgn} Let $(S_{i})_{i \geq 0}$ be a one-dimensional random walk with i.i.d.~increments $X_{i}$ of law $\mu$ on $ \mathbb{R}$. Suppose that 
$$ \frac{S_{n}}{n} \xrightarrow[n\to\infty]{a.s.} \mathcal{X},$$
for some  finite (a priori random) variable $  \mathcal{X} \in \mathbb{R}$. Then $\mu$ has a first moment and  $ \mathcal{X} = \mathbb{E}[X]$ a.s.
\end{lemma}
\noindent \textbf{Proof of the lemma.} Suppose that $ n^{-1}S_{n}$ converges almost surely as $n$ goes to infinity to an a priori random but finite variable $ \mathcal{X}$. In particular, we have the almost sure convergence $$ \frac{X_{n}}{n}  =  \frac{S_{n}}{n} - \frac{n-1}{n}\cdot \frac{S_{n-1}}{n-1} \xrightarrow[n\to\infty]{a.s.} \mathcal{X}- \mathcal{X} = 0.$$
We deduce that the event $\{ |X_{n}| \geq n\}$ happens only finitely many times a.s., and since those events are independent, by the second Borel--Cantelli lemma we deduce that 
$$ \infty > \sum_{n \geq 1} \mathbb{P}(|X_{n}| \geq n) \underset{ \mathrm{i.d.}}{=} \sum_{n \geq 1} \mathbb{P}(|X| \geq n) = \mathbb{E}[|X|],$$ where the last equality is a standard exercise using Fubini. We deduce that $ \mathbb{E}[|X|]< \infty$ and by the strong law of large numbers we have $  \mathcal{X} = \mathbb{E}[X]$ a.s. \qed 

\medskip 
Beware, the converse of the weak law of large number does not hold:
\begin{exo}[No converse to the weak law of large numbers]  \label{exo:nolfgn}Let $(S_{n})_{n \geq 0}$ be a one-dimensional random walk with symmetric step distribution $\mu_{k} = \mu_{-k}$ satisfying $$\mu_{k} \sim  \frac{1}{k^{2} \log k}, \quad \mbox{ as } k \to \infty.$$
Show that $ n^{-1}S_{n} \to 0$ in probability, but not almost surely, as $n \to \infty$.
\end{exo}
In words, if $(S)$ is random walk as in the previous exercise, then for most scales we have $S_{n} =o(n)$ whereas there exists exceptional scales where $|S_{n}| >>n$ due to an unlikely event of a large jump during this scale. In fact, an almost sure convergence can always be realized as a ``uniform'' convergence in probability in the following sense:
\begin{exo}[Almost sure convergence is a uniform convergence in probability] Let $ \mathcal{X}_{n}, \mathcal{X}$ be random variables taking values in a Polish space $(E, \mathrm{d})$. Show that 
  \begin{eqnarray*} \mathcal{X}_{n} \xrightarrow[n\to\infty]{( \mathbb{P})} \mathcal{X} &\iff& \mathrm{d}( \mathcal{X}_{n}, \mathcal{X}) \xrightarrow[n\to\infty]{( \mathbb{P})} 0\\
 \mathcal{X}_{n} \xrightarrow[n\to\infty]{a.s.} \mathcal{X} &\iff& \sup_{k \geq n} \mathrm{d}( \mathcal{X}_{k}, \mathcal{X}) \xrightarrow[n\to\infty]{( \mathbb{P})} 0.  \end{eqnarray*}
\end{exo}

\begin{exo}[Independence is crucial!] \label{exo:22} Construct a random walk $(S_{n} : n \geq 0)$ whose increments have the same law $ \frac{2}{3}\delta_{1} + \frac{1}{3} \delta_{-1}$ (but not independent) and so that $(S)$ is recurrent.
\end{exo}

\section{Heavy tailed random walks}
We will see below (Theorem \ref{thm:recfourier1}) a powerful recurrence criterion based on the Fourier transform, but let us use a probabilistic argument to construct transient yet oscillating (symmetric) random walks with heavy tails.

\begin{theorem} \label{prop:discretetransient}Let $\mu$ be a symmetric step distribution (i.e.~$\mu_k = \mu_{-k}$ for $k \in \mathbb{Z}$) satisfying 
  \begin{eqnarray} \mu_k \sim  \mathrm{c} \ k^{-\alpha}, \quad \mbox{ as }k \to \infty,  \label{eq:tailsymmetric}\end{eqnarray}
for some $ \mathrm{c} >0$ and with $\alpha \in (1,\infty)$. Then the walk is recurrent if and only if $\alpha \geq 2$.
 \end{theorem}
 \begin{remark} Notice that since $\mu$ is symmetric, the walk $(S)$ automatically oscillates. The case $\alpha=2$ is critical as already hinted in Exercise \ref{exo:nolfgn}.
 \end{remark}
 
 \noindent \textbf{Proof.} \textsc{When $\alpha >2$} the increments have a finite mean and $ \mathbb{E}[X_1]=0$ by symmetry so the result follows from Theorem \ref{thm:rec0}. \\
 \textsc{ Let us treat the case $\alpha \in (1,2)$} and show that $(S)$ is transient i.e.
  \begin{eqnarray} \label{eq:goalsumtrans}\sum_{n \geq 0 }\mathbb{P}( S_{n} = 0) < \infty.  \end{eqnarray}
The idea is to use the randomness produced by a \textit{single big jump} of the walk to produce a upper bound on $ \mathbb{P}(S_n =0)$. More precisely, let us introduce the stopping time  $ \tau_{n} = \inf \{  i \geq 1:  |X_{i}| > n^{1 + \varepsilon} \}$ where $ \varepsilon>0$ will be chosen small enough later on. We can write 
  \begin{eqnarray*} \mathbb{P}( S_{n} = 0) &\leq& \mathbb{P}( \tau_{n} > n) + \mathbb{P}(S_{n}=0  \mbox{ and } \tau_{n} \leq n)\\
   & \leq & \mathbb{P}( \tau_{n} > n) + \mathbb{P}(S_{n}=0  \mid \tau_{n} \leq n)   \end{eqnarray*}
The first term of the right-hand side is easy to evaluate:
$$ \mathbb{P}(\tau_{n} > n) = \big (1 - \mathbb{P}(|X|\geq n^{1 + \varepsilon})\big)^{n} \underset{ \eqref{eq:tailsymmetric}}{=} \exp\left( -  \frac{\mathrm{c}}{\alpha-1} \cdot n \cdot n^{(1+ \varepsilon)(1-\alpha)} (1+ o(1))\right)  \leq \exp(-n^{\delta}),$$
for some $\delta >0$ provided that $   \varepsilon < \frac{2- \alpha}{\alpha-1}$ is small enough. On the other hand, conditionally on $\{ \tau_{n} \leq n\}$, the increment $X_{\tau_{n}}$ is independent of $\tau_{n}$ and of the increments $\{X_1, \dots , \widehat{X_{\tau_{n}}}, \dots , X_n\}$ (beware, those are not i.i.d.~anymore) and its law $\nu_n$ is the law  of  $X$ conditioned on being of absolute value larger than $ n^{1+ \varepsilon}$;   in particular 
$$ \forall k \in \mathbb{Z},\qquad  \mathbb{P}(X_{\tau_{n}} = k \mid \tau_{n} \leq n) =  \nu_{n}(k) = \mathbf{1}_{|k| > n^{1+ \varepsilon}} \frac{\mathbb{P}(X =k)}{ \mathbb{P}(|X| > n^{1 + \varepsilon} )} $$ so that by \eqref{eq:tailsymmetric} we have  \begin{eqnarray} \label{eq:suploi} \sup_{k \in \mathbb{Z}} \nu_{n}(k) \leq  \mathrm{C}\, n^{-1- \varepsilon},  \end{eqnarray} for some constant $ \mathrm{C}>0$ for all $ n \geq 1$. Hence we can write 
  \begin{eqnarray*}\mathbb{P}(S_{n}=0 \mid \tau_{n} \leq n) &=& \mathbb{P}(X_{\tau_{n}} = -(X_{1}+ \dots + \widehat{X_{\tau_{n}} }+ \dots + X_{n}) \mid \tau_{n} \leq n )\\ 
  & \underset{ \mathrm{indept}}{=}& \mathbb{E}\left[\nu_{n}(-(X_{1}+ \dots + \widehat{X_{\tau_{n}} }+ \dots + X_{n})) \mid \tau_{n} \leq n \right]\\
  & \leq& \sup_{k \in \mathbb{Z}} \nu_{n}(k) \underset{ \eqref{eq:suploi}}{\leq} C n^{-1 - \varepsilon}.  \end{eqnarray*}
Gathering-up the pieces, we deduced that $ \mathbb{P}(S_{n}=0) \leq \exp(-n^{\delta}) +  \mathrm{C} n^{-1- \varepsilon}$ for $\delta >0$ provided that $ \varepsilon>0$ is small enough. The implies summability of the series \eqref{eq:goalsumtrans} and ensures transience of the walk. \\
We now use the same idea to treat the borderline line \textsc{case $\alpha =2$} and show that $(S)$ is recurrent by providing the lower bound
\begin{eqnarray} \label{eq:goalsumrec}\mathbb{P}( S_{n} = 0) \geq \frac{c}{n},  \end{eqnarray} for some $c>0$ thus ensuring the divergence of the expected number of visits to the origin \eqref{eq:cnsrecurrence}. We use the same idea as above but with  $ \varepsilon=0$. Let us consider the good event 
$$ \mathcal{G}_{n} = \{ |X_{i}| \leq n  \mbox{ for all } 1 \leq i \leq n \mbox{ except for two values } \tau_{n}^{1} \mbox{ and } \tau_{n}^{2} \}.$$
The probability of $ \mathcal{G}_{n}$ is easily computed and we have 
  \begin{eqnarray} \label{eq:asymGn} \mathbb{P}( \mathcal{G}_{n}) = {n \choose 2}  \mathbb{P}(|X|>n) \cdot \left(1- \mathbb{P}(|X|>n)\right) ^{n-2} \xrightarrow[n\to\infty]{ \eqref{eq:tailsymmetric}}  \frac{ \mathrm{c}^{2}}{2} \mathrm{e}^{- \mathrm{c}}>0,   \end{eqnarray} where $ \mathrm{c}>0$ appears in \eqref{eq:tailsymmetric}. In particular, this event is of asymptotically positive probability. Conditionally on $ \mathcal{G}_{n}$, the two values $ X_{\tau_{n}^{1}}, X_{\tau_{n}^{2}}$ are independent of  $\{X_1, \dots , \widehat{X_{\tau^{1}_{n}}}, \dots ,\widehat{X_{\tau^{2}_{n}}}, \dots, X_n\}$ and their common law is $\nu_{n}$, the law of $X$ conditioned on $\{|X|>n\}$.  In particular, the variable $  \mathrm{J}_n := X_{\tau_{n}^{1}} +  X_{\tau_{n}^{2}}$ is independent of $ \tilde{S}_n := X_1+ \dots + \widehat{X_{\tau^{1}_{n}}} +  \dots  + \widehat{X_{\tau^{2}_{n}}} +  \dots + X_n$. If we denote by $\nu^{\otimes 2}_n$ the law $\mathrm{J}_n$, then we clearly have $ \inf_{k \in \mathbb{Z}} \nu^{ \otimes 2}_n (k) =0$, but a moment's though shows that for any $A >0$, there exists $c_A>0$ so that  
 \begin{eqnarray} \label{eq:infavecA} \inf_{\begin{subarray}{c} k \in \mathbb{Z} \\ |k| \leq A n \end{subarray}} \nu^{ \otimes 2}_n (k) \geq \frac{c_A}{n}.  \end{eqnarray}
 On the other hand, conditionally on $ \mathcal{G}_n$, the variables $X_i$ for $i \notin \{\tau_n^1, \tau_n^2\}$ are independent and have the law of $X$ conditioned on $\{|X| \leq n\}$. In particular, there are centered, and their variance is equal to 
 $$  \frac{1}{\mu([-n,n])} \sum_{k=-n}^n  k^2 \mu_k  \underset{\eqref{eq:tailsymmetric}}{\sim} n.$$
 We deduce that the variance of $ \tilde{S}_n$ is asymptotic to $n^2$ and by Markov's inequality that 
  \begin{eqnarray} \label{eq:markovA} \mathbb{P}( |\tilde{S}_n| \geq An \mid \mathcal{G}_n ) \leq \frac{ \mathrm{Var}( \tilde{S}_n \mid \mathcal{G}_n) }{(An)^2} \leq \frac{2}{A^2},  \end{eqnarray} eventually as $n \to \infty$. Taking $A>2$ we can thus write 
   \begin{eqnarray*} \mathbb{P}(S_n = 0) &\geq& \mathbb{P}(S_n =0 \ \&\  \mathcal{G}_{n} \ \&\  |\tilde{S}_n| \leq An) \\
   &=& \mathbb{P}( \mathcal{G}_n) \cdot \mathbb{E}\left[  \mathbb{P}( J_n = - \tilde{S}_n  \ \& \ |\tilde{S}_n| \leq An \mid \mathcal{G}_n) \right]\\
   & \underset{ \mathrm{indept}}{=}& \mathbb{P}( \mathcal{G}_n) \cdot \mathbb{E}\left[  \mathbb{E}\left[ \nu^{\otimes 2}(- \tilde{S}_n)  \mathbf{1}_{|\tilde{S}_n| \leq An} \mid \mathcal{G}_n\right] \right]\\
      &  \geq & \mathbb{P}( \mathcal{G}_n) \cdot  \big( \inf_{\begin{subarray}{c} k \in \mathbb{Z} \\ |k| \leq A n \end{subarray}} \nu^{ \otimes 2}_n (k) \big) \cdot   \mathbb{P}\left( |\tilde{S}_n| \leq An \mid \mathcal{G}_n\right) \\
   & \underset{ \eqref{eq:asymGn},\eqref{eq:infavecA}, \eqref{eq:markovA} }{\geq}& \frac{ \mathrm{c}^2 \mathrm{e}^{- \mathrm{c}}}{4}  \cdot \frac{c_A}{n} \cdot (1- \frac{2}{A^2}), 
       \end{eqnarray*}
       for $n$ large enough. This shows \eqref{eq:goalsumrec} and completes the proof.\qed 
      
   \medskip

 The above result is initially due to  Shepp\footnote{\raisebox{-5mm}{\includegraphics[width=1cm]{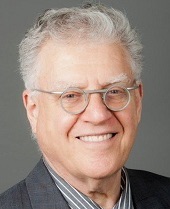}} Lawrence Alan Shepp (1936--2013), American} with a proof based on Theorem \ref{thm:recfourier1} below. He also showed the disturbing fact that there exist recurrent one-dimensional random walks with arbitrary fat tails (but not regularly varying):
\begin{theorem}[Shepp (1964)]For any positive function $ \epsilon(x) \in (0,1)$ tending to $0$ as $x \to \infty$, there exists a symmetric step distribution $\mu$ such that $\mu( \mathbb{R} \backslash [-x,x]) \geq \epsilon(x)$ for any $x \geq 0$ and such that the associated random walk $(S)$ is recurrent.
\end{theorem}



\section{Fourier transform} \label{sec:fourier}
In this section, we use the Fourier transform to give a recurrence criterion as well as a local version of the central limit theorem. Recall that if $\mu$ is the step distribution of a random walk $(S)$ on $ \mathbb{Z}$, then the Fourier\footnote{\raisebox{-3mm}{\includegraphics[width=1cm]{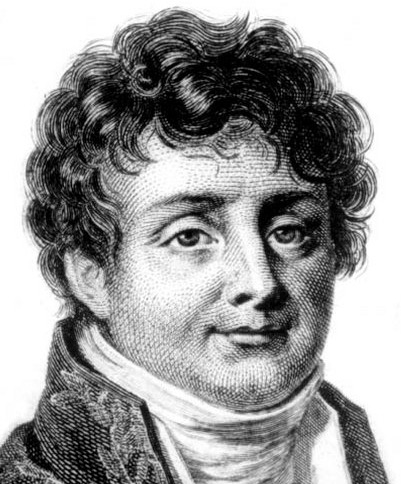}} Jean Baptiste Joseph Fourier (1768--1830), French} transform of the measure $\mu$ is defined by $$ \hat{\mu}(\xi) = \mathbb{E}[\mathrm{e}^{{ \mathrm{i}\xi  X_{1}}}]= \sum_{k \in \mathbb{Z}} \mathrm{e}^{ \mathrm{i} \xi k} \mu_k, \quad \mbox{ for } \xi \in \mathbb{R}.$$ To get information on the walk $(S)$ from $\hat{\mu}$, the main idea is of course to use Cauchy's formula to relate probabilities to integrals of powers of the Fourier transform:
 \begin{eqnarray} \label{eq:cauchy} \forall x \in  \mathbb{Z}, \quad  \mathbb{P}( S_{n}=x) &=& \sum_{k \in \mathbb{Z}} \frac{1}{2\pi} \int_{-\pi}^{\pi} \mathrm{d}\xi \,  \mathrm{e}^{- \mathrm{i} \xi  x}  \mathrm{e}^{{\mathrm{i} \xi  k}} \mathbb{P}(S_{n}=k)\\	
  &=&  \frac{1}{2\pi} \int_{-\pi}^{\pi} \mathrm{d}\xi \,  \mathrm{e}^{- \mathrm{i} \xi  x} \mathbb{E}[\mathrm{e}^{{\mathrm{i} \xi  S_{n}}}] = \frac{1}{2\pi} \int_{-\pi}^{\pi} \mathrm{d}\xi \,  \mathrm{e}^{- \mathrm{i}  \xi  x} \left( \hat{\mu}(\xi) \right)^{n}, \nonumber \end{eqnarray} where we used the fact that $ \mathbb{E}[\mathrm{e}^{\mathrm{i}\xi S_{n}}] = (\hat{\mu}(\xi))^{n}$ by independence of the increments and where the interchange of series and integral is easily justified by dominated convergence. 
\subsection{Chung-Fuchs}
  In this section, we give a criterion for recurrence of a one-dimensional random walk based on its Fourier transform. The criterion is valid mutatis mutandis for more general  random walks with values in $ \mathbb{R}^d$.
 \label{sec:fourierrec1}

\begin{theorem}[Easy version of Chung--Fuchs]\label{thm:recfourier1}\noindent The one-dimensional walk $(S)$ is recurrent if and only if we have
$$ \lim_{r \uparrow 1} \int_{-\pi}^{\pi}  \mathrm{d}\xi \  \mathfrak{Re} \left( \frac{1}{1- r  \hat{\mu}(\xi)} \right) = \infty. $$
\end{theorem}

\noindent \textbf{Proof.} By \eqref{eq:cnsrecurrence}, the walk $(S)$ is recurrent if and only if the series $\sum_{n \geq 0} \mathbb{P}(S_{n}=0)$ diverges. Recall from \eqref{eq:cauchy} that 
$$ \mathbb{P}( S_{n}=0) =  \frac{1}{2\pi} \int_{-\pi}^{\pi} \mathrm{d}t \, \mathbb{E}[\mathrm{e}^{{ \mathrm{i}tS_{n}}}] = \frac{1}{2\pi} \int_{-\pi}^{\pi} \mathrm{d}t  \left( \hat{\mu}(t) \right)^{n}.$$ We are lead to sum the last equality for $n\geq 0$, but before that we first multiply by $r^{n}$ for some $r \in [0,1)$ in order to be sure that we can exchange series, expectation and integral. One gets
$$ \sum_{n \geq 0} r^{n} \mathbb{P}( S_{n}=0) =  \frac{1}{2\pi} \int_{-\pi}^{\pi} \mathrm{d}t \, \sum_{n \geq 0} r^{n}  \left( \hat{\mu}(t) \right)^{n} = \frac{1}{2\pi}  \int_{-\pi}^{\pi}  \frac{\mathrm{d}t}{1-r \hat{\mu}(t)}.$$
Since the left-hand side is real, one can take the real part in the integral. Letting $r \uparrow 1$, the first series diverges if and only if $\sum_{n \geq 0} \mathbb{P}(S_{n}=0)= \infty$. This completes the proof of the theorem. \qed \medskip 

 In fact, there is a stronger version of Theorem \ref{thm:recfourier1} which is obtained by formally interchanging the limit and the integral in the last theorem: the random walk $(S)$ is transient or recurrent according as to whether the real part of $(1- \hat{\mu}(t))^{{-1}}$ is integrable or not near $0$ (we do not give the proof).  This can easily be proved when the law $\mu$ is \textbf{symmetric} (i.e.\,$\mu_k = \mu_{-k}$): In this case, $\hat{\mu}$ is real valued and notice that when $\hat{\mu}(t) \geq 0$ the function $r \mapsto (1- r \hat{\mu}(t))^{-1}$ is increasing, whereas if $\hat{\mu}(t) \leq 0$ we have $(1- r \hat{\mu}(t))^{-1} \leq 1$. Splitting the integral according to the sign of $\hat{\mu}(t)$ and using monotone and dominated convergence theorems on the respective parts shows that 
 $$ \lim_{r \uparrow 1} \int_{-\pi}^{\pi}  \mathrm{d}\xi \  \frac{1}{1- r  \hat{\mu}(\xi)} = \int_{-\pi}^{\pi}  \mathrm{d}\xi \  \frac{1}{1- \hat{\mu}(\xi)}.$$
 
%
%
%
%

This criterion can be used to give Fourier-proofs of some of the preceding results:

\begin{exo} \label{exo:25} Suppose $\mu$ is centered. 
\begin{enumerate}[(i)]
\item Show that $\hat{\mu}(\xi) = 1 + o(\xi)$ as $\xi \to 0$.
\item Give another proof of Theorem \ref{thm:rec0} (ii) using Theorem \ref{thm:recfourier1}.
\item Give another proof to Theorem \ref{prop:discretetransient}.
\end{enumerate}
\end{exo}

\begin{exo}[Sums of random walks] \label{exo:edouard}Let $(S_{n})_{n \geq 0}$ and $(S'_{n})_{ n\geq 0}$ be two independent one-dimensional random walks with independent increments of law $\mu$ and $\mu'$ on $ \mathbb{R}$.
\begin{enumerate}[(i)]
\item Give an example where $(S)$ and $(S')$ are transient and yet $(S+S')$ is recurrent.
\item We suppose that $\mu$ and $\mu'$ are both symmetric. Show that as soon as $(S)$ or $(S')$ is transient then so is $(S+S')$.
\item Give an example where $(S)$ is recurrent,  $(S')$ transient and yet $(S+S')$ is recurrent. 
\item (*) Can we have both $(S)$ and  $(S')$ recurrent and $(S+S')$ transient? \end{enumerate}
\end{exo}

\subsection{Local central limit theorem}
The central limit theorem is one of the most important theorems in probability theory and says in our context that the rescaled random walk $S_{n}/ \sqrt{n}$ converges in distribution towards a normal law provided that $\mu$ is centered and has finite variance. There are many proofs of this result, the most standard being through the use of Fourier transform and Lévy's criterion for convergence in law\footnote{ Here are a couple of other proofs: Lindeberg swapping trick, method of moments, Stein method, Skorokhod embedding theorem, approximation by discrete variables and de Moivre-Laplace, contraction method and Zolotarev metric... See the beautiful page by Terence Tao on this subject: https://terrytao.wordpress.com/2010/01/05/254a-notes-2-the-central-limit-theorem/ or the recent note \cite{CCW21CLT}}. We will see below that the central limit theorem can be ``disintegrated'' to get a more powerful \textbf{local} version of it.  The proof is again based on \eqref{eq:cauchy}.\medskip

When a one-dimensional random walk with mean $m$ and variance $\sigma^{2}$ satisfies a central limit theorem we mean that for any $a<b$ we have 
$$ \mathbb{P}\left(\frac{S_{n}- n m}{ \sqrt{n  }} \in  [a,b]\right) \xrightarrow[n\to\infty]{} \int_{a}^{b} \frac{ \mathrm{d}x}{ \sqrt{2 \pi \sigma^{2}}} \mathrm{e}^{{- x^{2}/(2 \sigma^{2})}}.$$
We say that we have a local central limit theorem if we can reduce the interval $[a,b]$ as a function of $n$ until it contains just \textbf{one point} of the lattice, that is if for $x \in \mathbb{Z}$ we have 
$$  \mathbb{P}( S_{n} = x) =\mathbb{P}\left(\frac{S_{n}- n m}{ \sqrt{n  }} \in  \left[ \frac{x-nm}{ \sqrt{n}},\frac{x-nm +1}{ \sqrt{n}}\right)\right) \approx \frac{1}{\sqrt{2 \pi \sigma^{2}}} \mathrm{e}^{{-  \frac{(x-nm)^{2}}{2 n \sigma^{2}}}} \frac{1}{ \sqrt{n}}.$$
It turns out that aperiodicity and finite variance are already sufficient to get the local central limit theorem (the result extends to higher dimensions and to the case of random walks converging towards stable Lévy process with {mutatis mutandis} the same proof):

\begin{theorem}[Local central limit theorem, Gnedenko]\label{thm:local}Let $\mu$ be a distribution supported on $ \mathbb{Z}$, aperiodic, with mean $m \in \mathbb{R}$ and with a finite variance $\sigma^{2}>0$.  If we denote by $\gamma_{\sigma}(x) = \frac{1}{ \sqrt{2\pi \sigma^2}} \mathrm{e}^{{-x^{2}}/(2 \sigma^{2})}$  the density  of the centered normal law of variance $\sigma^{2}$ then we have 
$$ \lim_{n \to \infty} \sup_{ x \in \mathbb{Z}}  n^{{1/2}} \left|  \mathbb{P}(S_{n}=x) - n^{-{1/2}} \gamma_{\sigma}\left(\frac{x-nm}{\sqrt{n}}\right)\right| = 0.$$
\end{theorem}

The usual central limit theorem follows from its local version. Indeed, if we consider the random variable $ \tilde{S}_{n}=S_{n} + U_{n} $ where $U_{n}$ is uniform over $[0,1]$ and independent of $S_{n}$, then the local central limit theorem shows that the law of $(\tilde{S}_{n} - nm)/ \sqrt{n}$ is absolutely continuous with respect to the Lebesgue measure on $ \mathbb{R}$ whose density $f_{n}$ converges pointwise towards the density of $\gamma_{\sigma}$. Scheffé's\footnote{\raisebox{-5mm}{\includegraphics[width=1cm]{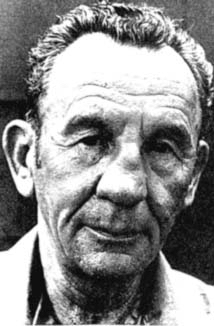}} Henry Scheff\'e (1907--1977), American} lemma (see Exercise \ref{exo:scheffe} below) then implies that $(\tilde{S}_{n}-nm)/ \sqrt{n}$ converges in law towards $ \gamma_{ \sigma} ( \mathrm{d}x)$ and similarly after removing the tilde. 
  \begin{exo}[Scheff\'e lemma]  \label{exo:scheffe}Let $X_{n},X$ be random variables taking values in a Polish space $(E,d)$ and whose distributions have densities $f_{n},f$ with respect to a background measure $ \mathrm{m}$ on $E$. We suppose that $f_{n} \to f$ pointwise $ \mathrm{m}$-almost everywhere. Prove that 
  \begin{enumerate}[(i)]
  \item $f_{n} \to f$ in $ \mathbb{L}^{1}(  \mathrm{m})$.
  \item $ \mathrm{d_{TV}}( f_{n} \mathrm{d} \mathrm{m} , f \mathrm{d} \mathrm{m}) \to 0$ where $ \mathrm{d_{TV}}$ is the total variation distance as $n \to \infty$,
  \item deduce that $X_{n} \to X$ in distribution as $n \to \infty$.
\end{enumerate}
\end{exo}

Before moving to the proof of the local CLT, let us translate the aperiodicity condition on the Fourier transform:
\begin{lemma} \label{ref:lem<1} When $\mu$ is aperiodic, we have \label{eq:<1}
  \begin{eqnarray*}  |\hat{\mu}(\xi)| < 1, \quad \mbox{ for }\xi \in (0,2\pi).  \end{eqnarray*}
  \end{lemma}
  \noindent \textbf{Proof.}
  Indeed, if we have $ |\hat{\mu}(\xi)| = |\mathbb{E}[\mathrm{e}^{{ \mathrm{i}\xi  X_{1}}}]|=1$ we have also $|\big(\hat{\mu}(\xi)\big)^{n}| = |\mathbb{E}[\mathrm{e}^{{ \mathrm{i}\xi   S_{n}}}]|=1$. This  implies by the equality case in the triangle inequality that all $ \mathrm{e}^{{ \mathrm{i}\xi  x}}$ for $x \in \mathrm{Supp}( \mathcal{L}(S_{n}))$ are (positively) aligned. Using the aperiodicity assumption, one can choose $n$ large enough so that the support of the law of $S_{n}$ contains $0$ and $1$.  This shows that $\xi \equiv 0$ modulo $[2 \pi]$. \qed 

\medskip

\noindent \textbf{Proof of the local central limit theorem.} The starting point is again Cauchy formula's relating probabilities to Fourier transform:
$$ \mathbb{P}(S_{n}=x) = \frac{1}{2\pi} \int_{-\pi}^{\pi} \mathrm{d}t \, \mathrm{e}^{-ixt} \mathbb{E}[\mathrm{e}^{ \mathrm{i} S_{n}t}] = \frac{1}{2\pi} \int_{-\pi}^{\pi} \mathrm{d}t \, \mathrm{e}^{-\mathrm{i}xt}( \hat{\mu}(t))^{n}.$$
Since $|\hat \mu (t)| <1$ when $t \ne 0$ the main contribution of the integral comes from the integration near $0$: we shall apply Laplace's method. Since we want to use the series expansion of the Fourier transform near $0$ it is natural to introduce the image measure $\nu$  of $\mu$ after translation by $-m$ so that $\nu$ is centered and has finite variance: we can write $\hat{\nu}(t) = 1 - \frac{\sigma^{2}}{2} t^{2} + o(t^{2})$ for $t$ small. The last display then becomes

 \begin{eqnarray*} \mathbb{P}(S_{n}=x) &=& \frac{1}{2\pi} \int_{-\pi}^{\pi} \mathrm{d}t \, \mathrm{e}^{- \mathrm{i}xt} \mathrm{e}^{{ \mathrm{i} nm t}}( \hat{\nu}(t))^{n} \\
 &=& \frac{1}{2\pi} \int_{-\pi}^{\pi} \mathrm{d}t \, \mathrm{e}^{-\mathrm{i}xt} \mathrm{e}^{{ \mathrm{i} nmt}}\left( 1 - \frac{\sigma^{2}}{2}t^{2} + o(t^{2})\right)^{n}\\
  &\underset{u =t \sqrt{n}}{=}&   \frac{1}{ \sqrt{n}}\frac{1}{2\pi} \int_{-\pi \sqrt{n}}^{\pi \sqrt{n}} \mathrm{d}u \, \mathrm{e}^{-\mathrm{i} u x/ \sqrt{n}} \mathrm{e}^{{ \mathrm{i}  \sqrt{n} m u}}\underbrace{\left( 1 - \frac{\sigma^{2}}{2n}u^{2} + o(u^{2}/n)\right)^{n}}_{\approx  \frac{ \sqrt{2\pi}}{\sigma}\gamma_{1/\sigma}(u)}.  \end{eqnarray*}
  We can then approximate the last integral by
  $$ \frac{ \sqrt{2\pi}}{\sigma}\frac{1}{ \sqrt{n}}\frac{1}{2\pi} \int_{-\infty}^{\infty} \mathrm{d}u \, \mathrm{e}^{-\mathrm{i} u x/ \sqrt{n}} \mathrm{e}^{{ \mathrm{i}  \sqrt{n} m u}} \gamma_{1/\sigma}(u) =  \frac{1}{ \sigma \sqrt{2 \pi n}}\mathbb{E}\left[\exp\left( \mathrm{i} \left(\sqrt{n}m - \frac{x}{ \sqrt{n}}\right)  \frac{ \mathcal{N}}{ \sigma}\right)\right],$$ where $ \mathcal{N}$ denotes a standard normal variable. Using the identity $ \mathbb{E}[\mathrm{e}^{{ \mathrm{i}t \mathcal{N}}}] = \mathrm{e}^{{-t^{2}/2}}$ the last display is indeed equal to $ \gamma_{\sigma} \left( \frac{x- nm}{ \sqrt{n}}\right) / \sqrt{n}$ as desired. It remains to quantify the last approximation. The error made in the approximation is clearly bounded above by the sum of the two terms:
  $$ A= \frac{1}{ \sqrt{n}}\frac{1}{2\pi} \int_{|u| > \pi \sqrt{n}} \mathrm{d}u \, \gamma_{1/\sigma}(u),$$
    $$ B= \frac{1}{ \sqrt{n}}\frac{1}{2\pi} \int_{|u| < \pi \sqrt{n} } \mathrm{d}u \,  \left| 	 \gamma_{1/\sigma}(u) - \left(\hat{\nu}\left( \frac{u}{ \sqrt{n}}\right)\right)^{n}\right|.$$
    The first term $A$ causes no problem since it is exponentially small (of the order of $\mathrm{e}^{-n}$) hence negligible in front of  $ 1/ \sqrt{n}$. The second term may be further bounded above by the sum of three terms
 \begin{eqnarray*} B &\leq& \frac{1}{ \sqrt{n}} \int_{|u| < n^{{1/4}}} \mathrm{d}u \,  \left| 	 \gamma_{1/\sigma}(u) - \left(\hat{\nu}\left( \frac{u}{ \sqrt{n}}\right)\right)^{n}\right| \\ & +&  \int_{n^{{1/4}}<|u| <  \pi \sqrt{n}} \mathrm{d}u \, 	 \gamma_{1/\sigma}(u) \\ &+& \int_{n^{{1/4}}<|u| <  \pi \sqrt{n}} \mathrm{d}u \, \left|\hat{\nu}\left( \frac{u}{ \sqrt{n}}\right)\right|^{n}.  \end{eqnarray*}
    The first of these terms is shown to be $o( n^{{-1/2}})$ using dominated convergence: in the region considered for $u$, the integrand converges pointwise to $0$; for the domination we may use the fact  for $|u| < \varepsilon \sqrt{n}$ we have by the expansion of $\hat{\nu}$ that  $\left|\hat{\nu}\left( \frac{u}{ \sqrt{n}}\right)\right|^{n} \leq ( 1 - \frac{ \sigma^{2} u^{2}}{4n})^{n} \leq  \mathrm{e}^{- \sigma^{2}u^{2}/4}$. The second term of the sum is handled as above and seen to be of order $\mathrm{e}^{- \sqrt{n}}$. For the third term, we bound the integrand by $\left|\hat{\nu}\left( \frac{u}{ \sqrt{n}} \right)\right|^{n} \leq  \mathrm{e}^{- \sigma^{2}u^{2}/4}$ for $|u| < \varepsilon \sqrt{n}$, as for $  \varepsilon \sqrt{n} < |u|< \pi \sqrt{n}$ we use the fact that $|\hat{\mu}(x)|<c<1$ for all $x \in [\varepsilon, \pi]$ by aperiodicity. The sum of the three terms is then of negligible order compared to $n^{{-1/2}}$ as desired.  \qed  \medskip

%
%

%
%
%
%

\noindent \textbf{Bibliographical notes.} The material in this chapter is standard and can be found in many textbooks, see e.g.~ \cite[Chapters I,II]{Spi76}, \cite[Chapter 8]{Chung74} or \cite{Kal07,LalCours}. The proof of Theorem \ref{prop:discretetransient} is due to Yuval Peres (personal communication) while Shepp's original proof  \cite{shepp1962symmetric} is based on the Fourier transform. Theorem \ref{thm:shepp} can be found in \cite{Shepp64}. The Fourier transform is a remarkable tool (whose efficiency is sometimes a bit mysterious) to study random walks  with independent increments. The local central limit theorem is valid in the much broader context of random walks converging towards stable Lévy processes, see Gnedenko's local limit theorem  in \cite[Theorem 4.2.1]{IL71}, and can be sharpened when we have further moment assumptions, see \cite{LL10}. It also applies when the variables are not exactly i.i.d.~\cite{davis1995elementary}. 

\medskip 

\noindent \textbf{Hints for Exercises.}\ \\ 
Exercise \ref{exo:nolfgn}: Use the truncated increments $X_{i} \mathbf{1}_{|X_{i}| < A n/\log n}$ for $1 \leq i \leq n$ for some large $A>0$.\\
Exercise \ref{exo:22}: Couple the increments $X_{i}$ so that they are the same for $2^{k} \leq i < 2^{k+1}$.\\
Exercise \ref{exo:edouard}: (i) is easy. (ii) Use the Fourier criterion. (iii) Use the $1$-stable Cauchy distribution (or a discrete version thereof). (iv) Edouard Maurel-Segala solved it here: \\
\textit{https://mathoverflow.net/questions/314312/sum-of-independent-random-walks}

\chapter{Skip-free random walks} \label{chap:WH}

\hfill Two simple but powerful observations.

\bigskip
 
In this chapter we still consider a one-dimensional random walk $(S)$  based on i.i.d.\;increments of law $\mu$ (whose support is not contained in $ \mathbb{Z}_{\geq0}$ nor in $ \mathbb{Z}_{\leq 0}$). But compared to the previous chapter, we  furthermore suppose that the walk is \textbf{skip-free} which means that 
$$ \mathrm{Supp}( \mu) \subset \{-1, 0, 1, 2, 3, \dots \}.$$ In other words, the only negative steps of $(S)$ are steps of size $-1$. We shall see that some combinatorial magic happens for such walks. Let us start by drawing a consequence of the last chapter: the expectation $$m = \sum_{k \geq -1} k \mu_k$$ is always well-defined and belongs to $ [-1, \infty]$ and so by Theorem \ref{thm:rec0} the walk is recurrent if $m=0$ and drifts otherwise. We will now perform two simple combinatorial operations on paths  (reversal and cycle shift) and explore their distributional consequences.

\section{Duality lemma}
We begin with a simple but surprisingly important observation called \textbf{duality}. This is valid for any random walk, not necessarily skip-free and not necessarily integer-valued.
\subsection{Duality}
\begin{proposition}[Duality]  \label{lem:duality}For each fixed $n \geq 0$, we have the following equality in distribution 
$$ (0=S_{0},S_{1}, \dots , S_{n})  \overset{(d)}{=} (S_{n}-S_{n},S_{n}-S_{n-1}, S_{n}-S_{n-2}, \dots , S_{n}-S_{1}, S_{n}-S_0).$$
 \end{proposition}
 \begin{figure}[!h]
 \begin{center}
 \includegraphics[width=13cm]{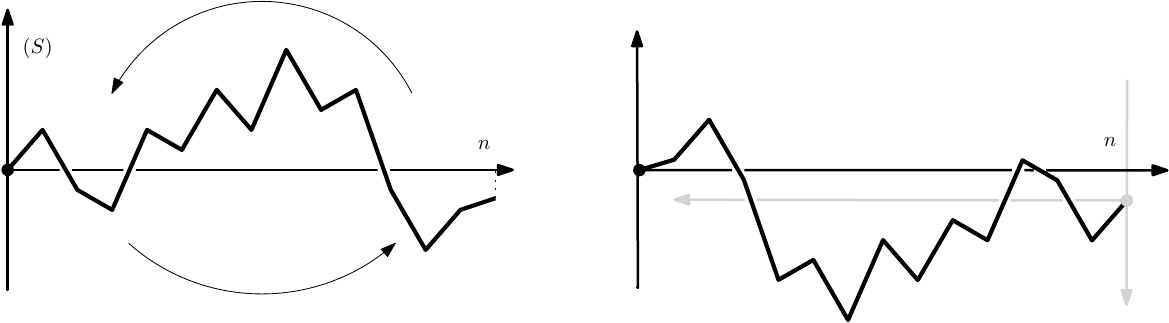}
 \caption{Geometric interpretation of the duality: the rotation by an angle $\pi$ of the first $n$ steps of the walk $(S)$ leaves its distribution invariant.}
 \end{center}
 \end{figure}
 \noindent \textbf{Proof.} It suffices to notice that the increments of the walk $(S_{n}-S_{n-1}, S_{n}-S_{n-2}, \dots , S_{n}-S_{1}, S_n-S_0)$ are just given by $(X_{n}, X_{n-1}, \dots , X_{1})$ which obviously has the same law as $(X_{1}, \dots ,X_{n})$ since the $(X_{i})_{i \geq 1}$ are i.i.d.\ hence exchangeable (i.e.~whose law is invariant under any fixed permutation). \qed \medskip
 
 Beware, the duality lemma can only be applied for $n$ fixed and not for all $n$ simultaneously, yet it can be useful to deduce asymptotic properties of the walk:
 
  \begin{corollary} \label{cor:duality} Suppose that $(S)$ is a one-dimensional (non necessarily skip-free nor integer valued). We denote by $ \overline{S}_{n} = \sup \{ 0 \leq k \leq n : S_{k}\}$ and $ \underline{S}_{n} = \inf \{ 0 \leq k \leq n : S_{k}\}$ the running supremum and infimum processes.  We suppose that  $(S)$ drifts towards $+\infty$ so that $ \min S = \underline{S}_{\infty} > -\infty$ a.s.. Then we have
  $$  \overline{S}_{n} - S_{n} \xrightarrow[n\to\infty]{(d)} - \underline{S}_{\infty} < \infty.$$
\end{corollary}
\noindent \textbf{Proof.} By duality we have for each $n \geq 0$
$$ \overline{S}_{n} - S_{n}  \quad \overset{(d)}{ =} \quad - \underline{S}_{n}$$
whereas since $(S)$ drifts towards $-\infty$ we have  $$- \underline{S}_{n} \quad  \xrightarrow[n\to\infty]{a.s.} \quad  - \underline{S}_{\infty} < \infty.$$\qed  \bigskip

\begin{figure}[!h]
 \begin{center}
 \includegraphics[width=14cm]{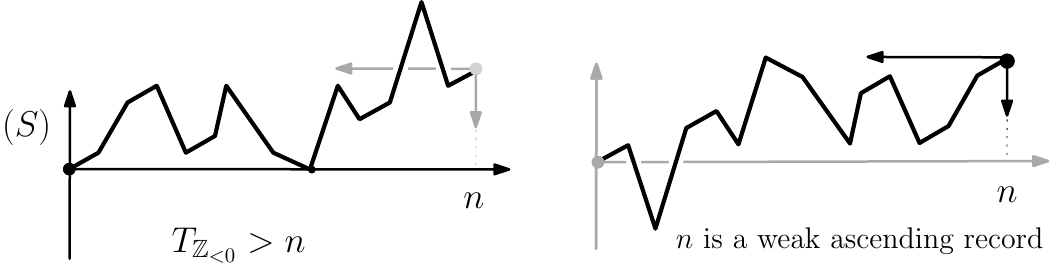}
 \caption{ \label{fig:dualitylln} Duality shows that $\mathbb{P}(T_{ \mathbb{Z}_{<0}}>n) = \mathbb{P}(n \mbox{ is a weak ascending record time})$.}
 \end{center}
 \end{figure}

One of the main application of duality is the following interpretation of hitting times of half-spaces. For $A \subset \mathbb{Z}$, we denote by $  {T}_{A} = \inf \{i \geq 0 : S_i \in A\}$ the hitting time of $A$ by the walk $(S)$. Then the previous proposition shows (see Figure \ref{fig:dualitylln}) that  for $n \geq 0$
  \begin{eqnarray*}
   \mathbb{P}( T_{ \mathbb{Z}_{<0}} > n) &=& \mathbb{P}( S_{0}=0, S_{1} \geq 0, \dots , S_{n}\geq 0)\\
   & \underset{ \mathrm{duality}}{=} &  \mathbb{P}( S_{n}-S_{n}= 0, S_{n}-S_{n-1} \geq 0, S_n-S_{n-2} \geq 0 \dots , S_{n} \geq 0)\\
   & =&  \mathbb{P}(S_{n} \geq S_{n-1}, S_{n}\geq S_{n-2}, \dots , S_{n} \geq S_{0})\\
    &=& \mathbb{P}(n \mbox{ is a new (weak) ascending  record time for the walk}),  \end{eqnarray*} where an \textbf{ ascending/descending (resp. weak) record time} is a time where the walk attains (or equals) a  new  maximum/minimum value so far i.e.~such that $S_{i} > \max\{S_{j} : 0 \leq j < i\}$ for strict ascending, $\geq$ for weak ascending, $\leq \min$ for weak descending and $< \min$ for strict descending. Summing over $n \geq 0$ we deduce that 
 \begin{eqnarray}  \sum_{n \geq 0}\mathbb{P}(T_{ \mathbb{Z}_{<0}}>n) = \mathbb{E}[T_{ \mathbb{Z}_{<0}}]= \mathbb{E}[ \#\mbox{ weak ascending record times}]. \label{eq:duality*}  \end{eqnarray}
 However, it is easy to see using the Markov property that the number of (weak) ascending record times is a geometric random variable which is finite almost surely if and only if the walk is bounded from above. In particular, we deduce that $T_{ \mathbb{Z}_{<0}}$ has finite expectation iff $m <0$ and since for skip-free random walk we have $ S_{T_{ \mathbb{Z}_{<0}}} = -1$ we deduce from Wald's identity that
 $$ \mathbb{E}[ T_{ \mathbb{Z}_{<0}}]  \underset{ \mathrm{Thm.} \ref{thm:wald}}{=} \frac{1}{|m|}.$$
 
 \begin{exo} \label{exo:legallcornell} Let $(S)$ be a centered skip-free random walk. Show using duality that for any $k \geq 1$ we have
 $$  \mathbb{P}(S_{T_{ \mathbb{Z}_{>0}}} = k) = \frac{1}{\mu_{-1}} \sum_{i \geq k} \mu_i.$$
\end{exo}

 \subsection{A proof of the law of large numbers}
 
To illustrate the power of the duality  lemma, let us use it to give a short proof of the law of large numbers. In this section only, let $X_{1}, X_{2}, \dots $ be i.i.d.~random variables not necessarily integer-valued with finite expectation and let $S_{n} = X_{1}+ \dots + X_{n}$ for $n \geq 0$ be the corresponding random walk. Kolmogorov's strong law of large numbers says that $ n^{-1}S_{n} \to  \mathbb{E}[X]$ almost surely as $n \to \infty$. Clearly, it is a consequence of the following lemma:
             
             \begin{lemma} \label{lem:LLN} Let $X_{1}, X_{2}, \dots $ be i.i.d.~r.v.~with $ \mathbb{E}[X]<0$. Then $ \sup_{n \geq 0} (X_{1} + \cdots + X_{n})$ is finite a.s.
             \end{lemma}
\noindent \textbf{Proof.}
 \textsc{Step 1. Bounding the increments from below.} Choose $C>0$ large enough so that by dominated convergence $ \mathbb{E}[X \mathbf{1}_{X>-C}] < 0$. We will show that the random walk $\tilde{S}_{n} = X_{1} \mathbf{1}_{X_{1}>-C} + \dots + X_{n} \mathbf{1}_{ X_{n}>-C}$ is a.s.~bounded from above which is sufficient to prove the lemma.\\
\textsc{Step 2. Duality.} Consider $\tilde{T}_{ \mathbb{Z}_{<0}} = \inf\{i \geq 0 : \tilde{S}_{i}<0\}$ and recall from \eqref{eq:duality*} that $$ \mathbb{E}[\tilde{T}_{ \mathbb{Z}_{<0}}] = \mathbb{E}[\# \mbox{ weak ascending record times of } \tilde{S}]$$ and the proof is complete if we prove that  $ \mathbb{E}[\tilde{T}_{ \mathbb{Z}_{<0}}]< \infty$ since this implies that almost surely there is a finite number of weak ascending records for $ \tilde{S}$, hence the walk is bounded from above a.s. \\
\textsc{Step 3. Optional sampling theorem.} To prove  $ \mathbb{E}[\tilde{T}_{ \mathbb{Z}_{<0}}]<\infty$, consider the same martingale as in the proof of Wald's identity (Theorem \ref{thm:wald}) namely $$ M_{n}~=~\tilde{S}_{n} - \mathbb{E}[X \mathbf{1}_{X>-C}]n, \quad \mbox{ for }n \geq0$$ (for the filtration generated by the $X_{i}$'s) and apply the optional sampling theorem to the stopping time $n \wedge \tilde{T}_{ \mathbb{Z}_{<0}}$ to deduce that 
$$ 0=\mathbb{E}[M_{n\wedge \tilde{T}_{ \mathbb{Z}_{<0}}}] \quad \mbox{ or in other words } \quad \mathbb{E}[X \mathbf{1}_{X>-C}] \mathbb{E}[n \wedge \tilde{T}_{ \mathbb{Z}_{<0}}] = \mathbb{E}[\tilde{S}_{n \wedge \tilde{T}_{ \mathbb{Z}_{<0}}}].$$
Since the increments of $\tilde{S}$ are bounded below by $-C$, the right-hand side of the last display is bounded from below by $-C$ as well. Recalling that $\mathbb{E}[X \mathbf{1}_{X>-C}] $ is negative we deduce that 
$$ \mathbb{E}[n \wedge \tilde{T}_{ \mathbb{Z}_{<0}}] =  \frac{\mathbb{E}[\tilde{S}_{n \wedge \tilde{T}_{ \mathbb{Z}_{<0}}}]}{\mathbb{E}[X \mathbf{1}_{X>-C}]} \leq \frac{C}{ |\mathbb{E}[X \mathbf{1}_{X>-C}]|} < \infty.$$
Letting $n \to \infty$, by monotone convergence we deduce that the expectation of $\tilde{T}_{ \mathbb{Z}_{<0}}$ is finite. Et voilà. \qed \medskip

\section{Cycle lemma} 
\label{sec:fellerskip}
 
 The following theorem has many names, equivalent forms and ramifications in the probabilistic and combinatorial  literature (Kemperman's formula, Otter--Dwass formula, Feller\footnote{\raisebox{-3mm}{\includegraphics[width=1cm]{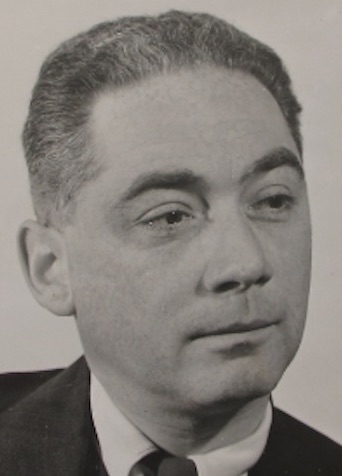}} William Feller (1906--1970), born Vilibald Sre\'cko Feller, Croatian and American} combinatorial lemma, D\'esir\'e Andr\'e cycle lemma, Lagrange inversion formula...). We shall start with the following deterministic statement:

Let $x_{1}, x_{2}, \dots , x_{n} \in \{-1, 0,1, \dots\}$ be integers which we consider as the increments of the skip-free walk $(s)$ defined by $$s_{0}=0, s_{1}=x_{1}, s_{2}=x_{1}+x_{2}, \quad  \dots, \quad s_{n}= x_{1}+ \dots + x_{n}.$$ If $ \ell \in \{ 0,1,2, \dots , n-1\}$ we consider $(s^{{(\ell)}})$ the $\ell$-th cyclic shift of the walk obtained by cyclically shifting its increments $\ell$ times, that is 
$$ s^{{(\ell)}}_{0}=0, s^{{(\ell)}}_{1}=x_{\ell+1}, \quad  \dots \quad,  s^{{(\ell)}}_{n-\ell}=x_{\ell+1}+ \dots + x_{n},\quad  \dots \quad, s_{n}^{{(\ell)}} = x_{\ell+1}+ \dots + x_{n}+ x_{1}+ \dots + x_{\ell}.$$
\begin{lemma}[Feller]\label{lem:feller}Suppose that $s_{n}=-k$ for $k \geq 1$. Then there are exactly $k$ cyclic shifts $(s^{{(\ell)}})$  with $ \ell \in \{ 0,1,2, \dots , n-1\}$ for which time $n$ is the hitting time of $-k$ by the walk $(s^{{(\ell)}})$.\end{lemma}
\noindent \textbf{Proof.} Let us first prove there is at least one such cycle shift. For this, consider the first time $\ell\in \{1,2,\dots , n\}$ such that the walk $(s)$ reaches its overall minimum $\min\{s_{i} : 0 \leq i \leq n\}$. Then clearly, after performing the cycle shift at that time, the new walk stays above $-k$ over $\{0,1, \dots , n-1\}$, see Figure \ref{fig:rerootmin} below.

\begin{figure}[!h]
 \begin{center}
 \includegraphics[width=13cm]{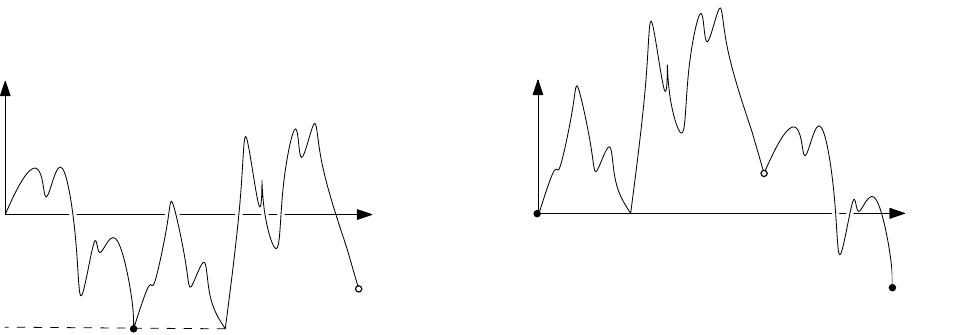}
 \caption{ \label{fig:rerootmin} Re-rooting the walk at the hitting time of the minimum gives a walk reaching its minimum at time $n$.}
 \end{center}
 \end{figure} 
 We can thus suppose without loss of generality that time $n$ is the hitting time of $-k$ by the walk.  It is now clear (again see Figure \ref{fig:reroot} below) that  the only possible cyclic shifts of the walk such that the resulting walk first hits $-k$ at time $n$ correspond to the hitting times of $ 0,-1, -2, \dots , -(k-1)$ by the walk (we use skip-free here to notice that those hitting times all appear).
 
 \begin{figure}[!h]
  \begin{center}
  \includegraphics[width=8cm]{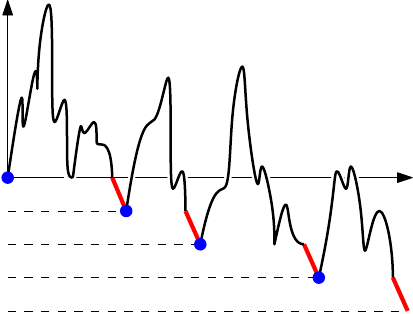}
  \caption{ \label{fig:reroot} The only possible cycle shifts (starting points in blue) for which the walk first hit $-k$ at time $n$ correspond to the hitting times of $0, -1, -2, \dots , -k+1$.}
  \end{center}
  \end{figure}  \qed \medskip

\begin{remark} \label{rem:shift}Beware, Feller's combinatorial lemma \textit{does not} say that the cyclic shifts $(s^{(\ell)})$ are distinct. Indeed, in the action of $ \mathbb{Z}/ n \mathbb{Z}$ on $\{ (s^{(\ell)}) : \ell \in \{0,1, \dots , n-1\}\}$ by cyclic shift, the size of the orbit is equal to $n / j$ where $j$ (which divides $n$) is the cardinal of the subgroup stabilizing $ (s^{(0)})$. In our case, it is easy to see that $j$ must also divide $k$ and in this case there are only $k/j$ distinct cyclic shifts having $n$ as the hitting time of $-k$. In particular, when  $k=1$ the $n$ cycle shifts are pairwise distinct.\end{remark}

\subsection{Kemperman's formula and applications}
Notice that Lemma \ref{lem:duality} does not require that the random walk has i.i.d.~increments: it holds as soon as the increments $( \mathcal{X}_k : k \geq 1)$ of a random process $ ( \mathcal{S}_k : k \geq 0)$ are invariant by time reversal i.e. $$ ( \mathcal{X}_1, \dots,  \mathcal{X}_{n_0}) \overset{(d)}{=} ( \mathcal{X}_{n_0}, \dots ,  \mathcal{X}_1),$$  for all $n_0 \geq 1$. In the application of Feller combinatorial lemma below, we shall use another invariance, by cycle shift, which amounts to ask 
$$ ( \mathcal{X}_1, \dots,  \mathcal{X}_{n_0}) \overset{(d)}{=} (\mathcal{X}_2, \mathcal{X}_3, \dots , \mathcal{X}_{n_0}, \mathcal{X}_1)  \overset{(d)}{=} (\mathcal{X}_3, \mathcal{X}_4, \dots , \mathcal{X}_{n_0}, \mathcal{X}_1,\mathcal{X}_2) = \dots.$$ for all $n_0 \geq 1$. Those properties are in particular satisfied as soon as the  increments $( \mathcal{X})$ are  exchangeable in the sense that 
$$ ( \mathcal{X}_1, \dots,  \mathcal{X}_{n_0}) \overset{(d)}{=} (\mathcal{X}_{\sigma(1)}, \mathcal{X}_{\sigma(2)}, \dots , \mathcal{X}_{\sigma(n_0)})$$ for any $n_0$ and any permutation $\sigma$ of $\{1,2, \dots , n_0\}$. For the connoisseur, De Finetti's theorem (not discussed in these pages) shows that those processes are mixture of random walks with i.i.d.~increments.

\subsection{Kemperman's formula}
As usual, for $k \in \mathbb{Z}$, we denote by $T_{k} = \inf \{ 0\leq i  \leq  n : S_i=k\}$ the hitting time of $k$ by the random walk $(S)$. An easy corollary of the cycle lemma is  Kemperman's\footnote{\raisebox{-3mm}{\includegraphics[width=1cm]{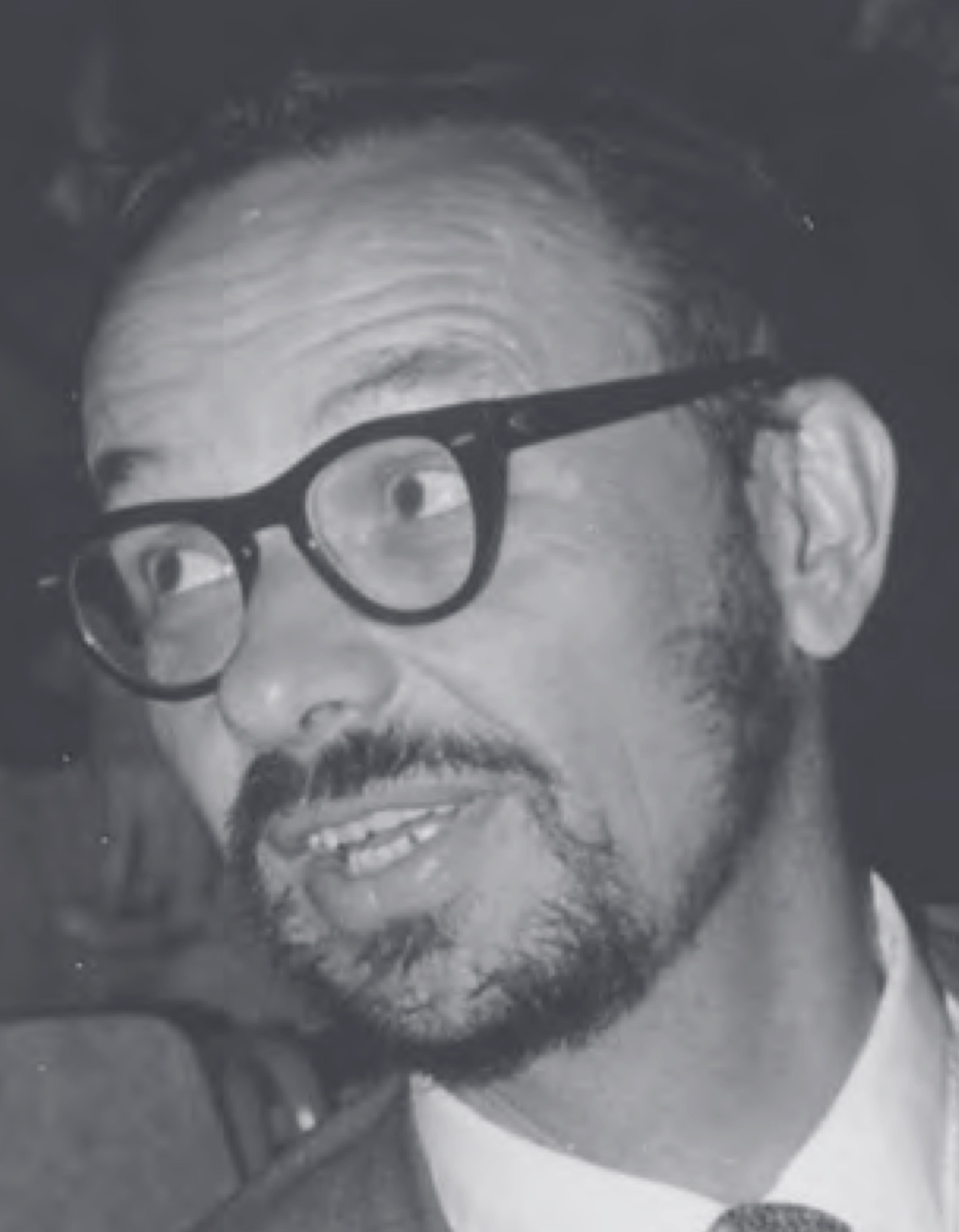}} Johannes Henricus Bernardus Kemperman (1924--2011), Dutch} formula:
\begin{proposition}[Kemperman's formula]\label{prop:kemperman} Let $(0=S_0,S_1, \dots , S_n)$ be a skip-free process with cyclically exchangeable increments. Then for every $n \geq 1$ and every $k \geq 1$ we have 
$$ \frac{1}{n}\mathbb{P}(S_{n}=-k) = \frac{1}{k} \mathbb{P}(T_{-k}=n).$$
\end{proposition}
\noindent \textbf{Proof}. Let us first re-write Lemma \ref{lem:feller} in a single equation
  \begin{eqnarray} \label{eq:equivfeller} \mathbf{1}_{s_{n}=-k} \quad =\quad  \frac{1}{k} \sum_{\ell=0}^{{n-1}} \mathbf{1}_{T_{-k}(s^{(\ell)})=n}.  \end{eqnarray}
Indeed, if the walk $(s)$ is such that $s_{n}=-k$ for $k \geq 1$, then there exists exactly $k$ shifts which do not annulate the indicator functions on the right-hand side. Since we divide by $k$ the total sum is one. We take expectation when $(s) = (S_{0}, S_{1}, \dots , S_{n})$ is the path made up of the first $n$ steps of our random walk.  Using exchangeability of the increments,  for all $0 \leq \ell \leq n-1$ we have $(S^{(\ell)}_{j})_{0 \leq j \leq n} = (S_{j})_{0 \leq j \leq n}$ in distribution.  We deduce Kemperman's formula. \qed \medskip

\begin{remark} \label{rek:kemp+local} Combining Kemperman's formula with the local central limit theorem (Theorem \ref{thm:local}), we deduce that if $(S)$ is an aperiodic skip-free random walk with centered increments having finite variance $\sigma^2$ then we have $$ \mathbb{P}(T_{-1}=n) \sim \frac{1}{\sqrt{2 \pi \sigma^{2}}} \cdot \frac{1}{n^{3/2}}, \quad \mbox{ as }n \to \infty.$$
\end{remark}

\begin{exo}\label{exo:shiftunif} Let $(S)$ be an integer-valued one-dimensional random walk, but non necessarily skip-free. For $n \geq 0$, let $K_{n} = \inf\{ 0 \leq k \leq n : S_{k}= \sup_{0 \leq i \leq n}S_{i}\}$ for the first time when the walk achieves its maximum over $\{0,1, \dots ,n\}$. Show that conditionally on $S_{n}=1$, the variable $K_{n}$ is uniformly distributed over $\{1,2, \dots , n\}$. Compare with Proposition \ref{prop:arcsine1st}.
\end{exo}

\subsection{Simple symmetric random walk} Let us give a first application of this formula in the case of the symmetric simple random walk whose step distribution is $ \mu= \frac{1}{2}( \delta_{1} + \delta_{-1})$. Due to parity reasons,  $T_{-1}$ must be odd, and by Kemperman's formula we have for $n \geq 1$
 \begin{eqnarray} \label{eq:T>ss} \mathbb{P}( T_{-1}=2n-1) &=& \frac{1}{2n-1} \mathbb{P}(S_{2n-1}=-1) = \frac{1}{2n-1} 2^{{-(2n-1)}}  { 2n-1 \choose n}  \\ &=& 2^{-2n+1} \frac{(2n-2)!}{n!(n-1)!} =  \frac{1}{2} \cdot {4^{-(n-1)}} \mathrm{Cat}(n-1),  \end{eqnarray} where for $n \geq 0$ we have put $ \mathrm{Cat}(n) = \frac{1}{n+1} {2n \choose n}$ for the $n$th Catalan\footnote{\raisebox{-3mm}{\includegraphics[width=1cm]{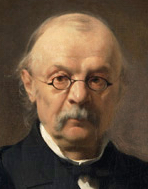}} Eugène Charles Catalan  (1814--1894), French and Belgian} number.
As an application of this formula, we can prove the famous arcsine law \footnote{there are at least three arcsine laws in the theory of random walk...}:

\begin{proposition}[1st Arcsine law] \label{prop:arcsine1st}Let $(S)$ be the simple symmetric random walk on $ \mathbb{Z}$. We put $K_{n} = \inf\{ 0 \leq k \leq n : S_{k} = \sup_{0 \leq i \leq n}S_{i}\}$ then 
 $$ \frac{K_{n}}{n} \xrightarrow[n\to\infty]{(d)} \frac{ \mathrm{d}x}{ \pi \sqrt{x (1-x)}} \mathbf{1}_{ x \in (0,1)}.$$
 \end{proposition}
  \begin{figure}[!h]
 \begin{center}
 \includegraphics[width=8cm]{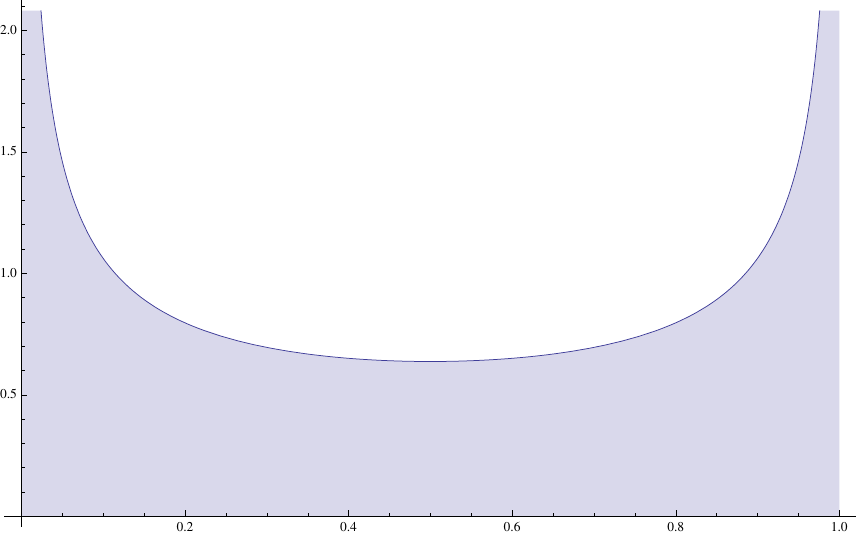}
 \caption{The arcsine distribution}
 \end{center}
 \end{figure}
 \begin{remark} The name arcsine comes from the cumulative distribution function of the right-hand side which is $\frac{2}{\pi} \mathrm{arcsin}(\sqrt{x})$. Quoting Feller “Contrary to intuition, the maximum accumulated gain is much more likely to occur towards the very beginning or the very end of a coin-tossing game than somewhere in the middle.” \end{remark}
 \noindent \textbf{Proof.} Putting $T_0^+ = \inf\{n >0 : S_n=0\}$ to be the first \textbf{return} time at $0$, using duality we can compute exactly for $k \in \{1,2, \dots,n\}$
  \begin{eqnarray*} \mathbb{P}(K_{n}=k) &=& \mathbb{P}(T_{-1} > n-k) \cdot  \mathbb{P}( {T}_0^+ > k \mbox{ and } S_{1} >0) \\
&\underset{ \mathrm{symm}}{=}&   \mathbb{P}(T_{-1} > n-k)  \cdot \frac{1}{2} \cdot \mathbb{P}( {T}_{-1} > k-1).  \end{eqnarray*}  For  $k=0$ we simply have $ \mathbb{P}(K_n=0)= \mathbb{P}(T_{-1}>n)$. Using \eqref{eq:T>ss} and Stirling's formula, the last display is shown to be equivalent to $\frac{1}{\pi} \frac{1}{ \sqrt{k(n-k)}}$   where the last asymptotic holds as $k$ and $n-k$ tend to $\infty$.  If we add a little blur to $K_{n}$ and consider $\tilde{K}_{n} = K_{n}+U_{n}$ where $U_{n}$ is independent of $K_{n}$ and uniformly distributed over $[0,1]$. Then clearly $ \tilde{K}_{n}/n$ has a density with respect to Lebesgue measure which converges pointwise towards the density of the arcsine law. It follows from Scheff\'e's lemma (Exercise \ref{exo:scheffe}) that $\tilde{K}_{n}/n$ converges in total variation towards the arcsine law and consequently $K_{n}/n$ converges in distribution towards the arcsine law since $U_{n}/n \to 0$ in probability. \qed


\subsection{Poisson random walk} \label{sec:poissonRW}
Another explicit application of Kemperman's formula is obtained by considering a random walk $(S)$ with step distribution given by the law of $  \mathfrak{P}({\alpha})-1$ where $ \mathfrak{P}({\alpha})$ is a Poisson random variable of parameter $\alpha$, namely 
$$ \mathbb{P}(X=k-1) = \mathrm{e}^{-\alpha} \frac{\alpha^k}{ k!}, \quad \mbox{ for } k \geq 0.$$
Clearly, if $\alpha > 1$ then the walk is transient and drifts towards $+\infty$. Using the additivity property of independent Poisson variables and  Kemperman's formula we have:
 \begin{eqnarray*} \mathbb{P}(T_{-1}=n) &\underset{ \mathrm{Kemperman}}{=} & \frac{1}{n} \mathbb{P}( S_n = -1)\\ 
 &\underset{ \sum \  \mathrm{id\ Poisson}}{=}& \frac{1}{n} \mathbb{P}( \mathfrak{P}({n \alpha}) = n-1) = \mathrm{e}^{-\alpha n} \frac{(\alpha n)^{n-1}}{n!}.  \end{eqnarray*}
 This law is named after Borel:\footnote{ \raisebox{-5mm}{\includegraphics[width=1cm]{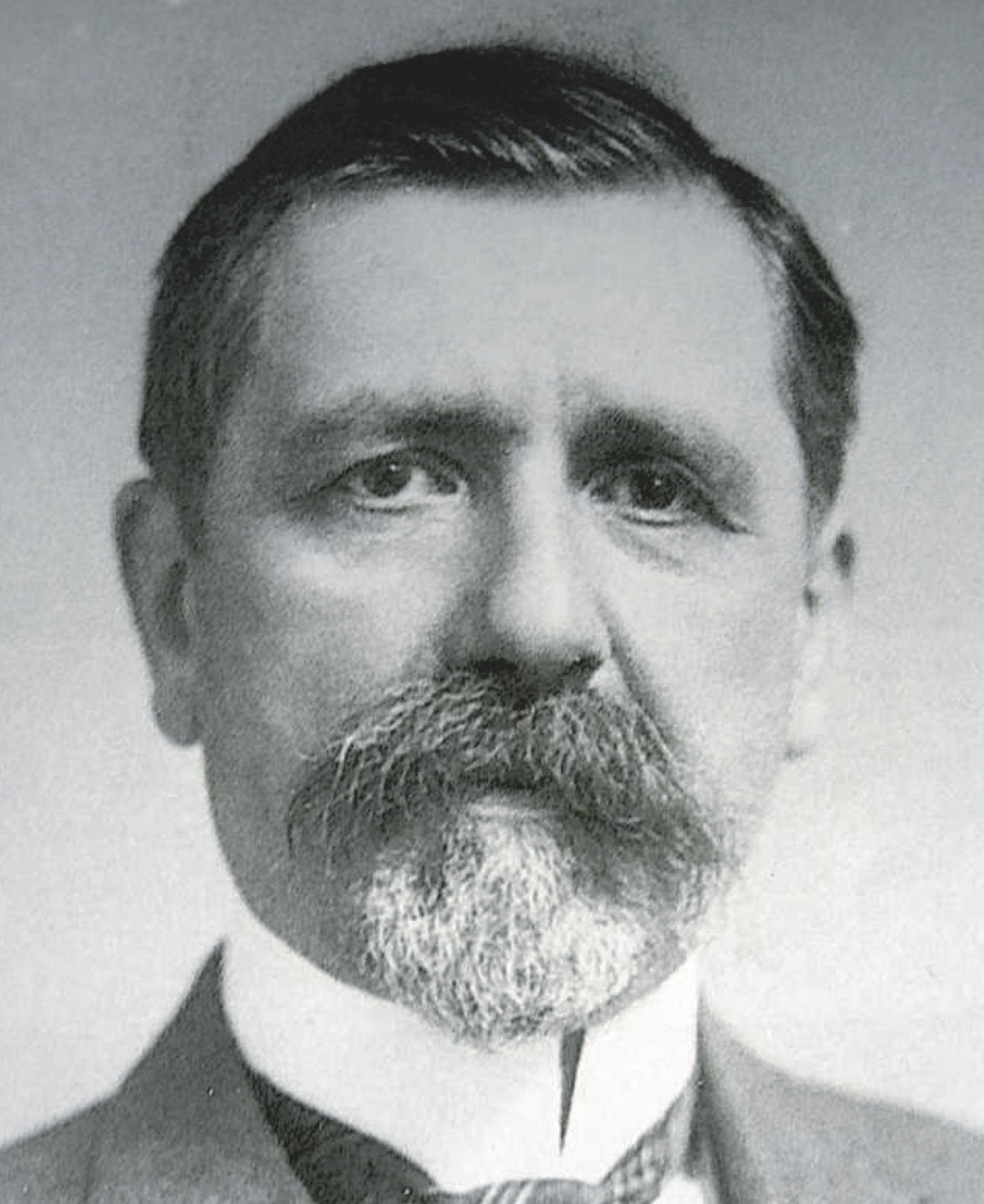}}Félix Édouard Justin Émile Borel (1871--1956), French.}
\begin{definition} \label{def:borel-tanner} For $\alpha \in [0,1]$, the Borel--Tanner distribution $\xi_\alpha$ is the law on $ \mathbb{Z}_{>0}$ given by 
$$ \xi_\alpha(n) = \mathrm{e}^{-\alpha n} \frac{(\alpha n)^{n-1}}{n!}, \quad \mbox{ for } n \geq 1.$$
\end{definition}
\begin{exo} Do you have an elementary way to see that the above display defines a probability distribution?
\end{exo}

\section{Ballot theorems}
Let us now turn our attention to ballot theorems when we require a positivity constraint on the walk. In the following we say that $(S)$ is \textbf{skip-free ascending} if $(-S)$ is a skip-free random process.
\subsection{Ballot theorem}
\begin{lemma} \label{lem:>0} Let $(S)$ be a skip-free ascending random walk. Then for every $n \geq 1$ and every $k \geq 1$ we have 
$$ \mathbb{P}( S_{i}>0, \forall 1 \leq i \leq n \mid S_{n}=k) = \frac{k}{n}.$$
\end{lemma}
\noindent \textbf{Proof.} Notice that the walk $(-S)$ is skip-free descending. So by time reversal (but not space reversal as in Lemma \ref{lem:duality})  and Kemperman's formula we have 
 \begin{eqnarray*} \mathbb{P}(S_{i}>0, \forall 1 \leq i \leq n \mbox{ and } S_{n}=k)  \underset{ \mathrm{time-rev.}}{=} \mathbb{P}(T_{-k}(-S) = n)
\underset{ \mathrm{Prop.}\  \ref{prop:kemperman}}{=} \frac{k}{n} \mathbb{P}(-S_{n}=-k)=  \frac{k}{n} \mathbb{P}(S_{n}=k).\end{eqnarray*} \qed \medskip 
 
Let us give an immediate application due to Bertrand\footnote{\raisebox{-5mm}{\includegraphics[width=1cm]{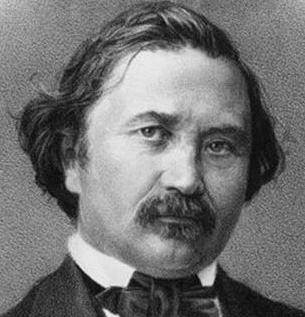}} Joseph Bertrand  (1822--1900), French} which is useful  during election days:
  \begin{theorem}[Ballot theorem]\label{the:ballot}During an election, candidates $A$ and $B$ respectively have $a >b$ votes. Suppose that votes are spread uniformly in the urn. What is the chance that during the counting of votes, candidate $A$ is always strictly ahead? 
  $$ \mbox{ answer: } \quad \frac{a-b}{a+b}.$$
  \end{theorem}
\noindent \textbf{Proof.} Let us model the scenario by a uniform path making only $+1$ or $-1$ steps which starts at $(0,0)$ and ends at $(a+b, a-b)$. The $+1$ steps correspond to votes for candidate $A$ and the $-1$ steps for votes for $B$. This path can be seen as the trajectory of a symmetric random walk conditioned to be equal to $a-b$ at time $a+b$.  The conclusion is given by the previous lemma. 
 \qed

\subsection{Staying positive forever}
\label{sec:stayposfacile}
Let $(S)$ be a one-dimensional random walk with integrable increments having positive mean.  Recall from Remark \ref{rek:>0pos>0} that the probability that walk stays positive after time $1$ is strictly positive. We compute below this probability in the case of skip-free ascending and skip-free (descending) walks:

\begin{corollary}  If $(S)$ is skip-free ascending such that $ \mathbb{E}[S_{1}] >0$ then we have 
$$ \mathbb{P}( S_{i}>0: \forall i \geq 1) = \mathbb{E}[S_{1}].$$
\end{corollary}
\noindent \textbf{Proof.} We have 
 \begin{eqnarray*} \mathbb{P}( S_{i}>0: \forall i \geq 1) &=& \lim_{n \to \infty} \mathbb{P}( S_{i}>0: \forall 1 \leq i \leq n)\\
 &=&\lim_{n \to \infty} \mathbb{E}\left[  \mathbb{P}\mathbb( S_{i}>0: \forall 1 \leq i \leq n \mid S_{n})\right]\\
  & \underset{ \mathrm{Lem.} \ref{lem:>0}}{=} &  \lim_{n \to \infty} \mathbb{E}\left[  \frac{S_{n}}{n} \mathbf{1}_{ S_{n} >0}\right]  = \mathbb{E}[S_{1}],  \end{eqnarray*}
  where for the last convergence we used the fact that  $\frac{S_n}{n} \mathbf{1}_{S_n>0} \to \mathbb{E}[S_{1}]$ almost surely by the strong law of large numbers together with the  fact that $|\frac{S_n}{n} \mathbf{1}_{S_n>0}| \leq 1$  (recall that $S_n \leq n$ since the walk is skip-free ascending) which enabled us to invoke dominated convergence. \qed 
  
\begin{proposition} \label{prop:GWdisguise}If $(S)$ is skip-free descending (with $\mu \ne \delta_{0}$) then  $ \mathbb{P}(S_{n} \geq 0, \forall n \geq 0)=1- \alpha$ where $\alpha$ is the smallest solution in $\alpha \in [0,1]$ to the equation:
  \begin{eqnarray} \label{eq:GWdisguise} \alpha = \sum_{k=-1}^{\infty} \mu_{k} \alpha^{k+1}.  \end{eqnarray}
\end{proposition}
\noindent \textbf{Proof.} Since $ \mu$ is supported by $\{-1,0,1,\dots\}$ its mean $m$ is well-defined and belongs to $[-1, \infty]$. We already know from the previous chapter that $ \mathbb{P}(S_{n} \geq 0, \forall n \geq 0)>0 $ if and only if $m >0$ (we use here the fact that the walk is not constant since $\mu \ne \delta_{0}$). We denote by $ T_{<0}$ the hitting time of $\{\dots,-3,-2,-1\}$ by the walk $(S)$. Since $(S)$ is skip free descending,if $T_{<0}$ is finite then necessarily $ S_{T_{<0}}=-1$. To get the equation of the proposition we perform one step of the random walk $S$: if $S_{1}=-1$ then $ T_{<0}< \infty$. Otherwise if $S_{1}\geq 0$ then consider the stopping times $$\theta_{0}=1,\quad  \theta_{1}= \inf\{ k\geq 1 : S_{k} = S_{1}-1\},\quad \theta_{2}= \inf\{ k\geq  \theta_{1} : S_{k} = S_{1}-2\}, \dots .$$ If $I = \sup\{ i \geq 0 : \theta_{i}< \infty\}$,  the strong Markov property shows that $(\theta_{i+1}-\theta_{i})_{0 \leq i \leq I}$ has the same law as  i.i.d.\ samplings of law $T_{ <0}$ until the first hit of $+\infty$. Furthermore on the event $S_{1}\geq 0$ we have 
$$ \{ T_{ <0} < \infty\} = \bigcap_{n=0}^{S_{1}} \{ \theta_{n+1}-\theta_{n} < \infty\}.$$
Taking expectation, we deduce that $ \mathbb{P}( T_{ < 0}< \infty)$ is indeed solution of  \eqref{eq:GWdisguise}. Now, notice that $F:\alpha \mapsto \sum_{k=-1}^{\infty} \mu_{k} \alpha^{k+1}$ is a convex function on $[0,1]$ which always admits $1$ as a fixed point. Since $F'(1) = m+1$ we deduce that $F$ admits two fixed points in the case $m>0$. But when $m>0$ we already know that $\alpha <1$ and so $\alpha$ must be equal to the smallest solution of \eqref{eq:GWdisguise}. \qed 
\begin{exo} \label{exo:32} Let $(S)$ be a skip-free descending random walk which drifts towards $+ \infty$. Compute the law of 
$\inf\{ S_{k} : k \geq 0\}$. Relate to Corollary \ref{cor:duality}. \end{exo}

\subsection{Parking on the line}
We finish our applications of the cycle lemma by giving a last, but nice, application of skip-free random walk to the parking problem on the line.
 Imagine an oriented discrete line with $n$ vertices, \textbf{the parking spots} (each vertex can only accommodate at most one car). The cars are labeled from $1$ up to $m$ and they arrive one after the other on some of the $n$ vertices. When arriving, they try to park at their arrival node, and,  if the parking spot is occupied, the cars drive towards the left of the line and park on the first available spot. If they do not find a free spot, they exit from the parking lot.
 
 \begin{figure}[!h]
  \begin{center}
  \includegraphics[width=15cm]{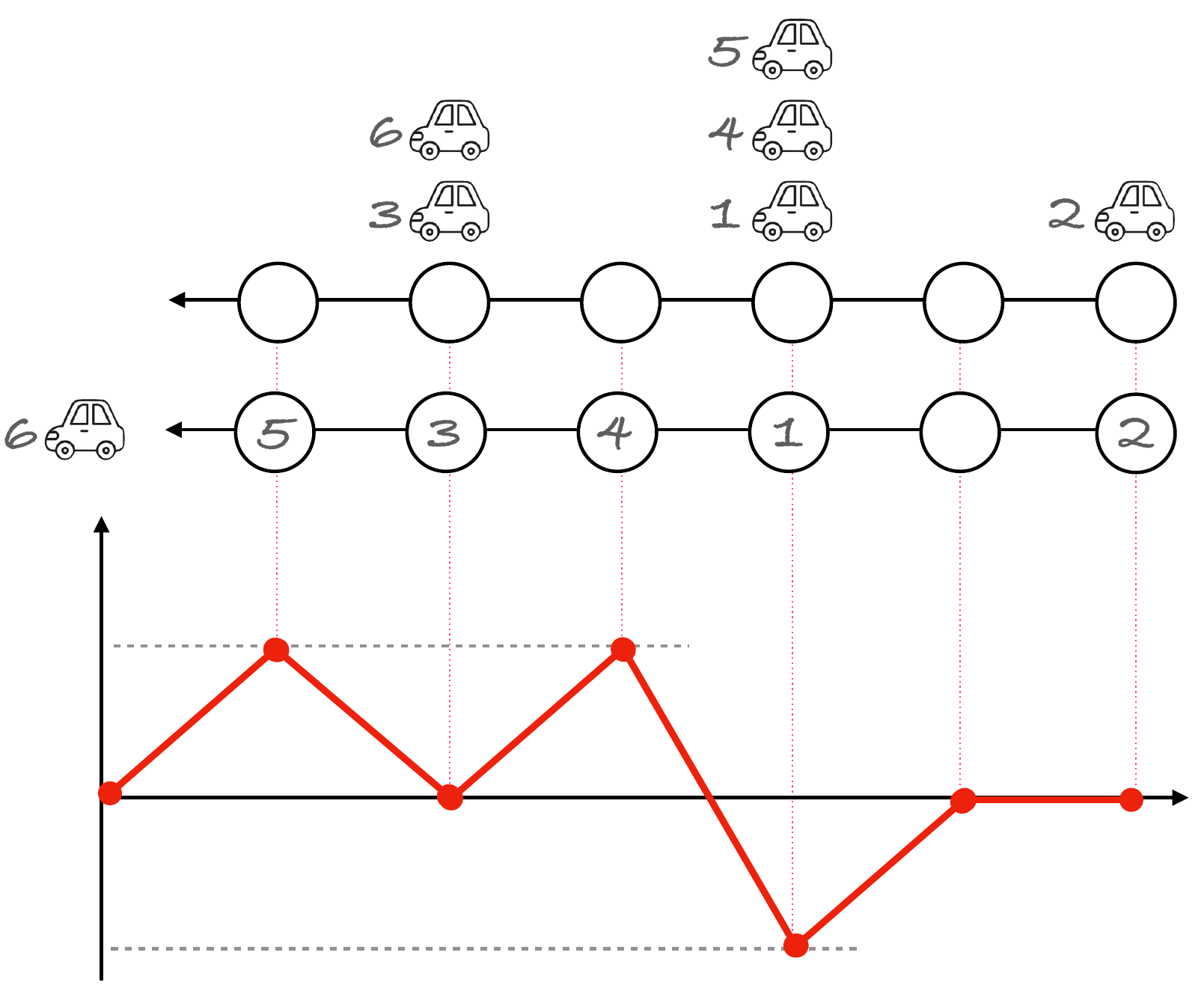}
  \caption{Illustration of the parking of $6$ cars on a line with $6$ spots. Notice that car number $6$ did not manage to find a spot and exited the parking lot. Below, the encoding of that the parking configuration by a skip-free ascending  walk. \label{fig:parking}}
  \end{center}
  \end{figure}
 An  Abelian property shows that the unlabeled final configuration as well as the number of cars exiting the parking lot does not depend on the order in which we try to park the cars (exercise!). In our random model, one shall imagine that the $m$ cars pick independently an arriving vertex $\in \{1,2, \dots , n\}$ uniformly at random. Of course, as long as $m >n$, it is impossible that all cars park.\bigskip

We shall prove the following theorem due to Konheim and Weiss:
\begin{theorem}[Konheim \& Weiss (1966)]\noindent Imagine that $m$ cars try to park uniformly and independently on $n$ vertices. The probability that they all manage to park is equal to 
$$  \frac{n+1-m}{n+1} \left( 1+\frac{1}{n}\right)^m.$$
In particular, if $m= [\alpha n]$ with $\alpha \in (0,1)$ the above probability converges to $(1-\alpha)  \mathrm{e}^{\alpha}$ as $n \to \infty$. \label{thm:parking}
\end{theorem}
\noindent \textbf{Proof}. The idea is to encode the parking situation by a walk $(S)$. Specifically, each vertex receiving $k$ cars corresponds to an increment of the walk of $1-k$, see Figure \ref{fig:parking}. The path we obtain this way is clearly skip-free ascending. By construction of the coding, for any $i \in \{1, \dots , n\}$ the value of the walk at time $i$ is equal to $i$ minus  the number of cars arriving on vertices on the left of it. It is easy to see that full parking for the cars corresponds to the fact that the walk stays non-negative until time $n$.  In our probabilistic model where the $m$ cars choose independently and uniformly their arrival vertices, the increments of the walk $(S)$ are not independent. However, we clearly have $S_n = n-m$ and the increments of this walk are exchangeable, we can thus apply Lemma \ref{lem:>0}. The slight problem is that Lemma \ref{lem:>0} evaluates the probability that the walk stays positive, and we need the probability that it stays \textit{non-negative}. To go around this problem, we imagine that we add a $(n+1)$th vertex at the left extremity of the line. Clearly, each successful parking configuration on $\{1, \dots , n\}$ corresponds to a single configuration where the $m$ cars choose to park on vertices in $\{1,\dots , n\}$ and that the vertex $n+1$ is empty at the end. In terms of the random walk, we precisely ask that it stays positive. Hence, by Lemma \ref{lem:>0}, the number of successful parking configurations with $m$ drivers and $n$ spots is equal to 
 $$ \frac{n+1-m}{n+1} \cdot (n+1)^m.$$
 The theorem follows immediately after dividing by the number of configurations in the initial model, i.e.~by $n^{m}$.
 \qed 
 
 \section{Wiener-Hopf factorization}
 In this section we extend the theory to the case of random walk with arbitrary step distribution $\mu$ which is non necessarily integer valued nor skip-free. We still denote $(S)$ a one-dimensional random walk starting from $0$ and with independent increments of law $\mu$ supported by $ \mathbb{R}$. We first need to introduce the so-called \textbf{ladder variables}. 
 
 \subsection{Ladder variables}
  
\label{sec:ladder}
  
  \begin{definition}[Ladder heights and epochs]  \label{def:prop:ladder} We define by induction $ T_{0}^{\footnotesize >}=T_{0}^{<}={T}_{0}^{\geq}= {T}_{0}^{\leq}=0$ as well as $ H_{0}^{>}=H_{0}^{<}={H}_{0}^{\geq}= {H}_{0}^{\leq}=0$ and for $i \geq 1$ we put 
 \begin{eqnarray*} 
 T_{i}^{>} &=& \inf\left\{ k > T_{i-1}^{>}: S_{k} > H_{i-1}^{>} \right\} \quad \mbox{ and } \quad H_{i}^{>} = S_{T_{i}^{>}}, \\
 {T}_{i}^{\geq} &=& \inf\left\{ k > {T}_{i-1}^{\geq}: S_{k} \geq {H}_{i-1}^{\geq} \right\} \quad \mbox{ and } \quad {H}_{i}^{\geq} = S_{{T}_{i}^{\geq}}, \\
 T_{i}^{<} &=& \inf\left\{ k > T_{i-1}^{<}: S_{k} < H_{i-1}^{<} \right\} \quad \mbox{ and } \quad H_{i}^{<} = S_{T_{i}^{<}}, \\
 {T}_{i}^{\leq} &=& \inf\left\{ k > {T}_{i-1}^{\leq}: S_{k} \leq {H}_{i-1}^{\leq} \right\} \quad \mbox{ and } \quad {H}_{i}^{\leq} = S_{{T}_{i}^{\leq}}.  \end{eqnarray*}
If $T_{i}^{*}$  is not defined (i.e.~we take the infimum over the empty set) then we put $T_{j}^{*}= H_{j}^{*} = \pm \infty$ for all $j \geq i$.  The variables $(T^{>}/{T}^{\geq})$ (resp. $(T^{<}/{T}^{\leq})$) are called the strict/weak ascending (resp.\ descending) ladder epochs. The associated $H$ process are called the (strict/weak ascending/descending) ladder heights.

  \end{definition}
\begin{figure}[!h]
 \begin{center}
 \includegraphics[width=16cm]{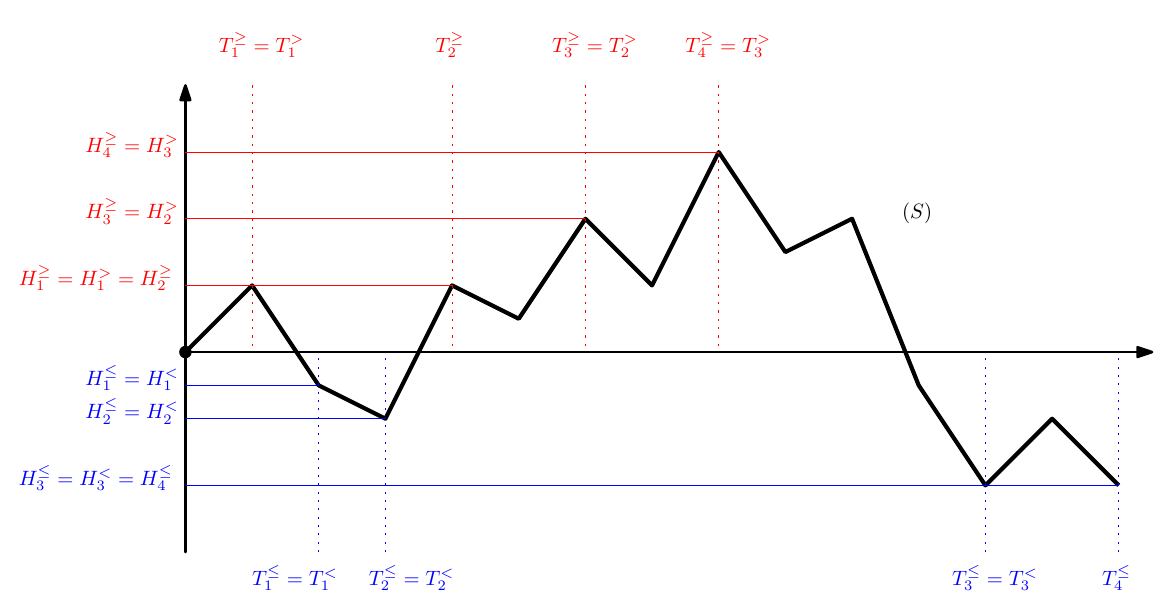}
 \caption{Illustration of the definition of the ladder heights and epochs.}
 \end{center}
 \end{figure}

When $\mu$ has no atoms, the walk $S$ does not take twice the same value a.s. so the weak and strict ladder variables are the same.
  In the following we write $H$ and $T$ generically for one of the four couples $(T^{\geq}, H^{\geq}),({T}^{>}, {H}^{>}),(T^{<}, H^{<}) \mbox{ or }(T^{\leq}, H^{\leq}).$ Since the ladder epochs are stopping times for the natural filtration generated by the walk, the strong Markov property then shows that $N=\inf\{ i \geq 0: T_{i} = \infty\}$  is a geometric random variable with distribution 
$$ \mathbb{P}( N = k ) = \mathbb{P}( T_{1} = \infty)\mathbb{P}( T_{1} < \infty)^{{k-1}},$$ and that conditionally on $N$ the random variables 
 $((T_{i}-T_{i-1}), (H_{i}-H_{i-1}))_{1 \leq i \leq N-1}$ are i.i.d.\ with law $(T_{1}, H_{1}) \mbox{ conditioned on } T_{1}< \infty$. In particular, $\limsup_{n \to \infty}S_n = +\infty$ a.s.~if and only if $ \mathbb{P}(T_1^> = \infty) = \mathbb{P}(T_1^\geq =\infty)=0$. \bigskip 
 
 One can now extend Feller's cycle lemma (Lemma \ref{lem:feller}) in this setup. The main difference is that when the walk is not skip-free, the number of records cannot be easily tightened to the value of the walk, that is why the ladders epoch and heights are needed. With the same notation as in Section \ref{sec:fellerskip}, we have the extension of \eqref{eq:equivfeller} (with mutatis mutandis the same proof): For every $n \geq 1$ and any measurable subset $A \subset \mathbb{R}_{+}^{*}$ we have 
  \begin{eqnarray*}  \mathbf{1}_{s_{n} \in A} = \sum_{i=0}^{{n-1}} \sum_{k=1}^{\infty} \frac{1}{k} \mathbf{1}_{T_{k}^{>}(s^{(i)})=n} \mathbf{1}_{H_{k}^{>}(s^{(i)})\in A}.  \end{eqnarray*}
Taking expectation and using the invariance of the walk by cycle shift we deduce the equality of measures generalizing Kemperman's formula:
  \begin{eqnarray} \label{eq:kempgen}\mathbf{1}_{x >0} \frac{ \mathbb{P}(S_{n} \in \mathrm{d}x)}{n}= \sum_{k=1}^{\infty} \frac{1}{k} \mathbb{P}( H_{k}^{>} \in \mathrm{d}x, T_{k}^{>} =n) \mathbf{1}_{x >0}. \end{eqnarray}

\subsection{Wiener--Hopf factorization}
The following result is an analytic translation of our findings.
\begin{theorem}[Spitzer--Baxter formula ; Wiener--Hopf factorization]\label{thm:WH}For $r \in [0,1)$ and $\mu \in \mathbb{C}$ so that $ \mathfrak{Re}(\mu)\geq0$ we have 
$$\left( 1 - \mathbb{E}\left[ r^{T_{1}^{>}} \mathrm{e}^{-\mu H_{1}^{>}} \right]\right) = \exp \left( - \sum_{n=1}^{\infty} \frac{r^{n}}{n} \mathbb{E}\left[ \mathrm{e}^{-\mu S_{n}} \mathbf{1}_{S_{n}>0} \right]	\right),$$
$$\left( 1 - \mathbb{E}\left[ r^{{T}_{1}^{\leq}} \mathrm{e}^{\mu {H}_{1}^{\leq}} \right]\right) = \exp \left( - \sum_{n=1}^{\infty} \frac{r^{n}}{n} \mathbb{E}\left[ \mathrm{e}^{\mu S_{n}} \mathbf{1}_{S_{n}\leq 0} \right]	\right).$$
\end{theorem}
\noindent \textbf{Proof.} First since $r \in [0,1)$ and $ \mathfrak{Re}(\mu) \geq0$ all the quantities in the last two displays are well defined. We only prove the first display since the calculation is similar for the second one. Let us start from the right hand side of the theorem and write
  \begin{eqnarray*}\exp \left( - \sum_{n=1}^{\infty} \frac{r^{n}}{n} \mathbb{E}\left[ \mathrm{e}^{-\mu S_{n}} \mathbf{1}_{S_{n}>0} \right]	\right) &\underset{  \eqref{eq:kempgen}}{=}& \exp \left( - \sum_{n=1}^{\infty} \frac{r^{n}}{n} \sum_{k=1}^{\infty} \frac{n}{k} \mathbb{E}\left[ \mathrm{e}^{-\mu H_{k}^{>}} \mathbf{1}_{T_{k}^{>}=n} \right]	\right) \\
  &=& \exp \Big( - \sum_{k=1}^{\infty} \frac{1}{k} \underbrace{\mathbb{E}\left[  \mathrm{e}^{-\mu H_{k}^{>}} r^{T_{k}^{>}} \right]}_{\left(\mathbb{E}\left[  \mathrm{e}^{-\mu H_{1}^{>}} r^{T_{1}^{>}} \right]\right)^{k} } \Big) \\
  &=& 1- \mathbb{E}\left[  \mathrm{e}^{-\mu H_{1}^{>}} r^{T_{1}^{>}} \right],  \end{eqnarray*}
  where in the last line we used the equality $\sum_{k =1}^{\infty}  \frac{x^{k}}{k} = -\log (1-x)$ valid for $ |x|<1$. Note that we implicitly used the fact that $r<1$ by putting  $r^{T_{k}^{>}}=0$ when $T_{k}^{>}= \infty$. This proves Spitzer's\footnote{\raisebox{-3mm}{\includegraphics[width=0.7cm]{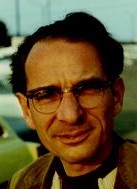}} Frank Ludvig Spitzer (1926--1992), Austrian \& American} formula\qed 

\begin{remark}[Explanation of the terminology of Wiener--Hopf factorization] If we write 
$$ \omega_{r}^{>}(\mu) = 	\exp \left( - \sum_{n=1}^{\infty} \frac{r^{n}}{n} \mathbb{E}\left[ \mathrm{e}^{-\mu S_{n}} \mathbf{1}_{S_{n}>0} \right]	\right) \quad \mbox{ and }\omega_{r}^{\leq}(\mu) = 	\exp \left( - \sum_{n=1}^{\infty} \frac{r^{n}}{n} \mathbb{E}\left[ \mathrm{e}^{-\mu S_{n}} \mathbf{1}_{S_{n}\leq0} \right]	\right),$$ then $\omega_{r}^{>}$ is analytic on the half-space $ \mathfrak{Re}(\mu) \geq 0$ whereas $\omega_{r}^{\leq}$ is analytic on $ \mathfrak{Re}(\mu) \leq 0$. On the imaginary line where the two functions are well defined we have 
 \begin{eqnarray} \omega_{r}^{>}(  it)\omega_{r}^{\leq}(  it) = 1-r \mathbb{E}[\mathrm{e}^{-it X_{1}}].   \label{eq:wiener-hopf}\end{eqnarray}
Hence, the characteristic function  of the increment of the walk (or a slight modification thereof) has been writing as a product of two analytic functions, each defined on a different half-space. The idea of writing a function on a line as a product of two functions defined on a half-space goes back to Wiener \& Hopf and is often useful since we can use the tools of complex analysis for each of the factors. \end{remark} 
\medskip 

There are many applications of the previous formula, we just mention two surprizing ones:
\begin{corollary}\label{cor:law return} Let $(S)$ be a one-dimensional random walk with symmetric and diffuse step distribution. Then the law of $T_{1}^{>}$ is given by 
$$ \mathbb{E}[r^{T_{1}^{>}}] = 1 - \sqrt{1-r}, \ r \in[0,1), \quad \mbox{ or equivalently } \quad \mathbb{P}(T_{1}^{>}= n) = \frac{(2n-2)!}{2^{2n-1} n! (n-1)!}, \ n \geq 1.$$
\end{corollary}
\noindent \textbf{Proof.} It suffices to take the first display of Theorem \ref{thm:WH} and to plug $\mu=0$. Since by symmetry of the increments and the lack of atoms we have $ \mathbb{P}(S_{n}>0)= \mathbb{P}(S_{n}\geq 0)= \frac{1}{2}$. It follows that
 \begin{eqnarray*}
 1 - \mathbb{E}[r^{T_{1}^{>}}] = \exp\left( - \sum_{n \geq1} \frac{r^{n}}{n} \mathbb{P}(S_{n}>0)\right) 
 = \exp\left( - \sum_{n \geq1} \frac{r^{n}}{n} \frac{1}{2}\right) = \exp(- \frac{1}{2} \log(1-r)) = \sqrt{1-r}. \end{eqnarray*}
 To get the exact values of $ \mathbb{P}(T_{1}^{>}= n)$ it suffices to develop $1 - \sqrt{1-r}$ in power series and to identify the coefficients. \qed

 \begin{corollary}[Back to the law of large numbers, again!] The random walk $(S)$ drifts towards $-\infty$ if and only if 
 $$ \log \mathbb{E}[{T}_{1}^{<}] = \sum_{n \geq 1} \frac{\mathbb{P}(S_{n}\geq0)}{n} < \infty.$$
 \end{corollary}
 \noindent \textbf{Proof.} From Theorem \ref{thm:WH} with $\mu=0$ we get for $r \in [0,1)$
 $$ 1- \mathbb{E}[r^{T_{1}^{\geq}}] = \exp \left( - \sum_{n\geq 1} \frac{r^{n}}{n} \mathbb{P}(S_{n}\geq0)\right).$$
 Letting $r \uparrow 1$ the left-hand side converges towards $1 - \mathbb{E}[ \mathbf{1}_{T_{1}^{\geq}<\infty}] = \mathbb{P}(T_{1}^{\geq}=\infty)$ whereas the right-hand side converges towards $\exp(-\sum_{n \geq 1} \frac{\mathbb{P}(S_{n}\geq0)}{n})$. But clearly $(S)$ drifts towards $-\infty$ if and only if $T_{1}^{\geq}$ may be infinite. In this case, recall that by \eqref{eq:duality*} and the fact that the increments of the ladders variables are independent that we have $ \mathbb{E}[T_{1}^{<}] = 1/ \mathbb{P}(T_{1}^{\geq}=\infty)$ which immediately implies the second claim. 
 \qed \bigskip 
 



    \noindent \textbf{Biliographical notes.} The study of skip-free random walk may be seen as a particular case of fluctuation theory for random walks, see e.g.~\cite{kyprianou2010wiener} for a more trajectorial approach. The combinatorial approach taken here and based on the cycle lemma is adapted from \cite[Chapter XII]{Fel71} and \cite[Section 8.4]{Chung74}; it has many ramifications in the combinatorial literature, see \cite{ABR08} for much more about Ballot theorems and \cite{diaconis2017probabilizing} for parking functions. Theorem \ref{thm:parking} can be found in \cite{konheim1966occupancy}. The proof of the law of large numbers based on duality is taken from \cite{CurLLN}. In general, path transformations are very useful tools in fluctuation theory for random walks (Spitzer-Baxter or Wiener-Hopf factorization). In particular, we mention the Sparre-Andersen identity relating the position of the maximum and the time spent on the positive half-line for a random walk of length $n$, see \cite[Chapter XII]{Fel71} for more details. More recent applications of fluctuation theory for random walks can be found e.g.~in \cite{alili2005fluctuation,marchal2001two,kwasnicki2020random}. \medskip 
    
    \noindent \textbf{Hints for exercises:}\ \\
    Exercise \ref{exo:legallcornell}: is  \cite[Lemma 1.9]{LG05} (but the proof there is different). This is a baby example of the Spitzer-Baxter or Wiener-Hopf factorization.\\
    Exercise \ref{exo:shiftunif}: the $n$ distinct cycle shifts  are such that $G_{n} = \{1, 2,\dots , n\}$.\\
    Exercise \ref{exo:32}: it is a geometric distribution.

\chapter{Bienaym\'e-Galton-Watson trees}

\label{chap:GW}

\hfill I will survive.
\bigskip 

In this chapter we use our knowledge on one-dimensional random walk to study random tree coding for the genealogy of a population where individuals reproduce independently of each other according to the same offspring distribution. These are the famous Bienaym\'e--Galton--Watson (BGW) trees.

            \begin{figure}[!h]
         \begin{center}
         \includegraphics[width=13cm]{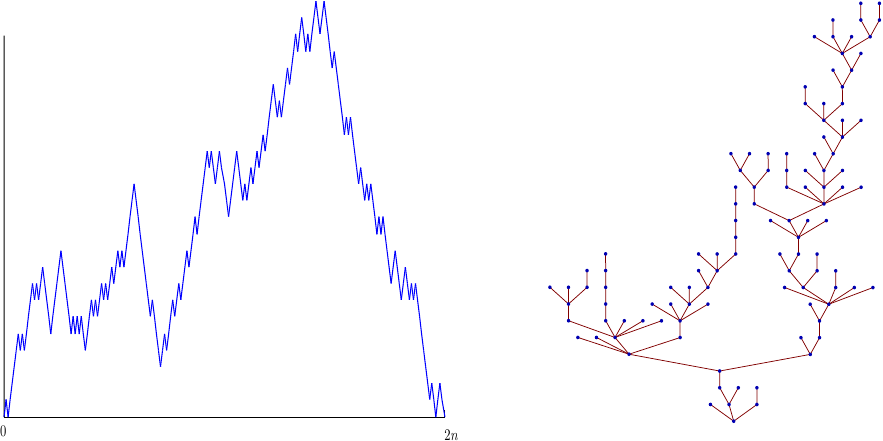}
         \caption{A large Bienaym\'e--Galton--Watson tree and its contour function}
         \end{center}
         \end{figure}
\section{Plane trees and Bienaym\'e--Galton--Watson processes}
\subsection{Plane trees}\label{sec:plane trees}

Throughout this chapter we will use the standard formalism for plane
trees as found in \cite{Nev86}. Let
 \begin{eqnarray*}
\mathcal{U}& = &\bigcup_{n=0}^{\infty} ( \mathbb{Z}_{>0})^n
 \end{eqnarray*}
where we recall that $\mathbb{Z}_{>0} = \{ 1,2, \ldots \}$ and $(\mathbb{Z}_{>0})^0 = \{
\varnothing \}$ by convention. An element $u$ of $\mathcal{U}$ is thus
a finite sequence of positive integers which we interpret as a \textbf{word} whose letters are positive integers. We let $|u|$ be the length of the word $u$.  If $u, v \in \mathcal{U}$,
$uv$ denotes the concatenation of $u$ and $v$. If $v$ is of the form
$uj$ with $j \in \mathbb{Z}_{>0}$, we say that $u$ is the \textbf{parent} of
$v$ or that $v$ is a \textbf{child} of $u$. More generally, if $v$ is of
the form $uw$, for $u,w \in \mathcal{U}$, we say that $u$ is an
\textbf{ancestor} of $v$ or that $v$ is a \textbf{descendant} of $u$. 

\begin{definition}\label{def:planetree} A
\textbf{plane tree} $\tau$ is a (finite or infinite) subset of
$\mathcal{U}$ such that
\begin{enumerate}
\item $\varnothing \in \tau$, the point $\varnothing$ is called the \textbf{root} of $\tau$,
\item if $v \in \tau$ and $v \neq \varnothing$ then the parent of $v$ also   belongs to $\tau$,
\item for every $u \in \tau$ there exists $k_u(\tau)\in \{ 0,1,2, \dots \} \cup \{ \infty\}$
  such that $uj \in \tau$ if and only if $j \leq k_u(\tau)$. The number $k_u( \tau)$ is then \textbf{the number of children} of $u$ in $\tau$.
\end{enumerate}
\end{definition}

\begin{figure}[!h]
 \begin{center}
 \includegraphics[width=5cm]{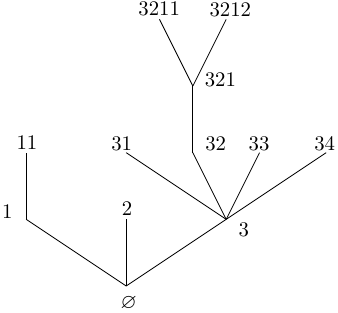}
 \caption{A (representation of a) finite plane tree.}
 \end{center}
 \end{figure}
 
 Since every $u \in \tau \backslash \{ \varnothing \}$  has a unique parent, we deduce that for finite plane trees $\tau$ we have 
  \begin{eqnarray} \label{eq:enfantparent}
   \#\tau -1= \sum_{u \in \tau}k_{u}(\tau).  \end{eqnarray}
A plane tree can be seen as a graph, in which an edge links
two vertices $u,v$ such that $u$ is the parent of $v$ or vice-versa. Notice that with our definition, vertices of infinite degree are allowed since $k_{u}(\tau)$ may be infinite. When all degrees are finite, the tree is said to be \textbf{locally finite}. In this case, this graph is of course a tree in the graph-theoretic sense (see Proposition \ref{prop:tree}), and we can draw it in the plane $ \mathbb{R}^2$ so that its edges are non-crossing and such that the edges from a vertex $u$
to its children $u\,1,\ldots,u\,k_u(\tau)$ and to its parent if $u \ne \varnothing$ are ordered in a clockwise fashion. Equivalently, a plane tree can be seen as a genealogical tree where the children of each vertex are ranked from the oldest to the youngest one. Unless explicitly mentioned, all the trees considered in this chapter are plane trees.

\begin{definition} \label{def:ulam} The set $ \mathcal{U}$ is a plane tree where $k_{u}( \mathcal{U})= \infty, \forall u \in \mathcal{U}$. It is called Ulam's tree. \end{definition}
 The integer $\#\tau$ denotes the number of vertices of $\tau$ and is called
the \textbf{size} of $\tau$. 
For any vertex $u \in \tau$, we denote the shifted tree at $u$ by $\sigma_u( \tau) :=\{v \in  \mathcal{U} : uv \in \tau\}$. The \textbf{height} of the tree $\tau$ is the maximal length of its words, 
$$ \mathrm{Height}(\tau) = \max \{ |u| : u \in \tau\} \in \{0,1,2,\dots\} \cup \{\infty\}.$$
The truncation at level $n$ of $\tau$ is denoted by $ [\tau]_n = \{ u \in \tau : |u| \leq n\}$ which is again a plane tree. Its boundary $\partial [\tau]_n$ is made of the individuals at generation exactly $n$ in the genealogical interpretation 
$$ \partial [\tau]_n = \{ u \in \tau : |u|=n \}.$$

\subsection{Bienaym\'e--Galton--Watson trees}

\label{sec:GWtrees}
Let $\mu$ be a distribution on $\{0,1,2,\dots \}$ which we usually suppose to be different from $\delta_{1}$. Informally speaking, a Bienaym\'e--Galton--Watson\footnote{\raisebox{-5mm}{\includegraphics[width=1cm]{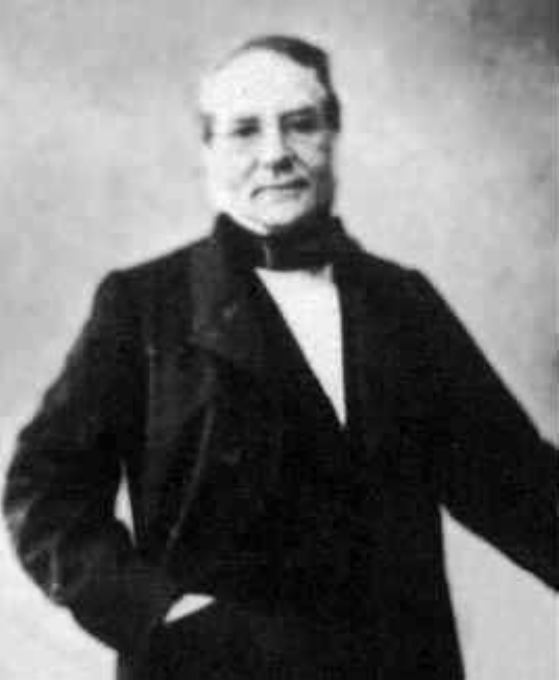}} \begin{tabular}{l}Irénée-Jules Bienaymé\\ (1796--1878), French \end{tabular} $\quad$  \raisebox{-5mm}{\includegraphics[width=1cm]{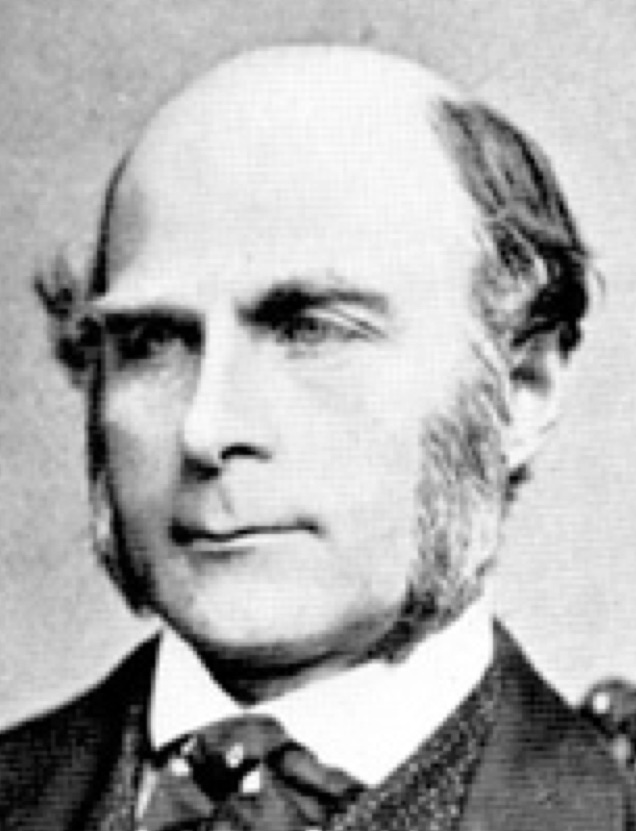}} \begin{tabular}{l} Francis Galton \\ (1822--1911), English \end{tabular} \ \ and \ \ \raisebox{-5mm}{\includegraphics[width=1cm]{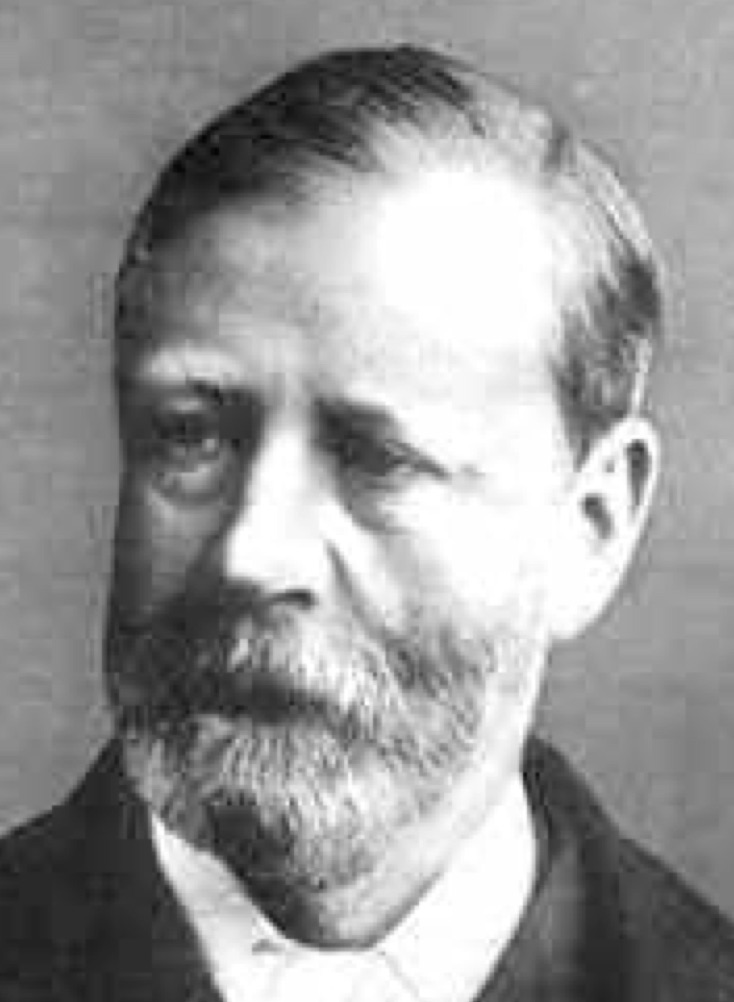}} \begin{tabular}{l}Henry William Watson\\ (1827--1903), English \end{tabular}} (BGW in short)
 tree with \textbf{offspring distribution} $\mu$ is a random (plane) tree coding the genealogy of a population starting with one individual and where all individuals reproduce independently of each other according to the distribution $\mu$. Here is the proper definition:
\begin{definition}[BGW tree] Let $(K_{u}:u \in \mathcal{U})$ be independent and identically distributed random variables of law $\mu$. We let $\mathcal{T}$ be the random plane tree made of all words $u = j_{1}j_{2}\dots j_{n} \in \mathcal{U}$ such that $j_{i} \leq K_{j_{1}\dots j_{i-1}}$ for all $1 \leq i \leq n$. In particular we have $k_u( \mathcal{T}) = K_u$ for all $u \in \mathcal{T}$. Then the law of $\mathcal{T}$ is the $\mu$-BGW distribution.
\end{definition}

Equivalently, the law of a $\mu$-BGW  tree $\mathcal{T}$ is characterized by the following \textbf{branching property}: Conditionally on the event $\{k_{ \varnothing}(\mathcal{T}) = \ell\}$ of probability $\mu_\ell$, then the $\ell$ random trees $\sigma_{i}(\mathcal{T})$ for $1 \leq i \leq \ell$ are independent and distributed as $\mathcal{T}$. Notice also that the $\mu$-BGW probability of a  \textit{finite} plane tree is explicit:
  \begin{eqnarray} \label{exo:GWprod}  \mathbb{P}(\mathcal{T} = \tau_{0}) = \prod_{u \in \tau_{0}} \mu_{k_{u}(\tau_{0})},  \end{eqnarray} but the previous display does not characterize the distribution since the random tree $ \mathcal{T}$ may very well be infinite.

We now link the BGW tree to the well-known BGW process. We first recall its construction. Let $(\xi_{i,j} : i \geq 0, j \geq 1)$ be i.i.d.\ random variables of law $\mu$. The $\mu$-Bienaym\'e--Galton--Watson process is defined by setting $Z_{0}=1$ and for $i \geq 0$ 
$$ Z_{i+1} = \sum_{j=1}^{Z_{i}} \xi_{i,j}.$$
It is then clear from the above construction that  if $\mathcal{T}$ is a $\mu$-Bienaym\'e-Galton--Watson tree, then the process $X_{n} = \# \{ u  \in \mathcal{T} : |u|=n\}$ has the law of a $\mu$-Bienaym\'e--Galton--Watson process.

\section{{\L}ukasiewicz walk and direct applications}
In this section we will encode (finite) trees via one-dimensional walks. This will enable us to get information on random BGW trees from our previous study of one-dimensional random walks. 
\subsection{{\L}ukasiewicz walk}
The \textbf{lexicographical or depth first} order $<$ on $ \mathcal{U}$ is defined as the reader may imagine: if $u=i_{1}i_{2}\dots i_{n}$ and $v= j_{1}j_{2}\dots j_{m}$ are two words then $u < v$ if $i_{\ell}< j_{\ell}$ where $\ell$ is the first index where $i_{ \ell} \ne j_{\ell}$, or if $n <m$ and $i_{1}i_{2}\dots i_{n} = j_{1}j_{2}\dots j_{n}$. 
The \textbf{breadth first} order on  $ \mathcal{U}$ is defined by $ u \prec v$ if $|u| < |v|$ and if the two words are of the same length then we require $u <v$ (for the lexicographical order).
\begin{definition} Let $\tau$ be a locally finite tree (i.e.~$k_u( \tau) < \infty$ for every $u \in \tau$). Write $u_{0}, u_{1}, \ldots$ for its vertices listed in the breadth first order. The 
  \textbf{{\L}ukasiewicz walk}  $ \mathcal{W}( \tau)= ( \mathcal{W}_n( \tau), 0 \leq n \leq  \#\tau)$ associated to $\tau$ is given by $ \mathcal{W}_0( \tau)=0$ and for $0 \leq n \leq \#\tau-1$:
$$ \mathcal{W}_ {n+1}( \tau)= \mathcal{W}_n( \tau)+k_ {u_{n}}( \tau)-1.$$
\end{definition}

\begin{figure}[!h]
 \begin{center}
 \includegraphics[width=15cm]{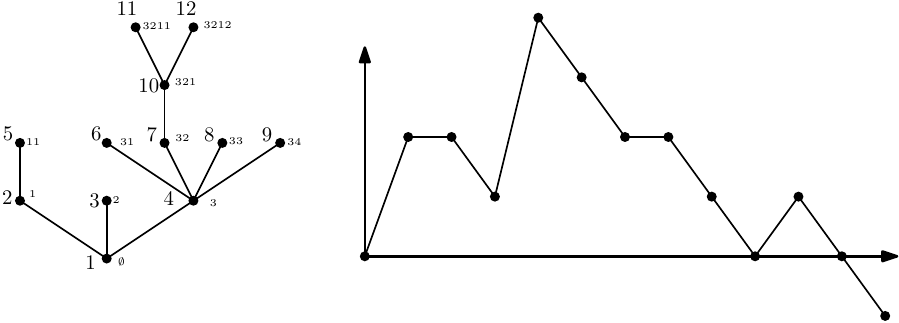}
 \caption{Left: a finite plane tree and its vertices listed in breadth-first order. Right: its associated {\L}ukasiewicz walk.}
 \end{center}
 \end{figure}
 In words, the {\L}ukasiewicz\footnote{\raisebox{-5mm}{\includegraphics[width=1cm]{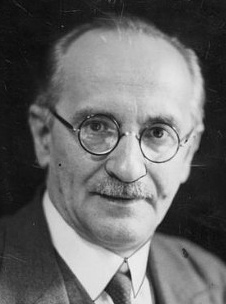}} Jan {\L}ukasiewicz  (1878--1956), Polish logician.}  walk consists in listing the vertices in breadth first order and making a stack by adding the number of children of each vertex and subtracting one (accounting for the exploration of the current vertex). In the case of a finite plane tree $\tau$, since the total number of children is equal to the number of vertices minus one, the following properties of $\mathcal{W}_{\cdot}( \tau)$ are easily checked:
 \begin{itemize}
 \item  the {\L}ukasiewicz walk starts at $0$, i.e.~$$\mathcal{W}_{0}( \tau)=0,$$
 \item it stays non-negative as long as all vertices have not been explored, i.e. $$\mathcal{W}_{i}( \tau)\geq 0 \quad \mbox{ for } 0 \leq i \leq \#\tau-1,$$
 \item it ends up at $-1$, i.e. $$W_{\#\tau}(\tau) =-1,$$
 \item the walk is skip-free in the sense of Chapter \ref{chap:WH}, i.e.$$  \mathcal{W}_{i+1}( \tau) - \mathcal{W}_{i}(\tau) \geq -1, \quad \mbox{  for any } 0 \leq i \leq \#\tau-1.$$
 \end{itemize}
  When the tree is infinite but locally finite,  every vertex of the tree will appear in the 
breadth first ordering\footnote{this is not true if we had chosen to explore the tree in the lexicographical (i.e.~depth first) order.} and the {\L}ukasiewicz path stays non-negative for ever. We leave the proof of the following as an exercise for the reader:
\begin{proposition} \label{exo:bijection}Let $  \mathbf{T}_{\ell oc }$ the set of all finite or infinite but locally finite plane trees. Let $ \mathbf{W}_{\ell oc }$ the set of all finite or infinite paths $ (w_{0},w_{1}, \dots , w_{n})$ with $n \in \{1,2, \dots \} \cup \{ \infty\}$  which starts at $w_{0}=0$ and ends at $w_{n}=-1$ and such that $w_{i+1}-w_{i} \geq -1$ as well as $w_{i} \geq 0$ for any $0 \leq i \leq n-1$. Then taking the {\L}ukasiewicz walk creates a bijection between $ \mathbf{T}_{\ell oc }$ and $ \mathbf{W}_{\ell oc }$.
\end{proposition}

\begin{remark}[Different types of exploration] The {\L}ukasiewicz path encodes the information when we discover a tree using the breadth first search. Although we shall only use this exploration in these notes, one can similarly discover the tree using the depth first search (i.e.~using the lexicographical total order to enumerate the vertices of a plane tree) or using more exotic type of exploration. In particular, the exploration of the Erd{\H{o}}s--R\'enyi graph (Chapter \ref{chap:poissonER}) will be based on a depth-first exploration. This flexibility in the exploration algorithm is at the core of many nice results in random tree theory, see e.g.~\cite{broutin2016new,CurStFlour,krivelevich2013phase}. See also the next chapters where the idea of discovering the underlying geometry with a given algorithm plays a key role.
\end{remark}

\subsection{{\L}ukasiewicz walk of a Bienaym\'e--Galton--Watson tree}
\begin{figure}[!h]
 \begin{center}
 \includegraphics[width=12cm]{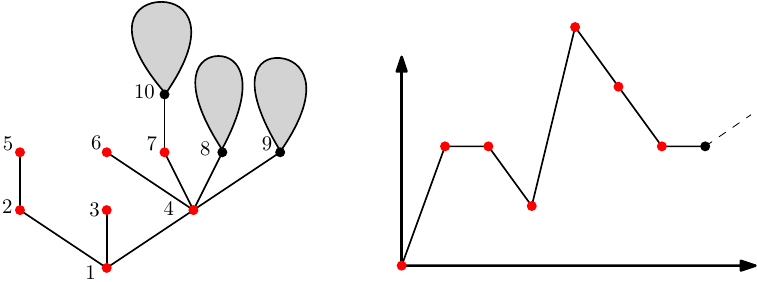}
 \caption{ \label{fig:partialuka} Fixing the first $n$ vertices  explored  (in red) during the breadth first exploration of a BGW tree. The black vertices and their subtrees (in gray) have not been explored yet.}
 \end{center}
 \end{figure}

As it turns out, the {\L}ukasiewicz walk associated to a $\mu$-BGW tree is roughly speaking a random walk. Recall that the offspring distribution $\mu$ is supported by $\{0,1,2, \dots \}$ so a $\mu$-BGW tree is locally finite a.s.
 \begin{proposition}\label{prop:luka}Let $ \mathcal{T}$ be a $\mu$-BGW tree, and let $(S_{n})_{n \geq0}$ be a random walk with i.i.d.\ increments of law $ \mathbb{P}(S_{1}=k) = \mu_{k+1}$ for $k \geq -1$. If $T_{-1}$ is the first hitting time of $-1$ by the walk $S$ (we may have $T_{-1}= \infty$) then we have 
 $$  \left( \W_{0}( \mathcal{T}), \W_{1}( \mathcal{T} ), \ldots, \W_{ \#\mathcal{T} }( \mathcal{T}) \right)  \quad \overset{(d)}{=} \qquad \ (S_{0},S_{1}, \ldots, S_{ T_{-1}}). $$
 \end{proposition}
\noindent \textbf{Proof.} Let $ (\omega^{0}_i: 0 \leq i \leq n)$ be the first $n$ steps of a skip-free random walk so that $n$ is less than or equal to the hitting time of $-1$ by this walk. By reversing the {\L}ukasiewicz construction we see that in order that the first $n$ steps of the {\L}ukasiewicz walk of the tree $ \mathcal{T}$ matches with $ (\omega^{0}_i: 0 \leq i \leq n)$ then the subtree $\tau_0$ of the first $n$ vertices of $ \mathcal{T}$ in breadth first order as well as their number of children are fixed by  $(\omega^{0}_i: 0 \leq i \leq n)$, see Figure \ref{fig:partialuka}.

 The probability under the $\mu$-BGW to see this event is  equal to
$$\prod_{u \in \tau_{0}} \mu_{k_{u}( \mathcal{T})} = \prod_{i=0}^{n-1} \mu_{\omega^{0}_{i+1}-\omega^0_i+1} = \mathbb{P}\left( (S_{i})_{0 \leq i \leq n} = (\omega^{0}_i)_{0 \leq i \leq n}\right).$$ The proposition follows. \qed \medskip

Combining the previous proposition with Remark \ref{rek:kemp+local} we deduce that if the offspring distribution is critical, aperiodic and has finite variance we have 
$$ \mathbb{P}( \#\mathcal{T}=n)\sim \frac{1}{\sqrt{2 \pi \sigma^{2}}} \cdot \frac{1}{n^{3/2}}, \quad \mbox{ as }n \to \infty.$$

%
%

\paragraph{Extinction probability.}
As a direct application of the previous proposition let us give a random walk proof of the following well-known criterion for survival of a Bienaym\'e--Galton--Watson process:
\begin{theorem}[Extinction probability]\label{thm:extinction}Let $\mu$ be an offspring distribution of mean $m \geq 0$ such that $\mu \ne \delta_{1}$. The probability that $ \mathcal{T}$ is finite is equal to the smallest solution $\alpha \in [0,1]$ to the equation 
 \begin{eqnarray} \alpha = \sum_{k\geq 0} \mu_{k} \alpha^k,   \label{eq:GWfixed}\end{eqnarray}
  in particular it is equal to $1$ if $m \leq 1$.
\end{theorem}
\noindent \textbf{Proof.}  With the same notation as in  Proposition \ref{prop:luka} we have that $ \mathbb{P}( \# \mathcal{T} = \infty) = \mathbb{P}( T_{-1} = \infty)$. Since the walk $S$ is non trivial (i.e.~non constant) and skip-free,   Proposition \ref{prop:GWdisguise} yields the statement. \qed \bigskip

Let us also recall the more ``standard'' proof of the previous theorem which is useful in Exercise \ref{exo:dekking}.  Let $ g(z) = \sum_{ k \geq 0} \mu_{k} z^{k}$ be the generating function of the offspring distribution $\mu$. In particular, if $ \mathcal{T}$ is a $\mu$-BGW tree then $g$ is the generating function of $ \#\partial [ \mathcal{T}]_1$. More generally, if $g_n$ is the generating function of $ \# \partial [ \mathcal{T}]_n$, then by the branching property of BGW trees and  standard operation on generating functions we have that $g_{n+1} = g \circ g_n$ for $n \geq 1$ so that 
$$g_n =  g \circ g \circ \cdots \circ g,$$ ($n$-fold composition). Specifying at $z=0$ we deduce that $ u_{n} = \mathbb{P}( \mathrm{Height}( \mathcal{T}) \leq n)$ follows the recurrence relation $ u_{0}=0$ and $u_{n+1} = g(u_{n})$. This recursive system is easily studied and $u_{n}$ converges towards the first fixed point of $g$ in $[0,1]$ which is strictly less than $1$ if and only if $g'(1) >1$ by convexity of $g$. We conclude using the fact that $ \{\mathcal{T} = \infty\}$ is the decreasing limit of the events $\{ \mathrm{Height}( \mathcal{T})\geq n\}$ as $n \to \infty$.

\begin{figure}[!h]
 \begin{center}
 \includegraphics[width=14cm]{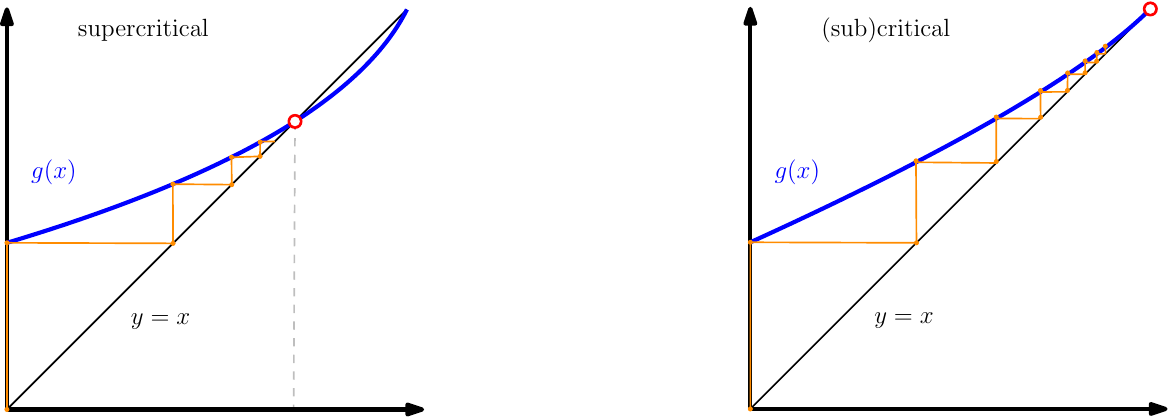}
 \caption{Illustration of the ``standard'' proof Theorem \ref{thm:extinction}. The extinction probability is computed as the limit of the recursive system defined by $u_0 =0$ and $u_{n+1} = g(u_n)$. \label{fig:GWclassic}}
 \end{center}
 \end{figure}

\begin{exo}[A theorem of Dekking \cite{Dek91b} and a \textbf{discontinuous phase transition}] \label{exo:dekking}
We say that an infinite tree $\tau$ contains an infinite binary tree (starting at the root) if it is possible to find a subset $S$ of vertices of $\tau$ containing the origin $\varnothing$ and such that each vertex in $S$ has exactly two children in $ S$. Let $ g(z) = \sum_{ k \geq 0} \mu_{k} z^{k}$ be the generating function of the offspring distribution $\mu$.
\begin{enumerate}
\item Show that the probability that a $\mu$-BGW tree $ \mathcal{T}$ contains no infinite binary tree (starting at the root) is the smallest solution $z \in [0,1]$ to 
$$ z = g(z) + (1-z) g'(z).$$
\item Application: in the case $p_{1}= (1-p)$ and $p_{3}=p$  with $p \in [0,1]$ show that there is no infinite binary tree in $ \mathcal{T}$ if and only if $p < \frac{8}{9}$ and that in the critical case $p= \frac{8}{9}$ this probability is in fact positive (contrary to the above case for survival of the tree).
\begin{figure}[!h]
 \begin{center}
 \includegraphics[width=8cm]{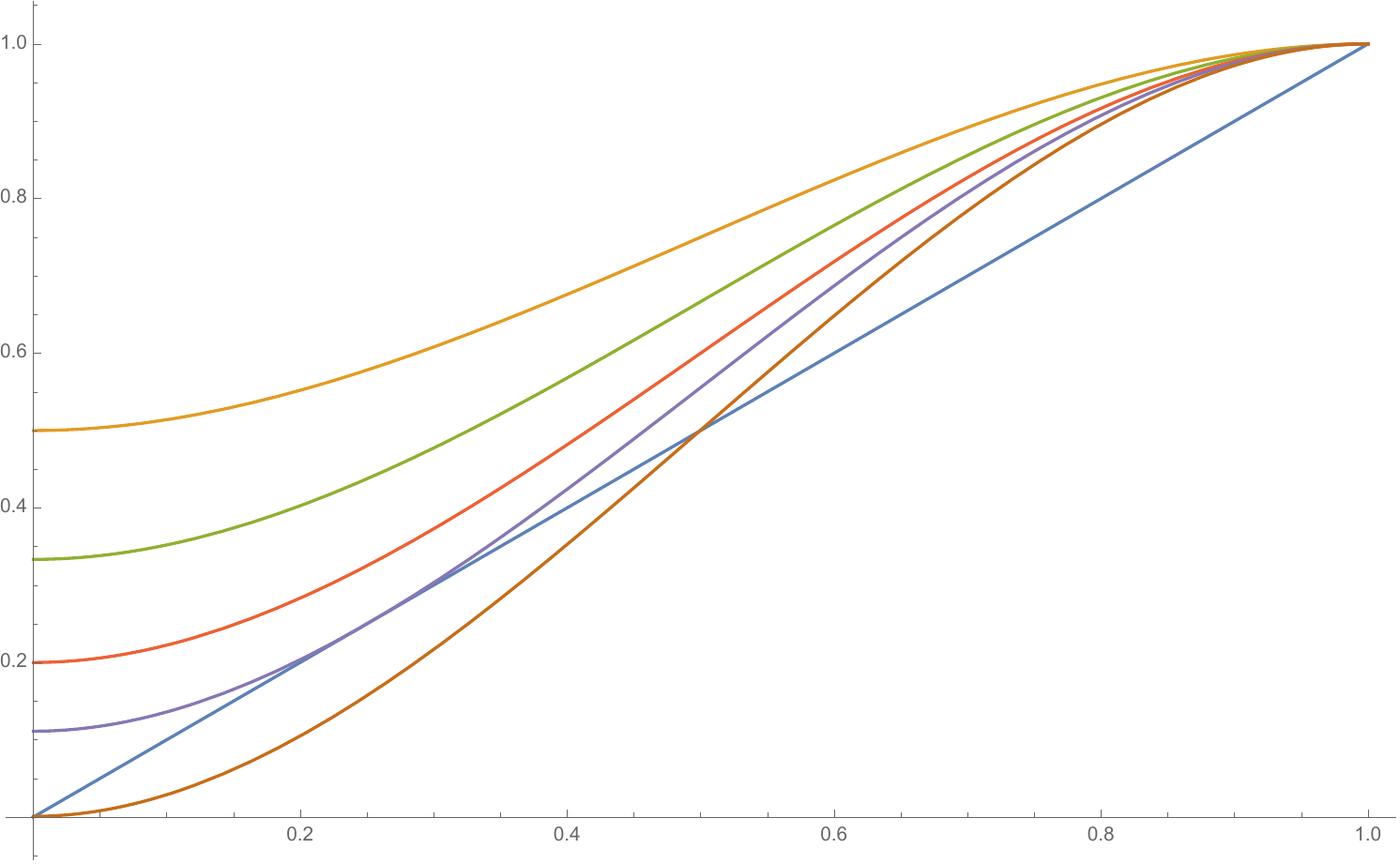}
 \caption{Plot of the function $g(z) + (1-z) g'(z)$ against the first bissector (in blue) where $g(z) = (1-p)z+p z^{3}$ for the values $ p = \frac{1}{2}, \frac{2}{3}, \frac{4}{5}$ in (yellow, green, red), the critical case $p= \frac{8}{9}$ in purple and $p = 1$ in brown.}
 \end{center}
 \end{figure}
\end{enumerate}
\end{exo}

         \begin{remark}[A historical remark] We usually attribute to Galton and Watson the introduction and study of the so-called Galton--Watson process in 1873 in order to study the survival of family names among British lords. However, in their initial paper devoted to the calculation of the extinction probability they concluded hastily that the extinction is almost sure whatever the offspring distribution! This is even more surprising since almost thirty years before, in 1845  Bienaymé considered the very same model and derived correctly the extinction probability. This is yet just another illustration of Stigler's law of eponymy!
         \end{remark}

\subsection{Lagrange inversion formula}
The Lagrange inversion is a closed formula for the coefficients of the reciprocal (composition inverse) of a power series. More precisely, imagine that $f(z) = \sum_{i \geq 0} f_{i} z^{i} \in \mathbb{C}[[z]]$ is a formal power series in the indeterminate $z$ (no convergence conditions are assumed) so that $f_{0}=0$ and $f_{1} \ne 0$. We recall the notation $[z^i]f(z) = f_i$. One  would like to invert $f$ i.e.~to find a power series $\phi \in \mathbb{C}[[z]]$ such that $ z = \phi( f(z)) = f(\phi(z))$. In combinatorics, the above equation is usually written in the ``Lagrange formulation'' by supposing that $f(z) = \frac{z}{R(z)}$ with $ R(z) \in \mathbb{C}[[z]]$ with $R(0) \ne 0$ so that the equation becomes 
 \begin{eqnarray} \phi(z) = z \cdot R( \phi( z)).   \label{eq:lagrange} \end{eqnarray}

\begin{theorem}[Lagrange inversion formula]  \label{thm:lagrange}Let $ R \in \mathbb{C}[[z]]$ be a formal power series in $z$ such that $[z^{0}]R \ne 0$. Then there exists a unique formal power series $\phi$ satisfying \eqref{eq:lagrange} and we have for all $k\geq0$ and all $n \geq 1$
$$ [z^{n}] \big(\phi(z)\big)^{k} = \frac{k}{n}[z^{n-1}] \left(z^{k-1} R(z)^{n}\right),$$
where $[z^{n}] f(z)$ in the coefficient in front of $z^{n}$ in the formal power series $f \in \mathbb{C}[[z]]$.
\end{theorem}
\noindent \textbf{Proof.} The idea is to interpret combinatorially the weights in the formal expansion $z \cdot R( \phi( z))$, where $R(z) = \sum_{i \geq 0} r_{i}z^{i}$. Indeed, using \eqref{eq:lagrange}, it easy to prove by induction on $n \geq 1$ that the coefficient  in front of $z^{n}$ in $\phi$ can be interpreted as a sum over all plane trees with $n$ vertices where the weight of a tree $\tau$ is given by 
 \begin{eqnarray*} \mathrm{w}( \tau) &=& \prod_{u \in \tau} r_{ k_{u}(\tau)}.   \label{eq:weightwalk}\end{eqnarray*}
 This is true for $n=1$ and for $ n \geq 1$ using \eqref{eq:lagrange} writing $\phi_{n} = [z^{n}] \phi(z)$ and $r_{n} = [z^{n}] R(z)$ we find 
  \begin{eqnarray*} \phi_{n+1}&=& r_{1} [z^{n}] \phi(z) + r_{2} [z^{n}] \phi^{2}(z) + r_{3} [z^{n}] \phi^{3}(z) + \dots\\
  &=& \sum_{\ell \geq 1} r_{ \ell} \sum_{k_{1} + \dots + k_{\ell} = n} \prod_{i=1}^{\ell} \phi_{k_{i}}\\
&\underset{ \mathrm{Induc.}}{=}&   \sum_{\ell \geq 1} r_{ \ell} \sum_{k_{1} + \dots + k_{\ell} = n} \prod_{i=1}^{\ell} \left( \sum_{ \tau \mathrm{ \ plane\  tree\  size\  } k_{i} }  \mathrm{w}( \tau)\right)\\
&=& \sum_{ \tau \mbox{ plane tree size } n+1 }  \mathrm{w}( \tau),  \end{eqnarray*}
since the latter equality just comes from the decomposition of a plane tree of size $n+1$ at its root vertex.
 
 \begin{figure}[!h]
  \begin{center}
  \includegraphics[width=10cm]{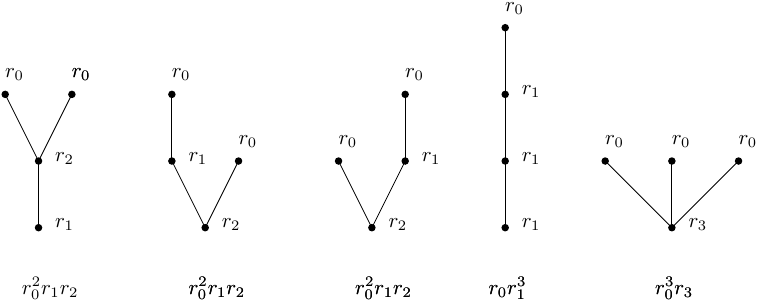}
  \caption{Interpretation of $[z^{4}]\phi(z)$ in diagrammatic form.}
  \end{center}
  \end{figure}
Similarly for $k \geq 1$, the coefficient of $z^{n}$ in $\phi^{k}$ is the total weight of  forests of $k$ trees having $n$ vertices in total. Now, using the {\L}ukasiewicz encoding, such  a forest can be encoded by a skip-free descending path $(S)$ with $n$ steps and reaching $-k$ for the first time at time $n$ where  the weight of such a path becomes $\mathrm{w}(S) = \prod_{i = 0}^{n-1} r_{S_{i+1}-S_{i} +1}$.  By Feller's combinatorial lemma, for a skip-free descending walk $(S)$ of length $n$ such that $S_{n}=-k$ there are exactly $k$ cyclic shifts so that $n$ is the $k$-th strict descending ladder time. So if we partition the set of all walks of length $n$ so that $S_{n}=-k$ using the cyclic shift as an equivalence relation, we know that in each equivalence class, the proportion of walks so that $T_{-k}=n$ is $ \frac{k}{n}$ (most of the classes actually have $n$ elements in it, but it could be the case that the subgroup of cyclic shifts fixing the walk is non-trivial and has order $\ell | k$, in which case there are $n/\ell$ elements in the orbit and $k/\ell$ are such that $T_{-k}=n$). Since the weight $ \mathrm{w}(\cdot)$  is constant over all equivalence classes we deduce that:

$$ \sum_{ \begin{subarray}{c}(S) \mathrm{\ walks \ of \ length\  }n \\ S_{n}=-k \end{subarray}} \mathrm{w}(S) = \frac{n}{k} \sum_{ \begin{subarray}{c}(S) \mathrm{\ walks \ of \  length \ }n \\ S_{n}=-k \mathrm{\  and \  }T_{-k}=n \end{subarray}}  \mathrm{w}(S).$$
 It remains to notice that 
$$ [z^{n-1}]  \left(z^{k-1} R(z)^{n}\right),$$
is exactly the weight of all paths $(S)$ of length $n$ such that $S_{n} = -k$. \qed
\medskip 

Here are two recreative (but surprising) applications of Lagrange inversion formula taken from the post ``What is Lagrange inversion formula good for?'' in \textit{Mathoverflow}:

\begin{exo} Let $F(x)$ be the be the unique power series such that for all $n \geq 0$ the coefficient of $x^n$ in $F^{n+1}(x)$ is equal to $1$. Show that $F(x) = \frac{x}{1-  \mathrm{e}^{-x}}$.\end{exo}

\begin{exo} \label{exo:lagrangebis} For $a \in (0,1/2)$ show that the positive solution $x=x(a)$ near $0$ of $x^5-x-a=0$ can be written as 
$$ x = - \sum_{k\geq 0} {5k \choose k} \frac{a^{4k+1}}{4k+1},$$
i.e.~we can ``solve" quintic equations (any quintic equation can be put into this form, see ``Bring radical" or ``Bring–Jerrard" on Wikipedia).
 \end{exo}

\section{Probabilistic counting of trees}
In this section we illustrate how to enumerate certain classes of trees using our knowledge on (random) walks. One underlying idea is  to design a random variable which is uniformly distributed on the set we wish to count.
\subsection{Prescribed degrees}
\begin{theorem}[Harary \& Prins \& Tutte (1964)]\label{thm:prescribed} \noindent The number of plane trees with $d_{i}$ vertices with $i \geq 0$ children, and with $n = 1+\sum i d_{i}= \sum d_{i}$ vertices is equal to 
 $$ \frac{(n-1)!}{d_{0}! d_{1}! \cdots  d_{i}! \cdots } = \frac{1}{n} {n \choose d_{0}, d_{1}, d_{2}, \dots }.$$
 \end{theorem}
\noindent \textbf{Proof.} Fix $d_{i}, k$ and $n$ as in the theorem. Notice that from \eqref{eq:enfantparent} we must have $n= 1 + \sum i d_{i} = \sum d_{i}$. By the encoding of plane trees into their {\L}ukaciewicz path it suffices to enumerate the number of paths starting from $0$, ending at $-1$ at $n$ and with $d_{i}$ steps of $i-1$ and which stay non-negative until time $n-1$. Clearly, if one removes the last assumption there are 
$$ {n \choose d_{0}, \dots, d_{i}, \dots} = \frac{n!}{ d_{0}! d_{1}!\cdots}$$
such paths. If we partition those paths according to the cyclic shift equivalence relation, then by Lemma \ref{lem:feller} (see also Remark \ref{rem:shift}) we know that each equivalence class has cardinal $n$ and has a unique element which stays non-negative until time $n-1$. Hence the quantity we wanted to enumerate is equal to 
$$ \frac{1}{n} {n \choose d_{0}, \dots, d_{i}, \dots} = {(n-1)!} \prod_{i} \frac{1}{d_{i}!}.$$ \qed

\begin{corollary}[Catalan's counting]\label{cor:catalan} For $n \in \{0,1,2,\dots\}$ we have $$ \# \big\{ \mathrm{plane\  trees \ with \ }n \mbox{ edges}\big\} = \# \big\{ \mathrm{plane\  trees \ with \ }n+1 \mbox{ vertices}\big\}= \frac{1}{n+1} {2n \choose n}.$$
\end{corollary}

\noindent \textbf{Proof.}  With the same notation as in the preceding theorem, the number of trees with $n \geq 1$ vertices is equal to 
$$  \sum_{\begin{subarray}{c}d_{0},d_{1}, d_{2}, \dots \\ 1+\sum i d_{i} =n = \sum d_{i} \end{subarray}} \frac{(n-1)!}{d_{0}!d_{1}! \cdots} =  \frac{1}{n} \sum_{\begin{subarray}{c}d_{0},d_{1}, d_{2}, \dots \\ 1+\sum i d_{i} =n = \sum d_{i} \end{subarray}} { n \choose d_{0}, d_{1}, \dots } =  \frac{1}{n} [z^{n-1}] \left( 1 + z+z^{2}+z^{3}+ \cdots \right)^{n}. $$
Using Lagrange inversion formula (Theorem \ref{thm:lagrange}) the last quantity can be expressed as $[z^{n}] \phi(z)$ where $\phi(z)$ is the formal power series solution to $\phi(z) = \frac{z}{1-\phi(z)}$ (i.e.~with $R(z) = \frac{1}{1-z}$). Solving explicitly we get $ \phi(z) = \frac{1}{2}(1- \sqrt{1-4z})$ and a coefficient extraction yields the desired formula. Alternatively, if we put $\phi(z) = z + z \psi(z)$, we find that $\psi$ satisfies the Lagrange equation $\psi(z) = z (1+ \psi(z))^2$ so that $\psi$ is amenable to an easy  Lagrange inversion:  we get that the number of plane trees with $n+1$ vertices is  $$ [z^{n+1}] \phi(z) = [z^n]\psi(z) = \frac{1}{n} [z^{n-1}]\left((1+z)^{2}\right)^n = \frac{1}{n} {2n \choose n-1}=\frac{1}{n+1} {2n \choose n}. $$  \qed 
\subsection{Uniform geometric BGW plane trees} \label{sec:uniform}
We denote by $ \mathbf{T}_{n}$ the set of all plane trees with $n$ edges and by $ \mathcal{T}_{n}$ a uniform plane tree taken in $ \mathbf{T}_{n}$. As we shall see $ \mathcal{T}_{n}$ can be interpreted as a conditioned version of a BGW tree:
\begin{proposition} Let $ \mathcal{T}$ be a Bienaym\'e--Galton--Watson tree with  geometric offspring distribution of parameter $1/2$, i.e.~$\mu_{k} = \left(\frac{1}{2}\right)^{k+1}$ for $k \geq 0$.  Then $ \mathcal{T}_{n}$ has the law of $ \mathcal{T}$ conditioned on having $n$ edges.
\end{proposition}
\noindent \textbf{Proof.} Let $\tau_{0}$ be a tree with $n$ edges. Then by Exercise \ref{exo:GWprod} we have 
$$ \mathbb{P}( \mathcal{T} = \tau_{0}) = \prod_{u \in \tau_{0}} 2^{-k_{u}(\tau)-1}.$$
However, from \eqref{eq:enfantparent} we have $\sum_{u \in \tau_{0}} k_{u}(\tau_{0}) = \#\tau_{0}-1 = n$ so that the last display is equal to $ \frac{1}{2}4^{-n}$. The point is that this probability does not depend on $\tau_{0}$ as long as it has $n$ edges. Hence, the conditional law  of $ \mathcal{T}$ on $ \mathbf{T}_{n}$ is the uniform law. \qed \medskip

Notice that the above proposition and its proof hold for any non trivial parameter of the geometric offspring distribution. However, we chose $1/2$ because in this case the offspring distribution is critical, i.e.\ it has mean $1$. We can give another proof of Corollary \ref{cor:catalan}:\medskip

\noindent \textbf{Proof of Corollary \ref{cor:catalan} (bis).} Combining the previous proposition with Proposition \ref{prop:luka} and Kemperman formula yields
$$ \mathbb{P}( \# \mathcal{T} =n+1) = \mathbb{P}( T_{-1} = n+1) \underset{ \mathrm{Prop.} \ref{prop:kemperman}}{=} \frac{1}{n+1} \mathbb{P}(S_{n+1}=-1), 
$$
where $(S)$ is the random walk whose increments are distributed as $ \mathbb{P}(S_{1}= k) = 2^{-k-2}$ for $k \geq -1$ or equivalently as $G-1$ where $G$ is the geometric offspring distribution of parameter $1/2$. Recall that $G$ is also the number of failures before the first success in a series of independent coin flips: this is the negative Binomial distribution with parameter $(1,1/2)$. Hence 
$ \mathbb{P}(S_{n+1}=-1) = \mathbb{P}( \mathrm{Binneg}(n+1,1/2) = n)$ where $ \mathrm{Binneg}(n,p)$ is the negative Binomial distribution with parameter $(n,p)$ --the discrete analog of the Gamma laws. This distribution is explicit and we have  $\mathbb{P}( \mathrm{Binneg}(n,p) = k) = {n+k-1 \choose n-1} p^{n} (1-p)^{k}$ which is our case reduces to 
$$ \frac{1}{n+1} \mathbb{P}(S_{n+1}=-1) = \frac{1}{n+1}\mathbb{P}( \mathrm{Binneg}(n+1,1/2) = n) =  \frac{1}{2} 4^{-n}\frac{1}{n+1} { 2n \choose n}.$$
By the previous proposition (and its proof) we have on the other hand 
$$ \mathbb{P}( \# \mathcal{T} =n+1) = \# \{ \mathrm{plane\  trees \ with \ }n+1 \mbox{ vertices}\} \cdot \frac{1}{2}4^{-n}.$$
The result follows by comparing the previous two displays. \qed \medskip 

\begin{exo}[Enumeration of plane forests] \label{exo:forest} Extend the above proof to show that the number of forests of $f \geq 1$ trees (i.e.~ordered sequence of $f$ trees) whose total number of edges is $n$ is equal to 
$$  \frac{f}{2n+f} {2n +f \choose n}.$$
 Give another proof of the last display using Lagrange inversion formula (Theorem \ref{thm:lagrange}).
\end{exo}

The above exercise is useful to show that the typical height of $ \mathcal{T}_{n}$ converge in law towards the Rayleigh\footnote{\raisebox{-3mm}{\includegraphics[width=1cm]{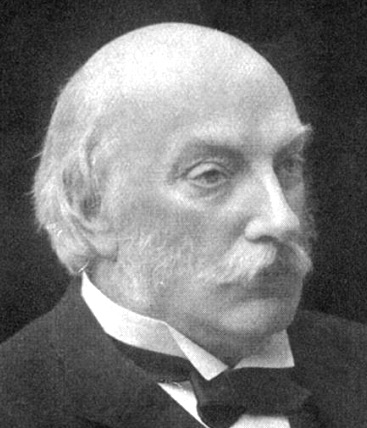}} John William Strutt, 3rd Baron Rayleigh (1842–1919), English} distribution $ \mathcal{R}$ which is the law of the norm of a two-dimensional normal vector:
 \begin{eqnarray} \label{def:rayleigh} \mathcal{R} \sim r \exp(-r^{2}) \mathbf{1}_{r>0} \mathrm{d}r.  \end{eqnarray}

\begin{corollary}[Typical height of uniform plane trees] \label{cor:heightplane} Let $ \mathcal{T}_{n}$ be a uniform plane tree with $n$ edges. Conditionally on $ \mathcal{T}_{n}$,  let $\delta_{n}$ be a uniformly chosen vertex of $ \mathcal{T}_{n}$ and denote  its height by $H_{n}$. Then we have 
$$ \mathbb{P}(H_{n} =h ) = \frac{2h+1}{2n+1} \frac{ {2n+1 \choose n-h}}{{2n \choose n}}.$$
In particular, we have the following convergence in distribution towards a scaled Rayleigh distribution  $$\frac{H_{n}}{ \sqrt{n}} \xrightarrow[n\to\infty]{(d)}  \frac{ \mathcal{R}}{ \sqrt{2}}.$$
\end{corollary}
\noindent \textbf{Proof.} We compute exactly the probability that the point $\delta_{n}$ is located at height $h\geq 0$. \\ 

\begin{minipage}{4cm}
\includegraphics[width=2.5cm]{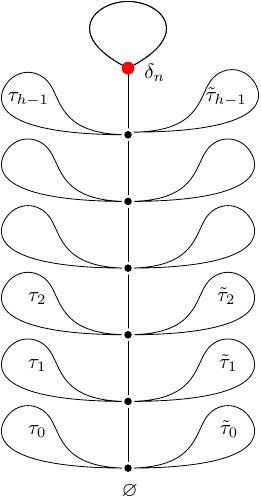}
\end{minipage}
\begin{minipage}{11.5cm} If so, the tree $ \mathcal{T}_{n}$ is obtained from the line joining $\varnothing$ to $\delta_{n}$ by grafting $h$ plane trees on its left, $h$ plane trees on its right and one on $\delta_{n}$, see the figure on the left. Obviously, the total number of edges of these trees must be equal to $n-h$. Using Exercise \ref{exo:forest} we deduce that 

$$ \mathbb{P}(H_{n}= h) = \frac{ \frac{2h+1}{2n+1} { 2n +1 \choose n-h}}{ {2n \choose n}}. $$
The second item of the theorem follows after applying Stirling formula and using Exercise \ref{exo:scheffe}. \qed
\end{minipage} 

\medskip

 \begin{exo} For any $p \geq 2$ a $p$-tree is a plane tree such that the number of children of each vertex is either $0$ or $p$. When $p=2$ we speak of binary trees. In particular, the number of edges of a $p$-tree must be a multiple of $p$. Show that for any $k \geq 1$ we have 
 $$ \# \{ p- \mbox{trees with }kp  \mbox{ edges }\} = \frac{1}{kp+1} {k p+1 \choose k},$$
 in three ways: using a direct application of Theorem \ref{thm:prescribed}, using a probabilistic approach via a certain class of random BGW trees, or via Lagrange inversion's formula Theorem \ref{thm:lagrange}.

 \end{exo}

\subsection{Cayley and Poisson BGW trees}
\begin{figure}[h!]
\begin{center}
\includegraphics[height=3cm]{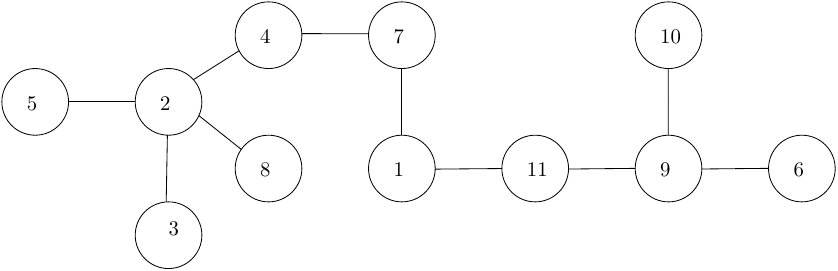}
\caption{ \label{fig:cayleytree}	A Cayley tree over $\{1,2,3,4,\dots, 11\}$.}
\end{center}
\end{figure}

In this section we focus on a different type of tree first studied by Cayley:\footnote{\raisebox{-5mm}{\includegraphics[width=1cm]{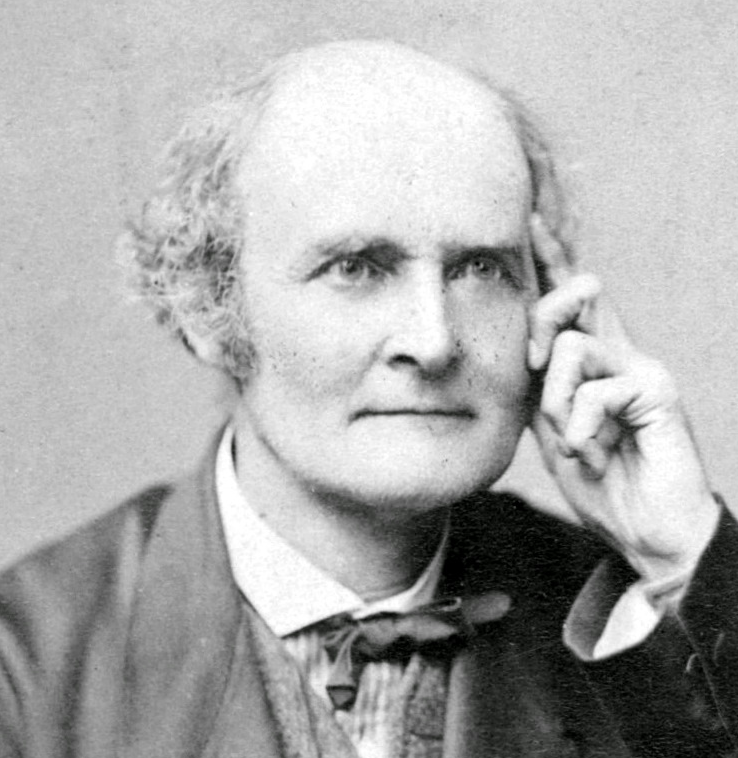}} Arthur Cayley  (1821--1895) receiving a phone call, English}
\begin{definition}
A \textbf{Cayley tree} of size $n$  is a tree over the  $n$ vertices $\{1,2, \dots , n\}$ without any orientation nor distinguished point. In other words, it is a spanning tree on $\mathbb{K}_{n}$, the complete graph over $n$ vertices. See Figure~\ref{fig:cayleytree}. 
\end{definition}
 Let $ \mathcal{T}$ be a  BGW (plane) tree with Poisson offspring distribution of parameter $1$ (in particular, the mean number of children is $1$ and we are in the critical case). As in the previous subsection (but with vertices instead of edges) we denote by $ \mathcal{T}_{n}$ the random tree $ \mathcal{T}$ conditioned on having $n$ vertices. 
 
 \begin{proposition} \label{prop:cayleyGW}
 Consider $ \mathcal{T}_{n}$ and assign the labels $\{1, \dots ,n\}$ uniformly at random to the vertices of $ \mathcal{T}_{n}$. After forgetting the plane ordering $ \mathcal{T}_{n}$ this produces a Cayley tree which we denote by $ \mathscr{T}_{n}$. Then $ \mathscr{T}_{n}$ is uniformly distributed over all Cayley trees with size $n$. 
 \end{proposition}
 \noindent \textbf{Proof.} Let us first compute the probability that $ \mathcal{T}$ has $n$ vertices. Using the {\L}ukasiewicz walk and the cyclic lemma we get that 
 $ \mathbb{P}(\# \mathcal{T} =n) = \frac{1}{n} \mathbb{P}(S_{n}=-1),$ where $S$ is the random walk whose increments are centered and distributed according to $ \mathrm{Poisson}(1)-1$ i.i.d.~variables. Recalling Section \ref{sec:poissonRW}, it follows that
 $$ \mathbb{P}(\# \mathcal{T} =n) = \frac{1}{n} \mathbb{P}( \mathrm{Poisson}(n) = n-1)  = \mathrm{e}^{-n} \frac{n^{n-2}}{(n-1)!}.$$
Fix a Cayley tree $ \mathfrak{t}$ and let us study the possible ways to obtain $ \mathfrak{t}$ by the above process. We first choose the root of the tree among the $n$ possibles vertices and obtain a rooted Cayley tree $ \mathfrak{t}^{\bullet}$. Once the origin is distinguished, there are $\prod_{ u \in \mathfrak{t}^{\bullet}} k_{u}( \mathfrak{t}^{\bullet}) !$ possible ways to give a planar orientation to the tree, where $k_{u} ( \mathfrak{t}^{\bullet})$ is the number of children of the vertex $u$ in $ \mathfrak{t}^{\bullet}$ (for this we only need the ancestor vertex, not the planar ordering). After these operations, each of the labeled, rooted, plane trees $(\tau, \ell)$ obtained appears with a probability (under the Poisson(1)-BGW measure) equal to 
$$ \frac{1}{n!} \mathrm{e}^{-n} \prod_{ u \in \tau} \frac{1}{k_{u}(\tau)!} = \frac{1}{n!} \mathrm{e}^{-n} \prod_{ u \in \tau} \frac{1}{k_{u}( \mathfrak{t}^{\bullet})!}.$$
Performing the summation, the symmetry factors involving the $k_{u}!$ conveniently disappear and we get
$$ \mathbb{P}( \mathcal{T}_{n} \to \mathfrak{t}) = n \times \frac{\mathrm{e}^{-n}}{n!} \left(\mathrm{e}^{-n} \frac{n^{n-2}}{(n-1)!}\right)^{-1} = \frac{1}{n^{n-2}}.$$
Since the result of the last display does not depend on the shape of $ \mathfrak{t}$, the induced law is indeed uniform over all Cayley trees and we have even proved:

 \begin{corollary}[Cayley's formula]\label{cayley}The  number of Cayley trees of size $n$ is $ n^{n-2}$.\end{corollary}
As a short application of the above corollary, we propose:
\begin{exo}[Pick a tree - any tree, \cite{chin2015pick}] \label{exo:pickatree} Let $T_{n}$ be a random labeled subtree (no planar ordering) of the complete graph $ \mathbb{K}_{n}$ over the vertices $\{1,2, \dots , n\}$. Show that
$$ \lim_{n \to \infty}\mathbb{P}(T_{n} \mathrm{ \ spans\  all\  vertices \ of \ } \mathbb{K}_{n}) =  \mathrm{e}^{- \mathrm{e}^{{-1}}}.$$
\end{exo}

\begin{exo}[Lagrange meets Cayley] \label{exo:lagrangecayley} Let $T(z)$ be the (formal) exponential generating series of Cayley trees with a distinguished vertex, i.e.
$$ T(z) = \sum_{k \geq 1}  n \cdot  \# \{ \mathrm{Cayley \ trees \ size\ }n\} \cdot \frac{z^n}{n!}.$$
Show using a recursive decomposition at the root that $T(z) = z \mathrm{e}^{T(z)}$. Apply Lagrange inversion formula (Theorem \ref{thm:lagrange}) to recover Corollary \ref{cayley}.
\end{exo}

We have the following generalization similar to Exercise \ref{exo:forest}:
\begin{exo}[Cayley forests] \label{exo:cayleyforest} Show that the number of (non-plane) forests on $\{1,2, \dots , n\}$ with $k$ trees with roots $1,2, \dots, k$ is counted by 
$$ \mathcal{F}(k,n) =  \frac{k}{n} n^{n-k}. $$	
\end{exo}

The previous exercise can be used to prove the same Rayleigh limit (recall \eqref{def:rayleigh} and Corollary \ref{cor:heightplane}) for the typical height in a large uniform Cayley tree:
\begin{corollary}[Typical height of uniform Cayley trees]  \label{cor:heightcayley} Let $ \mathscr{T}_{n}$ be a uniform Cayley tree of size $n$. Conditionally on $ \mathscr{T}_{n}$, let $\delta_{n}$ be a uniform vertex of $\{1,2, \dots , n\}$. Then the distance $D_{n}$ in $ \mathscr{T}_{n}$ between the vertices $1$ and $\delta_{n}$ has the following distribution 
$$ \mathbb{P}(  D_{n} = k-1) = \left( 1- \frac{1}{n}\right)\left( 1- \frac{2}{n}\right) \cdots \left( 1- \frac{k-1}{n}\right) \frac{k}{n}, \quad \mbox{ for }1 \leq k \leq n.$$
In particular we have 
$$\frac{ D_{n}}{ \sqrt{n}} \xrightarrow[n\to\infty]{(d)}  \mathcal{R}. $$
\end{corollary}
\noindent \textbf{Proof.} By symmetry, $D_{n}$ has the same law as the distance between two uniform vertices $U_{n},V_{n}$ of $ \mathscr{T}_{n}$ (possibly confounded). For $k=1$, the probability that $D_n=0$ is the probability that $U_n = V_n$ which is indeed $1/n$. Otherwise, for $k \geq 2$, the event $\{D_n =k-1\}$ happens if $ \mathscr{T}_n$ is obtained from an ordered line of $k$ vertices on which we graft a forest of $k$ Cayley trees with prescribed roots, and so that the selected vertices are on endpoints of this line. Diving by the obvious symmetry factors, the previous exercise shows that this probability is given by 
$$   \frac{1}{2} \times \frac{2}{n^2} \cdot n(n-1)\dots (n-(k-1)) \cdot \frac{\mathcal{F}(k,n)}{ n^{n-2}} =  \frac{n(n-1)\dots (n-(k-1))}{n^k} \frac{k}{n},$$ as desired. We recognize the law of the \textbf{first collision in the birthday paradox} on a year with $n$ days. In particular, for $k_n = [x \sqrt{n}]$ with $ x >0$ we have 
$$ \mathbb{P}(D_n \geq k_n-1)  = \prod_{i=1}^{k_n} (1- \frac{i}{n}) \sim \exp \left( -  \frac{1}{n} \sum_{i=1}^{k_n} i\right) \xrightarrow[n\to\infty]{} \mathrm{e}^{-x^2/2},$$ entailing the convergence to the Rayleigh distribution.
\qed

 \begin{exo}[Random mapping]  \label{exo:mapping} Let $M_{n} = \{1,2, \dots, n \} \to \{1, 2, \dots , n\}$ be a mapping  chosen uniformly at random among the  $n^{n}$ possibilities. We represent $M_{n}$ as an oriented graph where an arrow goes from $i$  to $M_n(i)$, see Fig.~\ref{fig:mapping}.  We denote by $ \mathcal{C}_{n} \subset \{1,2, \dots , n \}$ the cyclic points i.e.~the integers $i$ such that there exists  $m \geq 1$ with $(M_{n})^{\circ m}(i) = i$.
\begin{figure}[h!]
 \begin{center}
 \includegraphics[width=11cm]{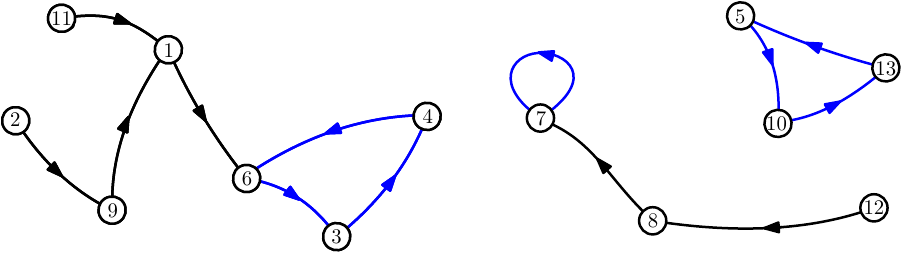}
 \caption{Illustration of the graph of the mapping $1 \to 6, 2\to 9, 3 \to 4, 4 \to 6, 5 \to 10, 6 \to 3, 7 \to 7, 8 \to 7, 9 \to 1, 10 \to 13, 11 \to 1, 12 \to 8, 13 \to 5$. \label{fig:mapping}}
 \end{center}
 \end{figure}
\begin{enumerate}
 \item Prove that $ \# \mathcal{C}_{n}$ has the same law as $D_n$ in Corollary \ref{cor:heightcayley}.
\item  Show that $$ \mathbb{P}(\mbox{the  (unoriented) graph of }M_{n} \mbox{ is connected}) =  \frac{1}{n^{n-2}}\sum_{k=1}^{n} {n \choose k} (k-1)! \cdot \frac{k}{n} n^{n-k}$$ and give its asymptotic when $n \to \infty$.
 \end{enumerate}
 \end{exo}

 \subsection{Contour function} \label{sec:contourfunction}
 We finish this section by mentioning another more geometrical encoding of plane trees which is probabilistically less convenient in the general BGW case but very useful in the case of geometric BGW trees.\medskip

\begin{figure}[!h]
 \begin{center}
 \includegraphics[width=15cm]{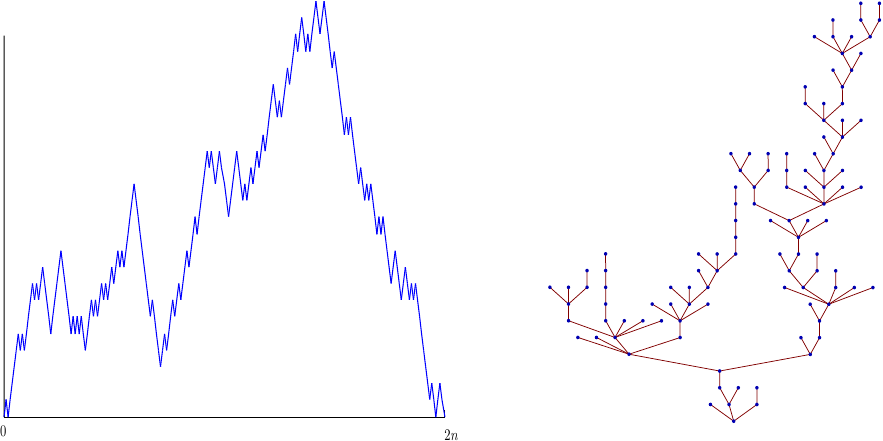}
 \caption{The contour function associated with a plane tree. \label{fig:contour}}
 \end{center}
 \end{figure}

Let $\tau$ be a finite plane tree. The contour function $ \mathcal{C}_{\tau}$ associated with $\tau$ is heuristically obtained by recording the height of a particle that climbs the tree and makes its contour at unit speed. More formally, to define it properly one needs the definition of a corner: We view
$\tau$ as embedded in the plane, then  a \textbf{corner} of a vertex in $\tau$
is an angular sector formed by two consecutive edges in clockwise
order around this vertex. Note that a vertex of degree $k$ in $\tau$
has exactly $k$ corners. If $c$ is a corner of $\tau$, $\mathrm{Ver}(c)$
denotes the vertex incident to $c$, see Figure \ref{fig:contour}.

The corners are ordered clockwise cyclically around the tree in the so-called {\em contour order}. If $\tau$ has $n \geq 2$ vertices we index the corners by letting $(c_{0},c_{1},c_{2}, \ldots,c_{2n-3})$ be the sequence of
corners visited during the contour process of $\tau$, starting from
the corner $c_0$ incident to $\varnothing$ that is located to the left
of the oriented edge going from $\varnothing$ to $1$ in $\tau$.
\begin{definition} Let $\tau$ be a finite plane tree with $n \geq 2$ vertices and let $(c_{0},c_{1},c_{2}, \ldots,c_{2n-3})$ be the sequence of corners visited during the contour process of $\tau$. We put $c_{2n-2}=c_{0}$ for notational convenience. The contour function of $\tau$ is the walk defined by 
$$ \mathcal{C}_{\tau}(i) = \# \mathrm{Ver}(c_{i}), \quad \mbox{ for }0 \leq i \leq 2n-2.$$
\end{definition}

Clearly, the contour function of a finite plane tree is a finite non-negative walk  of length $2 (\#\tau-1)$ which only makes $\pm 1$ jumps. Here as well, the encoding of a tree into its contour function is invertible:
\begin{exo} Show that taking the contour function creates a bijection between the set of all finite plane trees and the set of all non-negative finite walks with $\pm 1$ steps which start and end at $0$.
\end{exo}

Now, we give a probabilistic description of the law of the contour function of $ \mathcal{T}$ when $ \mathcal{T}$ is distributed as a  geometric(1/2)-BGW tree (i.e.~has the same law as in Section \ref{sec:uniform}).
\begin{proposition}[Contour function of Catalan trees] \label{prop:contourgeo} Let $ \mathcal{T}$ as above. Then its contour function $ \mathcal{C}_{ \mathcal{T}}$ has the same law as 
$$ (S_{0}, S_{1}, \dots , S_{T_{-1}}),$$
where $(S)$ is a simple symmetric random walk and $T_{-1}$ is the first hitting time of $-1$.
\end{proposition}
\noindent \textbf{Proof.} Notice first that $ \mathcal{T}$ is almost surely finite by Theorem \ref{thm:extinction} and so all the objects considered above are well defined. Let $\tau_{0}$ be a plane tree with $n$ edges. We have seen in the previous proposition that $ \mathbb{P}( \mathcal{T} = \tau_{0}) = \frac{1}{2} 4^{-n}$. On the other hand, the contour function of $\tau_{0}$ has length $2n$ and the probability that the first $2n$ steps of $(S)$ coincide with this function and that $T_{-1}= 2n+1$ is equal to 
$ 2^{-2n} \cdot \frac{1}{2} = \frac{1}{2}4^{-n}$. This concludes the proof.\qed \medskip 

\begin{exo} \label{exo:contourcat}Give a new proof of Corollary \ref{cor:catalan} using the contour function.
\end{exo}
\begin{exo} \label{exo:height} Let $ \mathcal{T}$ be a BGW tree with geometric(1/2) offspring distribution. The height of $ \mathcal{T}$ is the maximal length of one of its vertices. Prove that 
$$ \mathbb{P}(  \mathrm{Height}(\mathcal{T}) \geq n) = \frac{1}{n+1}.$$
\end{exo}

 \section{The Brownian continuum random tree}
 
The reader might be puzzled by the appearance of the Rayleigh distribution as the typical height in both uniform plane trees (Corollary \ref{cor:heightplane}) and uniform Cayley trees (Corollary \ref{cor:heightcayley}) of large size. This is only the tip of a much larger iceberg: many classes of random trees converge in the scaling limit towards a universal Continuum Random Tree (CRT) called the Brownian CRT. We briefly describe this fascinating object. We first describe a way to control globally the geometry of a random graph. 
\subsection{Gromov--Hausdorff topology}
The idea is to see a finite graph once endowed with its graph distance as a finite metric space, i.e.~a point in the space $$ \mathbb{K} = \left\{ \begin{array}{c} \mbox{isometry classes of compact metric spaces} \end{array}\right\},$$ since from the geometric point of view, it is impossible to distinguish two isometric metric spaces (in particular, in the following when we  speak of a metric space, the reader should think of its isometry class).  One might think that this set is monstrous and that its very definition could pose a problem. In reality, thanks to the compactness condition imposed on its points (i.e.~on the isometry classes of metric spaces), the space $ \mathbb{K}$ is quite ``small''; for example, any compact metric space can be seen as a closed subset of  $\ell^{\infty}( \mathbb{R})$.

 We will now equip $ \mathbb{K}$ with a distance, known as the Gromov--Hausdorff distance and denoted $ \mathrm{d_{GH}}$. Let $(E, \mathrm{d_E})$ and $(F, \mathrm{d_F})$ be two points of $\mathbb{K}$, i.e.~two (isometry classes of) compact metric spaces, then the Gromov--Hausdorff distance between $E$ and $F$ is 
 $$ \mathrm{d_{GH}}((E, \mathrm{d_E}),(F, \mathrm{d_F})) =  \inf\{ \mathrm{d}_{ \mathrm{Haus},G}( E',F')\}$$
where $\mathrm{d}_{ \mathrm{Haus},G}(E',F')$ is the Hausdorff distance between $E' \subset G$ and $F' \subset G$ two compacts of the same ambient space $G$ that are respectively isometric to $E$ and $F$.

\begin{figure}[!h]
 \begin{center}
 \includegraphics[width=12cm]{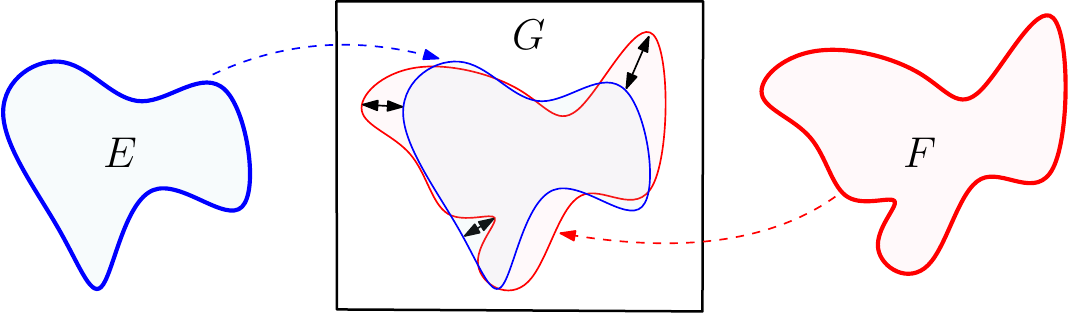}
 \caption{Illustration of the Gromov--Hausdorff distance: to compare two metric spaces, first embed them in a common metric space and use the Hausdorff distance.}
 \end{center}
 \end{figure}

\begin{theorem}
The space $( \mathbb{K},d)$ is a Polish metric space (i.e.~separable and complete).
\end{theorem}

We refer the reader to \cite[Chapter 7]{BBI01} for details concerning this space.  This formalism is very convenient and allows us to define the Brownian continuous tree as the ``scaling limit'' of renormalized random discrete trees. Indeed, if $\mu$ is a critical aperiodic offspring distribution with finite variance $\sigma^2 \in (0, \infty)$, one can consider $\mathcal{T}_{n}$ a $\mu$-BGW tree to have $n$ edges endowed with its graph distance as a random metric space. We have the following invariance principle: 

\begin{theorem}[Reformulation of Aldous by Le Gall ]
We have the following  convergence in distribution for the Gromov--Hausdorff topology
$$ \left( \mathcal{T}_{n} , \frac{1}{ \sqrt{n}} \mathrm{d_{gr}} \right) \quad \mathop{\longrightarrow}^{(d)}_{n \rightarrow \infty} \quad \left(  \mathfrak{T}, \frac{2}{ \sigma}\mathrm{d}\right),$$
where $( \mathfrak{T}, \mathrm{d})$ is a random compact continuous tree, called the Brownian continuum random tree, whose distribution does not depend on $\mu$. \label{thm:aldousCRT}
\end{theorem}
See Figure \ref{fig:CRT} for (an approximation of) a sampling of $ \mathfrak{T}$. The Brownian continuum random tree $ \mathfrak{T}$, frequently called CRT (for ``continuum random tree'') in the literature, is therefore a random metric space (for example, its diameter is random) but it has ``almost certain'' properties, i.e.~true with probability $1$:
\begin{itemize}
\item $ \mathfrak{T}$ is a.s.~a continuous tree, i.e.~a compact metric space, geodesic (in which any two points are connected by a single geodesic) and cycle-free.
\item for any $x \in \mathfrak{T}$, the space $ \mathfrak{T}\backslash \{x\}$ has at most $3$ connected components.
\item the fractal dimension of $ \mathfrak{T}$ is equal to $2$.
\end{itemize}

\subsection{Brownian excursion as continuous contour function}
\begin{figure}[!h]
 \begin{center}
 \includegraphics[width=6cm]{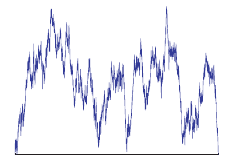}
 \caption{\label{fig:broexc} A simulation of a Brownian excursion.}
 \end{center}
 \end{figure}

At first glance, there is not much Brownian about the definition of $ \mathfrak{T}$. To understand where the name comes from, let us take a look at the contour function $C(\mathcal{T}_{n})=(C_{s}(\mathcal{T}_{n}))_{0 \leq s \leq 2n}$ of the conditioned BGW trees. As a proxy in the proof of the previous theorem, one usually shows the following convergence:
$$ \left( \frac{C_{ns}(\mathcal{T}_{n})}{ \sqrt{n}} : 0 \leq s\leq 1\right) \quad \mathop{\longrightarrow}^{(d)}_{n \rightarrow \infty} \quad \left( \frac{2}{\sigma}\mathbf{e}_{t} : 0 \leq t \leq 1\right)$$ for the uniform topology on $( \mathcal{C}([0,1]), \| \cdot \|_\infty)$ and where $\mathbf{e}$ is a random continuous function, called the \emph{Brownian excursion}, and which can informally be seen as a Brownian motion that starts from $0$ at time $0$, remains positive over the time interval $(0,1)$ and returns to $0$ at time $1$ (see Figure \ref{fig:broexc} for a simulation). 

The reason why Brownian motion appears is that, although in general the contour function is not a random walk (except in the case of Catalan trees, i.e.\ when the reproduction law is geometric, see Proposition  \ref{prop:contourgeo}), it can nevertheless be approximated by a random walk, so that the above convergence is an application (in a conditional setting) of Donsker's theorem, according to which suitably renormalized random walks converge to Brownian motion (this is the functional extension of the central limit theorem). In particular, the renormalization factor $\sqrt{n}$ is the same as in the central limit theorem, thanks to the finite variance assumption.
 It is then natural to expect the Brownian excursion to encode, in some sense, the Brownian continuous tree. This intuition was formalized by Duquesne \& Le Gall \cite{DLG05}, who mimicked the construction of a tree from its contour function in the discrete setting. More precisely, to any continuous function $f : [0,1] \to \mathbb{R}_{+}$ such that $f(0)=f(1)=0$, we associate a pseudo-distance on $[0,1]$, denoted $ \mathrm{d_{f}}$, and defined by 
$$\mathrm{d_{f}}(s,t) = f(s)+f(t) - 2 \min_{u \in [\min(s,t), \max(s,t)]} f(u).$$
It is easy to check that $ \mathrm{d_{f}}$ is a pseudo-distance and that the points $s,t \in [0,1]$ with zero distance are those that face each other under the graph of $f$. We can then consider the equivalence relation on $[0,1]$ obtained by putting $ s \sim t$ if $ \mathrm{d_{f}}(s,t)=0$. On the quotient space $[0,1]/\sim$ the (projection of) pseudo-distance $ \mathrm{d_{f}}$ is now a distance and $([0,1]/\sim, \mathrm{d_{f}})$ is a compact metric space, denoted $ \mathfrak{T}_{f}$, which is a continuous tree. When the previous construction is performed starting from the Brownian excursion, the random tree $ \mathfrak{T}_{ \mathbf{e}}$ is the continuous Brownian tree $( \mathfrak{T}, \mathrm{d})$ that appears in  Theorem \ref{thm:aldousCRT}. \bigskip


\noindent \textbf{Bibliographical notes.} The material about Bienaym\'e--Galton--Watson tree is rather classical. The coding of trees and the formalism for plane trees (the so-called Neveu's notation \cite{Nev86}) can be found in \cite{LG05}. The lecture notes of Igor Kortchemski \cite{Kor16} are a very good introduction accessible to the first years of undergraduate studies in math. The interested reader can also consult \cite{AD15}. Beware some authors prefer to take the lexicographical order rather than the breadth first order to define the {\L}ukasiewicz walk (in the finite case this causes no problem but this is not a bijection if the trees can be infinite).  The two exercices illustrating Lagrange inversion formula are taken from the MathOverFlow post ``What is Lagrange inversion good for?''. Exercise \ref{exo:pickatree} is taken from \cite{chin2015pick}.  Theorem \ref{thm:prescribed} is proved in \cite{harary1964number}. The idea of Gromov--Hausdorff topology was first discovered in theoretical physics by Edwards \cite{edwards1975structure} and later popularized by Gromov \cite{Gro07} in geometry. It was brought to the probability community mainly by Evans \cite{Eva08} and Le Gall \cite{LG06}. We refer to \cite{BBI01} for background. The theory of scaling limits of random trees is by now one of the pillar of random geometry. The pioneer papers of Aldous \cite{Ald91a,Ald91,Ald91c} are still the best references for background on the Brownian Continuum Random Tree. We refer to  \cite{LG05} for a nice introductory course and to \cite{LGMi10} for its applications in the theory of random planar maps.
\bigskip

\noindent \textbf{Hints for Exercises:}\\
Exercise \ref{exo:dekking} is taken from \cite{Dek91b}.\\
Exercise \ref{exo:lagrangebis}: Put $\tilde{x} = ax-1$ to recover a Lagrangian formulation.\\
Exercise \ref{exo:forest}: After concatenating their {\L}ukasiewicz paths, such forests are coded by a  skip-free walk of $n+f$ steps starting at $0$ and reaching $-f$ for the first time at $n$.\\
Exercise \ref{exo:lagrangecayley}: If $C_n$ is the number of Cayley trees on $[[1,n]]$ with a distinguished vertex, prove that  for $n \geq 2$ we have
$$ C_n = n \cdot \sum_{k \geq 1}  \frac{1}{k!}\sum_{\ell_1+\dots+ \ell_k = n-1} {n-1 \choose  \ell_1, \dots , \ell_k} C_{\ell_1} C_{\ell_2}\dots C_{\ell_k}.$$
Exercise \ref{exo:mapping}: Once the cyclic points have been chosen, the rest of the graph is obtained by grafting Cayley trees, then use Exercise \ref{exo:forest}.\\
Exercise \ref{exo:contourcat}: Using the contour, the number of plane trees with $n$ edges is also the number of $\pm 1$ paths going from $0$ to $0$ in $2n$ steps while staying non-negative.\\
Exercise \ref{exo:height}: Using the contour function, the probability that the height is larger than $n$ is the probability that a simple random walk started at $0$ reaches $n$ before $-1$.\\

\part[Erd{\H{o}}s-R\'enyi random graph]{Erd\"os-R\'enyi random graph 
                             \\ \\ 
\label{part:ER}
  \begin{center}
                     \begin{minipage}[l]{15cm}
       \normalsize 
       In this part we study the famous model of random graph due to Erd{\H{o}}s and R\'enyi: 
\begin{definition}[$ G(n,p)$ model] \label{def:erdosrenyi}The Erd{\H{o}}s--R\'enyi random graph $ G(n,p)$ with parameters $n~\geq1$ and $p \in [0,1]$ is the (distribution of a) random graph whose vertex set is $\{1,2, \dots , n\}$ and where for each pair $ i \ne j$ the edge $ i \leftrightarrow j$ is present with probability $p$ independently of all the other pairs.
\end{definition}
This is the most natural random graph model since conditionally on its number of edges $m$, the variable $ G(n, p)$ is uniformly distributed over the set $ \mathbb{G}_{n,m}$ of all labeled simple graphs on $\{1,2, \dots, n\}$ with $m$ edges. For convenience we shall use $ \mathbb{G}_n = \bigcup_{m\geq 1} \mathbb{G}_{n,m}$ for the set of all simple graphs on the vertex set $\{1,2, \dots , n\}$. For a fixed $n \geq 1$, we shall consider all Erd{\H{o}}s--R\'enyi graphs $(G(n,p ) : p \in [0,1])$ as coupled as in Section \ref{sec:percoabstrait}: for each $n$ we sample i.i.d.~uniform variables $(U_{\{i,j\}} :  1 \leq i \ne j \leq n)$ and declare that $\{i,j\}$ is present in $G(n,p)$ if $U_{\{i,j\}} \leq p$. 
                     \end{minipage}
                  \end{center}
                  \vspace{0cm}
                   \begin{figure}[!h]
 \begin{center}
 \includegraphics[width=0.68\linewidth]{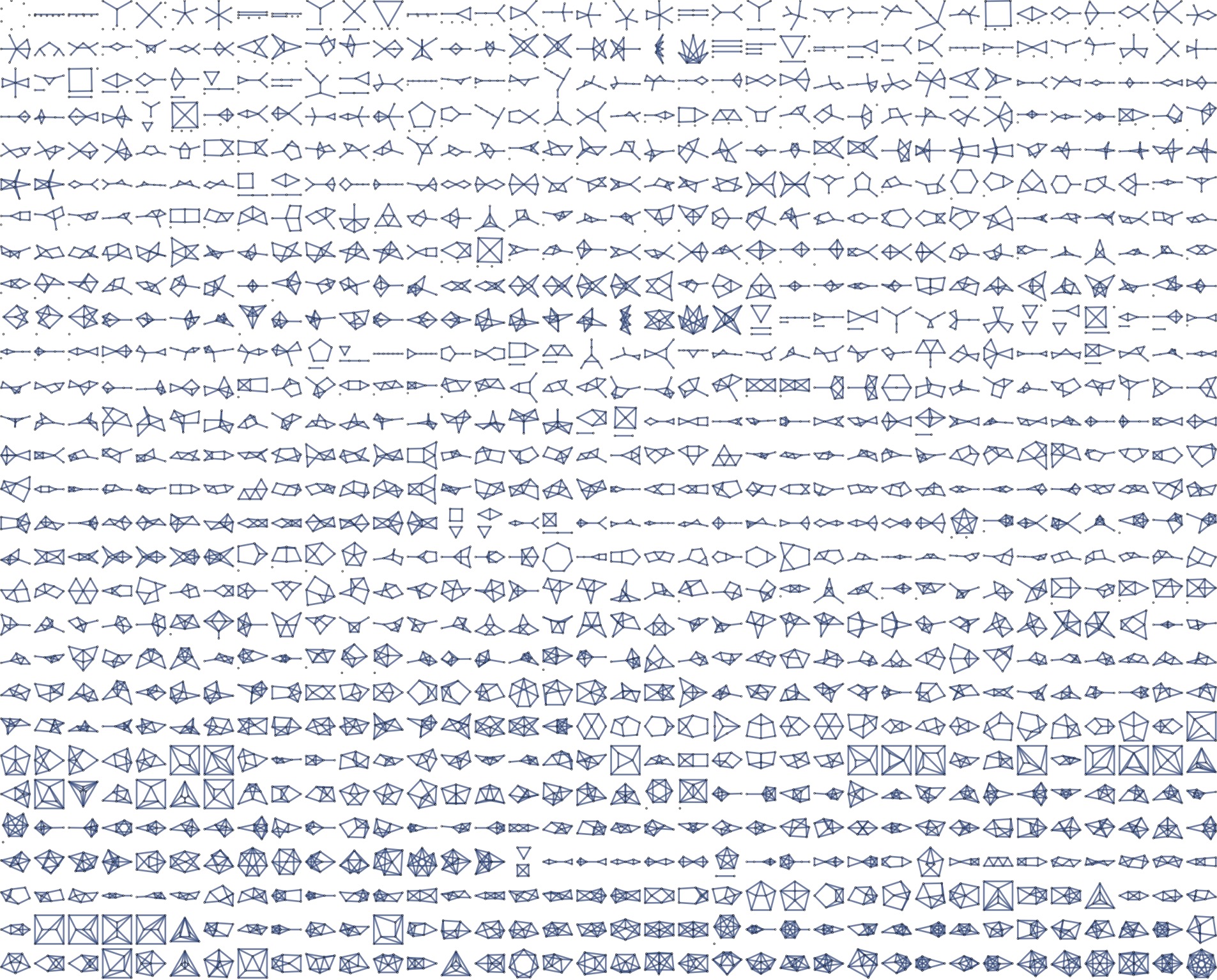}
 \caption{A list of the 1044 simple graphs on 7 vertices up to isomorphism.} 
 \end{center}
 \end{figure}
                }

\chapter{Local properties}
\hfill A tribute to the first and second moment method. \hfill 

\bigskip

In this chapter, we study ``local properties'' of $G(n,p)$ focusing mostly on the presence of certain subgraph in $G(n,p)$ as $p$ varies with $n$. We shall see that the presence of some subgraph sometimes satisfies a \textbf{phase transition} and governs interesting global properties of the graph such as the connectedness, or the spectral measure. Many proofs are based on the first and second method together with combinatorics and basic analysis.
\bigskip 

 In the following, a  \textbf{graph property} $A_{n}$ is just a subset of all simple graphs on $\{1,2, \dots , n \}$. We say that $A_{n}$ is \textbf{increasing} (resp. \textbf{decreasing}) if for any $ \mathfrak{g}, \mathfrak{g}' \in  \mathbb{G}_{n}$ satisfying $ \mathfrak{g} \subgraph \mathfrak{g}'$, we have $ \mathfrak{g} \in A_n \Rightarrow \mathfrak{g}' \in A_n$ (resp. with $ \mathfrak{g}' \subgraph \mathfrak{g}$). In words, adding (resp. removing) edges only help satisfying $A_{n}$. 
 
 \begin{example}[Appearance of subgraph or graph induced] Fix a graph $ \mathfrak{g}_{0}$ then the following graph properties
 $\{ \mathfrak{g} \in \mathbb{G}_n:\exists \mathfrak{g}'  \subgraph \mathfrak{g} \mbox{ with } \mathfrak{g'} \simeq \mathfrak{g}_{0}\}$ is an increasing graph property, whereas 
  $\{ \mathrm{Cayley\  trees \ of \  } \mathbb{G}_n\}$ is not an increasing graph property as soon as $n \geq 2$.
  \end{example}
  \begin{exo}[Erd\"os--R\'enyi is Cayley] For which $p\equiv p_{n}$ is $ \mathbb{P}(G( n, p) \mbox{ is a Cayley tree})$ maximal?  \label{exo:ercayley}\end{exo}
 If $(A_{n} : n \geq 1)$ is a sequence of graph properties and if the edge density $p\equiv p_{n}$ may depend on $n$, we say that $A_{n}$ holds for $G(n, p_{n})$ with \textbf{with high probability} (abbreviated by w.h.p) if
 $$\mathbb{P}( G(n,p_{n}) \in A_{n}) \xrightarrow[n\to\infty]{} 1.$$ When we are in presence of properties $A_{n}$ for which $ G(n,p) \in A_{n}$ or not depends on $p_{n}$ in a drastic way (as $n \to \infty$), we speak of sharp threshold phenomena. In what follows we shall only focus on increasing graph properties:
\begin{definition}[Sharp thresholds for graph properties]  \label{def:sharpthres}Let $(A_{n} : n \geq 1)$ be a sequence of increasing properties of $ \mathbb{G}_n$. We say that $(A_{n})_{n\geq 1}$ has a \textbf{sharp threshold transition} for $(G(n,p))_{ p \in [0,1]}$ at $p\equiv p_{n}$ if for every $  \varepsilon>0$ we have 
$$ \mathbb{P}( G(n, (1- \varepsilon) p_{n})\in A_n) \to 0 \qquad \mbox{ whereas } \qquad \mathbb{P}( G(n,(1+ \varepsilon) p_{n})\in A_n) \to 1 \quad \mbox{ as }n \to \infty.$$
\end{definition}
Notice that the location of the edge density threshold $p_{n}$ is unique up to asymptotic equivalence. An alternative ``dynamic'' way of speaking of sharp threshold is to consider the Erd{\H{o}}s--R\'enyi graphs $(G(n,p ) : p \in [0,1])$ as naturally coupled via uniform labelings on the edges as in Section \ref{sec:percoabstrait}; then we write 
$$ \tau_{n} := \inf\{ p >0 : G(n,p ) \in A_{n}\}.$$ For increasing graph properties, if $p < \tau_n$ then $ G(n, p) \notin A_n$ whereas if $p > \tau_n$ then $G(n,p) \in A_n$.
Definition \ref{def:sharpthres} is then equivalent to the following concentration
$$ \frac{\tau_{n}}{p_{n}} \xrightarrow[n\to\infty]{( \mathbb{P})}1.$$

\begin{figure}[!h]
 \begin{center}
 \includegraphics[width=10cm]{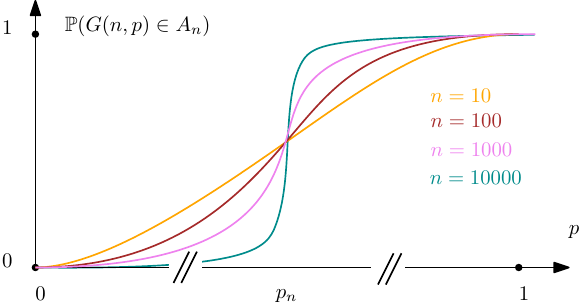}
 \caption{Illustration of the sharp threshold transition for an increasing graph property: the functions $x  \mapsto \mathbb{P}( G(n, x \cdot p_{n}) \in A_{n})$ converge pointwise on $[0,\infty) \backslash\{1\}$ towards the step function $ \mathbf{1}_{x < 1}$.}
 \end{center}
 \end{figure}

\begin{exo} \label{exo:weakthreshold} Let $A_n$ be a non-empty increasing graph property. Show that there exists $p_{n} \in (0,1)$ so that $A_{n}$ has a \textbf{weak threshold} at $p_{n}$ in the sense that for any sequence $\alpha_n \to \infty$ we have
$$\mathbb{P}\left( G(n, p_{n}/\alpha_{n}) \in A_n\right)) \to 0 \qquad \mbox{ whereas } \qquad \mathbb{P}( G(n,\alpha_{n} \cdot p_{n})\in A_n) \to 1 \quad \mbox{ as }n \to \infty.$$ 
\end{exo}

%
%

%
%
%

\section{Connectivity}
\label{sec:connectivity}

Probably the most natural question is to ask when the graph $G(n,p)$ becomes connected, i.e.~to consider the increasing graph property 
$$  \mathrm{Connected}_n = \{ \mathfrak{g} \in \mathbb{G}_n : \mathfrak{g} \mbox{ is connected}\}.$$
As will we see, with high probability this \textit{global} property is  in fact ruled by \textit{local} properties, namely the degrees of the vertices in $G(n,p)$.  As far as one given vertex $i$ is concerned, the situation is quite trivial since for every $i \in \{1, 2, \dots , n \}$ fixed, we have $$ \mathrm{deg}_{G(n,p)}(i) \overset{(d)}{=} \mathrm{Binomial}(n-1,p),$$ and so the expected degree of a given vertex is $(n-1)p$. But these degrees  are not independent! 

\subsection{Isolated vertices} 
We shall focus on  \textbf{isolated} vertices (i.e.\ of degree $0$), in the Erd{\H{o}}s--R\'enyi random graph. Consider the following increasing graph property
$$  \mathrm{NoIso}_n = \{ \mathfrak{g} \in \mathbb{G}_n: \forall 1 \leq i \leq n, \mathrm{deg}_ \mathfrak{g}(i)>0\},$$ and notice that we trivially have $ \mathrm{Connected}_n \subset \mathrm{NoIso}_n$ for every $n \geq 1$.

\begin{proposition} \label{prop:isolated}\noindent \noindent The sequence $( \mathrm{NoIso}_n)_{n\geq 1}$ has  a sharp threshold transition for  $(G(n, p))_{p \in [0,1]}$ at $$p_{n} = \frac{\log n}{n}.$$
\end{proposition}

\noindent \textbf{Proof.} We use the method of first and second moment. Since the degree of any single vertex in $ G(n, p_{n})$ follows a $ \mathrm{Bin}(n-1,p_{n})$ distribution, the first moment method shows that if $X(n,p)$ is the number of isolated vertices in $ G(n,p)$ then 
 \begin{eqnarray*} \mathbb{P}( G(n, p_{n}) \mbox{ has an isolated vertex})  &=& \mathbb{P}( X(n, p_{n})>0)\\ &\leq& \mathbb{E}[X(n,p_{n})]\\  &=& n \mathbb{P}( \mathrm{Bin}(n-1,p_{n}) = 0) = n (1-p_{n})^{n-1}.  \end{eqnarray*}
If $ p_{n} \geq (1 + \varepsilon)  \frac{\log n}{n}$ then the right-hand size clearly tends to $0$ as $n \to \infty$ and this shows that $ G(n, p_{n})$ has no isolated vertices w.h.p.\,in this regime. If now $ p_{n} \leq (1 - \varepsilon)  \frac{\log n}{n}$, we deduce from the last display that the expected number of isolated vertices diverges. To guarantee that  $\mathbb{P}(G(n,p_{n}) \mbox{ has an isolated vertex}) \to 1$, we use second moment method (Lemma \ref{def:second-moment}) and compute  
 \begin{eqnarray*}\mathbb{E}[X(n,p_{n})^{2}] &=& \sum_{1 \leq i,j \leq n} \mathbb{P}\left(i \mbox{ and }j \mbox{ are isolated in } G(n,p_{n}) \right)\\
 &=& n \mathbb{P}(1 \mbox{ is isolated}) + n(n-1) \mathbb{P}(1 \mbox{ and } 2 \mbox{ are isolated})\\
 & = & n(1-p_{n})^{n-1} + n(n-1)(1-p_{n})^{2n-3}.  \end{eqnarray*}
Notice that factor in the last display is not exactly equal to $\mathbb{P}(1 \mbox{ is isolated})^{2}$ since we only count  the  edge between the vertices $1$ and $2$ \textit{once}.  However, in the regime $p_{n} \leq (1- \varepsilon) \frac{\log n}{n}$ it is easy to see that $ \mathbb{E}[X(n,p_{n})]^{2} \sim \mathbb{E}[X(n,p_{n})^{2}]$ and so by Lemma \ref{def:second-moment} there are isolated vertices with high probability. \qed 

\subsection{Hitting time theorem}

Perhaps surprisingly, as soon as the graph $G(n,p)$ has no isolated vertices, it becomes instantaneously connected with high probability:
 \begin{theorem}[Erd{\H{o}}s--R\'enyi (1959)]\label{thm:connectedness}\noindent The sequence $( \mathrm{Connected}_n)_{n\geq 1}$ has  a sharp threshold transition for  $(G(n, p))_{p\in [0,1]}$ at $p_{n} = \frac{\log n}{n}.$ More precisely, in the coupled version of the Erd{\H{o}}s--R\'enyi random graphs, if we set 
 $$ \tau_{n} = \inf \{ p>0 : G(n,p) \mbox{ is connected}\}\quad \mbox{ and }\quad \theta_{n} = \inf \{ p>0 : G(n,p) \mbox{ has no isolated vertices}\}$$
 then we have $ \mathbb{P}(\tau_{n} = \theta_{n}) \to 1$ as $n \to \infty$. \label{thm:connectedfin}
 \end{theorem}
 We will see another proof of the first part of this result in the next chapter (Proposition \ref{prop:connectedness}).  \\

  Given Proposition \ref{prop:isolated}, it remains to understand whether the graph $G(n,p)$ can have several components which are not made of isolated vertices.  We say that a graph $ \mathfrak{g}$ has the \textbf{core property} if it is made of a (usually large, but possibly small) connected component of size at least $2$ together with isolated vertices. We denote by $ \mathrm{Core}_n$ the associated set of  graphs of $ \mathbb{G}_n$ satisfying the core property. Notice that this property \textit{is not increasing} (nor decreasing) and so some care is needed. 
  
\begin{figure}[!h]
 \begin{center}
 \includegraphics[width=5cm]{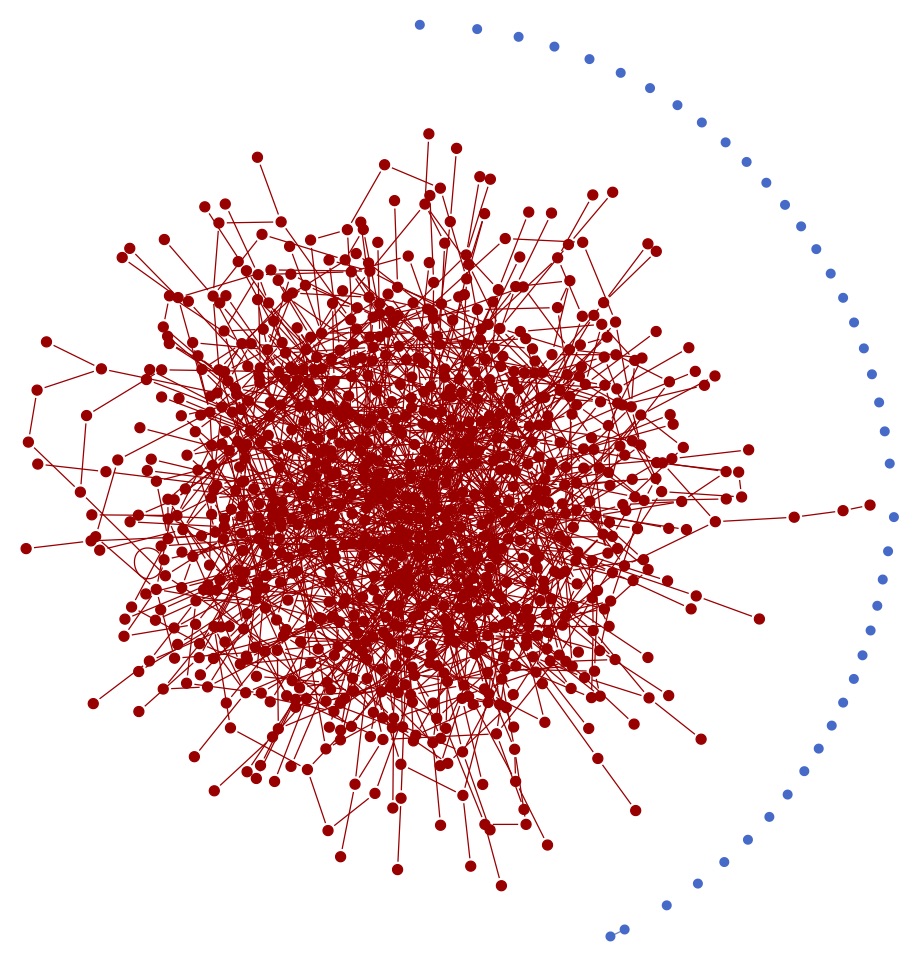}
 \caption{A graph having the core property: a single component (in red in the figure) together with isolated vertices (in blue). This property is however not stable by addition of edges.}
 \end{center}
 \end{figure}
 \begin{lemma}[Core property]  \label{lem:core} For any $ \frac{2}{3} \frac{\log n}{n} \leq p_n \leq 2 \frac{\log n}{n}$ we have  $G(n, p_{n}) \in \mathrm{Core}_n$  with high probability. \end{lemma}  \noindent \textbf{Proof of the lemma.} Actually, the proof will hint to the fact that the (sharp) phase transition for this property appears at $ \frac{1}{2} \frac{\log n}{n}$ (we put $2/3$ to be on safe ground). Let us denote $  \mathrm{Cut} \equiv \mathrm{Cut}(n,p)$ the number of ways to partition the vertices $\{1,2, \dots , n \}$  in two subsets $A \coprod B = \{1,2, \dots , n\}$ such that in $G(n,p)$ we have  $$ \left\{ \begin{array}{l} 2 \leq \# A \leq  \# B, \\ A \mbox{ is connected},\\  \mbox{there is no edge between } A  \mbox{ and }B.\end{array} \right.$$ Notice that if $G(n, p)$ does not have the core property, we can find two disjoint clusters of size at least $2$ and by adding components we can split the graph into two subsets $A,B$ as above (taking for $A$ the smallest subset). By the first moment method (applied twice) we have 
  \begin{eqnarray*} \mathbb{P}(G(n,p_n) \mbox{ has no Core}) &\leq& \mathbb{P}( \mathrm{Cut}(n,p_n) \geq 1) \\
  &\leq& \mathbb{E}[ \mathrm{Cut}(n,p_n)] \\ &=& \sum_{k=2}^{ \lfloor n/2 \rfloor}{n \choose k} \left((1-p_{n})^{n-k}\right)^{k}  \mathbb{P}( G(k, p_{n}) \mbox{ is connected})\\
  &=& \sum_{k=2}^{ \lfloor n/2 \rfloor}{n \choose k} \left((1-p_{n})^{n-k}\right)^{k}  \mathbb{P}(  \exists \mbox{ spanning tree in } G(k, p_{n}))\\
  &\leq&   \sum_{k=2}^{ \lfloor n/2 \rfloor}{n \choose k} \left((1-p_{n})^{n-k}\right)^{k}  \mathbb{E}\left[  \# \mbox{ spanning trees in } G(k, p_{n})\right]\\
  &\underset{ \mathrm{Cor}.\  \ref{cayley}}{=}& \sum_{k=2}^{ \lfloor n/2 \rfloor}{n \choose k} \left((1-p_{n})^{n-k}\right)^{k} k^{k-2} (p_{n})^{k-1},  \end{eqnarray*}
where the factor $ (1-p_{n})^{k(n-k)}$ counts for the probability that no edge is present between a subset of size $k$ and its complement in $ G(n,p_{n})$.
Take $ \varepsilon>0$ small enough so that if $k  \leq \varepsilon n$ we have $ \frac{2}{3}k(n-k) \geq \frac{7}{12} n k$. For those small $k$ we use the bound 
  \begin{eqnarray*} {n \choose k} \left((1-p_{n})^{n-k}\right)^{k} k^{k-2} (p_{n})^{k-1}  &\underset{\frac{2}{3} \frac{\log n}{n} \leq p_n \leq 2 \frac{\log n}{n}}{\leq}& \mathrm{Cst}^k\  \frac{n^{k}}{k!} \exp\left( -  \frac{7}{12} k n \frac{\log n}{n}\right) k^{k} \left( \frac{\log n}{n} \right)^{k-1}\\
   &\underset{  k^{k}/ k! \leq  \mathrm{Cst} \cdot \mathrm{e}^{k}}{\leq}&  \left( \mathrm{Cst} \cdot n \cdot n^{-7/12} \cdot  \mathrm{e} \cdot  \frac{\log n}{n} \right)^{k} \frac{n}{\log n}\\
   & \leq  &   \left(\mathrm{Cst} \cdot n^{-69/120} \right)^{k} n,  \end{eqnarray*} where $ \mathrm{Cst}$ is a universal constant that may vary from line to line.   Since $k \geq 2$, and $2 \times \frac{69}{120}>1$, those bounds are summable in $2 \leq k \leq n/2$ and are in fact dominated by the term $k=2$  which tends to $0$. For the large $k$, since $ n-k \geq n/2$ and $k \geq  \varepsilon n $, we have
$$ {n \choose k} \left((1-p_{n})^{n-k}\right)^{k}  \leq 2^{n} \exp\left( - \varepsilon n \log n / 100\right)  \leq ( 2 n^{ - \varepsilon/100})^{n},$$ and this bound can also by summed over the possible values of $k$ to get a vanishing quantity. The lemma is proved. \qed \medskip 
 
\noindent \textbf{Proof of Theorem \ref{thm:connectedness}.} The combination of Proposition \ref{prop:isolated} with Lemma \ref{lem:core} already shows that connectedness has a sharp threshold at $p_n = \frac{\log n}{n}$: w.h.p  there are still isolated vertices (and so the graph is not connected) at $ p = (1- \varepsilon)\frac{\log n}{n}$ by Proposition \ref{prop:isolated} whereas there are no isolated vertex at $ p_{n} = (1+ \varepsilon)\frac{\log n}{n}$ and by Lemma \ref{lem:core} the graph has the core property at this value of  $p_{n}$: it must be connected w.h.p. We cannot directly derive that $\tau_n = \theta_n$ because we cannot apply Lemma \ref{lem:core} to the random time $\theta_n$. However, the following strengthening of the lemma holds and enables us to conclude that $\tau_n = \theta_n$ w.h.p.\ :
  \begin{eqnarray} \mathbb{P}\left( G(n,p) \in \mathrm{Core}_n \mbox{\  simultaneously  for all }\frac{2}{3} \frac{\log n}{n} \leq p \leq 2 \frac{\log n}{n} \right) \xrightarrow[n\to\infty]{}1.  \label{eq:coreuniform}\end{eqnarray}
To see this, let us start from $ p_{n} = \frac{2}{3}  \frac{\log n}{n}$ where we know that the graph has the core property with high probability by Lemma \ref{lem:core}. Denote by $  \mathcal{C}_{n}$ its core and by $\ell_{1}, \dots , \ell_{X(n,p_n)}$ its isolated vertices. Recall from the proof of Proposition \ref{prop:isolated} that we have 
$$ \mathbb{E}[X(n,p_n)] = n (1-p_{n})^{{n-1}} \sim n \mathrm{e}^{- \frac{2}{3} \log n} = n^{1/3},$$ so that by Markov's inequality the event $ {H}_{n} = \{ X(n,p_n) \leq n^{5/12}\}$ happens with high probability as $n \to \infty$. Conditionally on $G(n,p_{n})$, consider for each isolated $\ell_{i}$ the next edge $(\ell_{i} \leftrightarrow x_{i})$ adjacent to $\ell_{i}$ to be added to $G(n,p_{n})$. These edges are not independent, but conditionally on $G(n,p_{n})$, for each $1 \leq i \leq X(n,p_n)$, the other extremity $x_{i}$ is uniform on $\{1,2, \dots , n\} \backslash \{ \ell_{i}\}$. In particular, the probability that $x_{i}$ does not belong to the core of $G(n,\frac{2}{3}  \frac{\log n}{n})$ is 
$$ \mathbb{P}(x_{i} \notin \mathcal{C}_{n}) = 1-\frac{ \# \mathcal{C}_{n}}{n} = \frac{X(n,p_n)-1}{n} \underset{ \mathrm{ on \ } H_{n}}{\leq} n^{-7/12}.$$
Conditionally on $G(n,p_{n})$ and on $H_{n}$, the expected number of $x_{i} \notin \mathcal{C}_{n}$ is bounded above by $n^{5/12}~\cdot~n^{-7/12}~=n^{-1/6}$ and by the first moment method we deduce that with high probability, for all $1 \leq i \leq X(n,p_n)$ the first edge connected to each isolated vertex $\ell_{i}$ after time $p_{n}$ will link it to the core $ \mathcal{C}_{n}$. In particular, no isolated vertices of $G(n,\frac{2}{3}  \frac{\log n}{n})$ get connected together and this entails \eqref{eq:coreuniform}. 
 \qed \medskip 
 
 We saw above that the variation of the individual degrees in $ G(n,p_{n})$ rules some large scale geometric properties. This variation disappears when $p_{n}\gg \frac{\log n}{n}$ and we leave the following as an exercise for the reader (after having given a look at Lemma \ref{lem:LDpoisson}):
\begin{exo}[Range of degrees]  \label{exo:asympdegree} For $\beta>0$ set $p_{n} = \beta \log n/n$. In $G(n, p_{n})$ let $m_{n}$ and $M_{n}$ respectively be the minimal and maximal vertex degrees. Show that we have 
$$ \frac{m_{n}}{\beta \log n} \xrightarrow[n\to\infty]{( \mathbb{P})} h(\beta), \quad \mbox{ and } \quad  \frac{M_{n}}{\beta \log n} \xrightarrow[n\to\infty]{( \mathbb{P})} H(\beta),$$ where $h(\beta)$ and $H(\beta)$ are the two solutions to the equation 
$$ \beta \cdot I(x) = 1, \quad \mbox{ with } I(a) := a \log a -(a-1),$$ where $h(\beta) = 0$ for $\beta < 1$ (there is only one solution). In particular $0=h(1)<H(1)= \mathrm{e}$ and $\beta^{-1}(H( \beta)- h( \beta)) \to 0$ as $\beta \to \infty$.
\end{exo}
 
 \section{Other thresholds via first and second moments}
 We present a few other sharp thresholds for appearance or disappearance of certain (induced) subgraphs in $G(n,p)$ whose proofs are also based on the first and second moment method.

\subsection{Diameter} \label{sec:diamER}
In this section, let us focus on the diameter of $G(n, p_{n})$, that is the maximal (graph) distance between any pairs of points in the graph. Of course, by the results of the last section, the diameter is finite only in the regime  $p_{n} \geq \frac{\log n}{n}(1+o(1))$. In a graph with maximum degree $d \geq 3$, by the same argument as in the proof of Proposition \ref{prop:lowerperco}, the number of vertices at distance less than or equal to $r$ from a given origin vertex is at most $1+ d + d(d-1) + d(d-1)^{2} + \cdots + d (d-1)^{r-1} = 1+ \frac{d}{d-2} ((d-1)^{r}-1)$. If the graph is connected and has $n$ vertices,  maximal degree $d$ and diameter $r$, we deduce the crude bound
$$  1+ \frac{d}{d-2}(d-1)^{r} \geq n.$$ Combining this with the rule of thumb ``$ \mbox{degrees} \approx \mbox{mean degree} \approx np_{n}$'' (valid as soon as $p_{n} \gg \frac{\log n}{n}$ by Exercise \ref{exo:asympdegree}) leads us to postulate that the diameter of $G(n,  p_{n})$ is roughly $ \frac{\log n}{\log np_{n}}$. Let us prove this fact in details for the case $ \mathrm{Diameter} =2$ which is the smallest non-trivial diameter. We denote by 
$$ \mathrm{Diam}^{2}_{n} = \{ \mathfrak{g} \in \mathbb{G}_{n} \mbox{ with diameter } \leq 2\}$$
the associated increasing graph property.

 \begin{proposition} The sequence $( \mathrm{Diam}^{2}_{n})_{n \geq 1}$ has a sharp  threshold  transition for $(G(n , p))_{p \in [0,1]}$ at $$p_{n} = \sqrt{ \frac{2 \log n}{n}}.$$
\end{proposition}
\noindent \textbf{Proof.} We again use the first and second moment method. The number $ D_n$ of pairs $\{i,j\}$ with $i \ne j \in \{1,2, \dots ,n\}$ so that $ \mathrm{d}_{ G(n,p_{n})}(i,j) > 2$ is easily computed since this just means that $i$ and $j$ are not neighbors nor share a neighbor:
  \begin{eqnarray} \mathbb{E}[D_n] &=& {n \choose 2} \mathbb{P}(1 \mbox{ and }2 \mbox{ are not neighbors and do not share a neighbor})\\
  &=& {n \choose 2} (1-p_{n})\cdot \mathbb{P}( \mbox{3 is not connected to both 1 and 2})^{n-2}\\ & =& {n \choose 2}(1-p_{n})( 1-p_{n}^{2})^{n-2} \sim \frac{n^{2}}{2}  \exp(- n p_{n}^{2}),  \end{eqnarray} when $p_{n} \to 0$. Hence if $ p_{n} \geq (1 + \varepsilon) \sqrt{2 \log n /n}$ the expected number of vertices at distance strictly larger than $2$ vanishes as $n \to \infty$. By the first moment method, this implies that w.h.p.~the diameter of $G(n, p_{n})$ is  less  than or equal to $2$ in this regime (it is then equal to $2$ unless $p_{n} = 1 - o(n^{-2})$; the diameter being $1$ if and only if all edges are present).\\
  
  We now suppose that $p_{n} = (1- \varepsilon) \sqrt{2 \log n/n}$ so that $n^2  \mathrm{e}^{-n p_n^2} \to \infty$ and the expectation of $D_n$ diverges. To prove that w.h.p.~there are vertices which do not share a neighbor we compute the second moment of $D_n$:
  $$ \mathbb{E}[D_n^{2}] = \sum_{i,j,k,l} \mathbb{P}\left(\begin{array}{c} i \mbox{ and }j \mbox{ have no common neighbor}\\
  i \mbox{ and }j \mbox{ are not neighbors}\\
  k \mbox{ and }l \mbox{ have no common neighbor}\\
    k \mbox{ and }l \mbox{ are not neighbors}\end{array}\right).$$
In the case when $i,j,k,l$ are all distinct, the possibilities of the induced subgraph $G(n,p_{n})[\{i,j,k,l\}]$ on $i,j,k,l$ are displayed in the following figure which account for a probability $(1-p_{n})^{6}+4p_{n}(1-p_{n})^{5}+2p_{n}^{2}(1-p_{n})^{4}$ which tends to $1$ as $n \to \infty$, and then the status of all other $n-4$ vertices contribute to a probability equal to $(1-p_{n}^{2})^{2(n-4)}$. 

\begin{figure}[!h]
 \begin{center}
 \includegraphics[width=13cm]{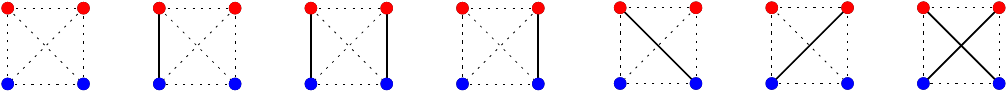}
 \caption{The possibilities for the induced subgraph on two pairs of vertices (here in red and blue) so that the distance between elements of each pair is at least $3$.}
 \end{center}
 \end{figure}
 The contribution of this case to the sum is then asymptotic to 
$$ \frac{n^4}{4}  \mathrm{e}^{-2np_n^2} \sim \mathbb{E}[D_n]^2.$$
Similarly, the contribution of overlapping pairs, i.e.\, the case when there are only three (resp.~two)  vertices among $i,j,k,l$ is negligible in front of the last term (we leave the details to the fearless reader). Since $ \mathbb{E}[D_n]$ diverges, only the first case prevails and $ \mathbb{E}[D_n^{2}] \sim  \mathbb{E}[D_n]^{2}$. By the second moment method (Lemma \ref{def:second-moment}) we conclude that $ \mathbb{P}(D_n>0) \to 1$ as $n \to \infty$. \qed 

\subsection{Ramsey theory,  max-clique, $ \mathrm{P=NP}$ and alien invasion}
A consequence of Ramsey's\footnote{\raisebox{-5mm}{\includegraphics[width=1cm]{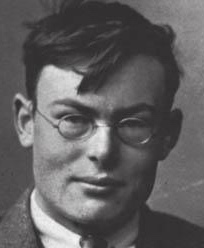}} Frank Ramsey  (1903--1930), English} theorem (for two colors) is that for any $k \geq 1$, there exists a number $ \mathcal{R}(k)$ such that every (simple) graph with more than $\mathcal{R}(k)$ vertices contains either a \textbf{clique} (an induced subgraph equal to the complete graph) of size $k$ or an \textbf{independent set} (induced subgraph with only isolated vertices)  with size $k$. Perhaps surprisingly, even the value of the Ramsey number $ \mathcal{R}(5)$ is unknown although it must lies in $\{43,44,45, \dots, 48\}$ (see wikipedia).  The bounds on $\mathcal{R}(k)$  are exponential in $k$ and actually random Erd{\H{o}}s--R\'enyi graphs achieve almost the best possible:

\begin{proposition} \label{prop:clique} Let $K_{n}$ and $ I_{n}$  respectively be the maximal size of a clique and of an independent set in $G(n,{\textstyle \frac{1}{2}})$. Then we have 
 $$ \frac{K_{n}}{ \log_{2}n} \xrightarrow[n\to\infty]{( \mathbb{P})} 2, \qquad  \frac{I_{n}}{ \log_{2}n} \xrightarrow[n\to\infty]{( \mathbb{P})} 2.$$ 
 \end{proposition}  

\noindent \textbf{Proof.}   Notice that $I_{n} = K_{n}$ in distribution since $ G(n,{\textstyle \frac{1}{2}})$ is self-dual in the sense that if we switch the status of all edges we obtain the same (law of) random graph. As the reader may have foreseen, we first compute the expected value of $ \chi_k$, the number of induced $k$-cliques in our random graph:   \begin{eqnarray} \label{eq:probaclique}  \mathbb{E}[ \chi_k]=\mathbb{E}\left[\# \left\{k-\mathrm{cliques\ } \mbox{ in }G(n,{\textstyle \frac{1}{2}})\right\}\right] = {n \choose k} 2^{- {k \choose 2}}.  \end{eqnarray}  It is easy to see that this tends to $0$ if $k\equiv k_n > (2 + \varepsilon) \log_{2} n$ as $ n \to \infty$. So by the first moment method, the size of the largest clique is less than $(2+ \varepsilon) \log_2 n$ w.h.p.

We now compute the second moment of the number of $k$-cliques and obtain  \begin{eqnarray*} \mathbb{E}[ \left(\chi_k \right)^2] &=& \sum_{\begin{subarray}{c}S,S' \subset \{1, \dots , n \}\\ \#S=\#S'=k \end{subarray}} \mathbb{P}( \mbox{the induced graphs on }S,S' \mbox{ are cliques})\\
 &=& {n \choose k} 2^{-{k \choose 2}} \sum_{\ell=0}^{k}{ k \choose \ell}{n-k \choose k-\ell} \cdot 2^{-{k \choose 2}+{\ell \choose 2}}.  \end{eqnarray*}
 To get the second line, we pick the first clique $S \subset \{1,2, \dots ,n \}$ and then partition according to the intersection $\ell = \#(S \cap S')$. If $S'$ shares some vertices with $S$, it is easier for it to be a clique since we only need to open ${ k \choose 2} -{\ell \choose 2}$ edges because the ${\ell \choose 2}$ edges in-between common vertices of $S$ and $S'$ are already present in (the induced subgraph of) $S$. Notice that when $\ell=0$ or $\ell=1$  the edges in-between vertices of $S$ and $S'$ are pairwise distinct. We then leave to the reader the tedious task of checking that the above sum is dominated by the term corresponding to $\ell=0$ when $k= k_n \leq (2-  \varepsilon) \log_2 n$, i.e.~that $ \mathbb{E}[\chi_{k_n}^2] \sim \mathbb{E}[\chi_{k_n}]^2$. By the second moment method (Lemma \ref{def:second-moment}), we deduce that indeed w.h.p, there are $k$-cliques for $k \leq (2- \varepsilon) \log_2 n$.  \qed \bigskip

\begin{remark}[Ramsey and Erd{\H{o}}s] Recall the definition of $ \mathcal{R}(k)$ as the smallest integer so that a graph with size larger than $ \mathcal{R}(k)$ must contain a clique or an independent set of size $k$. Proving that $ \mathcal{R}(k) < \infty$ is not  trivial and is in fact Ramsey's theorem. However, from \eqref{eq:probaclique} we deduce that if $ {n \choose k} 2^{-{k \choose 2}}<1/2$ then the expectation of $ \mathcal{X}_k + \mathcal{I}_{k}$ is less than $1$ where $ \mathcal{I}_{k}$ is the number of independent sets of size $k$, and this implies that $ \mathcal{X}_k+ \mathcal{I}_{k}$ is not almost surely larger than $1$ or equivalently that there \textbf{exists} a graph on $n$ vertices which has no clique nor independent set of size $k$. In other words, $$ \mathcal{R}(k) >n.$$ Although this reasoning  (one of the first instances of the \textbf{probabilistic method}) might appear simplistic, finding such a graph is a very difficult problem. Quoting Spencer \cite{spencer1994ten}: 
\begin{quote}``For the Ramsey function $ \mathcal{R}(k)$ no construction is known that gives nearly the lower bound that can be derived from the [above] proof... Erd{\H{o}}s asks us to imagine an alien force, vastly more powerful than us, landing on Earth and demanding the value of $ \mathcal{R}(5)$ or they will destroy our planet. In that case, he claims, we should marshall all our computers and all our mathematicians and attempt to find the value. But suppose, instead, that they ask for $ \mathcal{R}(6)$. In that case, he believes, we should attempt to destroy the aliens.'' \end{quote}
\end{remark}

\begin{remark}[Very sharp threshold] A careful inspection of the proof (and precise estimations) enable to reinforce Proposition \ref{prop:clique} as follows: There exists an integer $k_n \sim 2 \log_2 n$ so that we have 
$$ \mathbb{P}(K_n \in \{k_n, k_n+1\})  \xrightarrow[n\to\infty]{}1,$$ in other words, the maximal size of a clique is concentrated on only two values!
\end{remark}

\begin{remark}[Finding cliques] Although the previous result entails the existence of cliques of size $\approx 2 \log_2 n$ in $G(n, {\textstyle \frac{1}{2}})$, finding them is a very difficult task. Indeed, given a graph of size $n$, say by its adjacency matrix, an exhaustive search of a clique of size $\log_2 n$ costs
$$ {n \choose \log_2 n} \approx n^{\log_2 n} \mbox{ which is superpolynomial in }n.$$
It is known that finding the max clique in a (deterministic) graph is a NP-complete task, but more surprisingly it is open as of today whether we can find a clique of size $(1 + \varepsilon) \log_2 n$ in $G(n,{\textstyle \frac{1}{2}})$ -- hence a bit above half of the maximal size-- in a polynomial time! See the exercise below to find a clique of size approximatively $\log_2 n$.
\end{remark}   

\begin{exo}[Greedy construction of a clique] In $G(n, 1/2)$ , whose vertex set is $ \{1,2, \dots , n \}$, consider the following construction of a clique: Start with the vertex $X_0=1$. By induction, if $X_0<X_1<\dots < X_k$ have been constructed so that all edges $X_i \leftrightarrow X_j$ are present in $G(n, 1/2)$ for $0 \leq i < j \leq k$, let $X_{k+1}$ be the smallest vertex larger than $X_k$ which is connected to all $X_0, \dots , X_k$ in $G(n,1/2)$. If there is no such vertex the construction stops and output a complete induced subgraph, i.e.~a clique with $K_n$ vertices. Let $G_1, G_2, \dots $ be independent geometric variable with success parameter $2^{-i}$ i.e.
$$ \mathbb{P}(G_i = \ell) = 2^{-i}(1-2^{-i})^{\ell-1}.$$
\begin{enumerate}
\item Show that $ K_n = \min \{ k \geq 1 : 1+ \sum_{i=1}^k G_i > n\}$  in law.
\item Deduce that $$ \frac{K_n}{\log_2 n}  \xrightarrow[n\to\infty]{ ( \mathbb{P})}1.$$
\end{enumerate}
\end{exo}

\section{Higher moments}
So far, we have established sharp thresholds  for graph properties in $G(n,p)$ only using the first and second moment. When there is no sharp threshold or for  more refined probabilistic estimates such as convergence in distribution, we need to control higher moments. We shall exhibit two examples when we need to do so: the \textbf{Poisson paradigm} and the convergence of the spectral measure. Let us recall the classic method of moments:  
  
  \begin{lemma}[Method of moments]  \label{lem:metmom} Let $(\mu_{n} : n \geq 0)$ be probability measures on $ \mathbb{R}$ (resp.~random real variables $X_{n}$) such that for any $k \geq 0$, there exists $C_k \in \mathbb{R}$ such that we have
  $$   \int_{ \mathbb{R}} \mu_{n} ( \mathrm{d}x) \cdot x^k  \xrightarrow[n\to\infty]{} C_k, \qquad \left( \mbox{resp.} \quad \mathbb{E}[X_{n}^{k}] \xrightarrow[n\to\infty]{}C_{k} \right),$$ in particular the above moments all exist. We suppose furthermore that  for some $M >0$ we have $|C_k| \leq  M^k k!$ for all $k \geq 0$. Then there exists a probability measure  $\mu $ on $ \mathbb{R}$ (resp.~a random variable $X$) such that $\mu_{n} \to \mu$  in distribution as $n \to \infty$ (resp. $ X_{n} \to X$ in law).
  \end{lemma}
\noindent\textbf{Proof of the lemma:} Since $\mu_{n}$ have bounded first moment, $(\mu_n)_{n \geq 0}$ is tight and by dominated convergence its potential limits $\mu$ have the same moments $C_k$ for $k \geq 0$. However the growth condition $|C_k| \leq M^k k!$ implies that the moment generating function of (any possible limit) $\mu$ has a positive radius of convergence. By Fubini, the Laplace transform $ \mathcal{L}_\mu$ of $\mu$ also has a positive radius of convergence and, as every analytic function, $ \mathcal{L}_\mu$ is determined by its derivatives at $0$: it follows that $\mu$ is determined by its moments. In particular $\mu$ is unique, and $\mu_n \to  \mu$ weakly as $n \to \infty$. The translation in terms of random variables is straightforward.\qed \medskip

\begin{remark}[A trivial case: $ \sigma_n^2 \to 0$] When we have $ \mathbb{E}[X_n^1] \to C_1$ and $\mathbb{E}[X_n^2]\to C_1^2$ as $n \to \infty$ then we automatically have $ X_n \to C_1$ in distribution (and in probability). This was the case in most of the results of the previous section.
\end{remark}
 \subsection{The Poisson paradigm}
  \label{sec:poissonparadigm}
 In this section, we explain informally why the property 
 $$ \mathrm{Cycle}_{n} = \{   \mathfrak{g} \in \mathbb{G}_{n}: \mathfrak{g} \mbox{ contains a simple cycle as subgraph}\},$$
 actually has \textbf{no} {sharp threshold} transition for $G(n,p)$.  To fix ideas, let us look at the smallest non trivial cycle and let $\Delta(n, p)$ be the number of induced triangles in $G(n, p)$. One straightforwardly computes:
$$ \mathbb{E}[\Delta(n, p)]= {n \choose 3} p^{3},$$ and so when $p=p_{n}=o(1/n)$, by the first moment method, there is no triangle inside $ G(n, p)$ with high probability. When $p \gg 1/n$ then the expectation of the number of triangles blows up and we can show that the variance of the number of triangles is comparable to its squared mean (exercise!) so that by the second moment method, there is a triangle inside $G(n, p)$ with high probability. However when  $p = \frac{c}{n}$, the last expectation converge towards $c^{3}/6$ and actually $ \Delta(n, {\textstyle \frac{c}{n}})$ converges towards a Poisson variable of parameter $c^{3}/6$, this is the \textbf{Poisson paradigm}: the sum of many indicators of small probability which are roughly independent give a Poisson random variable in the limit. One way to prove it is to show that all moments of $ \Delta( n , {\textstyle \frac{c}{n}})$ converge towards the moments of $ \mathfrak{P}(c^{3}/6)$ and use the  method of moments (Lemma \ref{lem:metmom}). This requires a careful but not unbearable analysis which we will not do in this course.  \medskip 

In particular, existence of a triangle in $G(n, p)$ does not have a sharp threshold transition: its probability goes from $0$ to $1$ when $p$ ranges in the scale $1/n$ but does not jump from $0$ to $1$ abruptly ``in one scale''. Actually, the Poisson paradigm can be extended to consider cycles of length $3,4,5\dots$ simultaneously and actually they behave as independent Poisson variables with parameters $c^{k}/(2k)$ as $n \to \infty$ (each equivalence class of $k$ ordered points representing the same cycle has $2k$ members). In particular, for $c \in (0,1)$ we have 
$$ \mathbb{P}( G(n, {\textstyle \frac{c}{n}}) \mbox{ has no simple cycle}) \xrightarrow[n\to\infty]{} \exp\left( - \sum_{k \geq 3}  \frac{c^{k}}{2k} \right) = \sqrt{(1-c)}  \mathrm{e}^{ \frac{1}{4}c (c+2)}.$$

\begin{figure}[!h]
 \begin{center}
 \includegraphics[width=7cm]{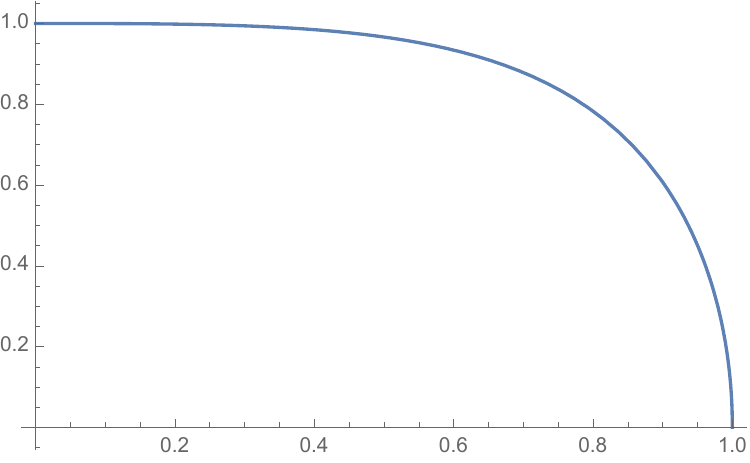}
 \caption{A plot of $\lim_{n \to \infty} \mathbb{P}( G(n, {\textstyle \frac{c}{n}}) \mbox{ has no simple cycle})$ for $ c \in [0,1]$. In particular the appearance of a simple cycle has no sharp threshold, but such a cycle should appear before $c=1$.}
 \end{center}
 \end{figure}

Recalling the first section of this chapter, although the presence of isolated vertices obeys a sharp threshold, if we refine the scale, the Poisson paradigm also appears and it is known that the number of isolated vertices in $ G\left(n, \frac{\log n +c}{n}\right)$ for $c \in \mathbb{R}$ converges in distribution towards a Poisson variable of parameter $   \mathrm{e}^{-c}$. In particular, we have the ``double exponential limit of Erd{\H{o}}s--R\'enyi''
  \begin{eqnarray} \label{eq:erconnectedpoisson} \mathbb{P}\left(  G\left(n, \frac{\log n}{n} +  \frac{c}{n}\right) \mbox{ has no isolated vertex}\right) \xrightarrow[n\to\infty]{} \mathrm{e}^{- \mathrm{e}^{-c}},  \end{eqnarray} see Theorem \ref{prop:connectedness} in Chapter \ref{chap:poissonER} for an proof of it.
 
 \begin{exo}[An application to random matrix]\label{exo:rangmatrix} Consider i.i.d.~vectors $X_1, \dots , X_k, \dots \in \{0,1\}^n$ such that $X_i$ has only zeros except at two positions chosen uniformly at random among the ${n \choose 2}$ possibilities. Evaluate $ \mathbb{P}(X_1, \dots, X_k \mbox{ are linearly independent over } \mathbb{Z}/2 \mathbb{Z})$ as a function of $k$ (for $n$ fixed but large). \label{exo:matrix2}
 \end{exo}

 \subsection{Spectrum} 
 In this section we shall study $G (n, p)$ from a spectral point of view. As the reader will see, this boils down to computing the (expected) number of (possibly backtracking) cycles in $G(n, p)$. In this section we focus on the case when 
\begin{center} $ \displaystyle p = \frac{c}{n}, \quad$ with $c >0$ fixed. \end{center}
 Let $A=A^{(n)}_c$ the adjacency matrix of $ G(n, {\textstyle \frac{c}{n}})$. More precisely, it is a symmetric square $n \times n$ matrix where $A_{i,j}=1$ if $i$ and $j$ are neighbors in $ G(n, p)$. Notice that the entries of $A$ are not centered and by convention we put $A_{i,i}= 0$. As any symmetric real matrix, $A^{(n)}_c$ has a spectral decomposition and $n$ eigenvalues 
 $$ \lambda_{1}^{(n)} \leq \lambda_{2}^{(n)} \leq \cdots \leq \lambda_{n}^{(n)}.$$
We  shall be interested in the empirical spectral measure 
 $$ \Lambda^{(n)}_c=   \frac{1}{n}\sum_{k=1}^{n} \delta_{\lambda_{i}^{(n)}},$$
 which is then a random probability measure (hence its law is an element of $ \mathcal{M}_{1}( \mathcal{M}_{1}( \mathbb{R}))$ where $ \mathcal{M}_{1}( \mathcal{X})$ is the set of all probability measures on $ \mathcal{X}$). The following theorem shows that this measure converges towards a deterministic measure:
 
 \begin{theorem}[Convergence of the spectral measure]\noindent With the above notation, for $c >0$ we have the following convergence in probability 
 $$ \Lambda^{(n)}_c \xrightarrow[n\to\infty]{( \mathbb{P})}  \mathcal{L}_{c},$$
 where $ \mathcal{L}_{c}$ is a (deterministic) probability measure on $ \mathbb{R}$.\end{theorem}
 The above convergence in probability just means that for any function $f \in \mathcal{C}_c( \mathbb{R})$  with compact support we have 
 $$  \int_{ \mathbb{R}} \Lambda^{(n)}_c( \mathrm{d}x) f(x) \xrightarrow[n\to\infty]{ (\mathbb{P})} \int_{ \mathbb{R}} \mathcal{L}_{c}( \mathrm{d}x) f(x).$$ It is possible to properly speak of convergence (in probability or in distribution) for random measures by defining a topology on probability distributions on $ \mathbb{R}$. We refer to the authoritative reference \cite{Kal86} for details about convergence of random measures. The limiting (deterministic) measure $ \mathcal{L}_{c}$ is poorly understood as of today (e.g.~the decomposition of $ \mathcal{L}_c$ in atomic and continuous part...).
 \begin{figure}[!h]
  \begin{center}
  \includegraphics[height=4cm]{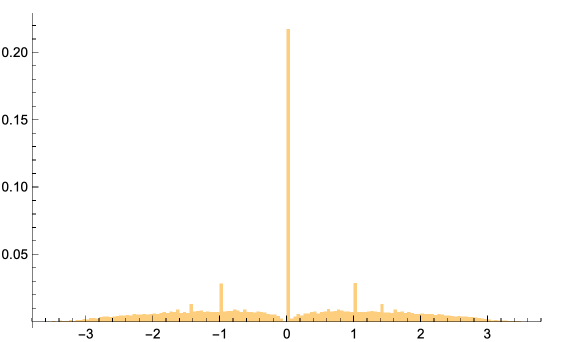}   \hspace{1cm} \includegraphics[height=4cm]{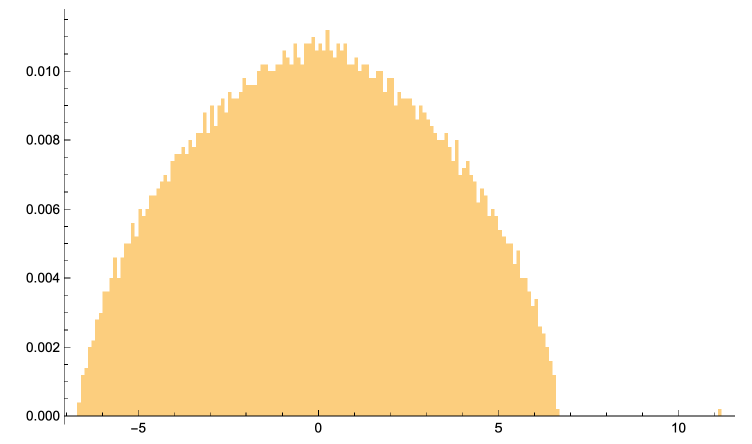}
  \caption{Simulations of $\Lambda^{(n)}_{c}$ for $c=2$ and $c=10$.}
  \end{center}
  \end{figure}

  \noindent \textbf{Partial proof.} We shall only prove a weak version of the theorem, namely that the \textbf{expected} empirical measure converges. To prove this, we shall use the method of moments (Lemma \ref{lem:metmom}) and prove convergence of the moments i.e.
  $$  \mathbb{E}\left[ \int \Lambda^{(n)}_c (\mathrm{d}x) \cdot x^k\right]  = \int  \underbrace{\mathbb{E}[\Lambda^{(n)}_c]}_{(*)} (\mathrm{d}x)  \cdot x^k \xrightarrow[n\to\infty]{}  \mathbb{E}\left[\int_{ \mathbb{R}}  \mathcal{L}_c ( \mathrm{d}x) \cdot x^k\right]=\int_{ \mathbb{R}} \underbrace{\mathbb{E}\left[ \mathcal{L}_c \right]}_{(**)}( \mathrm{d}x) \cdot x^k,$$ where $(*)$ and $(**)$ are the expected measures which are deterministic probability measures on $ \mathbb{R}$. The convergence in probability of the random measure $\Lambda^{(n)}_c$ is obtained by further establishing concentration of  the empirical moments (e.g.~by computing second moments), which we shall skip in these notes, see \cite{zakharevich2006generalization,khorunzhy2004eigenvalue} for details.  

\begin{figure}[!h]
   \begin{center}
   \includegraphics[width=14cm]{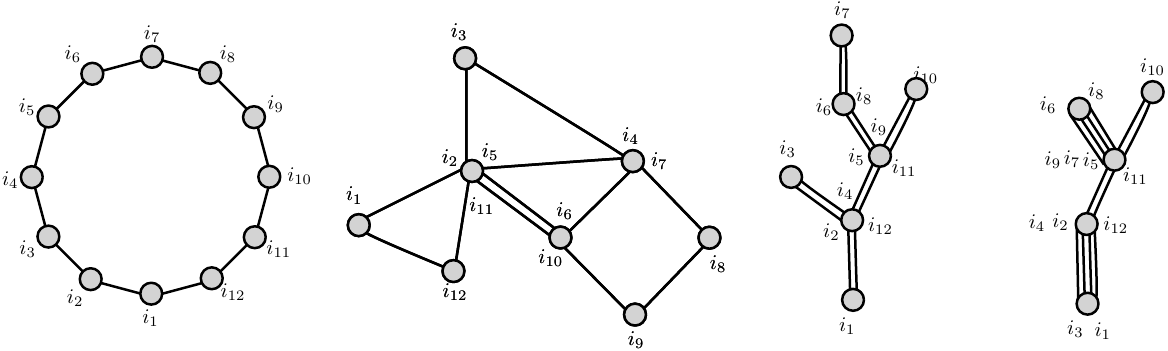}
   \caption{ \label{fig:graphdecomp} Expanding the expectation using a sum over Feynman diagrams. The only non zero asymptotic contribution comes from the trees with possibly several edges.}
   \end{center}
   \end{figure}

 Even the problem of the convergence of expectation of moments of $ \mathbb{E}[ \Lambda^{(n)}_c]$ might be complicated since the construction of the eigenvalues  of $ A^{(n)}_c$ is very intricate. The idea is to use the spectral decomposition and to take the expected trace of the powers of the matrix $A^{(n)}_c$: indeed by invariance of the trace under change of basis we have 
  $$n \cdot \int_{ \mathbb{R}} \Lambda^{(n)}_c( \mathrm{d}x) \cdot x^{k} = \sum_{k=1}^{n} (\lambda_{i}^{(n)}) ^{k} = \mathrm{Tr}\left(\left(A^{(n)}_c\right)^k\right)= \sum_{j=1}^{n}\sum_{j=i_{1}, \dots, i_{k}}A_{i_{1},i_{2}}A_{i_{2},i_{3}}\dots A_{i_{k},i_{1}}.$$
  After taking the expectation, all we need is to isolate the contribution of order $n$ in the above sum. We will  gather the terms in $\sum_{i_{1}, \dots, i_{k}} \mathbb{E}[A_{i_{1},i_{2}}A_{i_{2},i_{3}}\dots A_{i_{k},i_{1}}]$ which share the same the combinatorial structure for the cycle $i_1 \to i_2 \to \cdots \to i_k \to i_1$ and represent it  by a diagram. More precisely, we shall partition this sum according to the underlying multigraph $ \mathfrak{g}$ obtained by identifying in the ``free" cycle $i_{1} \to i_{2} \to \cdots \to i_{k} \to i_{1}$ the indices $i_j$ corresponding to the same vertex $\in \{1, \dots , n\}$. Those graphs are usually called Feynman's\footnote{\raisebox{-5mm}{\includegraphics[width=1cm]{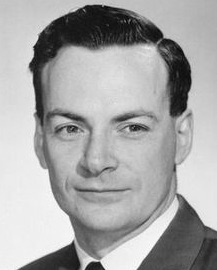}} Richard Feynman  (1918--1988), American} 	 diagrams in the physics literature. See Figure \ref{fig:graphdecomp} for examples.
  
   Once a multi-graph $ \mathfrak{g}$ with a rooted oriented spanning path is fixed, if $v$ is its number of vertices and $e$ its number of edges after collapsing the possible multi-edges, the corresponding contribution in the above sum is equal to 
   $$ {n \choose v} \left( \frac{c}{n} \right)^e \underset{n \to \infty}{\sim} \frac{c^e}{v!} n^{v-e}.$$
Hence, the main contribution to the above expectation is provided by Feynman diagrams for which $v-e$ is maximal: those are finite trees (for which we have $v-e=1$). More precisely, in this case $k= 2 \ell$ must be even and those objects are finite (non plane) trees with $e \leq \ell$ edges together with an image of the rooted $k$-cycle $1\to 2 \to \dots \to 2\ell \to 1$ which is surjective. If we denote by $ \mathrm{Fey}(2 \ell,e)$ the number of such combinatorial objects then we can summarize the discussion by 
 \begin{eqnarray*} \lim_{n \to \infty} \mathbb{E}\left[ \int \Lambda^{(n)}_c (\mathrm{d}x) \cdot x^k\right] &=& \lim_{n \to \infty}  \frac{1}{n} \sum_{j=1}^{n}\sum_{j=i_{1}, \dots, i_{k}} \mathbb{E}[A_{i_{1},i_{2}}A_{i_{2},i_{3}}\dots A_{i_{k},i_{1}}] \\ &=&  \mathbf{1}_{ k= 2 \ell \mathrm{\ is \ even}} \sum_{e=1}^{\ell} \mathrm{Fey}( 2 \ell ,e)  \frac{c^{e}}{(e+1)!} := C_{k}.  \end{eqnarray*}
 To check the growth condition on $C_{k}$ needed in Lemma \ref{lem:metmom}, notice that 
$$\sum_{e=1}^{\ell} \mathrm{Fey}( 2 \ell ,e)  \frac{c^{e}}{(e+1)!} \leq c^{\ell} \sum_{e=1}^{\ell} \mathrm{Fey}( 2 \ell, e) \leq c^{\ell} \underbrace{\#\{ \mbox{partitions of }\{1,2, \dots , 2 \ell \}\}}_{ \mbox{ Bell number}} \leq c^{\ell} (2 \ell)^{ 2 \ell} \leq M^{k} k!,$$ for some $M >0$. However, the number of diagrams corresponding to the moment of order $k$ grows quicker than exponentially: consider the case when the underlying tree is a star with $ \sqrt{n}$ vertices decorated by a walk of length $n$, it is easy to see that there are at least $( \sqrt{n})^{n- \sqrt{n}}$ diagrams and this in particular implies that the limiting measure $ \mathcal{L}_{c}$ has unbounded support. \qed

\paragraph{A warning to conclude: Thresholds and expectation thresholds.}   Before closing this chapter, let us warn the reader that the first (and second) moment method, although powerful, does not always yield the correct thresholds  for typical appearance of induced subgraph. Consider the following example of the ``pan graph''
 \begin{center}
 \includegraphics[width=4cm]{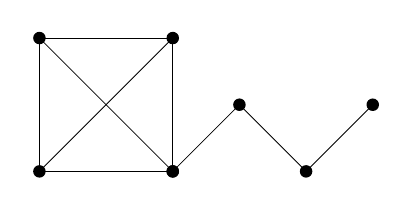}
 \end{center}
By the first moment, the mean number of such induced graphs in $G(n, p)$ is $ {n \choose 7} p^{9}$ and this blows up when $p \gg n^{-7/9}$. The naïve guess is then $p_{n} \approx n^{-7/9}$ for the (weak) threshold of appearance of the pan graph in $G(n, p_n)$. But if we focus on the appearance of the square with diagonals, we would guess a threshold of order $p_{n} \approx n^{-2/3}$. How come could the second threshold be larger than the first one since we consider a smaller subgraph? The reason is that the first guess $p_{n} \approx n^{-7/9}$ based on the first moment method is incorrect: with small probability, a square with diagonals may appear in $G(n, n^{-7/9 + \varepsilon})$ but then has many different ``tails'' causing a blow up of the expectation of such induced subgraphs...  The expectation threshold conjecture of Kahn \& Kalai  \cite{kahn2007thresholds} states that up to a multiplicative factor of $\log n$ the location of the (weak) threshold (see Exercise \ref{exo:weakthreshold}) is given by the above first moment method applied to all subgraphs and taking the maximum. This far-reaching conjecture was recently proved \cite{park2022proof}.  \bigskip

\noindent  \textbf{Bibliographical notes.} The Erd{\H{o}}s--R\'enyi model is probably the simplest and the most studied random graph model. It is a wonderful playground for combinatorics and probability. It has many variations and descendants such as the stochastic block model, the rank $1$ model, the configuration model...  which are more realistic models for real-life networks. The literature on this topic is vast, see e.g.~the recent monograph \cite{van2009random} or the classic books \cite{bollobas2001random,janson2011random}. There are also lecture notes available on the web  such as \cite{Bor15,blaszczyszyn2017lecture}  and \cite{wustatistical} for applications in statistics. Reading the original papers \cite{erdds1959random,erdHos1960evolution,erdos1963asymmetric} of Erd{\H{o}}s \& R\'enyi is still very inspiring. Theorem \ref{thm:connectedness} is proved in \cite{erdds1959random}. See the nice note \cite{cooper2019rank} for an application  of random graph theory to sparse random matrices (as in Exercise \ref{exo:matrix2}). \bigskip 

\noindent	\textbf{Hints for exercises.}\\
Exercise \ref{exo:ercayley}: By Cayley's formula the probability is equal to $n^{n - 2}p^{n - 1}(1 - p)^{{n \choose 2} - (n- 1)}$ and is maximal at $p = \frac{2}{n}$.\\
Exercise \ref{exo:weakthreshold}: Put $p_n$ such that $ \mathbb{P}(G(n,p_n) \in A_n) = \frac{1}{2}$. See   Bollobas \& Thomason  \cite{bollobas1987threshold} .\\
Exercise \ref{exo:asympdegree}: Use first and second moment method.\\
Exercise \ref{exo:rangmatrix}: Consider the graph whose vertices are the vectors and where there is an edge between two vectors if they share a non-zero coordinate in common. Then there is a non-trivial $ \mathbb{Z}/2 \mathbb{Z}$ relation iff the graph contains a cycle.

\chapter{Birth of a giant $1$, via $ \varepsilon$-cut}

\hfill Comment Gargantua nasquit en fa\c con bien estrange. (Rabelais)\bigskip

\vspace{2cm}

\begin{figure}[!h]
 \begin{center}
 \includegraphics[height=5cm]{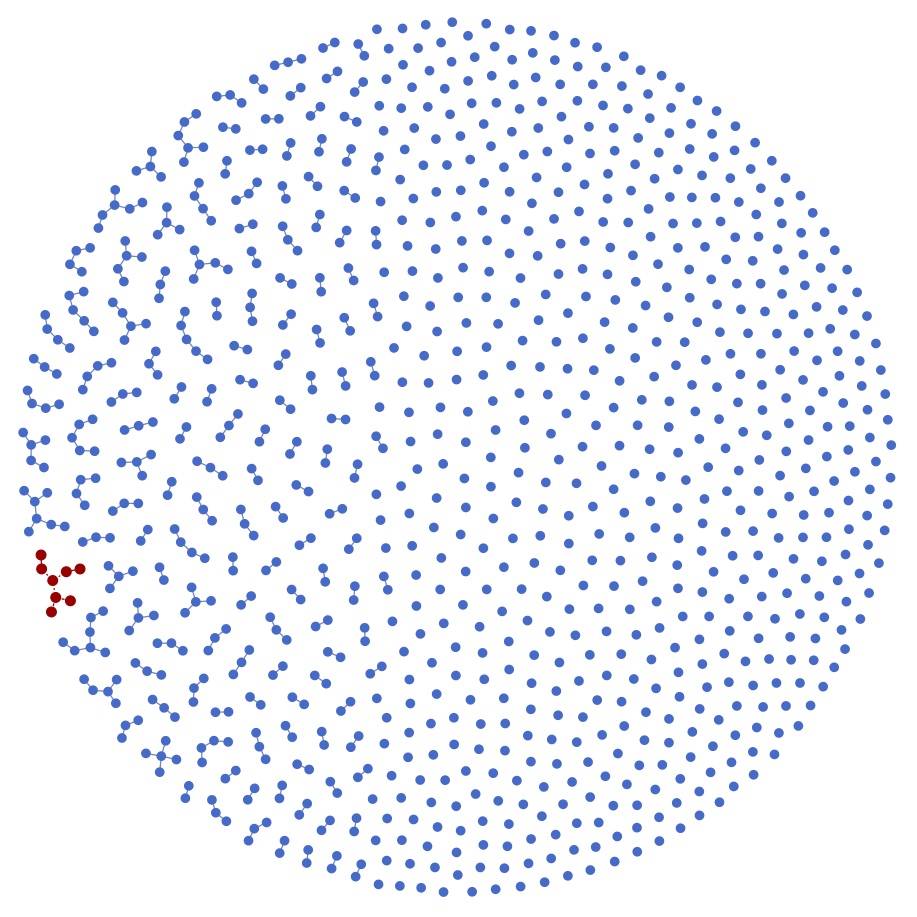}
  \includegraphics[height=5cm]{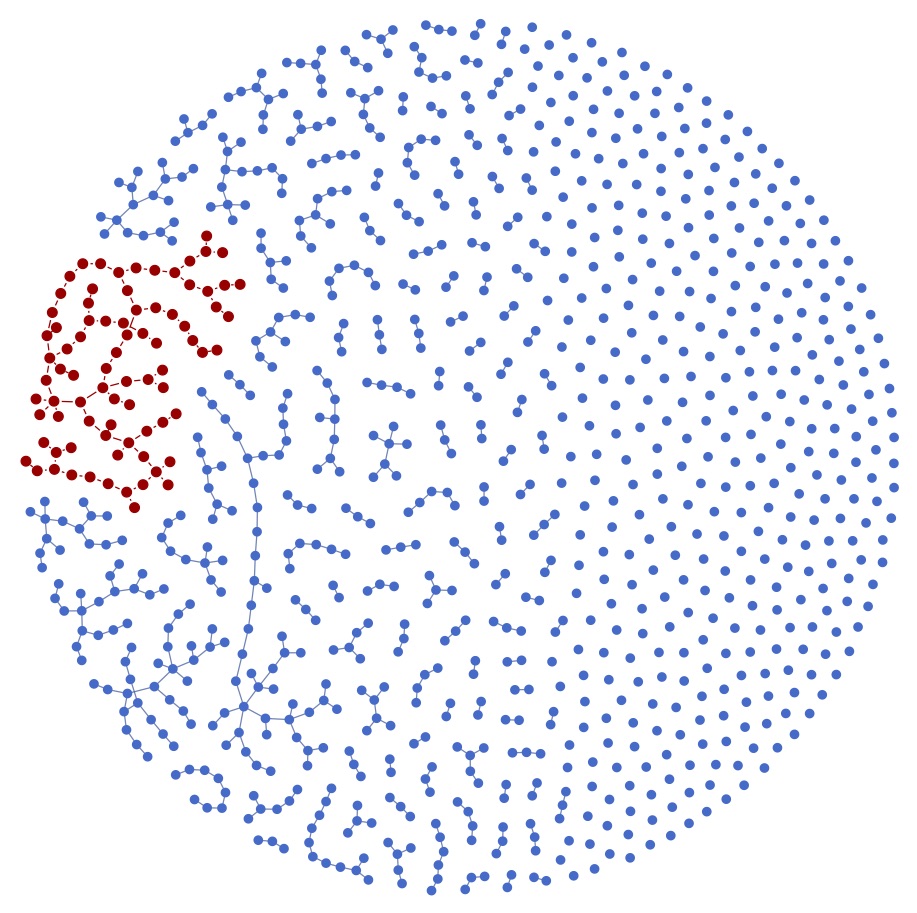}
   \includegraphics[height=5cm]{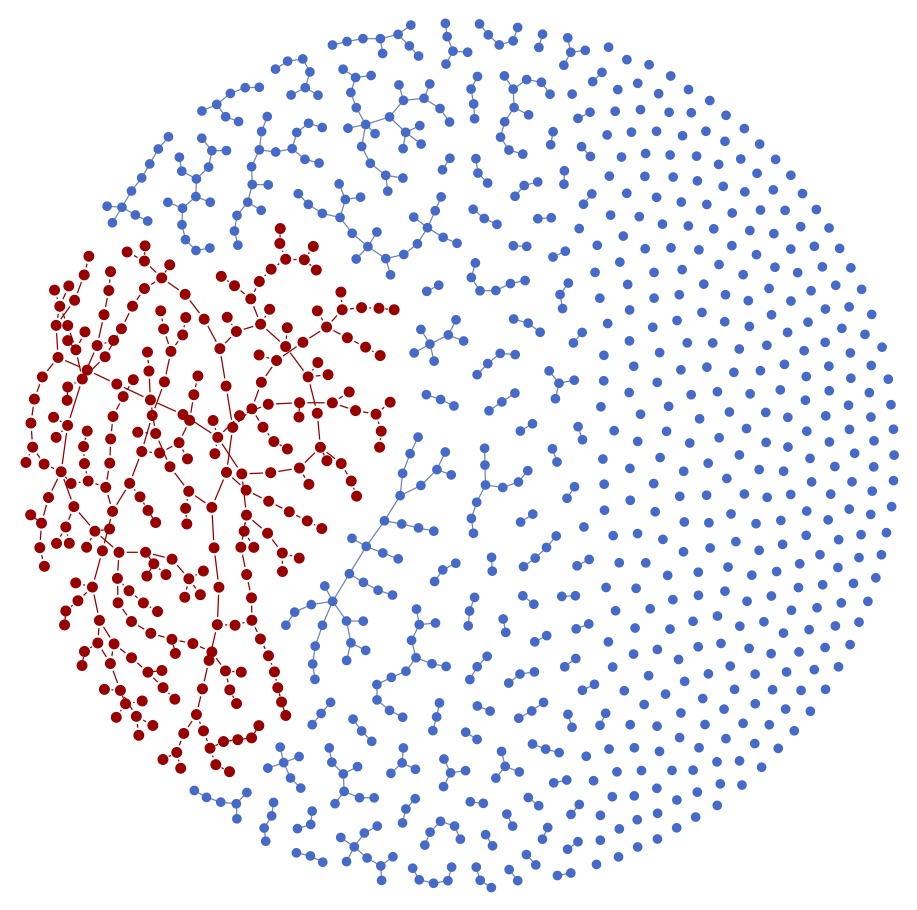}
      \includegraphics[height=5cm]{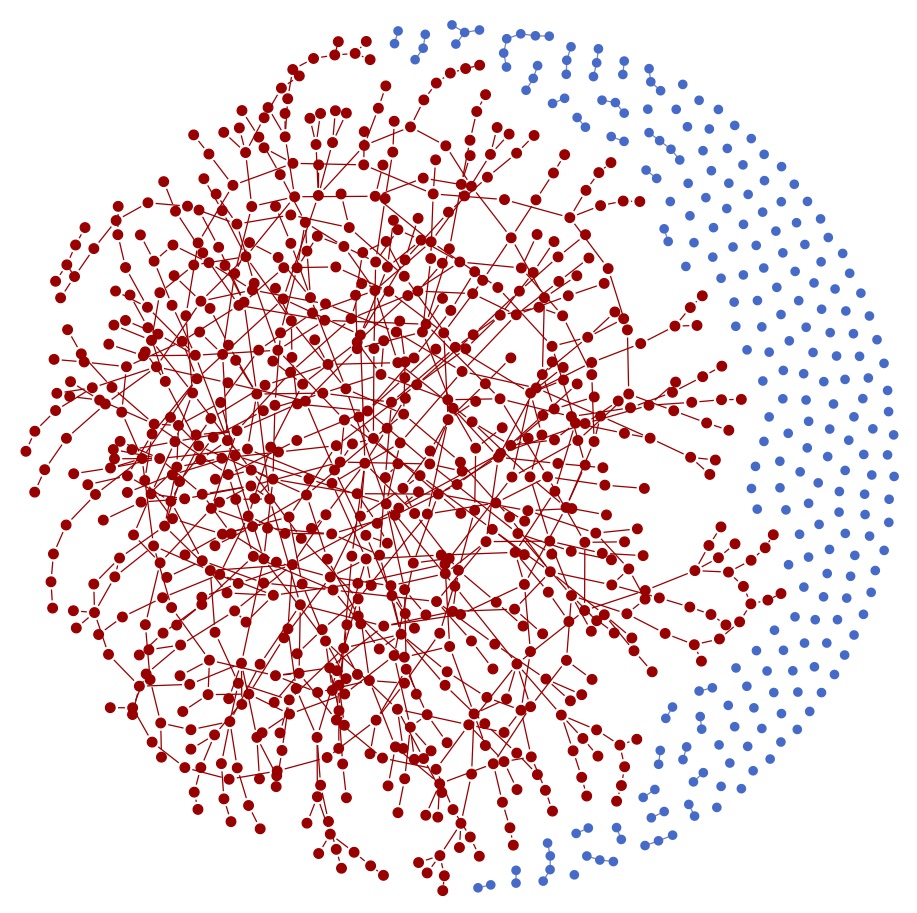}
         \includegraphics[height=5cm]{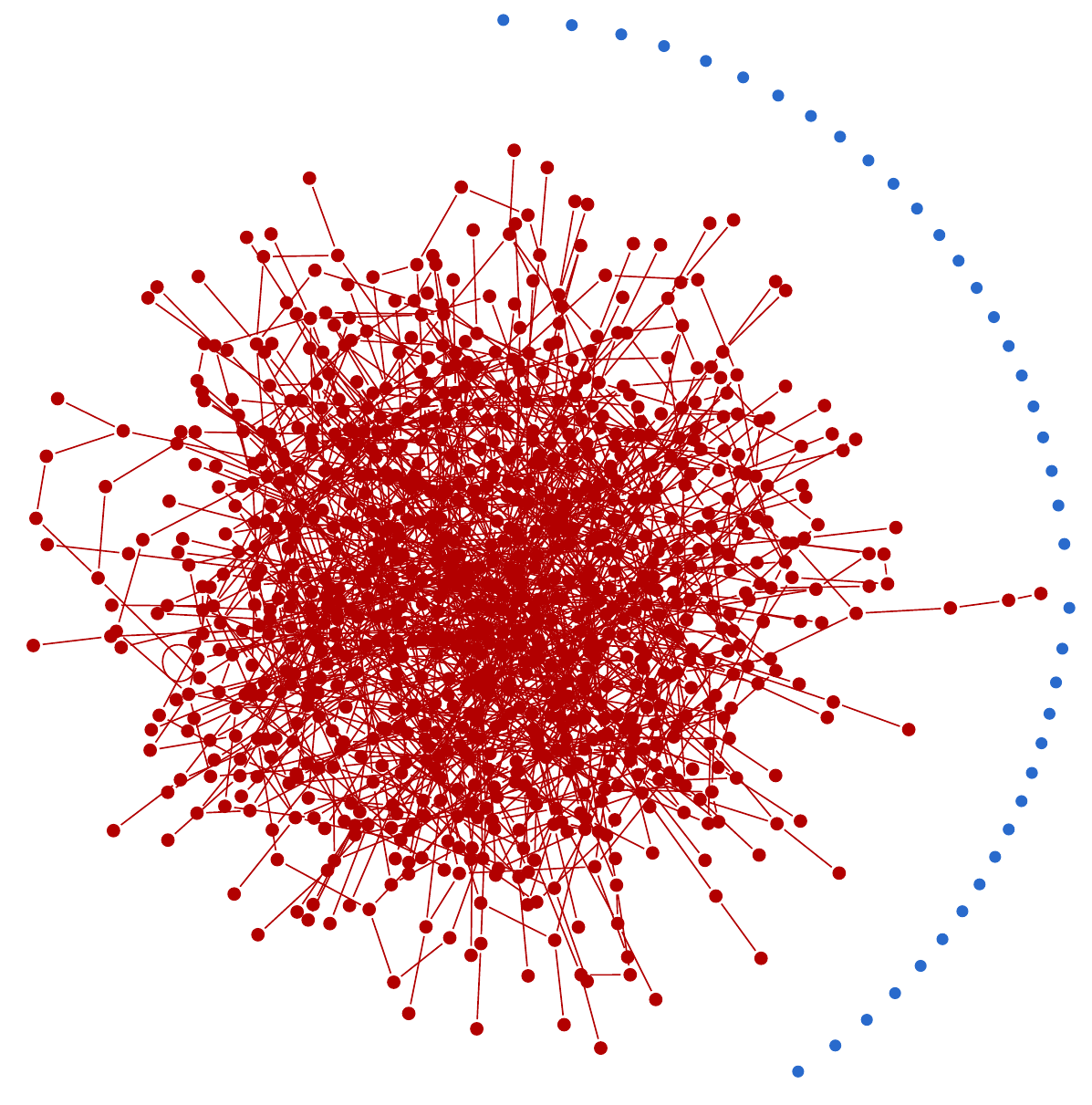}
            \includegraphics[height=5cm]{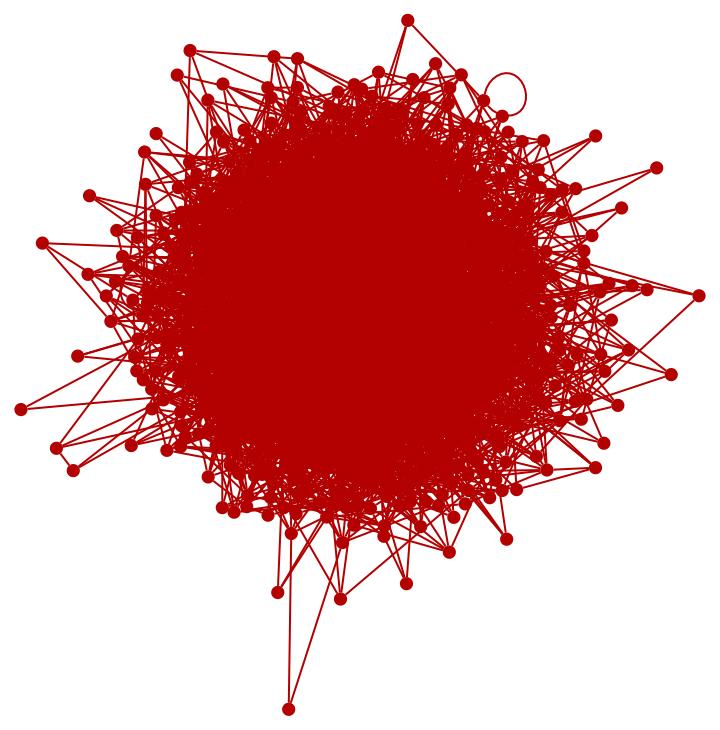}
 \caption{A large $G(n,{\textstyle \frac{c}{n}})$ graph with $n=1000$ and $c$ equals  (from left to right) to 
 $0.1\quad 0.5\quad 1 \quad 1.1 \quad  \frac{\log n}{2} \quad 2 \log n.$
 We  see the emergence of a giant component around $c \approx 1$ and that the graph becomes connected around $c \approx \log n$ (see Theorem \ref{thm:connectedfin}).}
 \end{center}
 \end{figure}
\clearpage 
We continue the study of the geometry of $G(n, p)$ as $p$ increases and prove a phase transition for the size of the clusters in $G(n, p)$:  suddenly a ``giant component'' carrying a positive proportion of the vertices appears around $p_{n} = \frac{1}{n}$.  More precisely, for a (finite) graph $ \mathfrak{g}$ and a vertex $v \in  \mathrm{V}(\mathfrak{g})$ we denote by $  \mathcal{C}^{v}( \mathfrak{g})$ the connected component of $v$ inside $ \mathfrak{g}$. We also denote by $ \mathrm{C}^{\max}_{1}( \mathfrak{g}),  \mathrm{C}^{{ \max}}_{2}( \mathfrak{g}),\dots$ the sizes (number of vertices) of the connected components of $ \mathfrak{g}$ in non-increasing order. We write $ \mathcal{C}_{\max}$ for a connected component of maximal size (with ties broken using the labelings of the vertices). Sometimes we drop the notation $ ( \mathfrak{g})$ when the underlying graph is clear from the context. The main theorem of this part is the following:
\begin{theorem}{Birth of a giant}\label{thm:erdosrenyi}\noindent There is a sharp threshold transition for the existence of a giant connected component at $p_{n} = \frac{1}{n}$. More precisely, if $p_{n} = \frac{c}{n}$ then  inside $G(n,  {\textstyle \frac{c}{n}})$: 
\begin{itemize}
\item \textbf{Subcritical: If  $ \mathbf{c<1}$} then there exists $A>0$ depending on $c>0$ such that w.h.p. we have $ \mathrm{C}^{\max}_{1}(G(n,{\textstyle \frac{c}{n}})) \leq A \log n$.
\item \textbf{Supercritical: If $ \mathbf{c>1}$} then there exists $A>0$ depending on $c>0$ such that w.h.p. we have $ \mathrm{C}^{\max}_{2}(G(n,{\textstyle \frac{c}{n}}))~\leq~A\log n$ whereas $ n^{-1} \mathrm{C}^{\max}_{1}(G(n,{\textstyle \frac{c}{n}})) \to  (1-\alpha(c))$ in probability where $\alpha(c)$ is the smallest solution in $(0,1]$ to the equation   \begin{eqnarray} \label{def:alphac} \alpha(c) = \mathrm{e}^{-c(1-\alpha(c))}.  \end{eqnarray}
\item \textbf{Critical: If $ \mathbf{c=1}$} then the vector $( n^{	-2/3}\mathrm{C}^{\max}_{i}(G(n,{\textstyle \frac{1}{n}})))_{i \geq 1}$ converges in law in the finite dimensional sense towards a positive infinite vector in $\ell^2$. This vector is in fact the ordered version of the length of the excursions of the function $t \mapsto B_{t} - \frac{t^{2}}{2}$ above its running infimum (sic!).
\end{itemize}
\end{theorem}

The goal of the following three chapters is to prove  the above result (multiple times). We will actually only prove points $(i)$ and $(ii)$ and just provide an upper bound for the size of the largest component for point $(iii)$ (see also Proposition \ref{prop:aldous} in a slightly different model). The intuition behind Theorem \ref{thm:erdosrenyi} is that the local neighborhood around a given vertex in $G(n, {\textstyle \frac{c}{n}})$ looks like a  BGW tree with offspring distribution $ \mathrm{Poisson}(c)$, see Proposition \ref{prop:cvlocpoisson}. When $c<1$ such a random tree dies out almost surely (in fact very quickly) and all the connected components in $G(n, {\textstyle \frac{c}{n}})$ are small. On the contrary, if $ c>1$ then the BGW process survives with positive probability equal to $1- \alpha(c)$ (see Proposition \ref{prop:GWdisguise} and Theorem \ref{thm:extinction}): in the finite setting this means that the component of the vertex in question is very large (a giant component). It turns out that this giant component is unique so that its density is asymptotically $1- \alpha(c)$, and the remaining components are small. We prove a weaker version of the above theorem using this sketch in this chapter then turn to a more modern proof using an exploration technique and estimates on random skip-free walks similar to those used in Chapter \ref{chap:GW}. This proof is shortened in Chapter \ref{chap:poissonER} by tricking a little the graph. 

\bigskip 

\label{sec:firstmomentgiant}
In this chapter we give a ``first moment'' proof of a weak version of Theorem \ref{thm:erdosrenyi} which is close in spirit to the historical proof of Erd{\H{o}}s \& R\'enyi \cite{erdHos1960evolution}. We first study the law of the connected components in $G(n,{\textstyle \frac{c}{n}})$ and show that they are indeed described by BGW trees with Poisson$(c)$ increments. We then use a first moment argument on $ \varepsilon$-cuts together with a ``sprinkling'' idea to deduce the existence of a giant component in the supercritical regime:

\begin{theorem}[The giant, weak version]  \noindent For $c>0$ denote by $\alpha(c)$  the smallest solution in $(0,1]$ to the equation  $\alpha(c) = \mathrm{e}^{-c(1-\alpha(c))}$ in particular $\alpha(c) =1$ when $c \leq 1$.  \label{thm:weak-giant} Then we have $$ \frac{\mathrm{C}^{\max}_1}{n} \xrightarrow[n\to\infty]{( \mathbb{P})} 1 - \alpha(c) \quad \mbox{ and } \quad  \frac{\mathrm{C}^{\max}_2}{n} \xrightarrow[n\to\infty]{( \mathbb{P})} 0.$$
  \end{theorem}

\section{The local limit}
We denote by   $ \mathcal{C}(n, p) \equiv \mathcal{C}^{1}(G(n,p))$ the cluster of the vertex $1$ in $G(n,p)$ which we see as a random labeled graph where its vertices have been relabeled in increasing order by $1,2, \dots , \#  \mathcal{C}(n,p)$. We write  $\mathrm{T}(c)$ for a Bienaymé--Galton--Watson plane tree with offspring distribution Poisson$(c)$ and denote by $ \mathcal{T}(c)$ the Cayley tree obtained by labeling its root by $1$ and the rest of its vertices by $2,3, \dots , \#  \mathrm{T}(c)$ uniformly at random. We put $ \mathcal{T}(c)= \dagger$ (a cemetery point) if $ \mathrm{T}(c)$ is infinite. 
\begin{proposition}  \label{prop:cvlocpoisson}Fix $c >0$ and suppose that $p_{n} \sim \frac{c}{n}$ as $n \to \infty$. For $k \geq 1$ and for any connected labeled graph $ \mathfrak{g} \in \mathbb{G}_k$ we have 
$$\lim_{n \to \infty }\mathbb{P}( \mathcal{C}(n,p_{n}) = \mathfrak{g})  =   \mathrm{e}^{-k \cdot c } \frac{c^{k-1}}{(k-1)!} \mathbf{1}_{  \displaystyle \mathfrak{g} \mbox{ is a Cayley tree}} = \mathbb{P}( \mathcal{T}(c) =  \mathfrak{g}).$$
\end{proposition}
\noindent \textbf{Proof.} If the connected labeled graph $ \mathfrak{g}$ with $k$ vertices and $\ell$ edges is fixed, we have 
$$\mathbb{P}( \mathcal{C}(n,{\textstyle \frac{c}{n}}) = \mathfrak{g}) \quad = \quad {n-1 \choose k-1} (1-p_{n})^{k(n-k) + {k \choose 2} -\ell} p_{n}^{\ell}$$
Taking limits as $n \to \infty$, the above display tends to $0$ if $\ell \geq k$ and towards $\mathrm{e}^{-k \cdot c } \frac{c^{k-1}}{(k-1)!}$ if $\ell = k-1$. In the latter case, this formula coincides with the probability that $ \mathcal{T}(c)$ lands on the tree $\mathfrak{g}$ as seen in the proof of  Proposition \ref{prop:cayleyGW}. \qed

\medskip 

In particular, when $c \leq 1$, since $ \mathrm{T}(c)$ is almost surely finite  we have that $$ \sum_{k \geq 1} k^{k-2} \mathrm{e}^{-k \cdot c } \frac{c^{k-1}}{(k-1)!} =1$$ so that $ \# \mathcal{C}(n,{\textstyle \frac{c}{n}}) \to \# \mathrm{T}(c)$ in distribution as $ n \to \infty$. We recall from Definition \ref{def:borel-tanner} that the latter is distributed according to the Borel--Tanner distribution with parameter $c$. As a quick corollary we can deduce that there is no giant component when $c \leq1$: If  $ \mathcal{C}_{\max}$ in a largest connected component in $G(n, {\textstyle \frac{c}{n}})$  we have
 \begin{eqnarray*}0 \xleftarrow[n\to\infty]{} \mathbb{P}(  \#\mathcal{C}(n,{\textstyle \frac{c}{n}}) \geq \varepsilon n) &=&  \frac{1}{n}\sum_{i=1}^{n} \mathbb{P}( \mathcal{C}^{i}(G(n, \frac{c}{n})) \varepsilon n)\\  &\geq&   \frac{1}{n} \sum_{i=1}^{n}\mathbb{E} \left[    \mathbf{1}_{ \mathrm{C}^{\max}_1 \geq \varepsilon n}  \mathbf{1}_{i \in \mathcal{C}^{\max}_{1}}\right]  \geq   \frac{ \varepsilon n}{n}  \mathbb{P}(  \mathrm{C}^{\max}_1 \geq  \varepsilon n).  \end{eqnarray*}
More precisely, we will see below that the proportion of vertices belonging to ``big clusters'' is concentrated. But before that let us state an easy lemma whose proof is straightfoward:
 \begin{lemma} \label{lem:condcluster}Conditionally on  $\mathcal{C}^{1}( G( n ,p ))$ the remaining graph\footnote{with vertices relabeled in increasing order} $ G(n,p) \backslash   	\mathcal{C}^{1} (G (n, p ))$ has law $G( n - \# \mathcal{C}^{1}( G( n ,p )))$.
 \end{lemma}

\begin{corollary} \label{cor:propotionfini} With  $\alpha (c)$   as defined in Theorem \ref{thm:weak-giant}, for all $ A \in \{0,1,2, \dots \}$ we have
$$  n^{-1} \, \#\left\{ 1 \leq i \leq n : \# \mathcal{C}^{i}\big({\textstyle G(n, \frac{c}{n})}\big) \leq A\right\} \xrightarrow[n\to\infty]{ ( 	\mathbb{P})} \mathbb{P}( \#\mathrm{T}(c) \leq A) \xrightarrow[A\to\infty]{} \alpha(c).$$
\end{corollary}
\noindent \textbf{Proof.} If $ N_{n}(A) = \#\left\{ 1 \leq i \leq n : \# \mathcal{C}^{i}\big({\textstyle G(n, \frac{c}{n})}\big) \leq A\right\}$, by the previous proposition, we  have the asymptotic of the expectation:
$$ \mathbb{E}[N_{n}(A)] =  \sum_{i=1}^{n} \mathbb{P}( \# \mathcal{C}^{i} \leq A) = n \mathbb{P}( \#\mathcal{C}(n, {\textstyle \frac{c}{n}}) \leq A) \underset{n \to \infty}{\sim}  n \cdot  \mathbb{P}( \# \mathrm{T}(c) \leq A),$$ and this asymptotic actually holds as soon as $p_{n} \sim \frac{c}{n}$. Using $ \mathcal{C}^{i} = \mathcal{C}^{i} (G (n , c/n))$ as a shorthand notation, the second moment is easily bounded:
 \begin{eqnarray*}  \mathbb{E}[N_{n}(A)^{2}]&=&\sum_{1 \leq i , j \leq n} \mathbb{P}( \# \mathcal{C}^{i} \leq A \mathrm{ \ and \  } \# \mathcal{C}^{j} \leq A) \\
 &=& \sum_{1 \leq i , j \leq n} \left(\mathbb{P}( \# \mathcal{C}^{i} \leq A \mbox{ and } j \in \mathcal{C}^{i}) + \mathbb{P}( \# \mathcal{C}^{i} \leq A \mbox{ and }  \# \mathcal{C}^{j} \leq A \mbox{ and } j \notin \mathcal{C}^{i})\right)\\
 &=& n \mathbb{E}\left[ \# \mathcal{C}^{1} \mathbf{1}_{\# \mathcal{C}^{1} \leq A}\right] + n(n-1) \mathbb{P}\left( \# \mathcal{C}^{1} \leq A \mbox{ and }  \# \mathcal{C}^{2} \leq A \mbox{ and } 2 \notin \mathcal{C}^{1}\right)\\
 & \leq &  n A + n(n-1)  \mathbb{E}\Big[ \mathbf{1}_{  \#\mathcal{C}^{1} \leq A}   \mathbb{E}\left[ \mathbf{1}_{\# \mathcal{C}^{2} \leq A \mbox{ and } 2 \notin \mathcal{C}^{1} } \mid \mathcal{C}^{1}\} \right] \Big] \\
 & \underset{ \mathrm{Lemma\ \ref{lem:condcluster}}}{=}& n A + n(n-1)  \mathbb{E}\left[ \mathbf{1}_{  \#\mathcal{C}^{1} \leq A}   \mathbb{E}\Big[ \mathbf{1}_{2 \notin \mathcal{C}_1} \underbrace{\mathbb{E}\big[ \mathbf{1}_{\# \mathcal{C}^{2} \leq A}  \big | \{ 2 \notin \mathcal{C}^{1} \} \mbox{ and } \mathcal{C}^{1} \big]}_{ \mathbb{P}(  \# \mathcal{C}( N, \frac{c}{n}) \leq A) } \Big | \mathcal{C}^{1} \Big] \right] \\
 & =&  n A + n(n-1) \mathbb{E}\left[  \frac{N}{n-1} \cdot \mathbb{P}\left(  \# \mathcal{C}\left( N,  \frac{c}{n}\right) \leq A\right)\right],
\end{eqnarray*}
 where we used Lemma \ref{lem:condcluster} to argue that once conditioned on the cluster $ \mathcal{C}^{1}$, the remaining graph is distributed as $G(N, {\textstyle \frac{c}{n}})$ where $N = n  - \# \mathcal{C}^{1}$, so that by symmetry $ \mathbb{P}( 2 \notin \mathcal{C}^1 \mid \mathcal{C}_1) = \frac{N}{n-1}$ and $\mathbb{E}\left[ \mathbf{1}_{\# \mathcal{C}^{2} \leq A}  \mid \{ 2 \notin \mathcal{C}^{1} \} \mbox{ and } \mathcal{C}^{1} \right] = \mathbb{P}(  \# \mathcal{C}( N, c/n) \leq A)$. Since $N \sim n$, we can apply the lemma above once more and deduce that the previous display is asymptotic to $n^{2} \mathbb{P}( \#\mathrm{T}(c) \leq A)^{2}$. We deduce that $ \mathbb{E}[n^{-1} N_{n}(A)] \to \mathbb{P}( \#\mathrm{T}(c) \leq A)$ and $ \mathrm{Var}( n^{-1} N_{n}(A)) \to 0$ as $n \to \infty$, which by Tchebytchev inequality entails the convergence in probability in the corollary. The convergence of $n^{-1} N_{n}(A)$ to $\alpha(c)$ as $A \to \infty$ follows from Theorem \ref{thm:extinction}. \qed

 \section{An easy giant via $ \varepsilon$-cut}
Fix a graph $ \mathfrak{g}$ with $n$ vertices and $ \varepsilon>0$. An \textbf{$ \varepsilon$-cut} is a partition of $\{1,2, \dots , n\}$ into two subsets $A$ and $B$ of size (number of vertices) $ \varepsilon n$ and $(1- \varepsilon)n$ so that there is no edge between $A$ and $B$ in $ \mathfrak{g}$. That notion was already used in the proof of Lemma \ref{lem:core}.  The following deterministic lemma relates the existence of $  \varepsilon$-cuts to the size of the largest component:
\begin{lemma}[Giants makes cutting difficult] Recall that $ \mathrm{C}^{\max}_{1}( \mathfrak{g})$ is the size of the largest component in $ \mathfrak{g}$. Then one can find an $ \varepsilon$-cut in $ \mathfrak{g}$ with 
$$ \left| \varepsilon - \frac{1}{2} \right| \leq \frac{1}{2} \frac{ \mathrm{C}^{\max}_{1}( \mathfrak{g})}{n}.$$
\end{lemma}
\noindent \textbf{Proof.} Let us put $x_i = n^{-1} \mathrm{C}^{\max}_i ( \mathfrak{g})$ for the renormalized cluster sizes in $ \mathfrak{g}$ ranked in non-increasing order so that $\sum_{i \geq 1} x_i =1$. Let $$\ell = \inf\left\{ j \geq 1 : \sum_{i=1}^j x_i \geq  \frac{1}{2}\right\}.$$ Since $x_\ell \leq x_1 = \frac{ \mathrm{C}^{\max}_1 ( \mathfrak{g})}{n}$ we deduce that either $\sum_{i=1}^{\ell-1} x_i$ or $\sum_{i=1}^{\ell} x_i$ is $ \frac{x_1}{2}$-close to $1/2$. Regrouping the vertices of the corresponding components,  we get the desired $ \varepsilon$-cut. \qed

\paragraph{An easy giant.} Let us use this lemma to quickly prove that there is a large component in $G(n,{\textstyle \frac{c}{n}})$ when $c > 4 \log 2$ (this is a much weaker statement compared to Theorem \ref{thm:erdosrenyi}). Indeed, one can upper-bound the expected number of $ \varepsilon$-cuts in $G(n,p)$ by 
  \begin{eqnarray} \mathbb{E}\left[\# \varepsilon-\mbox{cuts in }G(n,p) \right] \leq {n \choose \varepsilon n} (1-p)^{ \varepsilon(1- \varepsilon) n^2} \underset{p= \frac{c}{n}}{\leq} 2^n \mathrm{exp}(- c  \varepsilon(1- \varepsilon) n).   \label{eq:cutgiant}\end{eqnarray}
When $c > 4 \log 2$ the right-hand side tends exponentially fast to $0$ as soon as $ \varepsilon \in  ( \frac{1-\delta_c}{2},\frac{1+\delta_c}{2}  )$ where $(\delta_c)^2 =  {1 - \frac{4 \log 2}{c}}$. Summing over all  the at most $n$ possible values of $ \varepsilon$ in this range, we  deduce by the first moment method that for all $ \eta >0$ we have 
$$ \mathbb{P}\left( \exists\   \varepsilon- \mbox{cut in } G(n,{\textstyle \frac{c}{n}}) \mbox{ with } \left| \frac{1}{2} - \varepsilon \right | \leq \frac{\delta_c- \eta}{2} \right)  \xrightarrow[n\to\infty]{} 0,$$
hence by the above lemma $$    \quad \mathbb{P}\left( \exists \mbox{ a cluster of size at least } (\delta_c -\eta) \cdot n \mbox{ in } G(n,{\textstyle \frac{c}{n}})\right)  \xrightarrow[n\to\infty]{} 1.$$

The former reasoning becomes very useful when we already start from a graph having large clusters: Suppose that $ \mathfrak{g}\in \mathbb{G}_N$ is a graph having only clusters of size $A >0$ and denote by $ G(N,p) \cup \mathfrak{g}$ the graph obtained by superimposing it with an independent Erd{\H{o}}s--R\'enyi random graph (and deleting the possible multiple edges). Then we have:
\begin{lemma}[Sprinkling] \label{lem:sprinkling}Fix $ \delta, \varepsilon >0$ and $A\in \{1,2, \dots \}$.  The graph $G(N, \frac{\delta}{N}) \cup \mathfrak{g}$ has a giant component of size at least $(1 - 2\varepsilon)N$ with high probability as $N \to \infty$ as soon as we have  $$ - \delta \varepsilon(1- \varepsilon) +   \frac{\log 2 }{A} < 0.$$
\end{lemma}
\noindent \textbf{Proof.} As above, we compute the expected number of $ \varepsilon$-cuts in $G(N, \frac{\delta}{N}) \cup \mathfrak{g}$. Since those cuts have to be compatible with the initial structure of $ \mathfrak{g}$, there are at most $2^{K}$ choices where $K \leq N/A$ is the number of connected components of $ \mathfrak{g}$. Hence, the expected number of $  \varepsilon$-cuts is upper bounded by 
$$  \mathbb{E}\left[\# \varepsilon-\mbox{cuts in }G(N, {\textstyle \frac{\delta}{n}}) \cup \mathfrak{g} \right] \leq 2^{ N/A}  \left(1- \frac{\delta}{N}\right)^{ \varepsilon(1- \varepsilon) N^{2}},$$ and we conclude as above by the first moment method after summing over the at most $N$ possible values of $ \varepsilon$. \qed \medskip

\begin{exo} \label{exo:bissection} Suppose $n$ is even. Use Theorem \ref{thm:weak-giant} to prove that the existence of a $ \frac{1}{2}$-cut (i.e.~a partition of the vertices into two subsets of the same cardinality without edges between them) in $G(n, p)$ has a sharp threshold at $$p_{n} = \frac{\log 4}{n}.$$ 
\end{exo}

\section{Sprinkling}
We now gather Corollary \ref{cor:propotionfini} and Lemma \ref{lem:sprinkling} and prove Theorem \ref{thm:weak-giant}. The idea is to remark that the superimposition of two independent Erd{\H{o}}s--R\'enyi random graph is again an Erd{\H{o}}s--R\'enyi graph: for $c>1$ and $\delta>0$  we have that
  \begin{eqnarray} G\big(n, \frac{c}{n}\big) \underset{ \mathrm{indpt}}{\bigcup} G(n, \frac{\delta}{n}) \quad \overset{(d)}{=}\quad  G\Big(n, 1- \big( 1- \frac{c}{n}\big)\big( 1- \frac{\delta}{n}\big)\Big) \quad \subgraph \quad   G\big(n,  \frac{c + 2 \delta}{n}\big)   \label{eq:couplage2delta}\end{eqnarray} for $n$ large enough.  \medskip 
  

\noindent \textit{Proof of Theorem \ref{thm:weak-giant}.} Fix $c>1$, fix $ \varepsilon,\delta >0$ small  and $A>0$ large. Denote by $ \mathfrak{g}$ the subgraph of $G(n, { \textstyle \frac{c}{n}})$ spanned by the vertices in components of size larger than $A$. We know from  Corollary \ref{cor:propotionfini} that $\mathfrak{g}$ is of size $N  =  n \mathbb{P}( \#\mathrm{T}(c) \geq A) +o_{ \mathbb{P}}(n)$ and we assume that $A$  is large enough so that $$\mathbb{P}( \#\mathrm{T}(c) \geq A) \geq (1- \varepsilon) (1- \alpha(c)),$$ in particular we used here that $c>1$ so that $1-\alpha(c)>0$. Conditionally $G(n, { \textstyle \frac{c}{n}})$ and in particular on $ \mathfrak{g}$ and $N$, when $N \geq (1- \varepsilon)^{2}n$ we can apply Lemma \ref{lem:sprinkling} and deduce that in the graph $ G(n, {\textstyle \frac{\delta}{n}}) \cup \mathfrak{g}$ restricted to the vertices of $ \mathfrak{g}$, there is w.h.p.~a component of size at least $ (1- \varepsilon)N$ as soon as $$ \delta (1- \varepsilon)^{2} \varepsilon(1-  \varepsilon) - \log 2 /A >0.$$ Up to further increasing $A$, we can suppose that the former inequality is satisfied. We deduce that w.h.p.~there is a component of size $ (1-  \varepsilon)N \geq (1- \varepsilon)^{3}(1-\alpha(c)) + o_{ \mathbb{P}}(n)$ inside $G\big(n, \frac{c}{n}\big) {\bigcup} G(n, \frac{\delta}{n}) \subgraph    G\big(n,  \frac{c + 2 \delta}{n}\big)$ by \eqref{eq:couplage2delta}. Letting $A \to \infty$ while $ \varepsilon \to 0$, this shows the existence of a connected component of size \textit{at least} $(1-\alpha(c))n + o_{ \mathbb{P}}(n)$ in $ G(n, { \textstyle \frac{c'}{n}})$ for any $c'>c$. By continuity of $c \mapsto \alpha(c)$ we deduce the existence of a connected component of size \textit{at least} $(1 - \alpha(c))n + o_ \mathbb{P}(n)$ in $ G(n, { \textstyle \frac{c}{n}})$ whereas Corollary \ref{cor:propotionfini} entails that $ \alpha(c) n + o_{ \mathbb{P}}(n)$ of its vertices are in components of bounded size (irrespectively of $n$). This proves Theorem \ref{thm:weak-giant}.\qed \medskip 

\paragraph{Bibliographical notes.} The analysis of the phase transition for the emergence of the giant component is a classic in nowadays probability theory, see \cite{erdHos1960evolution} for the initial paper and \cite{janson1993birth} and \cite{aldous1997brownian} for a detailed analysis.  The proof of Section \ref{sec:firstmomentgiant} is directly inspired by the original proof of Erd{\H{o}}s and R\'enyi. The local limit paradigm is quite recent \cite{BS01} and has been a fruitful idea applied in the realm of random graphs, see \cite{AL07,BCstationary} for references. \medskip

\noindent{\textbf{Hints for Exercises.}}\ \\
Exercise \ref{exo:bissection}: At $p_{n} = \frac{\log 4}{n}$ the giant component in $G(n,p_{n})$ is of size $ n/2 + o_{ \mathbb{P}}(n)$. A sharper result is even proved in \cite{luczak2001bisecting}.\\


%

 \chapter{Birth of a giant $2$, exploration and fluid limit}

 We now turn to a more modern and powerful way of proving Theorem \ref{thm:erdosrenyi} based on  exploration techniques and stochastic analysis.  We define an exploration process of $G(n,p)$ which discovers its connected components one after the other in a Markovian way, by revealing its vertices one by one as $k=0,1,2, \dots , n$ and study  associated  $ \mathbb{R}$-valued Markov processes in the scaling limits. This will be the occasion to introduce the differential equation method, or fluid limit method whose applications are numerous.

 \section{Exploration process as a Markov chain}
\begin{figure}[!h]
 \begin{center}
 \includegraphics[width=15cm]{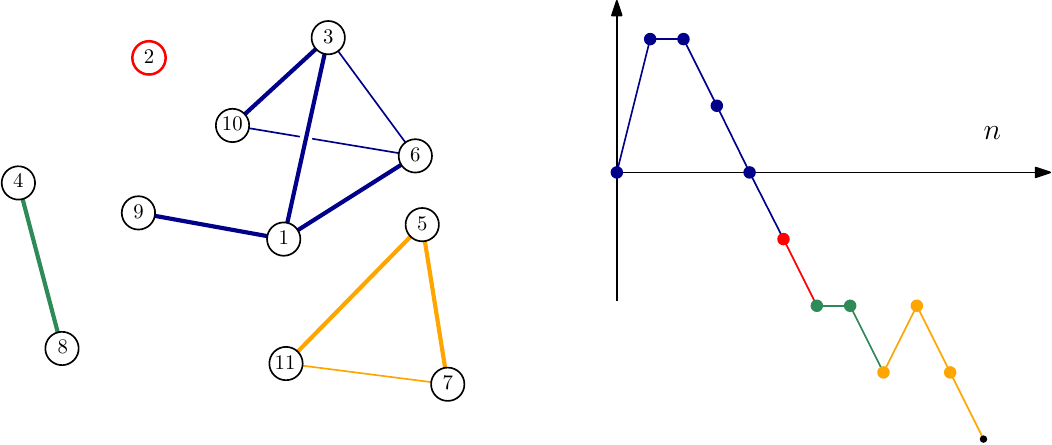}
 \caption{The {\L}ukasiewicz exploration of a random graph. The edges revealed during the exploration are in thick lines, they form spanning trees of each components. The concatenation (ordered by the minimal label of their component) of the {\L}ukasiewicz paths associated to those trees (explored by order of their labels)  is the {\L}ukasiewicz path of the graph.\label{fig:lukaER}}
 \end{center}
 \end{figure}

To properly define the exploration, we shall split the vertices $\{1,2, \dots ,n \}$ of $G(n,p)$ into three categories: the \textbf{untouched} vertices $ \mathcal{U}_{k}$, the  \textbf{explored}  vertices $ \mathcal{E}_{k}$ and the vertices in the current \textbf{stack} $ \mathcal{S}_{k}$ whose neighborhoods remain to be explored. The algorithm evolves as follows:
 \begin{itemize}
 \item at time $k=0$ we have $ \mathcal{E}_{0} = \varnothing$, the untouched vertices are $ \mathcal{U}_{0} = \{2,3,\dots\}$ and the only vertex in the stack is $ \mathcal{S}_{1}=\{1\}$.
 \item suppose $k \geq 0$ is given and such that $ \mathcal{S}_{k} \ne \varnothing$. We then select the vertex $x \in \mathcal{S}_{k}$ with minimal label (recall that the vertex set of $G(n,p)$ is $\{1,2, \dots, n\}$) and reveal all the neighbors $y_{1}, \dots , y_{j}$ of $x$ among $ \mathcal{U}_{k}$ (this could be an empty set!). We then put 
 $$ \mathcal{U}_{k+1} = \mathcal{U}_{k} \backslash \{ y_{1}, \dots , y_{j}\}, \quad \mathcal{S}_{k+1} = \big( \mathcal{S}_{k} \backslash \{x\} \big) \cup \{ y_{1}, \dots , y_{j}\}, \quad \mathcal{E}_{k+1} = \mathcal{E}_{k} \cup \{x\}.$$
 \item When the current stack is empty $ \mathcal{S}_{k} = \varnothing$ then the \textbf{first stage} of the  algorithm ends.
 \end{itemize}
It should be clear from the above exploration that at time $\tau_{1}$ when  the first stage ends, the set of explored vertices $ \mathcal{E}_{\tau_1}$ is precisely the connected component of $1$ in $ G(n,p)$. If the graph is not yet entirely discovered, we shall continue the exploration in the remaining graph starting from the vertex with minimal label and consider the \textbf{{\L}ukasiewicz path} $$ (\mathbb{S}_{k}: 0 \leq k \leq n)$$ obtained by starting from $0$ and whose increments are equal to the number of neighbors discovered in the untouched part minus $1$, see Figure \ref{fig:lukaER} for an illustration. In terms of the stack process $ \mathcal{S}_k$, this consists in immediately adding the vertex with minimal label yet untouched (as long as there are some untouched vertices left) when $ \mathcal{S}_k$ becomes empty (without performing a time step). In particular, the stack becomes empty if and only if the graph has been entirely explored.

 Note that the excursions above the running infimum of $ \mathbb{S}$ correspond to the explorations of the different connected components of the graph, and in particular, if we introduce the running infimum process $ \underline{\mathbb{S}}_{k} = \inf_{0 \leq j \leq k} \mathbb{S}_{j}$ then we can recover the size of the current stack $ \# \mathcal{S}_{k}$ as being   \begin{eqnarray} \label{eq:sizestack} \# \mathcal{S}_{k} = \mathbb{S}_{k} -  \underline{\mathbb{S}}_{k}+1, \quad \mbox{for } 0 \leq k \leq n-1.  \end{eqnarray}
   We shall denote by $ \mathcal{F}_k$ for $k=0,\dots , n$ the filtration generated by the first $k$ steps of this exploration. 
 
 \begin{proposition}[Markov property of the exploration] \label{prop:markov} For any $0 \leq k \leq n$, conditionally on $ (\mathcal{U}_{k}, \mathcal{E}_{k}, \mathcal{S}_k)$, each edge in $G(n,p)$ between $x$ and $y$ where $x, y \in \mathcal{U}_{k}$ or $x \in \mathcal{U}_{k}$ and $y \in \mathcal{S}_{k}$ is present independently with probability $p$. \end{proposition}
 \noindent \textbf{Proof.} Fix $k \geq 0$ and notice that given the status of the edges and vertices revealed by time $k$, one could deterministically change the status of all the edges between $ \mathcal{S}_{k}$ and $ \mathcal{U}_{k}$ or in-between vertices of $ \mathcal{U}_{k}$ and this would not have affected the exploration up to time $k$ (because these edges have not been explored by the algorithm). It is easy to see from this that those edges are indeed i.i.d.~present with probability $p$. \\
An alternative, and more ``algorithmic'' way to see this is to imagine that all the edges of the graph $ G(n,p)$ carry a question mark ``?'' which means that its status is currently unknown, present with probability $p$ and absent with probability $1-p$. When performing the exploration of the successive clusters, we reveal the status of certain edges (the question marks then disappear). The key point is to notice that since we are not allowed to use the randomness of unrevealed edges, at time $k$, conditionally on the past exploration, all the edges in question in the proposition still carry their  ``?'' and so they are i.i.d.~present with probability $p$ and absent otherwise.
\qed \bigskip

We deduce from the above that the process $  \mathbb{S}$ evolves in a Markovian fashion, if one also records its running infimum process:
 \begin{proposition} \label{prop:markovER}  For $0 \leq k \leq n-1$, conditionally on $ \mathbb{S}_{0}, \dots , \mathbb{S}_{k}$ the increment $ \Delta \mathbb{S}_{k}:=\mathbb{S}_{k+1}- \mathbb{S}_{k}$ is distributed as 
$$ \Delta \mathbb{S}_{k}  \overset{(d)}{=} \mathrm{Bin}(\#\mathcal{U}_{k}, p) -1 = \mathrm{Bin}(n-k- ( \mathbb{S}_{k} - \underline{\mathbb{S}}_{k} +1 ), p) -1.$$
\end{proposition}
\noindent \textbf{Proof.} This follows from the previous proposition, since the size of the stack $ \mathcal{S}_k =  \mathbb{S}_k - \underline{ \mathbb{S}}_k+1$ is given by \eqref{eq:sizestack}  and since the number of untouched vertices is $n-k$ minus the size of the current stack. \qed \medskip


\section{Differential equation method or fluid limit} \label{sec:fluidER}
\begin{figure}[!h]
 \begin{center}
 \includegraphics[width=12cm]{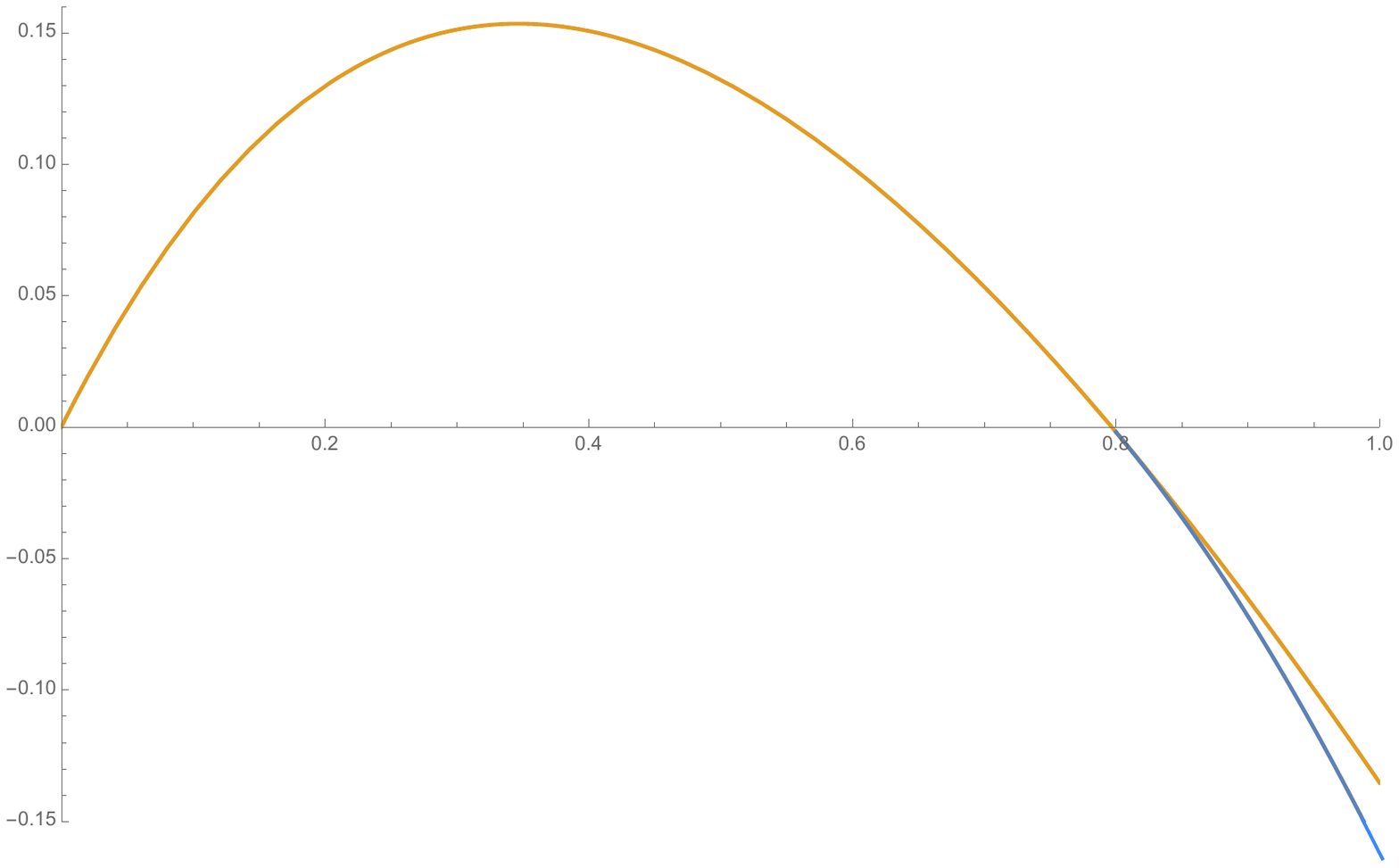}
 \caption{Plot of the function $f_2$: it follows the orange curve from $0$ to $1-\alpha(2) \approx  0.797$ and then the blue curve from $1-\alpha(2)$ to $1$. In particular, the function is not smooth at $t=1- \alpha(2)$.\label{fig:fc}}
 \end{center}
 \end{figure}

Fix $c >0$. In the rest of this section we take 
$$p = \frac{c}{n}.$$
 Taking expectations in Proposition \ref{prop:markovER},  according to a general principle that goes under the name of ``\textbf{fluid limit}'' or ``\textbf{differential equation method}'', we anticipate that the process $  \mathbb{S}$  behaves in the large scale limit as a deterministic function $f_{c}$  which satisfies the differential equation 
 \begin{eqnarray} \label{eq:diffinf} f'_{c}(t) = c\left(1-t - \big(f_c(t) - \underline{f_{c}}(t)\big)\right)-1,  \end{eqnarray}  and starts at $f_{c}(0)=0$, where we used the notation $\underline{g}(s) = \inf \{ g(u) : 0 \leq u \leq s\}$ for a continuous function $g$. This is not a standard differential equation due to the seemingly awkward dependence in the function $ \underline{f_{c}}$, but it is easy to convince oneself that the equation indeed has a unique solution and that this solution is either decreasing (when $c <1$) or  unimodal (when $c>1$). More precisely:\\
\noindent \textbf{For $\mathbf{c >1}$} the function $f_{c}$ is equivalently defined  as the solution to the differential equation 
$$ f_{c}'(t) = \left\{\begin{array}{lcl} c(1-t-f_{c}(t))-1 & \mbox{for } & 0 \leq t < \inf\{ s >0 : f_c(s)=0\} \\
c(1-t)-1 & \mbox{for } & \inf\{ s >0 : f_c(s)=0\} \leq t \leq 1. \end{array} \right.$$
In particular $f_{c}(t) =  1-\mathrm{e}^{-ct}-t$ until it comes back to $0$ at time $1- \alpha(c)$ where we recall from \eqref{def:alphac} that $\alpha\equiv \alpha(c)$ is the solution in $(0,1)$ to $ \alpha = \mathrm{e}^{-c(1-\alpha)}.$ For $t \geq (1- \alpha)$ the function follows the parabola  $$ f_c(t) = \frac{1}{2} (c (1 + \alpha - t)-2) (t-1+\alpha).$$ Although $f_c$ is $ \mathcal{C}^1$ over $[0,1]$, it is not $ \mathcal{C}^2$ at the point $1- \alpha(c)$ since its second derivative jumps from $-\alpha c^2$ to $-c$, see Figure \ref{fig:fc}.

\noindent\textbf{For $\mathbf{c\leq1}$}, since $f_c$ is decreasing we have $ \underline{f_{c}}(t) = f_c(t)$. It follows that  we have $\alpha=1$ and always have $f'_c(t)=c(1-t)-1$ so that  $$f_{c}(t) = \frac{t}{2}(2c-2-ct), \quad \forall t \in [0,1].$$

The above heuristic is indeed correct and we have:
\begin{theorem}[Fluid limit for the exploration]\noindent\label{thm:fluid}Fix $c>0$ and let $p= \frac{c}{n}$. Consider the {\L}ukasiewicz exploration $ (\mathbb{S}_{k} : 0 \leq k \leq n)$ of the random graph $G(n,p)$. Then we have the following convergence in probability 
$$ \big(n^{-1} \cdot \mathbb{S}_{\lfloor  nt \rfloor} : t \in [0,1] \big) \xrightarrow[n\to\infty]{} \big(f_{c}(t) : t \in [0,1] \big),$$
for the uniform norm.
\end{theorem}

\noindent \textbf{Proof.} Fix $c >0$.  The idea of the fluid limit theorem is to  argue that $ \mathbb{S}_k$ evolves  as a stochastic Euler's scheme based on the equation \eqref{eq:diffinf}. To be more precise, we shall compare $ \mathbb{S}_k $ with the discrete function $ \mathcal{L}_{k} = n f_c( \frac{k}{n})$ for $k=0,1,2, \dots , n$. We shall write $ \underline{ \mathbb{S}}_k:=\inf\{  \mathbb{S}_{j} : 0 \leq j \leq k \}$ and $ \underline{ \mathcal{L}}_{k} := \inf\{ \mathcal{L}_{j} : 0 \leq j \leq k \}$ for the running infimum processes  of $ \mathbb{S}$ and $  \mathcal{L}$ respectively. First, from  \eqref{eq:diffinf} and the fact that $f'_c$ is Lipschitz, it follows that we have the following Taylor approximation:
 \begin{eqnarray} \label{diffF}  \mathcal{L}_{k+1}- \mathcal{L}_{k} = \int_{k}^{k+1} \mathrm{d}s \, f'_c\left( \frac{s}{n}\right)  \underset{ \eqref{eq:diffinf}}{=}  c\left(1- \frac{k}{n} - \frac{ \mathcal{L}_{k}}{n} +  \frac{\underline{\mathcal{L}}_{k}}{n}\right)-1 +  \Theta(1/n),  \end{eqnarray}
where $\Theta(1/n)$ is a function bounded in absolute value by $ \mathrm{cst}/n$ independently of $0 \leq k \leq n$. We now analyse the process 
$$ X_k =  \mathcal{L}_{k}-\mathbb{S}_k.$$
Writing $ (\mathcal{F}_k : 0 \leq k \leq n)$ for the filtration generated by the exploration, we first compute the expected conditional increment of the process $X$:  \begin{eqnarray*}  \mathbb{E}[X_{k+1}-X_k \mid \mathcal{F}_k] &\underset{ \eqref{diffF} \ \& \  \mathrm{Prop.} \ref{prop:markovER}}{=}& c\left(1- \frac{k}{n} - \frac{ \mathcal{L}_{k}}{n} +  \frac{\underline{ \mathcal{L}}_{k}}{n}\right)-1 +  \Theta(1/n)\\
 & & - \left(   \mathbb{E}\left[ \mathrm{Bin}\left(n-k- \mathbb{S}_k +\underline{ \mathbb{S}}_k -1, \frac{c}{n}\right) \mid \mathcal{F}_{k}\right]-1\right)\\
 &=& c\left(\frac{\underline{ \mathcal{L}}_k- \underline{\mathbb{S}}_k}{n}-\frac{ \mathcal{L}_k- \mathbb{S}_k}{n}\right) + \Theta(1/n). 
 \end{eqnarray*} Remark that $| \underline{ \mathcal{L}}_{k}-\underline{\mathbb{S}}_k| \leq  \sup_{0 \leq i \leq k}| { \mathcal{L}}_{i}-{\mathbb{S}}_i|$ so that taking absolute values in the last display we deduce that  for all $ 0 \leq k \leq n-1$
 \begin{eqnarray*} \Big| \mathbb{E}[X_{k+1}-X_k \mid \mathcal{F}_k] \Big| & \leq &  \frac{C}{n}\left( 1 + \sup_{0 \leq i \leq k}|X_i|\right),  \end{eqnarray*}
 for some constant $C>0$. Furthermore, since the increments of $ \mathbb{S}$ are always stochastically dominated by $ \mathrm{Bin}(n, \frac{c}{n})$ it is plain to see that  up to increasing $C$ we have  
$$ \forall 0 \leq k \leq n-1, \quad \mathbb{E}[(X_{k+1}-X_k)^2] \leq C.$$
We are thus in position to apply the following ``stochastic'' version of Gronwall lemma to deduce that $ n^{-1}\sup_{0 \leq k \leq n}|X_k| \to 0$ in probability. This entails the theorem.
\subsection{Stochastic Gronwall lemma}

\begin{lemma}[Stochastic Gronwall lemma] Let $(X_k : 0 \leq k \leq n)$ be an adapted process with $X_0=0$. We define its supremum absolute value process $X_{k}^{*} = \sup \{ |X_{j}| : 0 \leq j \leq k \}$ for $0 \leq k \leq n$ and suppose that there exists $C >0$ satisfying for all $0 \leq k \leq n-1$
\begin{itemize}
\item $ |\mathbb{E}[X_{k+1}-X_k\mid \mathcal{F}_k]| \leq \frac{C}{n}\left(1+ X_k^{*}\right)$ almost surely,
\item $ \mathbb{E}[(X_{k+1}-X_k)^2] \leq C$.
\end{itemize}
Then we have 
$  n^{-1} \cdot X^{*}_{n} \to 0$ in probability as $n \to \infty$.
\end{lemma}
\noindent \textbf{Proof.} We decompose $X_k$ in its predictable and its martingale part by putting for $0 \leq k \leq n-1$
$$ X_{k+1}-X_k =  \underbrace{\mathbb{E}[X_{k+1}-X_k \mid \mathcal{F}_k]}_{=: \ D_k} + \underbrace{( (X_{k+1}-X_k) - \mathbb{E}[X_{k+1}-X_k \mid \mathcal{F}_k])}_{=: \ M_{k+1}-M_k},$$
so that if $M_0=0$ then $(M_k : 0 \leq k \leq n)$ is a martingale and    \begin{eqnarray} \label{eq:decompositionDoob}X_k = \sum_{i=0}^{k-1} D_i + M_k.  \end{eqnarray} Let us first take care of the martingale part: We have by (the conditional) Jensen's inequality
$$ \mathbb{E}[(M_{k+1}-M_k)^2] =  \mathbb{E}[\mathbb{E}[ \mathrm{Var}(X_{k+1}-X_k) \mid \mathcal{F}_k]] \leq   \mathbb{E}[\mathbb{E}[ (X_{k+1}-X_k)^2 \mid \mathcal{F}_k]] = \mathbb{E}[(X_{k+1}-X_k)^2].$$ Since the increments of a martingale are orthogonal  in $L^2$ by the above calculation we deduce that $ \mathbb{E}[M_n^2] \leq 4C n$.  By Doob's maximal inequality  we have $$ \mathbb{P}( \sup_{0 \leq k \leq n} |M_k| \geq A) \leq 4 \frac{\mathbb{E}[|M_n|^2]}{A^2}$$ and it follows that 
$$ 
\frac{\sup_{0 \leq k \leq n} |M_k|}{n} \xrightarrow[n\to\infty]{( \mathbb{P})} 0.$$ 
The rest of the argument is purely deterministic. By the hypothesis in the proposition we have $D_{i} \leq  \frac{C}{n}( X_{i}^{*}+1)$ and so  \eqref{eq:decompositionDoob} combined with the fact that the martingale part is negligible in front of $n$ yield that  for any $ \varepsilon >0$ on the event $\{ \sup_{0 \leq k \leq n}|M_{k}| \leq \varepsilon n\}$, whose probability tends to $1$ as $n \to \infty$, we have for all $ t \in [0,n]$
$$ X_{[t]}^{*} \underset{ \eqref{eq:decompositionDoob}}{\leq} \sum_{i=0}^{[t]-1}\frac{C}{n} X^{*}_{i} +  \left(\varepsilon + \frac{C}{n}\right) n \leq  2 \varepsilon n  + \frac{C}{n} \int_0^t \mathrm{d}s \, X^*_{[s]}.$$ On this event, by the usual (deterministic) Gr\"onwall\footnote{\raisebox{-5mm}{\includegraphics[width=1cm]{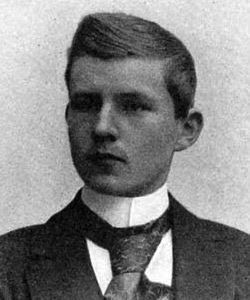}}Thomas Hakon Gr\"onwall  (1877--1932), Swedish} lemma we have  $ X^*_{[t]} \leq 2 \varepsilon n  \cdot \mathrm{exp}( \frac{C}{n} t),$ and in particular $ n^{-1} \cdot X^*_n \to 0$ as $n \to \infty$ in probability. \qed \bigskip 

The above strategy, called the ``differential equation method'' by Wormald \cite{wormald1995differential}, has found many applications in the realm of random graphs. Rather than giving an abstract convergence theorem, we propose to apply the strategy in the following exercise in order to estimate the size of the (random) greedy independent set on an Erd{\H{o}}s--R\'enyi random graph:
  \begin{exo}[Greedy independent set on $G(n,p)$] \label{exo:greedyER}Consider the graph $G(n,p)$ over the vertices $\{1,2, \dots , n\}$ for a parameter $p = \frac{c}{n}$ for some constant $c>0$. We will build inductively a random subset $  \mathbf{I}$ of $G(n,p)$ so that no vertices of $ \mathbf{I}$ are neighbors. To do this we put initially $ \mathcal{U}_{0}= \{1,2,3, \dots. , n\}$  (untouched) and $ \mathbf{I}_{0} = \varnothing$. Iteratively, for $i \geq 0$ as long as $ \mathcal{U}_{i} \ne \varnothing$ we select the vertex $x_{i}$ of smallest label in $  \mathcal{U}_{i}$ and denote its neighbors in $ G(n, p)$ by $\{y_{1}, \dots , y_{j}\}$ we then put 
  $$ \mathbf{I}_{i+1}= \mathbf{I}_{i} \cup \{ x_{i}\} \quad \mbox{ and }\quad \mathcal{U}_{i+1}= \mathcal{U}_{i} \backslash\{x_{i}, y_{1}, \dots , y_{j}\}.$$ We denote by $ \mathcal{F}_{n}$ the canonical filtration generated by this process and consider the stopping time $$\tau = \inf\{ k \geq 0 : \mathcal{U}_{k} = \varnothing\}.$$
\begin{enumerate}
\item Show that $ \mathbf{I}_{\tau}$ is an independent set, that is, no vertices of $ \mathbf{I}_{\tau}$ are neighbors.
\item Show that conditionally on $ \mathcal{F}_{k}$, the graph induced by $G(n,p)$ on $ \mathcal{U}_{k}$ is an Erd{\H{o}}s-R\'enyi random graph with parameter $p$. That is, all edges between vertices of $\mathcal{U}_{k}$ are independent and present with probability $p$.
\item Deduce that 
$$ \mathbb{E}[ \#\mathcal{U}_{k+1}- \# \mathcal{U}_{k} \mid \mathcal{F}_{k}] = - 1 - ( \# \mathcal{U}_{k}-1)p.$$
\item Recall that $p= p_{n} = \frac{c}{n}$. Use the differential equation method to prove that 
$$ (n^{-1} \# \mathcal{U}_{\lfloor nt \rfloor })_{t \in [0,1]} \to \left(\frac{(1+c-  \mathrm{e}^{{ct}}) \mathrm{e}^{{-ct}}}{c} \vee 0 : 0 \leq t \leq 1\right).$$
Hint : $f(t) = \frac{(1+c-  \mathrm{e}^{{ct}}) \mathrm{e}^{{-ct}}}{c}$ satisfies $f'(t) =-1-c f(t)$ and $f(0)=1$.
\item Deduce and explain why 
$$ n^{-1}\tau \xrightarrow[n\to\infty]{( \mathbb{P})}  \frac{\log(1+c)}{c}.$$
\end{enumerate}
\end{exo}

\section{Corollaries and refined estimates} \label{sec:fluidgiant1}
Let us deduce some geometric consequences of the convergence of the rescaled process $ \mathbb{S}$ towards $f_c$ (Theorem \ref{thm:fluid}). Recall from the definition of  the exploration of  the connected components of $G(n,p)$ that: 
\begin{enumerate}
\item The number of components in $G(n,p)$ is exactly $ -\underline{\mathbb{S}}_{n}$,
\item The sizes of the components in $G(n,p)$ correspond to the lengths of the excursions of $ \mathbb{S}$ above its running infimum $ \underline{\mathbb{S}}$. 
\end{enumerate}
As a direct corollary of Theorem \ref{thm:fluid} and the first item above we deduce:
\begin{corollary} The number  of components in $G(n,p)$ satisfies 
 \begin{eqnarray*} \frac{\# \mathrm{ConnComp}(G(n, {\textstyle \frac{c}{n}}))}{n} &\xrightarrow[n\to\infty]{ (\mathbb{P})}& f_c(1) = \frac{\alpha(c)( 2- c \alpha(c)) }{2}, \quad \mbox{ for }c>0.  \end{eqnarray*}
\end{corollary}
If $f : [0,1]  \to \mathbb{R}$ is a continuous function with $f(0)=0$ and running infimum process $\underline{f}(t) = \inf_{0 \leq s \leq t} f(s)$, we denote by $ \mathrm{Exc}(f)$ the (at most countably many) excursion intervals of $ f - \underline{f}$ away from $0$. We write $\|\mathrm{Exc}(f)\| \in \ell_1$ for the lengths of those excursions ranked in decreasing order.  We deduce the weak-giant property (Theorem \ref{thm:weak-giant}) from Exercise \ref{exo:continuityt} and the continuous mapping theorem. In particular, as in Section \ref{sec:firstmomentgiant}, we established the existence of the unique giant component but we did not give the logarithmic upper bounds for the second largest component stated in Theorem \ref{thm:erdosrenyi}. To prove it, we will use large deviations estimates in the next section.

\begin{exo}   \label{exo:continuityt} Consider the mapping 
$$ \mathcal{E} :  f \in \big(\mathcal{C}([0,1], \mathbb{R}), \|\cdot \|_\infty\big) \mapsto \| \mathrm{ Exc}(f)\| \in (\ell^1, \|\cdot \|_1).$$
\begin{enumerate}
\item Show that $ \mathcal{E}$ is not continuous in general.
\item However, show that $ \mathcal{E}$ is continuous at points $f$ where $f$ has no two-sided local minima. In particular, $ \mathcal{E}$ is continuous at $f= f_{c}$ for $ c \geq 0$.
\end{enumerate}
\end{exo}


In the rest of this section, we prove most of the refined estimates on the cluster size stated in Theorem \ref{thm:erdosrenyi}, especially the logarithmic upper bound in the subcritical case and for the second largest cluster in the supercritical case. We start with a stochastic domination of the typical cluster size coming from the exploration process. \medskip

By Proposition \ref{prop:markovER}, since we always have $ \# \mathcal{U}_{k} \leq n$, we deduce that the increments of $ \mathbb{S}$ are stochastically dominated by independent $ \mathrm{Bin}(n,{\textstyle \frac{c}{n}})$ variables minus $1$. We denote by $(S^{(n)}_{t} : t \geq 0)$  a random walk  starting from $0$ with i.i.d.~increments of law $ \mathrm{Bin}(n,{\textstyle \frac{c}{n}})-1$. In particular, the size of the cluster of $1$ in $G(n, {\textstyle \frac{c}{n}})$ is stochastically dominated by $T_{-1}( S^{(n)})$, the hitting time of $-1$ by that random walk. This can be evaluated via Kemperman's formula (Proposition \ref{prop:kemperman}) and an explicit computation:

 \begin{eqnarray} \mathbb{P}( T_{-1}( S^{(n)}) = k) &\underset{ \mathrm{Kemperman}}{=}& \frac{1}{k} \mathbb{P}( \mathrm{Bin}(k\cdot n,{\textstyle \frac{c}{n}}) = k-1) \nonumber \\
 &=& \frac{1}{k}{nk \choose k-1} \left(	 \frac{c}{n}\right)^{k-1} \left( 1-\frac{c}{n}\right)^{nk- (k-1)} \nonumber \\
 & \underset{ \mathrm{Stirling}}{\leq}& \mathrm{Cst} \cdot  \frac{n}{k} \frac{ \sqrt{nk}(nk)^{nk}}{ \sqrt{k} k^{k-1} \sqrt{nk} ((n-1)k)^{(n-1)k+1}}\left(\frac{c}{n} \left(1- \frac{c}{n}\right)^{n-1}\right)^k \nonumber \\
 &\leq &\mathrm{Cst} \cdot   \frac{1}{k^{3/2}} \cdot \frac{(nk)^{nk}}{k^k((n-1)k)^{(n-1)k}} \left(\frac{c}{n} \left(1- \frac{c}{n}\right)^{n-1}\right)^k\nonumber \\
  &\leq &\mathrm{Cst} \cdot   \frac{1}{k^{3/2}} \cdot \left(\frac{n^{n}}{(n-1)^{(n-1)}} \left(\frac{c}{n} \left(1- \frac{c}{n}\right)^{n-1}\right)\right)^k\nonumber \\
 & \leq&  \mathrm{Cst}\cdot  \frac{1}{k^{3/2}} \cdot  \left( c \left(\frac{n-c}{n-1}\right)^{n-1}\right)^k , \label{eq:tailexact}\end{eqnarray}
 for some  constant $ \mathrm{Cst}>0$ that may vary from line to line but which is independent of $k\geq 1$ and $n\geq 2$. When $c<1$ the term in the parenthesis tends to $ c \mathrm{e}^{1-c} < 1$ as $n \to \infty$, whereas for $c=1$ this term is equal to $1$.

\subsection{Subcritical case} \label{sec:logarithmicsubcritical}
We can now prove the logarithmic upper bound on the size of the clusters in the subcritical regime in Theorem \ref{thm:erdosrenyi}: Suppose that $c<1$ and recall that $ \mathcal{C}^i \equiv \mathcal{C}^i( G(n, \frac{c}{n}))$ is the cluster of the vertex $i \in \{1,2,3, \dots ,n\}$ inside $G(n,{\textstyle \frac{c}{n}})$. Using the above bound, we deduce that for any $A >0$, we have in $G(n,{\textstyle \frac{c}{n}})$:
$$ \mathbb{P}( \# \mathcal{C}^1 \geq A) \leq \mathbb{P}( T_{-1}( S^{(n)}) \geq A) \underset{ \eqref{eq:tailexact}}{\leq}  \mathrm{Cst} \cdot \eta^A,$$ where $\eta<1$ is independent of $A$. Taking $A = \frac{2}{\eta} \log n$, we deduce using the union bound that 
 \begin{eqnarray*} \mathbb{P}( \mathrm{C}_{1}^{\max}  \geq  \frac{2}{\eta} \log n) &\underset{ \mathrm{union \ bd}}{\leq}& n \mathbb{P}\left ( \#  \mathcal{C}^1  \geq \frac{2}{\eta} \log n\right) \\  &\leq& n \mathbb{P}( T_{-1}( S^{(n)}) \geq  \frac{2}{\eta} \log n) \leq  \mathrm{Cst}\cdot n \exp(- 2 \log n) \to 0.  \end{eqnarray*}
\subsection{Critical case} The same strategy can be used in the critical case $c=1$ (together with a little size-biasing trick). More precisely, imagine that we  pick (independently of $ G(n,{\textstyle \frac{1}{n}})$) a vertex $U_{n}$ uniformly in $\{1,2, \dots , n\}$. The size of the cluster of $U_{n}$ has the same law as that of the vertex $1$ and so is stochastically dominated by $T_{-1}(S^{(n)})$. We can thus write 
 \begin{eqnarray*}  \mathbb{P}( T_{-1}(S^{(n)}) \geq A) &\geq& \mathbb{P}( \# \mathcal{C}^{ U_{n}} \geq A)\\ & \geq &\mathbb{P}(  \#\mathcal{C}_{\max} \geq A \mbox{ and }U_{n} \in\mathcal{C}_{\max}) \geq  \frac{A}{n} \mathbb{P}(\mathrm{C}^{\max}_{1} \geq A).  \end{eqnarray*}
 Now taking $A = \lambda n^{2/3}$ and using  \eqref{eq:tailexact} where  the exponential factor disappears when $c=1$, we find 
 $$ \mathbb{P}( \mathrm{C}^{\max}_{1} \geq \lambda n^{2/3}) = O( \lambda^{-3/2}),$$
 which already gives the good order of magnitude of the largest cluster in $ G( n, {\textstyle \frac{1}{n}})$. Getting the full distributional convergence of $(n^{-2/3} \mathrm{C}^{\max}_{i})_{i \geq 1}$ requires to understand in much more details the exploration process. See the next chapter for such a result (Proposition \ref{prop:aldous}) in a slight variant of the Erd{\H{o}}s--R\'enyi random graph.
 
 \subsection{Supercritical case}
 For the supercritical case, we shall establish a common phenomenon in statistical physics: in the supercritical regime, the complement of the giant behaves as a \textbf{subcritical} system. See Exercise \ref{exo:superfini} for an instance of this phenomenon in BGW trees. \medskip

 Let $c >1$. By   Theorem \ref{thm:fluid}  we know that after the giant component of $G(n, {\textstyle \frac{c}{n}})$ has been explored we are left with a graph over $\approx \alpha(c)n$ vertices with edge density $ \frac{c}{n}$. This graph is close to being an Erd{\H{o}}s--R\'enyi:

 \begin{lemma} Conditionally on  $  \mathrm{C}^{\max}_{1}$ and on $ \mathrm{C}^{\max}_{1} >  \mathrm{C}^{\max}_{2}$,  the remaining graph\footnote{with vertices relabeled in increasing order} $ G(n,p) \backslash   	\mathcal{C}_{\max}$ has law $G( n - \mathrm{C}^{\max}_{1},p)$ conditioned on having clusters of size strictly less than $\mathrm{C}^{\max}_{1}$. \label{lem:condmax}
 \end{lemma}
 \noindent \textbf{Proof.} Fix a connected component $ \mathfrak{g}_{\max}$ on $\{1,2, \dots , n\}$ and a given graph $ \mathfrak{g}_{ \mathrm{rem}}$ on the remaining vertices so that no component of $ \mathfrak{g}_{ \mathrm{rem}}$ has a cluster  of size larger or equal to $  \# \mathrm{V}(\mathfrak{g}_{\max})$. Then we have (with a slight abuse of notation)
   \begin{eqnarray*} &&\mathbb{P}(\mathcal{C}^{\max} = \mathfrak{g}_{\max} \mbox{ and } G(n,p) \backslash  \mathcal{C}^{\max} = \mathfrak{g}_{ \mathrm{rem}})\\
    &=&  {p}^{ \# \mathrm{E}(\mathfrak{g}_{\max})} p^{\# \mathrm{E}(\mathfrak{g}_{ \mathrm{rem}})} (1-p)^{ {n \choose 2} - \# \mathrm{E}(\mathfrak{g}_{\max}) - \# \mathrm{E}(\mathfrak{g}_{ \mathrm{rem}})} \\
   &=& \mathbb{P}\Big(G\big(n - \# \mathrm{V}( \mathfrak{g}_{\max}),p\big) = \mathfrak{g}_{ \mathrm{rem}}\Big) \cdot   \mathrm{Cst}(  \mathfrak{g}_{\max}), \end{eqnarray*} where the constant $\mathrm{Cst}(  \mathfrak{g}_{\max})$ only depends on $ \mathfrak{g}_{\max}$ and not on $ \mathfrak{g}_{ \mathrm{rem}}$ as long as its components have size strictly less than $\# \mathrm{V}( \mathfrak{g}_{\max})$. This proves the lemma. \endproof \medskip

\begin{figure}[!h]
 \begin{center}
 \includegraphics[width=8cm]{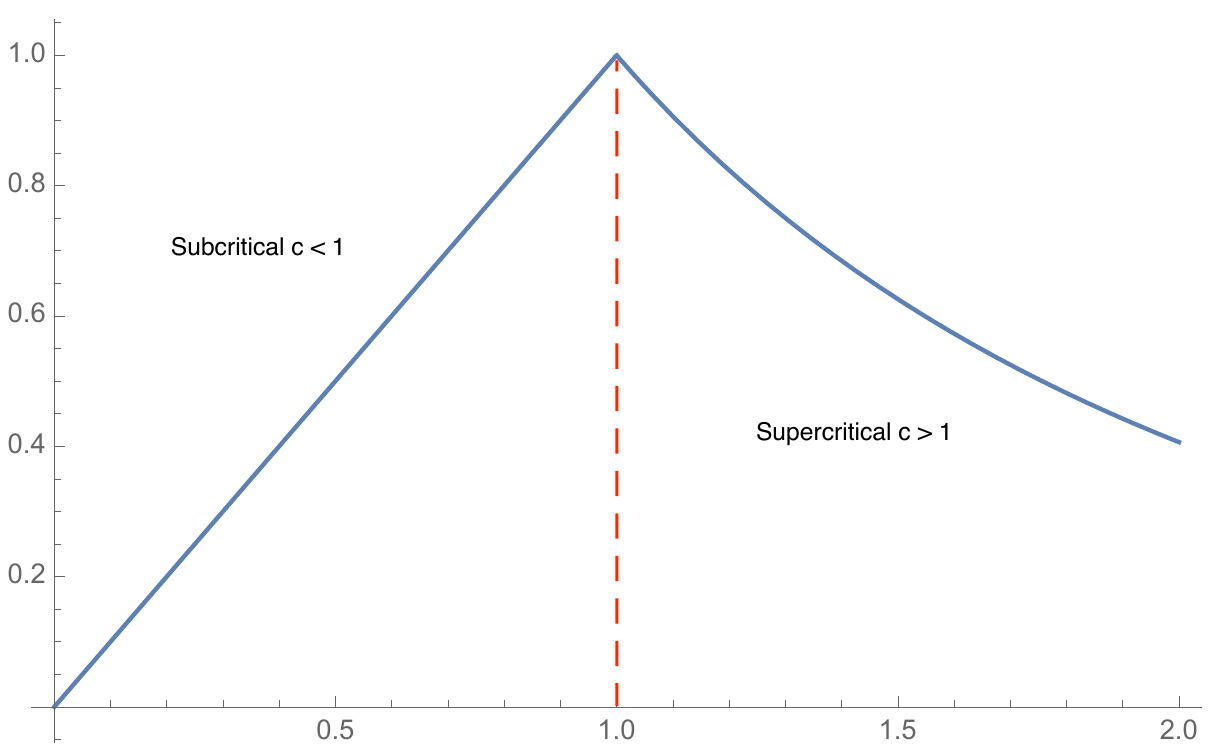}
 \caption{ \label{fig:calphac} Plot of $c \mapsto c\, \alpha(c)$ displaying the criticality in the remaining graph when the giant has been removed.}
 \end{center}
 \end{figure}

When $c >1$, notice that the complement of the giant component is \textbf{subcritical} since 
$$\forall c >0, \qquad  c \alpha(c)  <1,$$ see Figure \ref{fig:calphac}. 
 
We can thus prove point (ii) in Theorem \ref{thm:erdosrenyi}: Fix $c>1$ and $ \alpha (c)/2 > \varepsilon>0$. By Theorem \ref{thm:weak-giant}, the event $$ \{| \mathrm{C}^{\max}_{1} - (1-\alpha (c)) n| \leq \varepsilon n\} \cap \{  \mathrm{C}^{\max}_{2} < \varepsilon n \}$$ has a probability tending to $1$ and conditionally on it, the complement of the giant $G(n, \frac{c}{n}) \backslash  \mathcal{C}^{\max}$ is an Erd{\H{o}}s--R\'enyi with  $N \leq  ( \alpha(c) + \varepsilon)n$ vertices and edge density $ \frac{c}{n}$, conditioned on having no cluster of size larger than $ \mathrm{C}^{\max}_{1}$. If $ \varepsilon>0$ is small enough so that $c(\alpha(c) + \varepsilon)<1$, we know from Section \ref{sec:logarithmicsubcritical} that $G(N, \frac{c}{n})$ has no cluster of size larger than $A \log n$ for some $A >0$, so the previous conditioning does not affect its law asymptotically and we deduce (ii) in Theorem \ref{thm:erdosrenyi}.

  \begin{exo}[Supercritical BGW conditioned to be finite are subcritical BGW] \label{exo:superfini} Let $\mathcal{T}$ be a BGW tree with a \textit{supercritical} offspring distribution $\mu$ with generating function $g$. We denote by $\tilde{\mathcal{T}}$ the tree $ \mathcal{T}$ conditioned on the event $\{ \#\mathcal{T} < \infty\}$ whose probability is  equal to the unique solution $\alpha \in (0,1)$ to  $g(\alpha)= \alpha$, see Figure \ref{fig:GWclassic}. Show that $ \tilde{ \mathcal{T}}$ is a BGW with offspring distribution $\tilde{\mu}$ whose generating function $ \tilde{g}$ is given by 
  $$ \tilde{g}(z) = \frac{1}{\alpha} g (\alpha \cdot z), \quad \mbox{ for } z \in [0,1].$$
  \end{exo}

\paragraph{Bibliographical notes.} Although the exploration process of $G(n,p)$ is well known (see e.g.~\cite[(11.12)]{alon2016probabilistic} for Proposition \ref{prop:markovER}), the existence of the giant component using fluid limit for the exploration process seems to be new, although it is inspired by the much more precise analysis made in Aldous \cite{aldous1997brownian} and in \cite{nachmias2010critical}. More generally, the differential equation method has been used widely in random graph theory, see \cite{wormald1995differential}. Many formulations of the fluid limit paradigm, with increasing level of generality, can be found in the literature, see e.g.~ \cite{wormald1995differential,warnke2019wormald,darling2002fluid,darling2008differential}. Studying the emergence and the structure of the giant component in $G(n,p)$ is still a vivid subject in probability theory, see e.g.~\cite{ABBGM13} for very recent results relating the critical Erd{\H{o}}s-R\'enyi graph to the minimal spanning tree or \cite{schramm2005compositions} for a a connection with mixing time for the composition of transpositions on the symmetric group. We refer to \cite{van2009random} for extensions and more references.

\noindent{\textbf{Hints for Exercises.}}\ \\
Exercise \ref{exo:greedyER}: The result is first proved in \cite{pittel1982probable}.\\
Exercise \ref{exo:continuityt}: The mapping $ \mathcal{E}$ is not continuous at the function $f : x \mapsto |x||x-1/2||x-1|$ as limit of the functions $ f_{ \varepsilon} : x \mapsto |x|\big(|x-1/2| + \varepsilon\big) |x-1|$.\\
Exercise \ref{exo:superfini}. Use \eqref{exo:GWprod} and massage it. See Corollary 2.7 in \cite{abraham2015introduction} for a proof.

\chapter{Birth of giant $3$, Poissonized}

\label{chap:poissonER}

\hfill Pas frais mon poisson ? (Ordralfabétix)\bigskip

\begin{figure}[!h]
 \begin{center}
 \includegraphics[height=8cm]{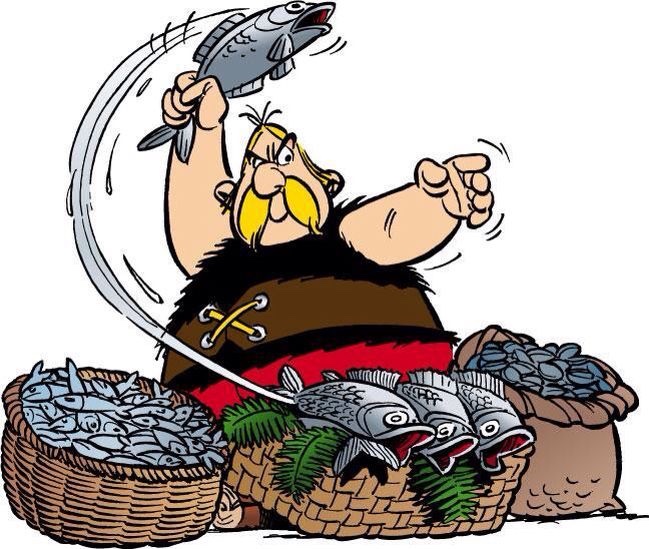}
 \caption{Ordralfabétix (\copyright\  Goscinny et Uderzo).}
 \end{center}
 \end{figure}

 We introduce a variant of the Erd{\H{o}}s--R\'enyi random graph where infinitely ``stack'' vertices are added on the side. A very simple Markov property of the  model entails that the {\L}ukasiewicz exploration is made of simple increments related to the repartition function of i.i.d.~uniforms. Using the standard Glivenko--Cantelli theorem,  this enables us to give very short proofs of classical results such as the phase transition for the giant component (Theorem \ref{thm:erdosrenyi}) or the connectedness for the standard Erd{\H{o}}s--R\'enyi model (Theorem \ref{thm:connectedness}).

\section{The stacked model and its exploration}
\begin{figure}[!h]
  \begin{center}
  \includegraphics[width=16cm]{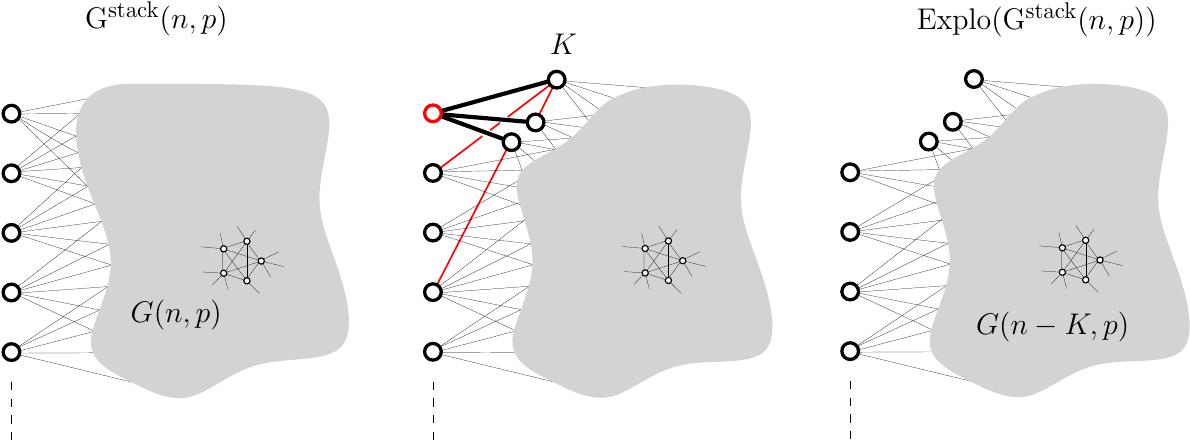}
  \caption{A stacked Erd{\H{o}}s--R\'enyi random graph and one step of exploration. The stack is made of the white vertices on the left part while the core is represented by the gray part.   After one step of exploration, the explored vertex (in red) is deleted as well as the edges linking the discovered vertices between each other or  to the former stack. Conditionally on the number $n'$ of vertices remaining in the core after this exploration, the resulting graph (after relabeling of its vertices) is distributed as $ G^{ \mathrm{stack}}(n',p)$. \label{fig:ERP}}
  \end{center}
  \end{figure}

We shall consider a variant of the Erd{\H{o}}s--R\'enyi model where we add infinitely many additional vertices ``in a stack on the side''. Formally, for $n \geq 1$ and $p$ fixed we consider the graph $ \mathrm{G}^{ \mathrm{stack}}(n,p)$ on the vertex set $\{1,2, \dots , n\} \cup \{1^{*},2^{*}, 3^{*}, \dots \}$, the vertices of $\{1,2, \dots , n \}$ form the \textbf{core} of the graph, whereas the vertices $\{1^{*}, 2^{*}, \dots \}$ form the \textbf{stack}. Then, each pair of core and stack vertices are connected by an edge with probability $p$ independently. \emph{There are no edges between vertices of the stack}. See Figure \ref{fig:ERP}. 

\paragraph{Markov property.} A step of \textbf{exploration} in $  \mathrm{G}^{ \mathrm{stack}}(n,p)$ is the following: Fix a vertex $\rho$ of the stack (independently of the core) and reveal its neighbors $ y_1, \dots , y_K$  with $K\geq 0$ inside the core. Then, see those vertices $y_1, \dots  , y_K$ as new vertices of the stack, in particular erase all possible edges between $y_1,\dots,y_K$ and between $y_1,\dots,y_K$ and other vertices of the stack. Denote by $ \mathrm{Explo}( \mathrm{G}^{\mathrm{stack}}(n,p))$ the resulting random graph whose vertices are relabeled by $\{1,2, \dots , n-K\}$ and $\{1^{*}, 2^{*}, \dots \}$ accordingly. The following is trivially verified:

\begin{lemma}[Markov property of $ \mathrm{G}^{\mathrm{stack}}(n,p)$] \label{lem:Markov}  \label{lem:markov} Let $K \geq 0$ be the number of neighbors in the core of $   \mathrm{G}^{\mathrm{stack}}(n,p)$ of the stack vertex $\rho$.  Then $ K \sim \mathrm{Bin}(n,p)$ and conditionally on $K$, we have the equality in law $$  \mathrm{Explo}( \mathrm{G}^{\mathrm{stack}}(n,p)) \quad \overset{(d)}{=}  \mathrm{G}^{\mathrm{stack}}(n-K,p).$$\end{lemma}

\begin{figure}[!h]
       \begin{center}
       \includegraphics[width=16cm]{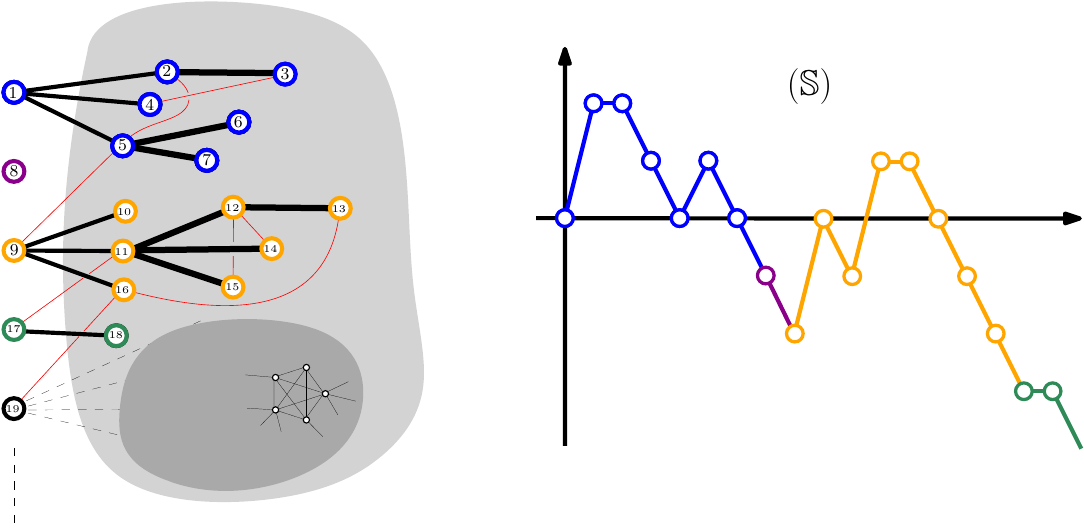}
       \caption{{\L}ukasiewicz exploration of the graph $ \mathrm{G^{stack}}(n,p)$: the numbering reflects the order in which the vertices have been explored. The thick edges are kept  whereas the thin red edges are discarded in the exploration. The thick (and very thick)  edges form $ \mathrm{F^{stack}}(n, p)$ and the very thick ones form $ \mathrm{F'^{stack}}(n, p)$ . The process $ \mathbb{S}^{(n, p)}$ on the right is  obtained by concatenating the successive number of neighbors $-1$.  \label{fig:lukaERpoi}}
       \end{center}
       \end{figure}

We shall now consider successive exploration steps and denote by $K \equiv K_{1}, K_{2}, \dots$ the number of vertices of the remaining core discovered at each step. In the rest of the chapter, we shall focus on a specific  exploration of the graph: we shall assume that iteratively, the discovered vertices are placed on top of the stack and that we successively explore the first vertex of the stack. We get the so-called \textbf{{\L}ukasiewicz exploration} of the graph $ \mathrm{G^{stack}}(n,p)$ similar to the one used in the previous chapter, see Figure \ref{fig:lukaERpoi}. We encode it in a process $$ (\mathbb{S}_{k}^{(n, p)} : k \geq 0)$$ or in short $(\mathbb{S}_{k} : k \geq 0)$, the {\L}ukasiewicz walk,  defined by $ \mathbb{S}^{(n, p)}_{0}=0$ and where $ \Delta \mathbb{S}^{(n, p)}_{i} = \mathbb{S}^{(n, p)}_{i}- \mathbb{S}^{(n, p)}_{i-1} = K_{i}-1$ is the number of neighbors discovered at step $i$ minus one.

              \paragraph{Relation to components.} Since $ \mathrm{G^{stack}}(n, p)$ has an infinite stack of vertices linked to each vertex of the core independently with probability $p$, as soon as $p >0$, the graph is a.s.~connected and in fact all vertices of the core have infinite degree almost surely. However, if we only consider the edges that are truly used in the {\L}ukasiewicz exploration (i.e.~not the edges between stack and revealed vertices, nor edges between revealed vertices) we obtain a spanning forest 
       $$ \mathrm{F^{ \mathrm{stack}}}(n,p) \subgraph \mathrm{G^{ \mathrm{stack}}}(n,p),$$ whose {\L}ukasiewciz walk is precisely $ \mathbb{S}$,  
       see Figure \ref{fig:lukaERpoi}. In particular, new minimal records of $ \mathbb{S}$ correspond to the discovery of a new tree component in $ \mathrm{F^{ \mathrm{stack}}}(n,p)$. If we further remove all vertices of the initial stack (together with the adjacent edges) we split $\mathrm{F^{ \mathrm{stack}}}(n, p)$ into a finer forest $\mathrm{F^{ \mathrm{stack},'}}(n,p)$ which spans the core and we can check the following graph inclusions          \begin{eqnarray} \label{eq:inclusion} \mathrm{F^{ \mathrm{stack},'}}(n,p) \subgraph \underbrace{G(n,p)}_{ \mathrm{Core}} \subgraph \mathrm{F^{ \mathrm{stack}}}(n,p) \subgraph  \mathrm{G^{ \mathrm{stack}}}(n,p).  \end{eqnarray}

\section{Law of the increments}
The advantage of the stacked version compared to the standard Erd{\H{o}}s--R\'enyi studied in the previous chapter is that the law of the increments of $ \mathbb{S}$ is simpler as it does not involved the running infimum process (compare Proposition \ref{prop:cid} with Proposition \ref{prop:markovER}). To make it even simpler, it is useful to randomize the size of the core. We first start with the description of $ (\Delta \mathbb{S}^{(n, p)}_{k} : k \geq 1)$ in the fixed-size case.

\subsection{Fixed size}
Consider the unit interval $[0,1)$ which is split in infinitely many subintervals
$$ [0,1) = \bigsqcup_{k \geq 1} \underbrace{\big[x^{(p)}_{k-1},x^{(p)}_{k}\big[}_{:= I^{(p)}_{k}}, \quad \mbox{ where } x^{(p)}_{k} = 1 - \left( 1-p\right)^{k}  \mbox{ for }k \geq 0,$$
so that for each $k \geq 1$, the length of $I^{(p)}_{k}$ is exactly $p$ times the total length of $I^{(p)}_{k},I^{(p)}_{k+1}, \dots$. We then throw $(U_{i} : 1 \leq i \leq n)$ independent identically distributed uniform r.v.~on $[0,1]$. The observation is:
\begin{lemma}  \label{prop:cid}The law of $(\Delta  \mathbb{S}^{(n, p)}_{k} +1 : k \geq 1)$ is equal to the law of $$ \big(\#\{ 1 \leq i \leq n : U_{i} \in I_{k}^{(p)}\}\big)_{k \geq 1}.$$
\end{lemma}
\noindent \textbf{Proof}. Denote by $ \tilde{K}_{j} = \#\{ 1 \leq i \leq n : U_{i} \in I_{j}^{(p)}\}$. Clearly $\tilde{K}_{1} \sim \mathrm{Bin}(n,p)$ in law. Furthermore, using the fact that the variables are uniform, we see that conditionally on $ \tilde{K}_{1}$, the sequence $ \tilde{K}_{2}, \tilde{K}_{3},\dots$ has the law of $(\tilde{K}_{1}, \tilde{K}_{2},\dots)$ where $n$ has been replaced by $n'=n- \tilde{K}_{1}$. Comparing with Lemma \ref{lem:markov} this suffices to prove equality of the laws recursively. \qed \medskip

If we write $F_{n}(x) = \# \{ 1\leq i \leq  n : U_{i} \leq x\}$ for the repartition function of the $n$ i.i.d.~uniforms, using the above proposition we can write \emph{simultaneously} for all $k \geq 0$
  \begin{eqnarray} \label{eq:walkcid}
   \mathbb{S}_{k} &=& F_{n}(x_{k}^{(p)}) -k.  \label{eq:couplagecid} \end{eqnarray}
   For our application, we recall the classical Glivenko-Cantelli
   \footnote{  \raisebox{-5mm}{\includegraphics[width=1cm]{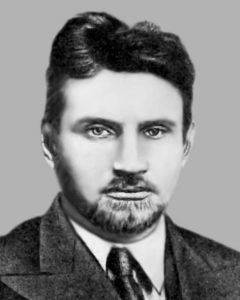}} Valery Ivanovich Glivenko (1897--1940), Ukrainian   \raisebox{-5mm}{\includegraphics[width=1cm]{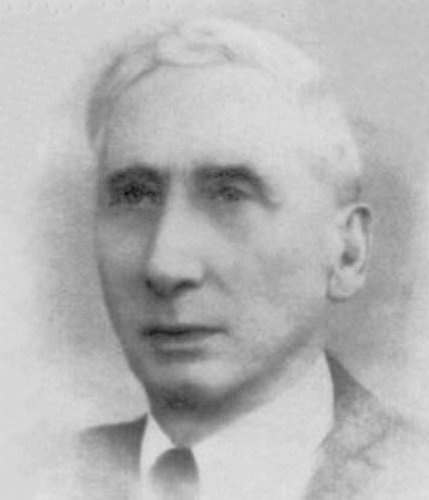}}Francesco Paolo Cantelli (1875 -- 1966)  Italian}
   theorem:
     \begin{eqnarray} \label{eq:glivenko} \left( \frac{F_{n}(x)}{n} : x \in [0,1]\right) & \xrightarrow[n\to\infty]{ (\mathbb{P})} & ( x : x \in [0,1]),  \end{eqnarray}
     for the $L^{\infty}$ metric. Before drawing probabilistic consequences of the above observations, let us consider the model where the size of the core is itself random which yield to further simplifications (and which gave the name to the chapter).

\subsection{Poissonized version} Fix $\alpha >0$ and suppose that $n\equiv N$ is first sampled at random, with law $\mathfrak{P}(\alpha)$ and conditionally on $N$ we perform the above construction. The resulting stacked graph will be denoted by $\mathrm{G_{Poi}^{stack}}(\alpha,p)$ and we denote the resulting {\L}ukasiewicz walk by $ \mathbb{S}^{[\alpha,p]}$.  By the classical Poisson thinning observation, in Lemma \ref{lem:markov} we then have $K \sim \mathrm{Bin} (N,p)  \sim \mathfrak{P}(\alpha p)$ and furthermore $K$ is independent of $N-K \sim \mathfrak{P}((1-p) \alpha)$. Iterating the above lemma, we deduce that in the Poissonized version the increments $ \Delta \mathbb{S}^{[\alpha,p]}_{k}+1$ of the {\L}ukasiewicz walk is now a sequence of \emph{independent} Poisson random variables with expectation $ \alpha p, \alpha p (1-p), \dots , \alpha p (1-p)^{k}, \dots$ whose total sum is just a Poisson variable of parameter $ \alpha p \sum_{i \geq 0} (1-p)^{i} = \alpha$, recovering the total number of vertices $N$ in the core as expected.\medskip 

As in \eqref{eq:walkcid} we can write in this case \emph{simultaneously} for all $k \geq 0$
  \begin{eqnarray}
   \mathbb{S}^{[\alpha,p]}_{k} &= &(\mathfrak{P}(\alpha p )-1) + ( \mathfrak{P}(\alpha p (1-p))-1) + \dots + ( \mathfrak{P}(\alpha p (1-p)^{k-1})-1)\nonumber \\ &=&   \mathfrak{P}\left( \alpha p \cdot \sum_{i=0}^{k-1} (1-p)^{i} \right)-k = \mathfrak{P}\left( \alpha (1-(1-p)^{k})\right) -k, \label{eq:lukapoisson} \end{eqnarray}
    where all the Poisson random variables written above are independent and where $ (\mathfrak{P}(t): t \geq 0)$ is a standard unit-rate Poisson counting process on $ \mathbb{R}_{+}$. We shall only use the following  standard limit theorems on the Poisson counting  process 
      \begin{eqnarray} \label{eq:lawll}
      \frac{ \mathfrak{P}(t)}{t} \xrightarrow[t\to\infty]{a.s.} 1, \quad \mbox{ and } \quad \frac{ ( \mathfrak{P}(tn)-tn)}{ \sqrt{n}} \xrightarrow[n\to\infty]{(d)} (B_{t} : t \geq 0),  \end{eqnarray} where $(B_{t} : t \geq 0)$ is a standard linear Brownian motion. The left-hand side follows from the law of large numbers and the right-hand side from Donsker's invariance principle.

\section{Phase transition for the giant}
Let us use the {\L}ukasiewicz exploration of the stacked version of the Erd{\H{o}}s--R\'enyi random graph to give a straightforward proof of Theorem \ref{thm:weak-giant}.

\subsection{Existence of the giant component}
Fix $c>0$. Let $p \equiv p_n = \frac{c}{n}$ and recall the notation   $ \mathbb{S}^{(n,p)}$ for the  {\L}ukasiewicz walk encoding the fixed size stacked Erd{\H{o}}s--R\'enyi random graph. Since we have $$x_{\lfloor nt \rfloor }^{(\frac{c}{n})} \sim \left(1- \frac{c}{n}\right)^{ \lfloor nt\rfloor } \to 1-\mathrm{e}^{-ct}, \quad \mbox{  as $n \to \infty$ uniformly over $ t \in \mathbb{R}_{+}$},$$ using  \eqref{eq:walkcid} and the Glivenko-Cantelli theorem \eqref{eq:glivenko}, we immediately deduce the analog of Theorem \ref{thm:fluid}:
\begin{proposition}[Fluid limit]\label{prop:fluid}We have the following convergences in probability
$$ \sup_{t \geq 0} \left\| \left(n^{-1} \cdot {\mathbb{S}^{(n, \frac{c}{n})}_{\lfloor nt\rfloor}}\right) -  \left( 1- \mathrm{e}^{{-ct}}-t\right)\right\| \xrightarrow[n\to\infty]{( \mathbb{P})}0.$$
\end{proposition}
  \begin{figure}[!h]
  \begin{center}
  \includegraphics[width=10cm]{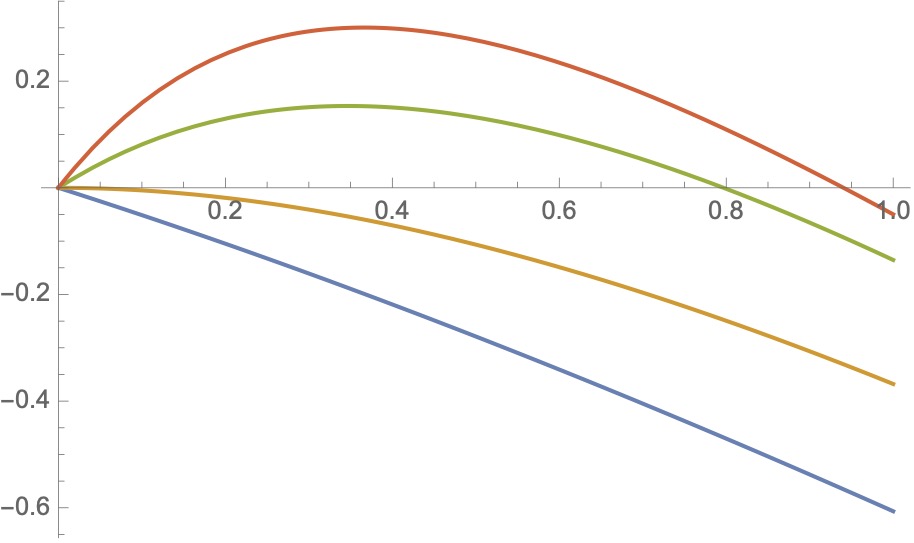}
  \caption{Graphs of the functions $( 1- \mathrm{e}^{{-ct}}-t)_{t \geq 0}$ for different of values of $c$: in blue $c=1/2$, in orange $c=1$, in green $c=2$ and in red $c=3$. Notice the root $1-\alpha(c)$ and compare with Figure \ref{fig:fc}.\label{fig:graphsfunctions}}
  \end{center}
  \end{figure}
  
Notice that the above convergence is not restricted to a compact time interval compared to Theorem \ref{thm:fluid}. However, when $c >1$, the function $ t \mapsto 1- \mathrm{e}^{-ct}-t$ coincides with the function $f_{c}$ defined in Section \ref{sec:fluidER} up to its first root at time $t = 1- \alpha(c)$, where we recall that $\alpha(c)$ is the smallest root to $ \alpha(c) = \mathrm{e}^{-c(1- \alpha(c))}$ and in particular $\alpha(c)=1$ if and only if $c \in [0,1]$. We give a  proof of the existence of the giant component in the Poissonized Erd{\H{o}}s--R\'enyi (that is Theorem \ref{thm:weak-giant}) using the same lines as in Section \ref{sec:fluidgiant1}:

 \begin{corollary}[Phase transition for $G(n, \frac{c}{n})$] \label{cor:giant} If $c <1$ then the largest connected components \textbf{in the core} of $ \mathrm{G^{stack}}(n,{\textstyle \frac{c}{n}})$ has size $o_{ \mathbb{P}}(n)$, whereas if $ c >1$ it contains a unique giant component of size $ (1-\alpha(c))n + o_ \mathbb{P}(n)$, and the second largest component has size  $o_{ \mathbb{P}}(n)$.  \end{corollary}

 \noindent \textbf{Proof.} Using the sandwiching of \eqref{eq:inclusion} it suffices to prove the similar statements for $ \mathrm{F^{stack}}$ and $  \mathrm{F^{stack,'}}$. The size of the connected components in $\mathrm{F^{stack}}(n, \frac{c}{n})$ are given by the lengths of the excursions of $ \mathbb{S}^{(n, \frac{c}{n})}$ above its running infimum process $$ \underline{ \mathbb{S}}^{(n, \frac{c}{n})}_k := \inf \{ \mathbb{S}^{(n,  \frac{c}{n})}_j : 0 \leq j \leq k\}.$$ We denote by $(L^{(n, \frac{c}{n})}_{i} : i \geq 1)$ those excursion lengths ranked in decreasing order. Notice that the excursion lengths above the running infimum of the function $ t \mapsto 1-  \mathrm{e}^{-ct}-t$ are given by $( 1-\alpha(c), 0, 0, \dots)$. Using Proposition \ref{prop:fluid} and (a variation on) Exercise \ref{exo:continuityt} shows that 
 $$\left( \frac{L^{(n, \frac{c}{n})}_{i}}{n} : i \geq 1\right) \xrightarrow[n\to\infty]{( \mathbb{P})} ( 1-\alpha(c), 0, 0, \dots)$$ for the $\ell^{\infty}$ norm. This proves the statement of the corollary for the random graph $\mathrm{F ^{ \mathrm{stack}}}(n, {\textstyle \frac{c}{n}})$. In the case $c \leq 1$, since $\mathrm{F^{stack,'}} \subgraph  \mathrm{F^{stack}}$ and $1-\alpha(c)=0$ there is nothing more to prove. However, when $c>1$ the removal of the initial stack vertices may split the giant component of $\mathrm{F^{stack}}(n, \frac{c}{n})$ of size $(1-\alpha(c))n +o_{ \mathbb{P}}(n)$ into several components but a moment's though using the {\L}ukasiewicz walk and  Proposition \ref{prop:fluid} again shows that one component of size $(1-\alpha(c))n + o_{ \mathbb{P}}(n)$ must remain.
\qed 

\subsection{Critical case}
In this section we turn to refined estimates on the cluster sizes  in the case $ \alpha=n$ and $p\equiv p_n =\frac{1}{n}$. For technical simplicity, we focus on the Poissonized version $ \mathrm{G^{stack}_{Poi}}$ for which we can use the Brownian limit in \eqref{eq:lawll}. This is an analog of point (iii) in Theorem \ref{thm:erdosrenyi} (where we take $\lambda=0$ below). Getting from those results the analogs for the fixed-size Erd{\H{o}}s--R\'enyi via depoissonization is doable, but is not covered in these pages.

\begin{proposition}[Near critical case] \label{prop:aldous} Fix $\lambda \in \mathbb{R}$. For $ p\equiv p_n = \frac{1}{n} + \frac{\lambda}{n^{{4/3}}}$ with $\lambda \in \mathbb{R}$, the {\L}ukasiewicz walk $ \mathbb{S}^{[n, \frac{1}{n} + \frac{\lambda}{n^{4/3}}]}$ of the Poissonized version satisfies 
$$  \left(n^{-1/3} \cdot {\mathbb{S}^{[n,\frac{1}{n} + \frac{\lambda}{n^{4/3}}]}_{\lfloor n^{2/3}t \rfloor }}\right)_{t \geq 0} \xrightarrow[n\to\infty]{(d)} \left( B_{t} + \lambda t - \frac{t^{2}}{2}\right)_{ t \geq 0},$$ where the convergence holds in distribution for the uniform norm over every compact of $ \mathbb{R}_{+}$.
\end{proposition}

\noindent \textbf{Proof.} Fix $A>0$. Putting $k = \lfloor n^{2/3}t\rfloor $ for $t \in [0,A]$ in the equation \eqref{eq:lukapoisson}, we have 
  \begin{eqnarray} \label{eq:dse}n \left(1-(1- \frac{1}{n}- \frac{\lambda}{n^{4/3}})^{\lfloor n^{2/3}t\rfloor}\right) =  tn^{2/3} + \lambda t n^{1/3} - \frac{t^{2}}{2} n^{1/3} + o(n^{1/3}),  \end{eqnarray} as $n \to \infty$ and where the little $o$ is uniform in $t \in [0,A]$. The second item of \eqref{eq:lawll} together with Skorokhod representation theorem show that on a common probability space we can build for each $m \geq 1$ a Poisson counting process $ \mathfrak{P}^{(m)}$ and a Brownian motion $B$ so that  we have the \emph{almost sure} convergence:
  \begin{eqnarray} \label{eq:skorokhod} \frac{ ( \mathfrak{P}^{(m)}(tm)-tm)}{ \sqrt{m}} \xrightarrow[m\to\infty]{a.s.} (B_{t} : t \geq 0)  \end{eqnarray}   for the uniform norm over every compact of $ \mathbb{R}_{+}$. Recalling  \eqref{eq:lukapoisson} those observations yield  for $m = \lfloor n^{2/3} \rfloor	$ 
 \begin{eqnarray*} \left(\frac{\mathbb{S}^{[n,p_{n}]}_{ \lfloor n^{2/3}t \rfloor }}{n^{1/3}} \right)_{0 \leq t \leq A}  & \overset{(d)}{\underset{ \mathrm{ for\  each\  }n}{=}}&  \left(\frac{\mathfrak{P}^{(m)}\left( n \left(1-\left(1- \frac{1}{n}- \frac{\lambda}{n^{4/3}}\right)^{ \lfloor n^{2/3}t \rfloor}\right)\right) - \lfloor n^{2/3}t \rfloor}{n^{1/3}}\right)_{0 \leq t \leq A}\\
& \underset{ \eqref{eq:dse}}{=}&  \left(\frac{\mathfrak{P}^{(m)}\left(  tm + \lambda t  \sqrt{m} - \frac{t^{2}}{2}  \sqrt{m} + o( \sqrt{m})) \right)- t m +o( \sqrt{m})}{ \sqrt{m} + o(1)}\right)_{0 \leq t \leq A}  \\
 &\underset{ \eqref{eq:skorokhod}}{\xrightarrow[n\to\infty]{a.s.}}&   \left( B_{t} + \lambda t - \frac{t^{2}}{2} \right)_{0 \leq t \leq A},  \end{eqnarray*} and this proves the proposition. \qed

  \section{Connectedness}
  As another application of our modification of the Erd{\H{o}}s--R\'enyi random graph,  let us give a short proof of the (very) sharp phase transition for connectedness in the fixed-size Erd{\H{o}}s--R\'enyi which is mentioned in \eqref{eq:erconnectedpoisson}:
  \begin{theorem}[Critical window for connectedness \cite{erdds1959random}]   For $c \in \mathbb{R}$ we have \label{prop:connectedness}
  $$ \mathbb{P}\left( G\left(n, \frac{\log n +c}{n}\right) \mbox{ is connected}\right) \xrightarrow[n\to\infty]{} \mathrm{e}^{- \mathrm{e}^{-c}}.$$
  \end{theorem}
  \noindent \textbf{Proof.} Let $p\equiv p_n = \frac{\log n +c}{n}$.  Connectedness of the core $ G(n,p_{n})$ is equivalent to the fact  that $ \mathrm{F^{stack}}(n,p_{n})$ has only one non-trivial component (the others being isolated vertices of the stack), or equivalently that the {\L}ukasiewicz walk $( \mathbb{S}^{(n,p_{n})})$ starts with a (large) excursion and once it has reached level $-1$, it makes only jumps of $-1$ forever. That is, $ \mathbb{S}^{(n,p_{n})}_{n+1} =-1$ and time $n+1$ is the first hitting time of $-1$. In particular, in the notation \eqref{eq:walkcid} we must have  $F_{n}( x^{(p_{n})}_{n+1}) = n$ or equivalently, that no uniform $U_{i}$ for $1 \leq i \leq n$ falls after the point 
  $$ x^{(p_{n})}_{n} = 1 - \left( 1- \frac{\log n +c}{n}\right)^{n} \sim \frac{ \mathrm{e}^{-c}}{n}, \quad \mbox{ as }n \to \infty.$$
  Computing this probability is routine and we have 
  $$ \mathbb{P}\left( \max_{1 \leq i \leq n} U_{i} \leq x^{(p_{n})}_{n}\right) \sim \left(1- \frac{ \mathrm{e}^{-c}}{n}\right)^{n} \xrightarrow[n\to\infty]{} \mathrm{e}^{- \mathrm{e}^{-c}}.$$
  To finish the proof, one shows that as long as this event is realized, then the core is connected with high probability. In term of the {\L}ukasiewicz walk this boils down to:
  \begin{lemma} For $p\equiv p_n = \frac{\log n +c}{n} $ we have 
 
  $$ \mathbb{P}\left( \mathbb{S}^{(n,p_{n})}_{k} \geq 0 : \forall 1 \leq k \leq n \mid F_{n}(x^{(p_{n})}_{n})=n\right)  \xrightarrow[n\to\infty]{}1.$$
  \end{lemma}
  \noindent \textbf{Proof of the lemma.} 
  Notice that the event on which we are conditioning is of asymptotically positive probability, so it suffices to shows that $\mathbb{P}( \exists 1 \leq k \leq n-1: \mathbb{S}^{(n,p_{n})}_{k}=0 \mbox{ and }\mathbb{S}^{(n,p_{n})}_{n} = 0)$ tends to $0$. We perform a union bound over all such $k's$ and compute
  \begin{eqnarray*} && \mathbb{P}( \exists 1 \leq k \leq n-1 : \mathbb{S}^{(n,p_{n})}_{k} = 0 \ \& \  \mathbb{S}^{(n,p_{n})}_{n} = 0) \\
  &\leq& \sum_{k=1}^{n-1}\mathbb{P}\left(  \begin{array}{c} \# \{ 1 \leq i \leq n : U_{i} \in I_{1}\cup I_{2} \cup \dots \cup I_{k}\} =k \\ \# \{ 1 \leq i \leq n : U_{i} \in I_{k+1}\cup I_{2} \cup \dots \cup I_{n}\} = n-k \end{array}\right) \\
&\leq &   \sum_{k=1}^{n/2}  \mathbb{P}( \mathrm{Bin}(n, x^{(p_{n})}_{k}) \leq k) + \sum_{k=n/2}^{n-1} \mathbb{P}( \mathrm{Bin}(n,x^{(p_{n})}_{n} - x^{(p_{n})}_{k}) \geq n- k). \end{eqnarray*}
For $ \varepsilon_{n}, \delta_{n}$ tending to $0$ such that $ n \varepsilon_{n} \to \infty$ as $n \to \infty$ we use the bound 
$$ \mathbb{P}( \mathrm{Bin}( n, \varepsilon_{n}) \leq \delta_{n} \varepsilon_{n}) \leq  \mathrm{e}^{{} - \mathrm{c} \ n \varepsilon_{n}} \quad \mbox{ and } \quad 
 \mathbb{P}( \mathrm{Bin}( n, \varepsilon_{n}) \geq n-\delta_{n} \varepsilon_{n}) \leq  \mathrm{e}^{- \mathrm{c} \ n \varepsilon_{n}},$$
 for some $c>0$. Since for $k \leq 10\frac{ n}{\log n}$ we have $k = o (n x^{(p_{n})}_{k})$ we can apply the above bound and get for some $ \mathrm{c}'>0$
 $$\sum_{k=1}^{10 n/\log n}  \mathbb{P}( \mathrm{Bin}(n, x^{(p_{n})}_{k}) \leq k) \leq \sum_{k =1}^{10 n/ \log n} \exp( - \mathrm{c}\ n x^{(p_{n})}_{k}) \leq \sum_{k =1}^{10 n/ \log n} \exp( - \mathrm{c'}\ k \log n) = o(1).$$ The case when $ 10\frac{ n}{\log n} \leq k \leq n/2$ is even easier since we have 
  \begin{eqnarray*}\mathbb{P}( \mathrm{Bin}(n, x^{(p_{n})}_{k}) \leq k) &\leq& \mathbb{P}\left( \mathrm{Bin}(n, x^{(p_{n})}_{10 \frac{n}{ \log n}}) \leq  \frac{n}{2}\right)\\
  & \leq & \mathbb{P}( \mathrm{Bin}(n, 1- 2\mathrm{e}^{-10}) \leq n/2) \leq \mathrm{e}^{-{c}'' n},  \end{eqnarray*}for some $c''>0$ by a large deviation estimate since $1- 2 \mathrm{e}^{-10}> 1/2$. Summing-up those estimates we deduce that 	$ \sum_{k=1}^{n/2}  \mathbb{P}( \mathrm{Bin}(n, x^{(p_{n})}_{k}) \leq k) \to 0$ as $n \to \infty$. A similar reasoning shows that $\sum_{k=n/2}^{n-1} \mathbb{P}( \mathrm{Bin}(n,x^{(p_{n})}_{n} - x^{(p_{n})}_{k}) \geq n- k) \to 0$ as well, and we leave the verification as an exercise for the reader. \qed 
  \bigskip

\noindent \textbf{Bibliographical notes.} The content of this chapter is adapted from the author's paper \cite{curien2022erd} and from the master's thesis of Damian Cid (promotion 2023-2024) who elegantly depoissonized the initial arguments. Various modifications of the Erd{\H{o}}s--R\'enyi random graph with nicer probabilistic properties have been used in the literature, see e.g.~the Poisson cloning model \cite{kim2007poisson}.

\part[Random tree growth]{Random tree growth 
                             \\ \\ 
  \begin{center}
                     \begin{minipage}[l]{15cm}
       \normalsize 
       In this part, we study several models of  random growing trees where vertices are attached to the preceding structure according to some rule. The prototype is the random recursive tree process $(T_{n} : n \geq 0)$ where $T_{n+1}$ is obtained  from $T_{n}$ by attaching a new vertex labeled $n+1$ onto a uniform vertex of $T_{n}$. We will study this process both from a static point of view (statistics of uniform random permutations), and from a dynamical process as $n$ increases (Polya urns and continuous time embedding).
                     \end{minipage}
                  \end{center}
                  \vspace{1cm}
                   \begin{figure}[!h]
 \begin{center}
 \includegraphics[height=5cm]{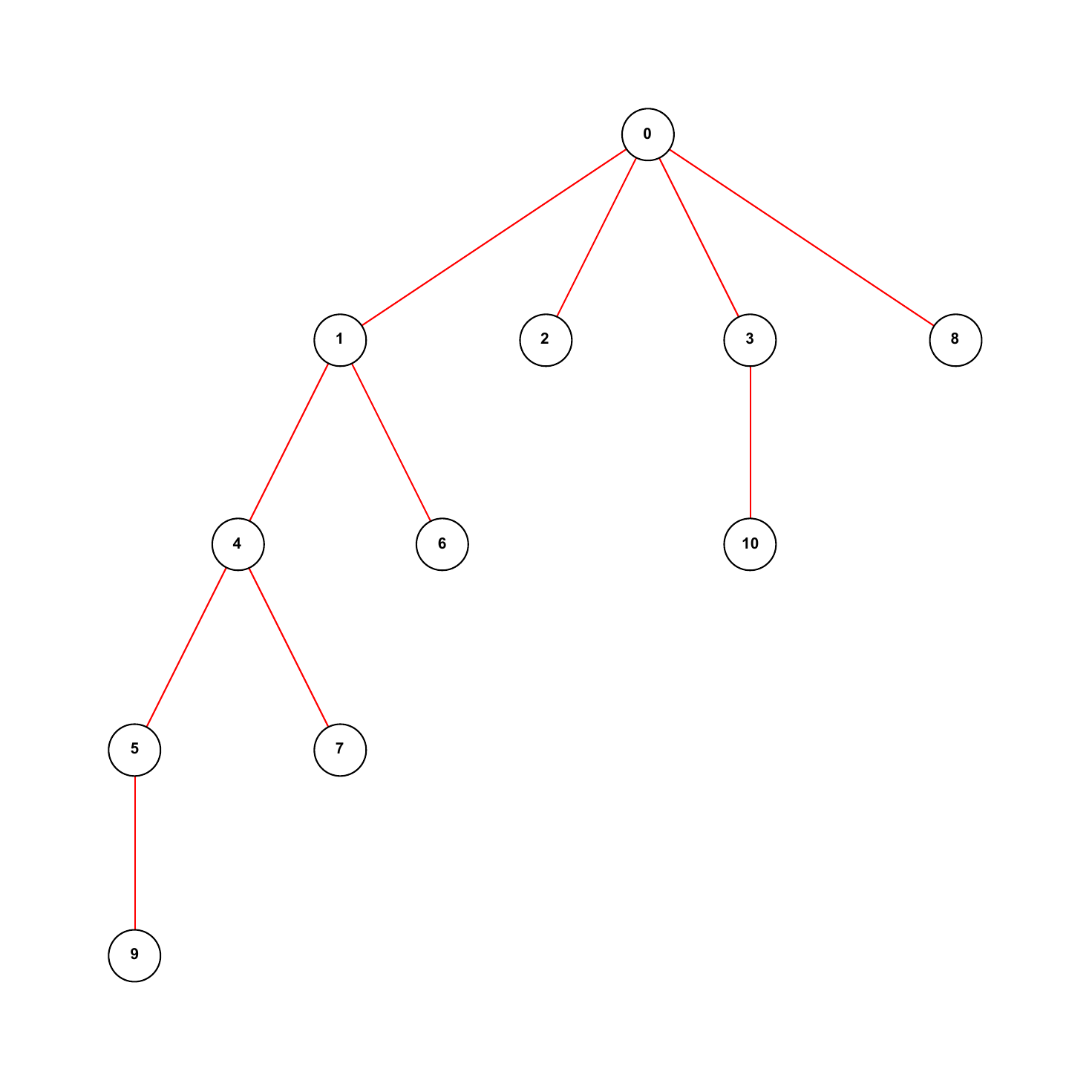}
  \includegraphics[height=5cm]{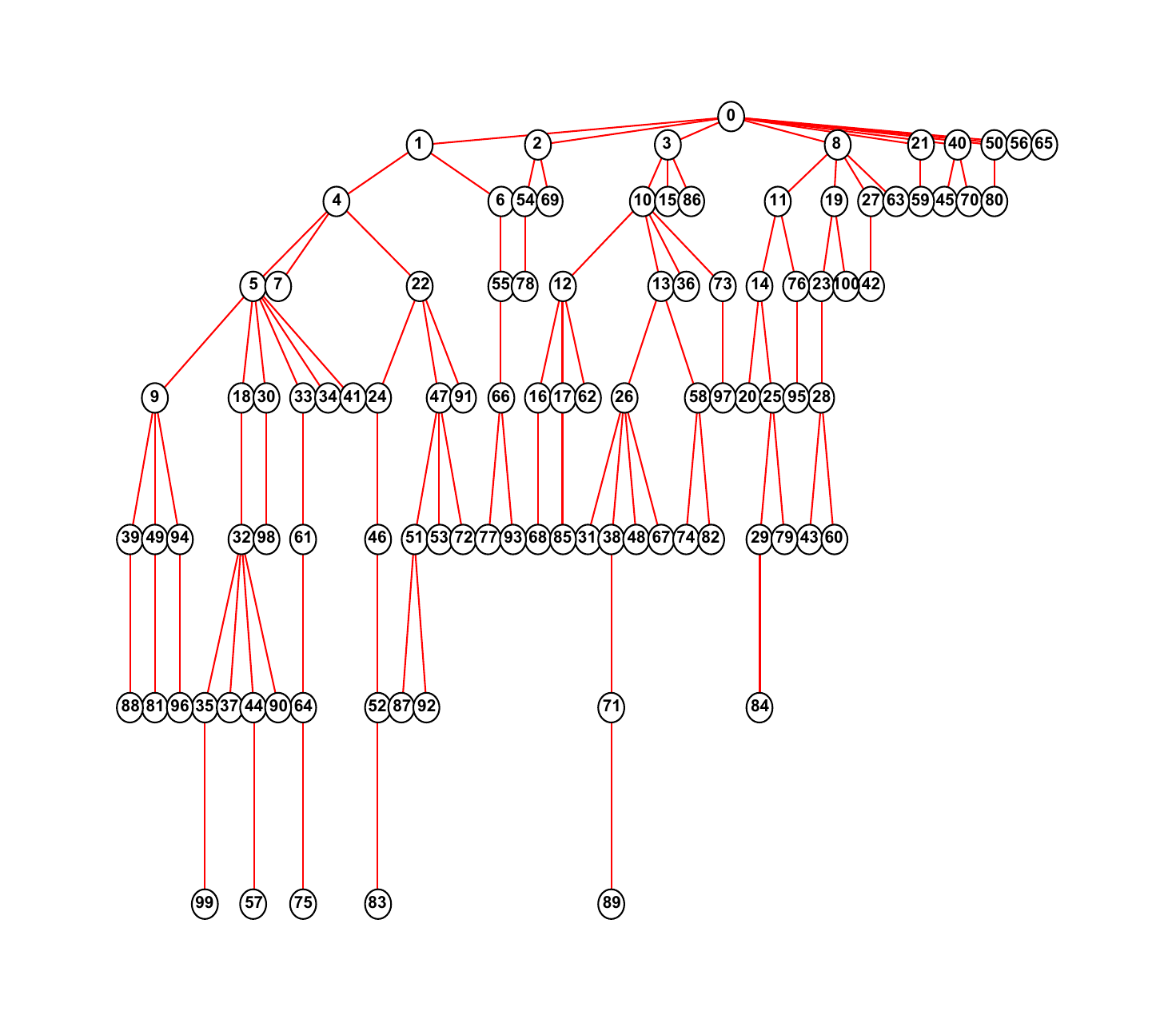}\\
  \includegraphics[width=14.5cm,height=5cm]{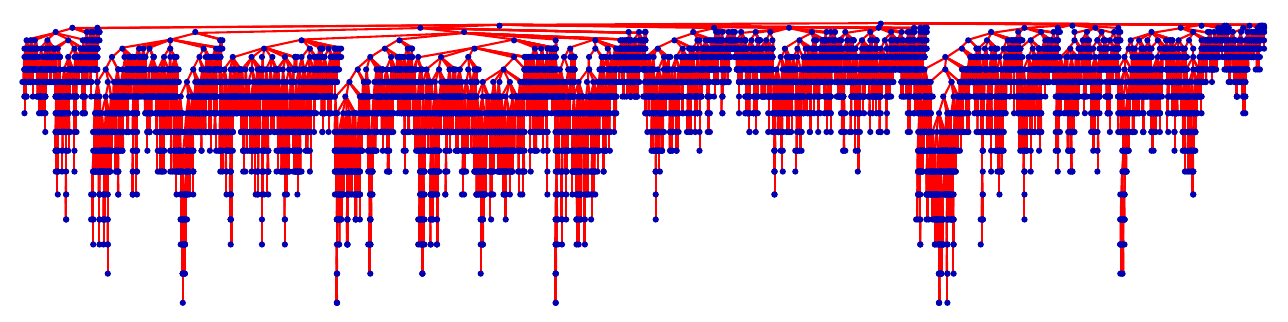}
 \caption{A random recursive tree at stages $10, 100$ and $10000$.}
 \end{center}
 \end{figure}
                }

\chapter{Random permutations}
\hfill Many points of view on $n!$\bigskip
\label{chap:permu}

In this chapter, we study the law of the cycle decomposition of a random  permutation $\boldsymbol{\sigma}_{n}$ chosen uniformly in the symmetric group $ \mathfrak{S}_{n}$ over $n$ elements $\{1,2, \dots , n\}$. In particular, we shall establish Poisson statistics for the number of shorts cycles and the Poisson--Dirichlet limit for the large cycles.  

\section{Feller coupling}
In 1945, Feller (the author of Lemma \ref{lem:feller}) introduced a coupling between the cycle structure of a uniform permutation $\boldsymbol{\sigma}_n$  and the spacings between successes in a sequence of $n$ independent Bernoulli variables of parameters $ \frac{1}{n}$, $ \frac{1}{n-1}$, \dots $  \frac{1}{2},1$. This will be the main tool used in this chapter. The key idea of this representation is to explore a given permutation along its cycles ordered by their minimal element. A concept which is sometimes called the \textbf{Foata\footnote{  \raisebox{-5mm}{\includegraphics[width=1cm]{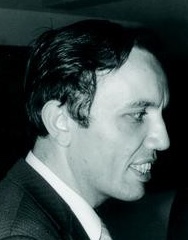}} Dominique Foata (1934--), French} correspondence}. 

\subsection{Foata correspondence}
\begin{figure}[!h]
 \begin{center}
 \includegraphics[width=13cm]{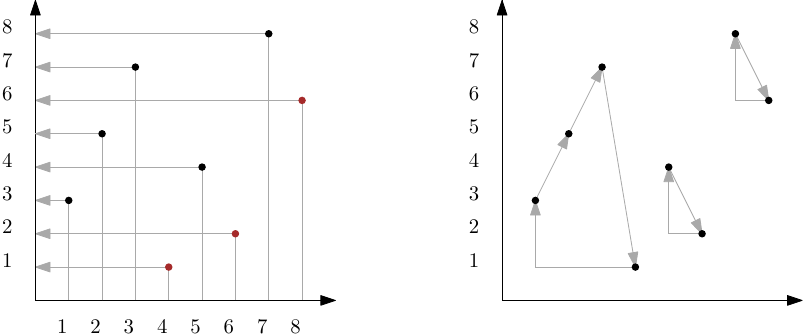}
 \caption{Foata correspondence: on the left a description of a permutation via its images, on the right the description of a permutation by exploration of its cycles ranked in increasing order of their minimal element. This bijection transforms the number of cycles into the number of minimal records (in red on the left). \label{fig:foata}}
 \end{center}
 \end{figure}

A permutation $\sigma_{n} \in \mathfrak{S}_{n}$ can obviously be described by a sequence  $(i_{1}, i_{2}, \dots , i_{n})$ representing the $n$ values $\{1,2, \dots,n\}$, the most obvious way is to prescribe the permutation by its values $ \sigma_{n}(1) = i_{1}, \sigma_{n}(2) = i_{2}, \dots , \sigma_{n}(n)= i_{n}$. Yet another way is to imagine that $(i_{1}, i_{2}, \dots , i_{n})$ is the sequence of values we discover when exploring the cycles of $\sigma_{n}$ ordered by their minimal values, see Figure \ref{fig:foata}. Specifically, let us denote by 
 \begin{eqnarray} \label{eq:foata} \big(a_1^{(1)}, \dots , a_{k_1}^{(1)}\big)\big(a_1^{(2)}, \dots , a_{k_2}^{(2)}\big) \cdots \big(a_1^{(\ell)}, \dots , a_{k_\ell}^{(\ell)}\big), \end{eqnarray} the decomposition of $\sigma_n$ into $\ell$ cycles with disjoint supports of length $k_1, \dots , k_\ell \geq 1$.  We suppose that those cycles are ranked according to their minimal element, which is placed at the end of each cycle in this representation:
$$ 1=a_{k_1}^{(1)} = \min_{1 \leq i \leq k_1} a_i^{(1)}  < a_{k_2}^{(2)} = \min_{1 \leq i \leq k_2} a_i^{(2)} < \cdots < a_{k_\ell}^{(\ell)} = \min_{1 \leq i \leq k_\ell} a_i^{(\ell)}.$$
Then, the Foata encoding of $\sigma_n$ is the permutation we obtain by reading the numbers in  \eqref{eq:foata} from left to right, namely 
$$ \mathrm{Foata}(\sigma_n) = \big(a_1^{(1)}, \dots , a_{k_1}^{(1)},a_1^{(2)}, \dots , a_{k_2}^{(2)}, \cdots ,a_1^{(\ell)}, \dots , a_{k_\ell}^{(\ell)}\big).$$ It is then clear that $ \mathrm{Foata} : \mathfrak{S}_n
 \to \mathfrak{S}_n$ is a bijection and furthermore that the number of cycles of $\sigma_n$ is equal to the number of minimal records of $ \mathrm{Foata}(\sigma_n)$, i.e.~the values $k$ such that $ \mathrm{Foata}(\sigma_n)_k = \min \{ \mathrm{Foata}(\sigma_n)_i : i \geq k \}$ and that the length of the cycles correspond to the spacing between those records.

 \begin{exo}[Law of a typical cycle]  \label{exo:taillecycle}Show using the Foata correspondence that the size of the cycle containing $1$ in a uniform permutation is uniform on $\{1,2, \dots , n \}$.
\end{exo}

 \subsection{Feller coupling} \label{sec:fellercoupling}
 Keeping in mind the Foata encoding of a permutation, we now present the famous result of Feller. We consider $n$ independent Bernoulli variables $ \mathrm{Ber}( \frac{1}{k})$ of success parameters  $$ \frac{1}{n}; \quad  \frac{1}{n-1}; \quad \dots \quad  \frac{1}{2}; \quad \frac{1}{1}.$$ Denote by $ n \geq I_1 > \dots > I_\ell = 1$ the indices (the reciprocal of the parameter) of the variables equal to $1$ and consider the $\ell$ spacings $ \mathcal{S}_{n} =  ((n+1)-I_1 , I_1-I_2, \dots , I_{\ell-1}-I_\ell)$ between the points $n+1 > I_1 > \dots > I_\ell$. The sum of those spacings is equal to  $n$. 
 
 \clearpage 
 
  \begin{figure}[h]
 \begin{center}
 \includegraphics[width=14cm]{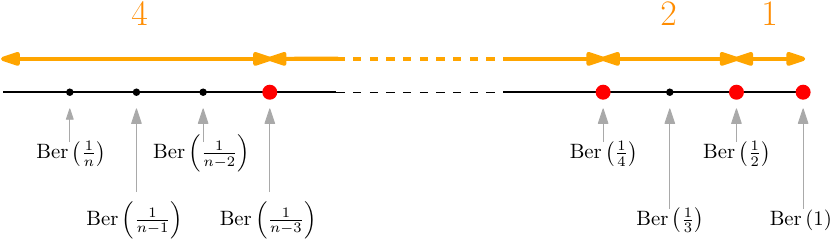}
 \caption{Constructing the law of the length of the cycles (in orange above) in a random permutation via the spacings in  Bernoulli trials with parameters $1/(n-i)$ for $i\in \{0,1,2, \dots , n-1\}$. The red dots correspond to successes. Notice that we start with a space of length $1$ just before the first trial of parameter $1/n$.}
 \end{center}
 \end{figure}

 \begin{theorem}[Feller]  The spacings $ \mathcal{S}_{n}$ between successes of the above Bernoulli variables have the same law as the cycle lengths of a uniform permutation $ \boldsymbol{\sigma}_{n} \in \mathfrak{S}_{n}$ when ordered as in the Foata construction \eqref{eq:foata}.\label{thm:fellercoupling}  \end{theorem}

\noindent \textbf{Proof.} Let us explore the cycle structure of $\boldsymbol{\sigma}_{n}$ step by step. Consider first the cycle containing $1$ in $ \boldsymbol{\sigma}_{n}$. Then, $1$ is a fixed point with probability $1/n$ --this corresponds to success of $ \mathrm{Ber}(1/n)$-- otherwise, it is sent via $\boldsymbol{\sigma}_{n}$ to a value $\boldsymbol{\sigma}_{n}(1)$ uniformly distributed over $\{2,3, \dots , n\}$. Conditionally on $\boldsymbol{\sigma}_n(1)\ne 1$, a simple calculation shows that we have $\boldsymbol{\sigma}_n( \boldsymbol{\sigma}_{n}(1))=1$  with probability $ \frac{1}{n-1}$ --this corresponds to the success of $ \mathrm{Ber}(1/(n-1))$-- or it is sent to a value $\boldsymbol{\sigma}_{n}^{2}(1) \notin\{ 1, \boldsymbol{\sigma}_{n}(1)\}$. Iteratively, if after $k \geq 2$ iterations, conditionally on $\boldsymbol{\sigma}_n^{j}(1) \ne 1$ for all $1 \leq j \leq k-1$, we have $\boldsymbol{\sigma}^{k}_{n}(1)=1$  with probability $ \frac{1}{n-k}$ --corresponding to the success of $ \mathrm{Ber}(1/(n-k))$-- otherwise the cycle continues. Hence, the length of the cycle containing $1$ indeed has the same law as the first spacing $(n+1)-I_1$ in the Bernoulli trials. Once the cycle of $1$ of length $k_1$, has been entirely explored, if $k_1 < n$ we can relabel the remaining values in increasing order by $\{1,2, \dots , n-k_1\}$ and it is easy to see that the permutation $\tilde{\boldsymbol{\sigma}}_{n-k_1}$ induced by $\boldsymbol{\sigma}_{n}$ on these values, is, conditionally on the exploration of the first cycle, uniform over $ \mathfrak{S}_{n-k_1}$ so that we can iterate the procedure. 
\qed \medskip
  
A direct consequence of the above theorem is that the law of the length of the cycle containing the point $i_0 \in \{1,2, \dots , n\}$ in the random permutation $\boldsymbol{\sigma}_n$ is a uniform variable over $\{1,2, \dots , n\}$ (see Exercise \ref{exo:taillecycle}). Also, the number of cycles $ \mathcal{C}_n$ of $ \boldsymbol{\sigma}_n$   can be expressed as $\sum_{1 \leq k \leq n} B_{k}$ where $B_{k} \sim\mathrm{Ber}(1/k)$ are independent, which is easily  handled:

 \begin{proposition}[Law of the number of cycles in a uniform permutation]\label{thm:numbercycles}\noindent\noindent For any $z \in \mathbb{C}$ we have
 $$ \mathbb{E}[z^{ \mathcal{C}_{n}}] = \prod_{j=0}^{n-1} \frac{z+j}{j+1}.$$ As a result, its expectation  and variance satisfy $ \mathbb{E}[ \mathcal{C}_{n}] = \sum_{k=1}^{n} \frac{1}{k} \sim \log n$ and $  \mathrm{Var}( \mathcal{C}_{n}) \sim \log n$ as $n \to \infty$, and we have a central limit theorem
$$\frac{ \mathcal{C}_{n}- \log n}{ \sqrt{\log n}} \xrightarrow[n\to\infty]{(d)} \mathcal{N}(0,1).$$
 \end{proposition}
 
 \noindent \textbf{Proof.} The formula for the generating function is easily proven using the equality in law $ \mathcal{C}_n =\sum_{1 \leq k \leq n} \mathrm{Ber}(1/k)$ where the Bernoulli variables are independent as in Theorem \ref{thm:fellercoupling}. Taking expectation  yields the harmonic sum, while taking variance yields to the sum of the variances which is $ \sum_{k=1}^{n} ( \tfrac1k - \tfrac1{k^{2}}) \sim \log n$. The central limit theorem can be proved  by evaluating the Fourier transform and using L\'evy's theorem (but we shall see another estimation-free route in Proposition \ref{prop:degreeroot}). \qed \medskip

\begin{exo} \label{exo:heiseiberg} Show that $ \mathbb{E}[2^{ \mathcal{C}_n}] = n+1$. Do you have a combinatorial interpretation?
\end{exo}

\section{Large cycles and Poisson--Dirichlet distribution}
In this section, we use Theorem \ref{thm:fellercoupling} to compute the law of the large cycles of a uniform permutation in the scaling limit. Perhaps surprisingly, the law of  the random partition of $1$ we obtain pops-up in other contexts such as in the factorization of large random  integers. 
\subsection{Stick breaking construction}
Let $U_1, U_2, \dots $ be a sequence of independent identically distributed uniform variables on $[0,1]$. We use these variables to perform a ``stick breaking'' of the interval $[0,1]$ by setting 
$$ X_1 = (1-U_1), \quad X_2 = U_1(1-U_2), \quad X_3 = U_1 U_2(1-U_3)\dots.$$
By the law of large numbers we have $$ \prod_{i \geq 1}^{n} U_i = \exp \Big( \underbrace{\sum_{i=1}^{n} \log U_{i}}_{ = -n + o_{a.s.}(n)} \Big) \xrightarrow[n\to\infty]{a.s.} 0.$$ and in particular we have  $ \sum X_i =1$ with probability one. 
\begin{definition}[Poisson--Dirichlet]  \label{def:PD} The Poisson--Dirichlet distribution is the law of the lengths $X_1, X_2, \dots$ in the above stick-breaking construction. 
\end{definition}

\paragraph{Ranked version.}
Although the variables $X_i$ are stochastically decreasing in $i$,  the sequence $(X_i : i \geq 1)$ is not decreasing in general. Sometimes the law of $(X_i : i \geq 1)$ is called the GEM (Griffiths, Engen, McCloskey) law and the Poisson--Dirichlet is its version $(X_i^\downarrow : i \geq 1)$ ranked in decreasing order. In these notes, we shall use the name Poisson--Dirichlet for both laws, the context making clear what we mean. The ranked version may seem more appropriate (at least to state convergence results), but actually the initial version is much more convenient from a probabilistic point of view.

\begin{figure}[!h]
 \begin{center}
 \includegraphics[height=3cm]{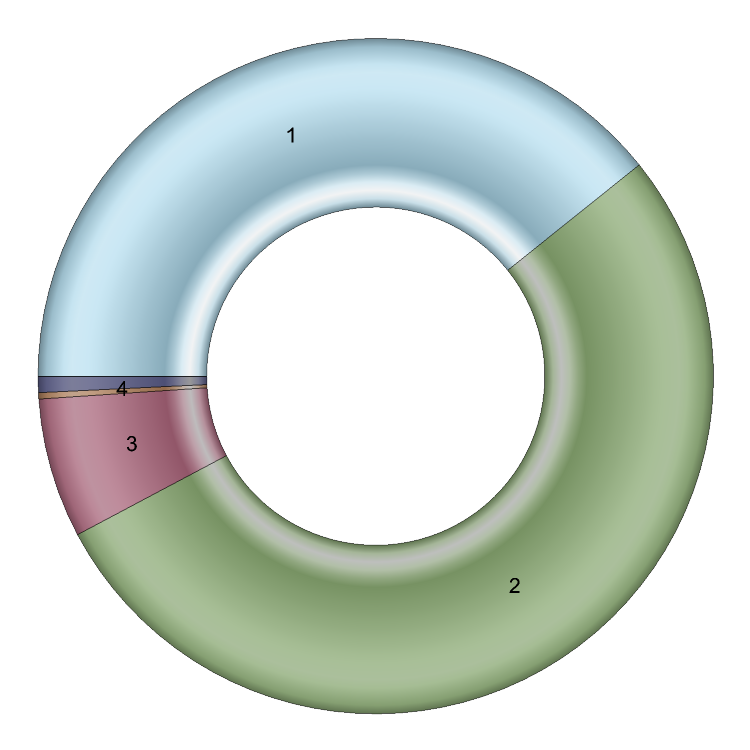}
  \includegraphics[height=3cm]{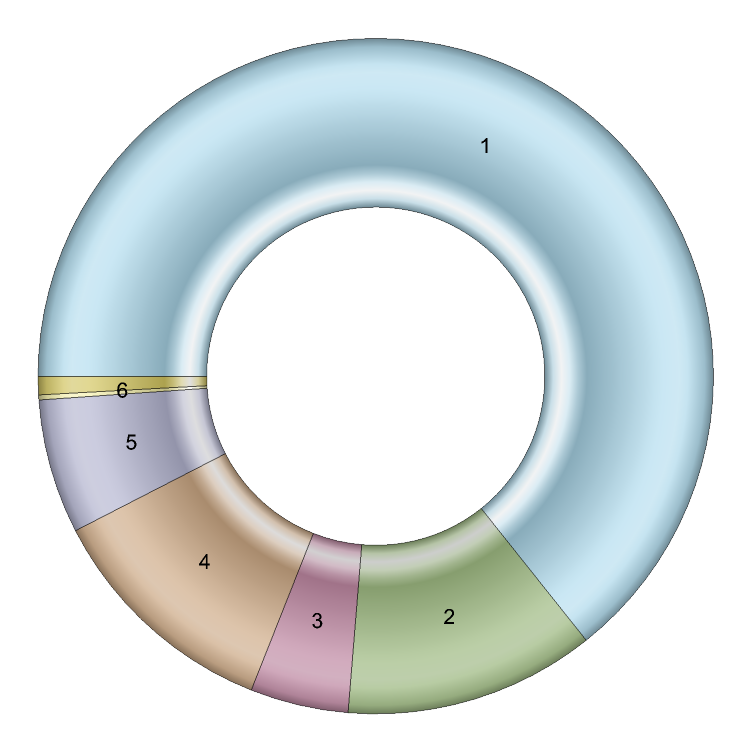}
   \includegraphics[height=3cm]{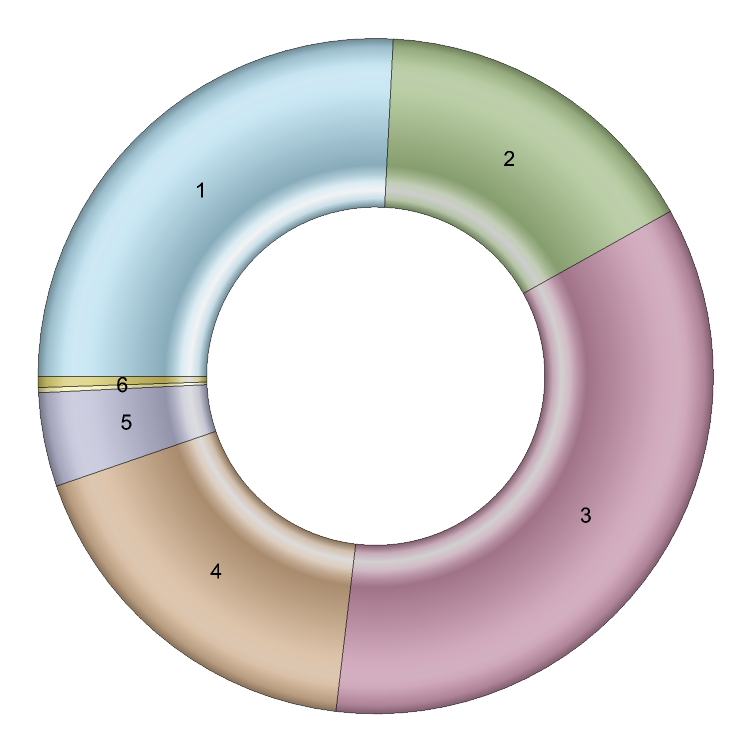}
      \includegraphics[height=3cm]{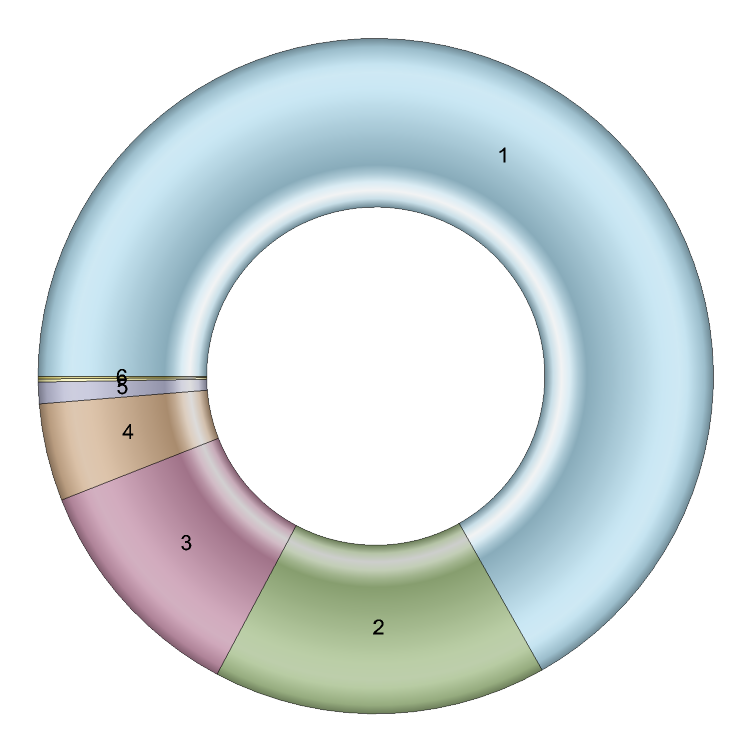}
            \includegraphics[height=3cm]{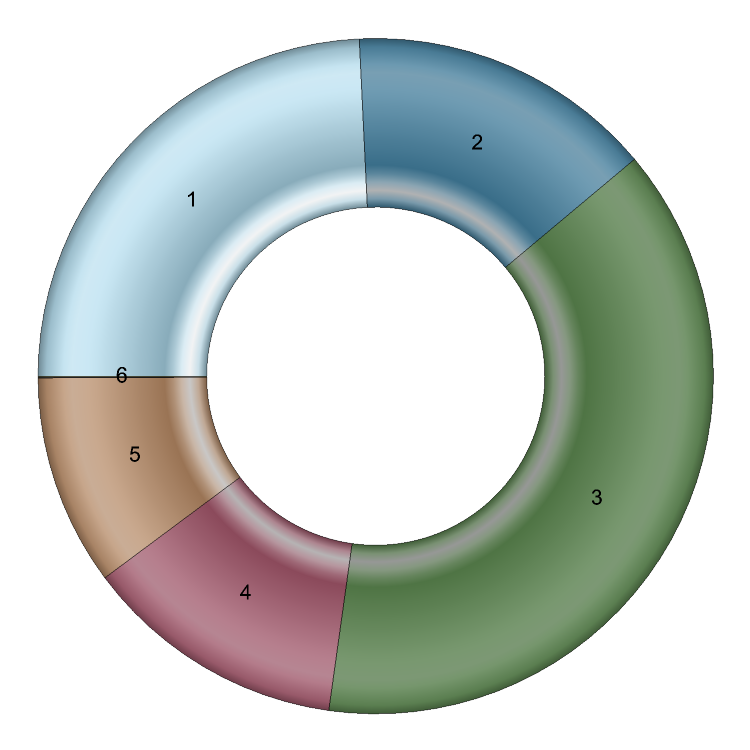}
 \caption{Five  simulations of the Poisson--Dirichlet (unranked) partition.}
 \end{center}
 \end{figure}

A corollary of Theorem \ref{thm:fellercoupling} is the following:

\begin{theorem}[Poisson--Dirichlet as limit of cycle length] \label{thm:PDforcycles}\noindent\noindent For $n \geq 0$ we denote by $K_{1}(  \boldsymbol{\sigma}_{n}), K_{2}( \boldsymbol{\sigma}_{n}), \dots$ the cycle lengths appearing in the Foata encoding of a uniform permutation $\boldsymbol{\sigma}_n \in \mathfrak{S}_n$ as in \eqref{eq:foata}. Then we have the following convergence in distribution 
 \begin{eqnarray}
  \label{eq:cvPD} \left(\frac{K_{i}(\boldsymbol{\sigma}_{n})}{n} : i \geq 1 \right ) \xrightarrow[n\to\infty]{(d)}   ( X_{i} : i \geq 1),  \end{eqnarray}
  for the $\ell_1$-distance on the space of sequences $  \ell_1^{(1)} = \{(x_i)_{i \geq 1} : x_i >0 \mbox{ and } \sum_i x_i =1 \}$. Consequently, if $K_1^\downarrow(\boldsymbol{\sigma}_n) \geq K_2^\downarrow(\boldsymbol{\sigma}_n) \geq \cdots$ are the cycle lengths of $\boldsymbol{\sigma}_n$ ranked in non-increasing order, then we have 
   \begin{eqnarray}
  \label{eq:cvPDbis} \left(\frac{K_{i}^\downarrow(\boldsymbol{\sigma}_{n})}{n} : i \geq 1 \right ) \xrightarrow[n\to\infty]{(d)}   ( X^\downarrow_{i} : i \geq 1),  \end{eqnarray}
\end{theorem}
\noindent \textbf{Proof.} In Feller's coupling,  it is straightforward to compute the law of the first spacing which is $ N = (n+1)- \sup\{ k \leq n : \mathrm{Ber}(1/k) = 1\}$. As already remarked (see Exercise \ref{exo:taillecycle}), this law is uniform over $\{1,2, \dots , n\}$ and conditionally on it, the remaining spacings have the law of $ \mathcal{S}_{n-N}$. It follows that if $ \mathcal{S}_{n}(1), \mathcal{S}_{n}(2)$ are the ordered spacings (when read from the parameter $1/n$ down to $1$) satisfy $ n^{-1} \cdot \mathcal{S}_{n}(1) \to U_{1}$ and recursively  
$$ \left(\frac{ \mathcal{S}_{n}(i)}{n} : i \geq 1\right) \xrightarrow[n\to\infty]{(d)} (X_{i} : i \geq 1),$$ in terms of finite-dimensional convergence. Actually, since we know that $ n^{-1} \mathcal{S}_n$ and $(X_i : i \geq 1)$ belong to $\ell_1^{(1)}$ (they sum-up to $1$) the finite dimensional convergence implies the $\ell_1$ convergence in law. The last convergence follows by the mapping theorem since  reordering of a sequence is a continuous operation on $\ell_1^{(1)}$.   \qed \medskip

\begin{remark}[Size-biasing and split merge dynamic] Let us give two distributional properties of the ranked Poisson--Dirichlet partition $(X_{i}^{\downarrow} : i \geq 1)$ which are not easy to prove in the continuous setting, but whose analogs in the discrete setting are obvious. 

Let us imagine $(X_{i}^{\downarrow} : i \geq 1)$ as a stick breaking of the interval $[0,1]$ into countably many intervals, and let $V\sim \mathrm{Unif}[0,1]$ be a uniform point chosen independently of this stick breaking. Then the size of the interval containing the point $V$ (there is almost surely no tie) is uniformly distribution on $[0,1]$. This can be shown by considering the cycle length of a uniform point $V_{n} \in \{1,2, \dots , n \}$ in $ \boldsymbol{\sigma}_{n}$. \\
Similarly, there is a natural dynamic on random permutations $\boldsymbol{\sigma}_{n}$ of $ \mathfrak{S}_{n}$ which preserves the uniform distribution: just compose (to the left or to the right) $\boldsymbol{\sigma}_{n}$ by a transposition $\tau_{i,j}$ where $i,j \in \{1,2, \dots , n\}$ are i.i.d.~uniform. In terms of the cycle structure, this gives rise to a split-merge transform. In the continuous setup,  this  boils down to sampling $V,V'$ independently of the stick breaking $(X_{i}^{\downarrow} : i \geq 1)$: if the two points fall into two distinct intervals, then those two pieces are merged. Otherwise, the interval containing both $V$ and $V'$ is split into two intervals uniformly. The Poisson--Dirichlet law is an invariant measure for this dynamic (and is in fact the only one, see \cite{diaconis2004}).
\end{remark}
\bigskip 

Perhaps surprisingly, the Poisson--Dirichlet law appears in many other ``logarithmic combinatorial structures'' such as factorization of random polynomials over finite fields or prime factorization of large random integers:

\begin{theorem}[Billingsley]\noindent Let $N \in \{1,2, \dots , n\}$ be a uniform integer less than or equal to $n$ and let $p^\downarrow_1(N) \geq p^{\downarrow}_2(N) \geq \dots$ its prime factors (with possible repetition). Then we have 
$$ \left( \frac{\log p^{\downarrow}_i(N)}{\log n} ; i \geq 1\right) \xrightarrow[n\to\infty]{(d)} (X_i^\downarrow: i \geq 1).$$
\end{theorem}
We refer to \cite{arratia2003logarithmic} for details.

\subsection{Dickman function and prisoners}

In this section, we present Dickman\footnote{  \raisebox{-5mm}{\includegraphics[width=1cm]{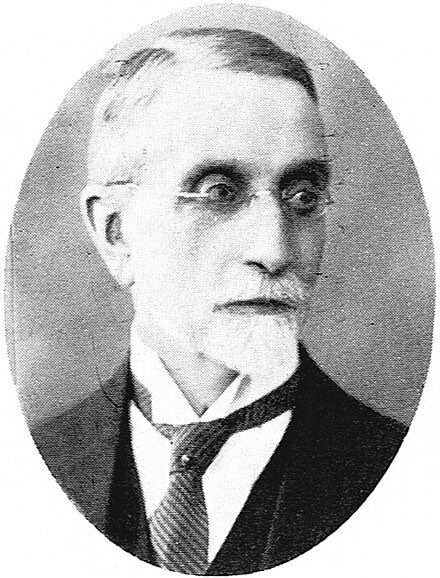}} Karl Dickman (1861--1947), Swedish. He was actuary and  published  only one article in mathematics \cite{dickman1930frequency} introducing this function when he was around $70$.} function which is essentially the tail distribution $X_1^\downarrow$, the scaling limit of the longest cycle in a random permutation. This function pops up in various places in analytic number theory and has intriguing properties.

\begin{proposition}[Dickman function] Consider $(X_i^\downarrow : i \geq 1)$ the ranked version of a Poisson--Dirichlet distribution. Then for $x \geq 0$ we have 
$$ \mathbb{P}( X_1^\downarrow \leq x) = \rho(1/x),$$
where $x \mapsto \rho(x)$ is Dickman's function defined by

$$\left\{ \begin{array}{lc} \rho(x)=1 & \mbox{ for } 0 \leq x \leq 1\\
x \rho'(x) = - \rho(x-1)& \mbox{ for } x \geq 1. \end{array}\right.$$

\begin{figure}[!h]
 \begin{center}
 \includegraphics[width=10cm]{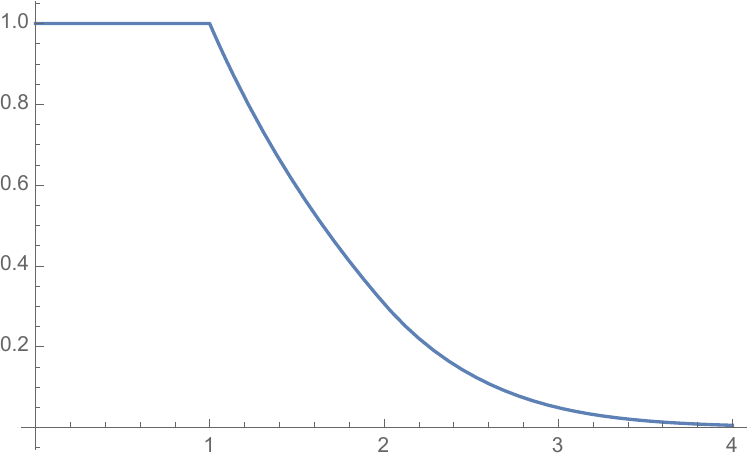}
 \caption{Dickman's function}
 \end{center}
 \end{figure}
\end{proposition}
\noindent \textbf{Proof.} We use the notation $\mathbb{P}( X_1^\downarrow \leq x) = \rho(1/x)$ extended to $\rho(u)=1$ for $u \in [0,1]$. In the unranked version $(X_{i} : i \geq 1)$ of the Poisson--Dirichlet partition we can write after conditioning on the first uniform variable $U_{1}$  \begin{eqnarray*} \mathbb{P}( X_1^\downarrow \leq x) &=& \mathbb{P}(\{X_{1} \leq x\} \cap \{X_{2}, \dots , X_{n}, \dots \leq x\})\\
&=& \mathbb{E}\left[ \mathbf{1}_{1-U_{1} \leq x} \mathbb{P}\left( \tilde{X}_{1}^{\downarrow} \leq  \frac{x}{U_{1}}\right)\right],  \end{eqnarray*} which give the following integral equation 
 \begin{eqnarray*} \rho(1/x) =\mathbb{P}( X_1^\downarrow \leq x) &=& \int_0^x \mathrm{d}u \, \mathbb{P}( X_1^\downarrow \leq  \frac{x}{1-u})\\  &=& \int_0^x \mathrm{d}u \, \rho\left( \frac{1-u}{x} \right)\\
  \rho(y) & \underset{\begin{subarray}{c}v=u/x\\ y=1/x \end{subarray}}{=}&  \frac{1}{y} \int_0^1  \mathrm{d}v \, \rho(y-v).  \end{eqnarray*}
 Differentiating the equality $y \rho(y) = \int_0^1 \mathrm{d}v \, \rho(y-v)$ with respect to $y$, we recover the delayed differential equation of the proposition. \qed
 
 \bigskip 
 
Related to the Dickman function, let us state a famous riddle:
\begin{quote}
The director of a prison offers 100 death row prisoners, who are numbered from 1 to 100, a last chance. A room contains a cupboard with 100 drawers. The director randomly puts one prisoner's number in each closed drawer. The prisoners enter the room, one after another. Each prisoner may open and look into 50 drawers in any order. The drawers are closed again afterwards. If, during this search, every prisoner finds their number in one of the drawers, all prisoners are pardoned. If even one prisoner does not find their number, all prisoners die. Before the first prisoner enters the room, the prisoners may discuss strategy — but may not communicate once the first prisoner enters to look in the drawers. What is the prisoners' best strategy?
\end{quote}
Opening $50$ drawers at random (independently for each prisoner) is a hopeless strategy since the probability that they all manage to find their numbers is $(1/2)^{100} \approx 0$. However, they can correlate their searchs if the $i$th prisoner starts with $i$th drawer, looks at the discovered label and successively follows the cycle of the underlying permutation of the labels. The probability of success is the probability that no cycle of the permutation of the labels has a length larger than $50$ which is approximately $ \mathbb{P}(X_1^\downarrow \leq 1/2) = 1- \log 2 \approx 30 \%$.

\paragraph{Formulas without words.}

  \begin{eqnarray*} \mbox{Euler's constant} & = & \int_1^\infty\left(  \frac{1}{\lfloor x \rfloor } -\frac{1}{x} \right) \mathrm{d}x = \log \int_0^\infty \mathrm{d}x\, \rho(x).\\
  \mbox{Golomb-Dickman constant} &=&  \int_0^1  \mathrm{d}x\,\mathrm{exp} \left( \int_0^x  \frac{\mathrm{d}t}{ \ln t} \right) = \mathbb{E}[X_1^\downarrow] = \int_0^\infty  \mathrm{d}t\frac{ \rho(t) }{(t+1)^2}.\\
  \sum_{k\geq 1} \prod_{j=1}^k U_j   &\overset{(d)}{=}&  \mathrm{e}^{-\gamma} \rho( x) \mathrm{d}x.\\
  \mbox{where $U_i$ are i.i.d.~uniforms on $[0,1]$} && 
  \end{eqnarray*}

\section{Poisson count for short cycles}
In the previous section, we saw that the Poisson--Dirichlet law is the limit law of the large cycles in a random uniform permutation. However, the information about the small cycles is lost in this limit and we will see below that they are ruled by the Poisson paradigm already encountered in Section \ref{sec:poissonparadigm}. 

\subsection{An abstract limit from the Feller coupling} 
Recall the setup of Theorem \ref{thm:fellercoupling}, let $( B_{k}  \sim \mathrm{Ber}(1/k) : k \geq 1)$ be independent Bernoulli variables of parameter $1/k$ and denote by $1 = \tilde{I}_{1} < \tilde{I}_{2} < \cdots$ the indices of the variables equal to $1$ (beware we see those variables as indexed ``in the other direction'' compared to the previous section).  In a sense, the spacings between $\mathcal{S}_{\infty} := ( \tilde{I}_k: k \geq 1)$ could be seen as the cycle structure of an ``infinite permutation''.  Down to earth, we have 
$$ \forall A \geq 1, \quad \sum_{k=1}^\infty \mathbb{P}(  \tilde{I}_{k+1}-\tilde{I}_{k} = A) =  \sum_{i=1}^\infty \sum_{k=1}^\infty \mathbb{P}(  \tilde{I}_k =i, \tilde{I}_{k+1} = i+A) \leq  \sum_{i=1}^\infty \frac{1}{i(i+A)} < \infty,$$ so that the Borel--Cantelli lemma shows that $ \tilde{I}_{k+1}- \tilde{I}_{k} \to \infty$ almost surely as $k \to \infty$. In particular, we can define the increasing rearrangement of the spacings between consecutive points in $( \tilde{I}_k: k \geq 1)$ and their count
$$   \mathcal{N}_{i} := \# \{ k \geq 1 :   \tilde{I}_{k+1}- \tilde{I}_{k} =i\} < \infty.$$

Below we write $N_i(\boldsymbol{\sigma}_n)$ for the \textbf{number of cycles} of length $i$ in the decomposition of the random uniform permutation $\boldsymbol{\sigma}_n$ into product of cycles with disjoint supports. Given Theorem \ref{thm:fellercoupling}, it is rather straightforward to show that $(N_{i}( \boldsymbol{\sigma}_{n}) : i \geq 1)$ converge in law as $n \to \infty$ :

\begin{proposition}  \label{prop:limitfellereasy} We have the convergence in law (in the sense of finite-dimensional marginals)
 \begin{eqnarray} ( N_{i}( \boldsymbol{\sigma}_{n}) : i \geq 1) \xrightarrow[n\to\infty]{(d)} ( \mathcal{N}_{i} : i \geq 1).   \label{eq:limitPoisson?}\end{eqnarray}
 \end{proposition}
 \noindent \textbf{Proof.} Feller's coupling (Theorem \ref{thm:fellercoupling}) provides a way to couple  uniform permutations $( \boldsymbol{\sigma}_n^{ \mathrm{fel}}: n \geq 1)$ on a common probability space so that $\boldsymbol{\sigma}_n^{ \mathrm{fel}} = \boldsymbol{\sigma}_n$ in law and such that the cycle structure of $ \boldsymbol{\sigma}_n^{ \mathrm{fel}}$ coincides with the spacings between the points $ 1=\tilde{I}_{1} < \tilde{I}_{2}  < \cdots < \tilde{I}_{\ell_n} < (n+1)$ where $\tilde{I}_{\ell_n}$ is the last index strictly before $(n+1)$. In this coupling we \textit{nearly} have the almost sure convergence $ N_i( \boldsymbol{\sigma}_n^{\mathrm{fel}}) \to \mathcal{N}_i$ as $n \to \infty$. The reason that the coupling falls short of proving this point-wise  convergence is that if $(n+1)$ is large and located precisely $i_0$ unit after a point of $ \mathcal{S}_\infty$ (with no other point in-between) then we have $N_{i_0}(\boldsymbol{\sigma}_n^{\mathrm{fel}}) = \mathcal{N}_{i_0} + 1$. However, for any positive function $f$ bounded by $C>0$ and any $k_0 \geq 1$ we have 
  \begin{eqnarray*} &&  \Big|\mathbb{E}[f( N_i( \boldsymbol{\sigma}_n^{ \mathrm{fel}}) : 1 \leq i \leq k_0)] - \mathbb{E}[f(  \mathcal{N}_i: 1 \leq i \leq k_0)]\Big| \\ &\leq& C \Big(\mathbb{P}( S_\infty \cap \{n-k_0,\dots, n-1,n\} \ne \varnothing) + \mathbb{P}( \exists  \tilde{I}_{\ell} \geq n \mbox{ with }  \tilde{I}_{\ell+1} - \tilde{I}_\ell \leq k_0)\Big)\\ & \xrightarrow[n\to\infty]{} &0.  \end{eqnarray*}
 The desired convergence in law follows. \qed 
 \bigskip 
 
We will see in Theorem \ref{thm:poissoncycle} below that  the law of  $ ( \mathcal{N}_{i} : i \geq 1)$ is actually super simple! A simple way to see this, is to take a small detour using Cauchy's formula and to randomize the permutation's length. This operation, usually called Poissonnization, will be made crystal clear in Chapter \ref{chap:AK}.

\subsection{Cauchy's formula and interpretation}

The starting point is a famous formula due to Cauchy giving the exact law of the cycle-counting function.  With the notation above we have:

\begin{proposition}[Cauchy] For any $c_1,c_2, \dots , c_n \in  \mathbb{Z}_{\geq 0}$ so that $\sum_{i=1}^n i c_i =n$ we have 
$$ \mathbb{P}(N_i(\boldsymbol{\sigma}_n) = c_i, \forall 1 \leq i \leq n) = \prod_{i=1}^n\frac{(1/i)^{c_i}}{(c_i) !}.$$
\end{proposition}
\noindent \textbf{Proof.} Once the cycle structure $(c_{i} : 1 \leq i \leq n)$ of the permutation $ \boldsymbol{\sigma}_{n}$ has been fixed (with the obvious constraint), the number of possible candidates is obtained by: 
\begin{itemize} 
\item distributing the $n$ numbers $1,2, \dots , n$ into the $\sum c_{i}$ boxes of sizes $1,2, \dots , i, \dots, n$: since the $c_{i}$ boxes of size $i$ are  indistinguishable, there are 
$$ \left( \frac{n!}{  \cdots  \underbrace{i! \cdot i!}_{c_{i} \mathrm{\ times}}\cdots }\right)\cdot  \prod_{i=1}^{n} \frac{1}{ c_{i} !} \quad \mbox{ such choices}.$$ 
\item then constructing an  $i$-cycle with the numbers in each box of size $i$: there are $(i-1)!$ possibilities each.
\end{itemize}
We deduce that the probability in the proposition is given by 
$$   \frac{1}{ n!} \cdot \frac{n!}{ \prod_{i=1}^{n} (i!)^{c_{i}}} \prod_{i=1}^{n} \frac{1}{ c_{i} !}  \prod_{i=1}^{n} \big((i-1)!\big)^{c_{i}} = \prod_{i=1}^n\frac{(1/i)^{c_i}}{(c_i) !}.$$ \qed \bigskip 

Let us put  our probabilist's glasses on and interpret the previous formula as follows:
 \begin{eqnarray*}
 \mathbb{P}(N_i(\boldsymbol{\sigma}_n) = c_i, \forall 1 \leq i \leq n) &=& \prod_{i=1}^n\frac{(1/i)^{c_i}}{(c_i) !}\\
 &=&   \mathrm{e}^{1+ \frac{1}{2} + \cdots + \frac{1}{n}} \cdot \prod_{i=1}^n \mathrm{e}^{-1/i} \frac{(1/i)^{c_i}}{(c_i) !} \\
 &=& \mathrm{e}^{ \mathrm{h}_n} \prod_{i=1}^n  \mathbb{P}(Z_i = c_i),\end{eqnarray*}
where $(Z_i : i \geq 1)$ are independent Poisson random variables with mean $1/i$, and where $ \mathrm{h}_n$ is the $n$th harmonic sum. In other words, the vector $ (N_i(\boldsymbol{\sigma}_n) : 1 \leq i \leq n)$ has the same law as $ (Z_i : 1 \leq i \leq n)$ conditioned on the event $\{ \sum_{1 \leq i \leq n} i Z_i =n \}$. This observation, due to Kolchin, can actually be pushed a little further as remarked by Lloyd \& Shepp. Denote by $\sigma_0$ the permutation with $0$ cycles so that $N_i(\sigma_0) =0$ for all $i \geq 1$. For $x \in (0,1)$, for any sequence $(c_i : i\geq 1)$ of integers which is eventually the null sequence, if we denote by $N = \sum i c_i$ then we have 
 \begin{eqnarray*}    (1-x) \sum_{n =0}^\infty x^n \mathbb{P}\big(N_i(\boldsymbol{\sigma}_n) = c_i: \forall 1 \leq i \leq n\big) &=&  (1-x) x^N  \prod_{i=1}^n \frac{(1/i)^{c_i}}{(c_i) !} \\
 &=&  \underbrace{(1-x) \mathrm{e}^{x+ \frac{x^2}{2} + \cdots + \frac{x^n}{n} +  \cdots}}_{1} \prod_{i=1}^\infty \mathrm{e}^{- x^i /i } \frac{(x^i/i)^{c_i}}{(c_i) !}
  \end{eqnarray*}
  This means:
  \begin{lemma}  If $ \mathbf{n}_x \in \{0,1,2, \dots\}$ is a geometric random variable with mean $ \frac{x}{1-x}$ and if, conditionally on $ \mathbf{n}_x$, we let $\boldsymbol{\sigma}_{\mathbf{n}_x}$ be a uniform permutation on $  \mathfrak{S}_{\mathbf{n}_x}$, then the cycle counts $(N_i( \boldsymbol{\sigma}_{\mathbf{n}_x}) : i \geq 1)$ has the same law as independent Poisson random variables with means $ \frac{x^i}{i}$ for $i \geq 1$. \label{lem:shepplloyd}
  \end{lemma}
  
  We will see in Chapter \ref{chap:AK} that the above lemma follows from combining the construction of the random recursive tree from a Yule process in continuous time and the Chinese restaurant process (sic!).
  
  \begin{exo}[Random $\zeta$-number]  \label{exo:zeta} For $s>1,$ consider $ \mathbf{N}_s \in \{1,2,\dots\}$ a random number sampled according to 
  $$ \mathbb{P}( \mathbf{N}_s = n) = \frac{1}{\zeta(s)} n^{-s}.$$
  Show that the $p$-valuations $(\nu_p(  \mathbf{N}_s) : p \in \mathcal{P})$ are independent geometric random variables with success parameters $(1/p)^s$ for all prime numbers $p \in \mathcal{P}$.
  \end{exo}

\subsection{Poisson limit}
We are now armed to prove the following:
\begin{theorem}[Goncharov, Kolchin]   Recall that  $N_i(\boldsymbol{\sigma}_n)$ is the number of cycles of length $i\geq 1$ in the decomposition of the uniform permutation $\boldsymbol{\sigma}_n$ into product of cycles with disjoint supports. Then we have the following convergence in law for the finite-dimensional marginals
$$ \big( N_i(\boldsymbol{\sigma}_n) : i \geq 1\big) \xrightarrow[n\to\infty]{(d)} \big( \mathrm{Poi}( 1/i) : i \geq 1 \big),$$ where the Poisson random variables on the right-hand side are independent and of mean $1/i$ for $i \geq 1$. \label{thm:poissoncycle} 
\end{theorem}
\begin{remark}[Derangements] We recover  the famous asymptotic of the number of derangements (permutations without fixed points) since the last theorem implies in particular that as $n \to \infty$ we have 
$$  \frac{\# \{ {\sigma_n} \in \mathfrak{S}_n : {\sigma}_n \mbox{ has no fixed points}\}}{n!} =  \mathbb{P}(N_1(\boldsymbol{\sigma}_n) =0)  \xrightarrow[n\to\infty]{}  \mathbb{P}( \mathrm{Poi}(1) =0) =  \mathrm{e}^{-1}.$$
In fact, the inclusion-exclusion principle shows that we have the explicit series representation $ \sum_{k=0}^{n} (-1)^k \frac{n!}{k!}$ for the number of derangements of $ \mathfrak{S}_n$.  
\end{remark}

\noindent \textbf{Proof.}  We already know from \eqref{eq:limitPoisson?} that $( N_i( \boldsymbol{\sigma}_n) : i \geq 1)$ converges in law towards some limiting vector $( \mathcal{N}_i : i \geq 1)$ as $n \to \infty$. On the other hand, if we let $ x \to 1$ in Lemma \ref{lem:shepplloyd} we deduce that $ \mathbf{n}_x \to \infty$ in probability. Since conditionally on $ \mathbf{n}_x$ the permutation $\boldsymbol{\sigma}_{ \mathbf{n}_x}$ is uniform, we deduce that 
$$ \begin{array}{rcl} ( N_i(\boldsymbol{\sigma}_{ \mathbf{n}_x}) : i \geq 1) & \xrightarrow[x\to1]{(d)} & ( \mathcal{N}_i : i \geq 1)\\
 \rotatebox{90}{=} \quad \mbox{in law}&& \quad \rotatebox{90}{=} \quad \mbox{in law} \\
 \big( \mathrm{Poi}( x^i/i) : i \geq 1\big) & \xrightarrow[x\to1]{(d)} & \big( \mathrm{Poi}( 1/i) : i \geq 1\big).  \end{array}$$ where all the Poisson variables are independent. \qed

\begin{remark}[Direct calculation] It can be seen directly that the variables in \eqref{eq:limitPoisson?} are independent Poisson variables with mean $1/i$ without referring to random permutations. In fact, once the limit has been re-interpreted as the spacings between records of i.i.d.~uniforms on $[0,1]$, it is a consequence of a more general theorem due to Ignatov on the Poissonnian structure of records values of a Markov process. We refer the interested reader to \cite{najnudel2020feller} and \cite{resnick2008extreme}  for details. \end{remark}

\paragraph{Bibliographical notes.} There are many references on the popular subject of random permutation, see e.g. the Saint-Flour lectures of Pitman \cite{Pit06} in particular Section 3.1 or the Bible in combinatorics  \cite{Flajolet:analytic}. Various sets of lecture notes are also available on the web such as \cite{feray2019random,gerinmini} and more recent results about  ``logarithmic combinatorial structures'' can be found in \cite{arratia2003logarithmic}. Feller's coupling is proved in   \cite{feller1945fundamental},  and the Poisson counting limit is due to Goncharov and Kolchin, but our proof based on Lemma \ref{lem:shepplloyd} is inspired from Lloyd and Shepp \cite{shepp1966ordered}.  For more about appearance of Dickman's function in probabilistic and analytic  number theory, see \cite{tenenbaum2015introduction} and \cite{chamayou1973probabilistic}. We also refer to \cite{bureaux2015methodes} for other applications of the randomization technique to random partitions. \bigskip

\noindent{\textbf{Hints for Exercises.}}\ \\
Exercise \ref{exo:taillecycle}: The size of the cycle containing $1$ in $\sigma_n$ is equal to the value of the pre-image of $1$ in $ \mathrm{Foata}(\sigma_n)$.\\
Exercise \ref{exo:heiseiberg}: A random permutation sampled according to $2^{\# \mathrm{ number \ of\  cycles}}$ appears in Toth's representation of the quantum Heisenberg ferromagnet on the complete graph (sic!), see \cite{toth1993improved}.\\
Exercise \ref{exo:zeta}: Re-interpret the Eulerian product formula $ \displaystyle \prod_{p \in \mathcal{P}} \left( \frac{1}{1- \frac{1}{p^s}}\right) = \sum_{n=1}^\infty \frac{1}{n^s}.$

\chapter{Random recursive tree}
\label{chap:RRT}


\hfill L’arbre, c’est cette puissance qui lentement épouse le ciel.
\medskip

\hfill A. de Saint-Exupéry
\bigskip

In this chapter we study the following random tree growth model:
\begin{definition}[RRT] The \textbf{random recursive tree} (RRT) is the Markov chain with values in the set of all unoriented  labeled trees such that $T_0 = $ \raisebox{.5pt}{\textcircled{\raisebox{-.9pt} {0}}} and so that for $n \geq 1$, conditionally on $T_{n-1}$, the labeled tree $T_{n}$ is obtained by attaching the new vertex  \raisebox{.5pt}{\textcircled{\raisebox{-.6pt} {$n$}}}  onto a uniform vertex of $T_{n-1}$.
\end{definition}

\begin{figure}[!h]
 \begin{center}
 \includegraphics[width=16cm]{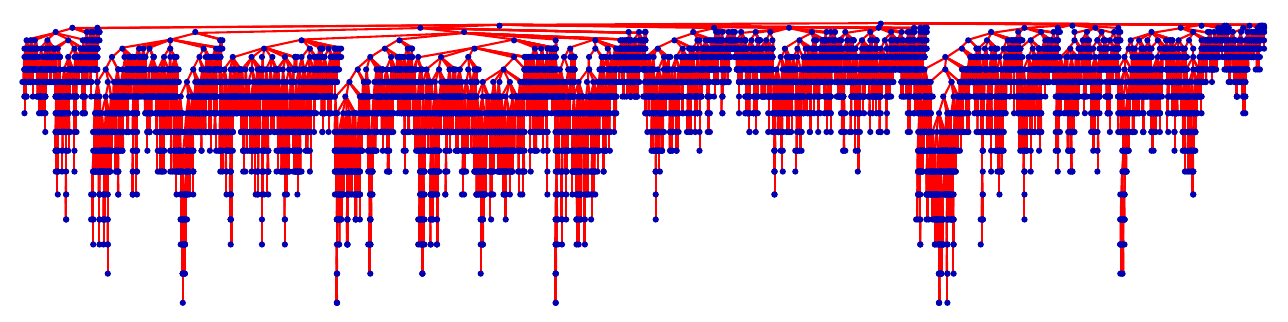}
 \caption{A simulation of $T_{10 000}$ where the root vertex $ \noeud{0}$ is placed at the top. Clearly, the random recursive tree seems ``short and fat''.}
 \end{center}
 \end{figure}

Obviously there are $n!$ possible values for $T_n$: these are all \textbf{increasing labeled trees with $n+1$ vertices} i.e.~unoriented trees labeled from $0$ up to $n$ and so that the labels along each branch starting from \raisebox{.5pt}{\textcircled{\raisebox{-.9pt} {0}}} are increasing.  For each $n\geq 0$, the RRT $T_{n}$ is a uniform random variable over this set. We shall start with a link between this model of random tree and random permutations of the symmetric group $ \mathfrak{S}_{n}$ over $n$ elements.

\section{Chinese Restaurant process}

\label{sec:RRTandpermut}

Since there are $n!$  equiprobable values for $T_n$, the RRT stopped at time $n$ can be seen as an encoding of  a uniform permutation $ \boldsymbol{ \sigma}_n$ of $ \mathfrak{S}_{n}$. Moreover, it is possible to couple these encodings in a particularly nice way so that it is coherent for all $n \geq 1$ simultaneously: this is the so-called \textit{Chinese restaurant process} (CRP). This coupling is different from Feller's coupling seen in the previous chapter. \medskip 

\subsection{Coupling CRP-RRT}
\label{sec:CRP-RRT}
Let $\sigma_{n} \in \mathfrak{S}_{n}$ be a permutation over $\{1,2, \dots , n\}$. If $ n \geq 2$, we can canonically associate with $\sigma_n$ a permutation $[\sigma_{n}]_{n-1} \in \mathfrak{S}_{n-1}$ as follows: it is the permutation  defined for $ k  \in \{1, \dots , n-1\}$ by  
$$ \left\{ \begin{array}{cl}
[\sigma_{n}]_{n-1}(k) = \sigma_{n}(k)  & \mbox{ if } \sigma_{n}(k) \ne n,\\

[\sigma_{n}]_{n-1}(k) = \sigma_{n}(n)  & \mbox{ if } \sigma_{n}(k) = n.\\
\end{array}\right.
$$
The effect of removing the value $n$ from $ \sigma_n$ is better understood on the cycle decomposition: the permutation $[\sigma_n]_{n-1}$ is obtained by removing the value $n$ in the cycle of $ \sigma_n$ which contains it. By extending the restriction step by step we can define $[\sigma_{n}]_{k}$ for all $k \leq n$ and it is easy to see that  if $\boldsymbol{\sigma}_{n}$ is uniform over $ \mathfrak{S}_{n}$ then $[ \boldsymbol{\sigma}_{n}]_{k}$ is also uniformly distributed over $ \mathfrak{S}_{k}$. 

Actually, it is easy to reverse the procedure and construct a  sequence of  random permutations $(\boldsymbol{\sigma}_{n}^{ \mathrm{cr}}  \in \mathfrak{S}_{n} : n \geq 1)$ as a Markov chain. Specifically, let $\boldsymbol{\sigma}^{ \mathrm{cr}}_{1}= (1)$ and for $n \geq 2$, conditionally on $\boldsymbol{\sigma}_{n-1}^{ \mathrm{cr}}$, the permutation $\boldsymbol{\sigma}_{n}^{ \mathrm{cr}}$ is obtained with probability $1/n$ by just declaring $\boldsymbol{\sigma}_n^{ \mathrm{cr}}(n)=n$ and with probability $1- \frac{1}{n}$ by picking a uniform integer $k \in \{1, \dots , n-1\}$ and declaring that  $$\boldsymbol{\sigma}^{ \mathrm{cr}}_{n}(k) = n\quad \mbox{and }\quad \boldsymbol{\sigma}^{ \mathrm{cr}}_{n}(n) = \boldsymbol{\sigma}^{ \mathrm{cr}}_{n-1}(k),$$ the others values being unchanged between $\sigma_{n}^{{ \mathrm{cr}}}$ and $\sigma_{n-1}^{{ \mathrm{cr}}}$.
With the above notation we have $[\boldsymbol{\sigma}^{ \mathrm{cr}}_{n}]_{k}= \boldsymbol{\sigma}_{k}^{ \mathrm{cr}}$ and this Markov chain produces a coupling of permutations uniformly distributed over $ \mathfrak{S}_{n}$ for each $n$. \medskip 

The evolution of the cycle structure of $ \boldsymbol{\sigma}^{ \mathrm{cr}}_{n}$ in the previous Markov chain is described by the following mechanism called the \textbf{Chinese restaurant process}\footnote{  \raisebox{-5mm}{\includegraphics[width=1cm]{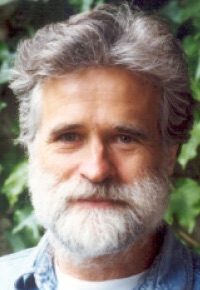}} Jim Pitman (1949--), Australian}: In this process, customers $1,2,3\dots$ arrive sequentially in an imaginary (Chinese) restaurant. At step $n=1$, the customer \raisebox{.5pt}{\textcircled{\raisebox{-.9pt} {$1$}}} arrives and sits at a new table. Inductively at step $n \geq 2$,  the customer \raisebox{.5pt}{\textcircled{\raisebox{-.6pt} {$n$}}} sits at the right of any of $n-1$ previous customers with probability $\tfrac1n$ or creates a new table with probability $\tfrac1n$. It should be clear from the above construction that the tables in the Chinese restaurant process describe the cycle structure of the growing sequence of permutations $(\boldsymbol{\sigma}^{ \mathrm{cr}}_n : n \geq 1)$. \medskip

\begin{figure}[!h]
 \begin{center}
 \includegraphics[width=16cm]{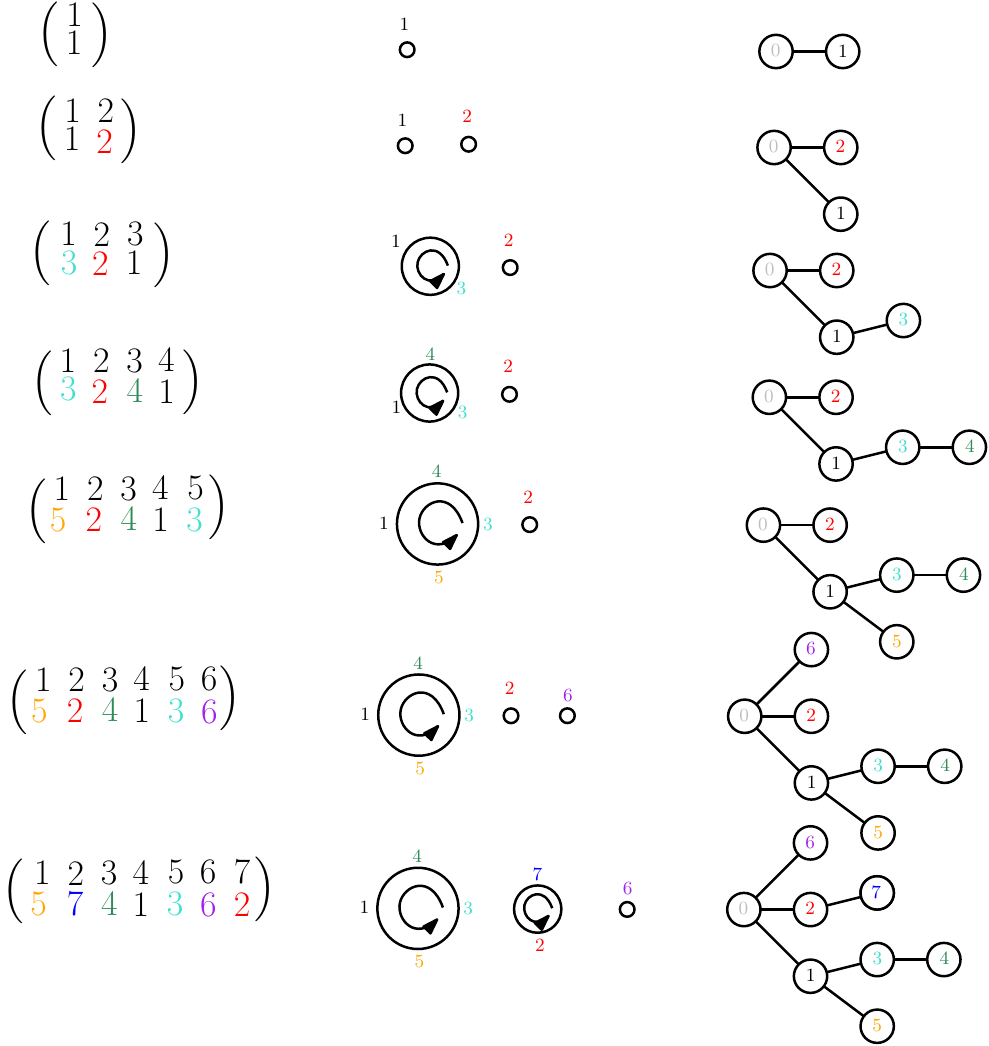}
 \caption{ \label{fig:CRP} Illustration of the coupling between growing permutations $( \boldsymbol{\sigma}_n^{ \mathrm{cr}} : n \geq 1)$ on the left, the Chinese restaurant process in the middle and the random recursive tree $(T_n : n \geq 1)$ starting with initial state $0$ on the right.}
 \end{center}
 \end{figure}

The Chinese restaurant process is canonically coupled with the RRT $(T_n : n \geq 0)$ by declaring that the new customer $n$ corresponds to the vertex \raisebox{.5pt}{\textcircled{\raisebox{-.9pt} {$n$}}} and it attaches in $T_n$ to the vertex corresponding to the customer on its left, or to the vertex \raisebox{.5pt}{\textcircled{\raisebox{-.9pt} {$0$}}} if this customer creates a new table. See Figure \ref{fig:CRP}. Thanks to this coupling, we deduce in particular that the degree of $\raisebox{.5pt}{\textcircled{\raisebox{-.9pt} {$0$}}}$ in $T_{n}$ is equal to the number of cycles in the cycle decomposition of $\boldsymbol{\sigma}_{n}^{ \mathrm{cr}}$.

\subsection{Pólya urn and almost sure convergence towards Poisson--Dirichlet} \label{sec:polyaurn}
The Chinese restaurant coupling $( \boldsymbol{\sigma}_n^{ \mathrm{cr}} : n \geq 1)$ is a different coupling compared to Feller's coupling $( \boldsymbol{\sigma}_n^{ \mathrm{fel}} : n \geq 1)$ used in the proof of Proposition \ref{prop:limitfellereasy}. Roughly speaking, in the Chinese restaurant coupling, the structure of large cycles converges almost surely (see below), whereas in Feller's coupling the structure of small cycles (nearly) converges almost surely. Recalling Theorem  \ref{thm:PDforcycles},  we have here:  

\begin{theorem}[Almost sure convergence of the Chinese restaurant process] Let $ (K_i(\boldsymbol{\sigma}_n^{ \mathrm{cr}}) : i \geq 1)$ be the cycle lengths of $\boldsymbol{\sigma}_n^{ \mathrm{cr}}$ in the Foata encoding \eqref{eq:foata} i.e.~the table sizes ranked by order of creation in the Chinese restaurant process. If $(X_i : i \geq 1)$ is the (unranked) Poisson--Dirichlet random partition of $[0,1]$ (see Definition \ref{def:PD}) then we have the following convergence in law in  $\ell_1^{(1)}$ \label{thm:CVPDps}
$$ \left( \frac{ K_i(\boldsymbol{\sigma}_n^{ \mathrm{cr}})}{n} : i \geq 1\right) \xrightarrow[n\to\infty]{a.s.} (X_i: i \geq 1).$$
\end{theorem}

To prove the theorem let us first focus on the behavior of the process 
$$ (R_n,B_n) := \left( K_1( \boldsymbol{\sigma}_n^{ \mathrm{cr}}) , 1+\sum_{i \geq 2} K_i(\boldsymbol{\sigma}_n^ {  \mathrm{cr}})\right),$$ for $n \geq 1$.
It is clear from the definition of the Chinese restaurant process that this is a Markov chain starting from $R_1=1, B_1=1$ and with transition probabilities given by 
  \begin{eqnarray} \mathbb{P}(R_{n+1} = R_n + 1 \mid R_n,B_n) = \frac{R_n}{R_n+B_n}, \quad \mathbb{P}(B_{n+1} = B_n + 1 \mid R_n,B_n) = \frac{B_n}{R_n+B_n}.   \label{eq:polya}\end{eqnarray}
We recognize here the (law of the) famous \textbf{Pólya\footnote{\raisebox{-5mm}{\includegraphics[width=1cm]{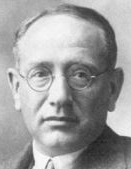}} George (György) Pólya (1887-1985), Hungarian} urn}, which is the stochastic system informally described as follows: initially at time $n=1$ an urn contains one red ball and one blue ball. At each step, a ball is drawn from the urn uniformly at random and is replaced in the urn \textbf{together} with a new ball of the same color (re-inforcement). Then the number  $(R_n,B_n)$ of red and blue balls at step $n$ is clearly a Markov chain with  transitions \eqref{eq:polya}. 

\begin{figure}[!h]
 \begin{center}
 \includegraphics[width=15cm]{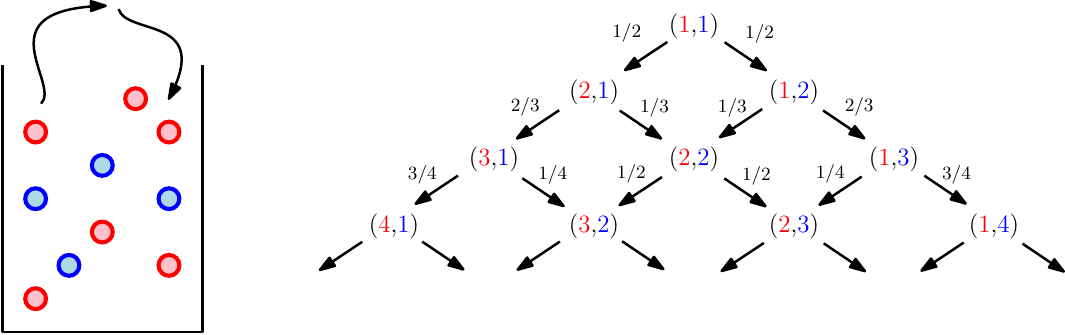}
 \caption{Mecanism of the standard Polya urn: a ball is drawn uniformly at random and is replaced together with a ball of the same color (reinforcement). On the right, the transitions for the number of balls of each color over the first steps of the process.}
 \end{center}
 \end{figure}

\begin{proposition}[Convergence of proportions] \label{prop:polyaunif} In the standard Polya urn started with 1 ball of each color, the proportion of red balls converges towards a uniform random variable on $[0,1]$.
\end{proposition}
\noindent \textbf{Proof.}  It is straightforward to check that $R_n/(R_n+B_n)$ is a bounded martingale (for the canonical filtration) which thus converges almost surely towards a limiting proportion $U  \in [0,1]$.  An easy induction on $n \geq 1$ shows that $R_{n}$ is uniformly distributed over $\{1,2, \dots , n\}$ and so
$$ \frac{R_{n}}{R_{n}+B_{n}} = \frac{R_{n}}{n+1} \xrightarrow[n\to\infty]{a.s.} U \sim \mathrm{Unif}[0,1].$$ In the next chapter, we will see another proof of this result based on continuous time techniques. \qed \bigskip 

\begin{exo}[Asymmetric starting configuration] \label{exo:polyaassy} Compute the law of the limiting proportion of red balls when the Polya urn starts from $R_{0}=a$ and $B_{0} = N_{0}-a$ balls.
\end{exo}

\noindent \textbf{Proof of Theorem \ref{thm:CVPDps}.} The above discussion, together with Proposition \ref{prop:polyaunif} translated in the framework of the theorem, shows the almost sure convergence $ n^{-1} \cdot K_{1}( \boldsymbol{\sigma}_{n}^{ \mathrm{cr}}) \to  U_{1}$ where $U_{1}$ is uniform over $[0,1]$. However, it is easy to see that conditionally on the values $$ u_{k} = \inf \left\{ t \geq 0 : 1+ \sum_{i\geq 2} K_{i}( \boldsymbol{\sigma}_{t}^{ \mathrm{cr}}) = k\right\},$$ the restricted process $( K_{i}( \boldsymbol{\sigma}_{u_{k}}^{ \mathrm{cr}}) : i \geq 2)_{k \geq 1}$ has the law of a Chinese restaurant process (thus  independent of $U_{1}$). By successive applications of the above reasoning we deduce that 
$$  \frac{K_{i}( \boldsymbol{\sigma}_{n})}{n} \xrightarrow[n\to\infty]{a.s.} (1-U_{1}) \cdots (1 - U_{i-1}) U_{i},$$ for independent random variables $U_{i}: i \geq 1$ uniformly distributed on $[0,1]$ as desired. \qed \medskip

\section{Degrees} \label{sec:degreeRRT}
In this section, we study the degrees of the vertices in $T_n$. More precisely, for $0 \leq i \leq n$ the \textbf{outdegree} (number of children) of  \raisebox{.5pt}{\textcircled{\raisebox{-.9pt} {$i$}}} in $T_{n}$ will be denoted by $$\mathrm{deg}^{+}_{T_{n}}( \raisebox{.5pt}{\textcircled{\raisebox{-.9pt} {$i$}}}) = \#\{ i < j \leq n : \noeud{i} \sim \noeud{j} \mbox{ in } T_{n}\}.$$

\begin{figure}[!h]
 \begin{center}
 \includegraphics[width=8cm]{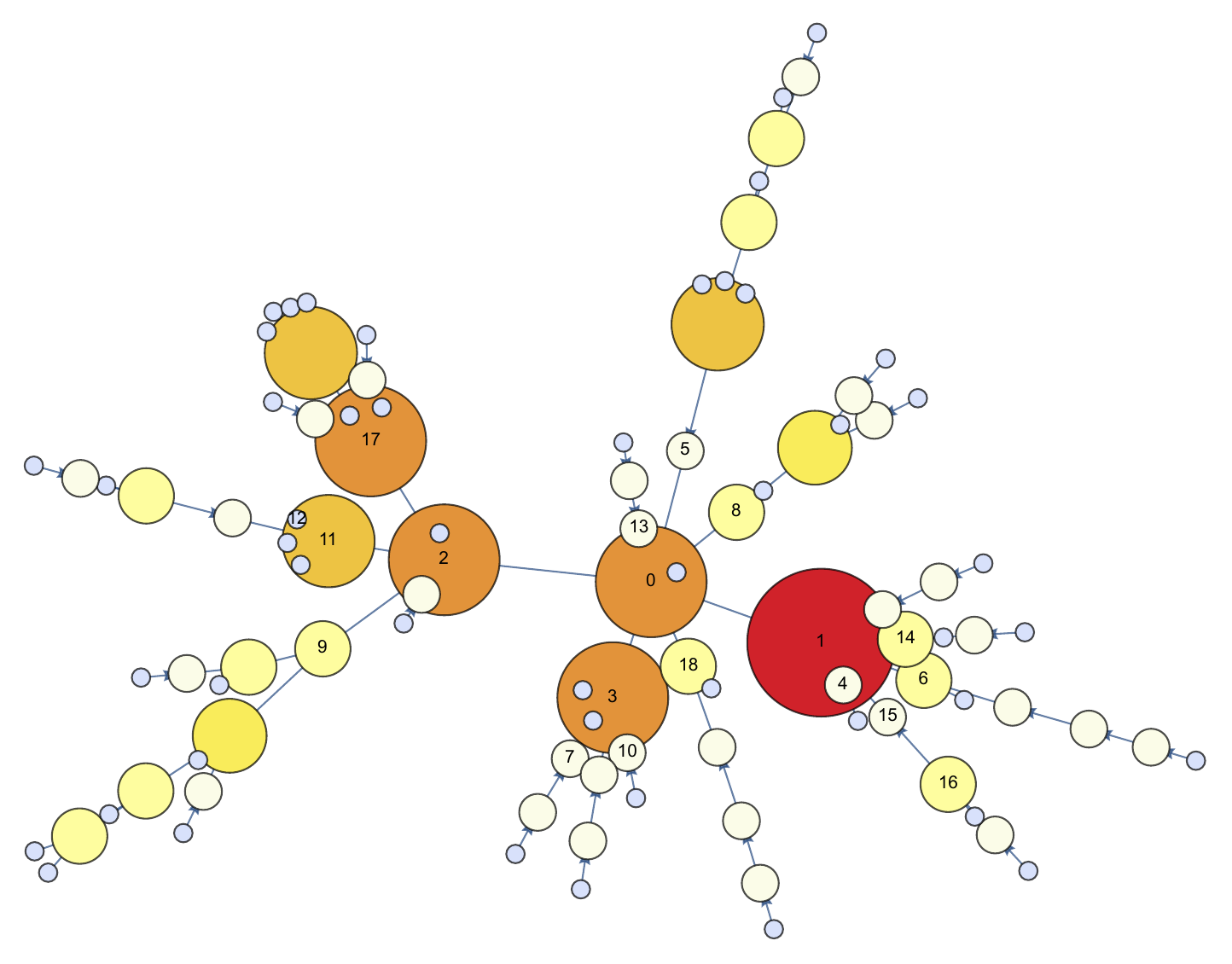}
 \caption{Simulation of $T_{100}$ where the size and color of vertices illustrate their degrees. The first $20$ vertices have their labels displayed.}
 \end{center}
 \end{figure}

\subsection{Degree of  fixed vertices}
 By construction, for any $i \geq 0$ fixed, we  have  \begin{eqnarray}\label{eq:sumbernoulli} \left(\mathrm{deg}^{+}_{T_{n}}( \raisebox{.5pt}{\textcircled{\raisebox{-.9pt} {$i$}}}) : n \geq  0 \right) =  \left( \sum_{k=i+1}^n B_{k} : n \geq 0 \right),  \end{eqnarray} where the Bernoulli random variables $ B_{k} \sim \mathrm{Ber}(1/k)$ are independent and of parameter $1/k$ for $k \geq 1$. Since $\sum_{k \geq 1} \frac{1}{k} = \infty$, the Borel--Cantelli lemma implies that the (out)degree of any vertex $\raisebox{.5pt}{\textcircled{\raisebox{-.9pt} {$i$}}}$ in $T_n$ tends to $\infty$ a.s.\ as $n \to \infty$. Also, by the coupling of the preceding section (or using Theorem \ref{thm:fellercoupling}) we deduce that for any $n \geq 1$ $$\mathrm{deg}^{+}_{T_{n}}( \raisebox{.5pt}{\textcircled{\raisebox{-.9pt} {$0$}}})  \overset{(d)}{=}  \mathcal{C}_n,$$ where we recall from Proposition \ref{thm:numbercycles} that $ \mathcal{C}_{n}$ is the law of the number of cycles in a random uniform permutation $  \boldsymbol{\sigma}_{n} \in \mathfrak{S}_{n}$ (with the CRP coupling, we have $\mathrm{deg}^{+}_{T_{n}}( \raisebox{.5pt}{\textcircled{\raisebox{-.9pt} {$0$}}}) = \# \mathrm{Cycles}(\boldsymbol{\sigma}_n^{ \mathrm{cr}})$). In particular, we deduce from Proposition \ref{thm:numbercycles} that  for each $i_0\geq 0$ fixed we have 
   \begin{eqnarray} \label{eq:smalldegreelog}  \frac{\mathrm{deg}^{+}_{T_{n}}( \raisebox{.5pt}{\textcircled{\raisebox{-.9pt} {$i_0$}}})}{\log n} \xrightarrow[n\to\infty]{( \mathbb{P})} 1,  \end{eqnarray}
and we will see later (Proposition \ref{prop:degreeroot}) that the convergence actually holds almost surely. 
\subsection{Empirical degree distribution}  \label{sec:empiricalRRT}
Let us now focus on the \textbf{empirical degree distribution} in $T_{n}$: We know from \eqref{eq:smalldegreelog} above that the vertices with small labels typically have a logarithmic degree, but as in any tree with $n+1$ vertices, the mean degree in $T_{n}$ is equal to $ \frac{2n}{n+1} \to 2$ as $n \to \infty$. So there must be (a lot) of vertices with small degrees. More precisely, we let $ \mu_{n} $ be the (random) empirical distribution of the out-degrees defined by 
$$ \mu_{n} =  \frac{1}{n+1}\sum_{i=0}^{n} \delta_{ \mathrm{deg}_{T_{n}}^{+}( {\footnotesize \noeud{i}})}.$$
It turns out that for large $n$'s the empirical degree distribution converges towards a deterministic distribution (a stronger version will be proved in Section \ref{sec:rrtfromyule}):
   \begin{proposition}[Convergence of the empirical degree distribution] \label{prop:empirical}The empirical distribution of the out-degrees in $T_{n}$ converges in probability towards the critical geometric distribution of parameter $1/2$, i.e.~for each $k_{0} \geq 0$ we have 
   $$ \mu_{n}(\{k_0\})  \xrightarrow[n\to\infty]{( \mathbb{P})} 2^{-k_{0}-1}.$$
   \end{proposition}
   \begin{exo} Prove the above proposition by computing the first and second moment of $\mu_{n}( \{k_{0}\})$.
   \end{exo}
   

\subsection{Maximal degree}
 By \eqref{eq:smalldegreelog}, the typical degree of vertices with fixed label is of order $ \log n$. Actually, the \textbf{largest} degree is much larger and is close to what would be the maximum of $n$ i.i.d.~critical geometric random variables, or in other words, as if we were sampling $n$ i.i.d.\ degrees distributed according to the limiting empirical degree distribution computed in Proposition \ref{prop:empirical}:

\begin{theorem}[Devroye \& Lu] Let $ \mathrm{MaxDegree}(T_n) = \max\{ \mathrm{deg}^{+}_{T_{n}}( \raisebox{.5pt}{\textcircled{\raisebox{-.9pt} {$i$}}}) : 1 \leq i \leq n\}$ be the largest vertex (out)-degree in $T_{n}$. \label{prop:maxdegreeRRT}  Then we have 
$$ \frac{\mathrm{MaxDegree}(T_n)}{ \log n} \xrightarrow[n\to\infty]{a.s.}  \frac{1}{\log 2} \approx 1.44\dots$$
\end{theorem}
\noindent \textbf{Teasing for the proof.} The convergence in probability can be approached using the first and second moment method, but the computations are really technical... A neat proof goes through a representation of the RRT in continuous time (a.k.a.~Rubbins/Athreya  construction) via a Yule process, see Chapter \ref{chap:AK}. \qed 

 \section{Height}
 We now turn to the study of   \textbf{heights} in $T_{n}$, i.e.~the distances of the vertices to the root $\noeud{0}$ in $T_{n}$.  More precisely, for $ n \geq 0$, we denote by $H_{n}$ the height (or generation) of the vertex \raisebox{.5pt}{\textcircled{\raisebox{-.9pt} {$n$}}} in the random recursive tree $T_{m}$ for $m \geq n$ (the definition does not depend on $m \geq n$ since the vertex $\noeud{n}$, once attached, is fixed in $T_{m}$ for $m \geq n$). 

 \subsection{Typical height}
Clearly, the height of the first few vertices $H_{1}, H_{2}, \dots$ are small and are given by the first stages in the construction of $(T_{n} : n \geq 0)$. We shall prove below the surprising fact that $H_{n}$ has the same law as $ \mathrm{deg}^{+}_{T_{n}}( \noeud{0})$, which is the law $ \mathcal{C}_{n}$ of the number of cycles in a uniform permutation $\boldsymbol{\sigma}_{n}$:
 \begin{proposition} For any $n \geq 0$ we have $H_{n} = \mathcal{C}_{n}$ in law.
\end{proposition}
\begin{remark} The above proposition shows that for \textbf{fixed} $n \geq 0$, we have $H_{n} = \mathcal{C}_{n} = \mathrm{deg}^{+}_{T_{n}}( \raisebox{.5pt}{\textcircled{\raisebox{-.9pt} {$0$}}})$ in law, but the previous equality does not hold in terms of process in $n \geq 0$:
$$ ( H_{n} : n \geq 0) \quad \ne \quad  \left(\mathrm{deg}^{+}_{T_{n}}( \raisebox{.5pt}{\textcircled{\raisebox{-.9pt} {$0$}}}) : n \geq 0\right).$$
Indeed, the process in the right-hand side is non-decreasing and tends to $+\infty$ a.s. (see  \eqref{eq:smalldegreelog}), while the first one does not: because the degree of \raisebox{.5pt}{\textcircled{\raisebox{-.9pt} {$0$}}} is unbounded as $n \to \infty$, there are infinitely many vertices grafted on \raisebox{.5pt}{\textcircled{\raisebox{-.9pt} {$0$}}} and so infinitely many values for which $H_{n} =1$.
\end{remark}

\begin{figure}[!h]
 \begin{center}
 \includegraphics[width=10cm]{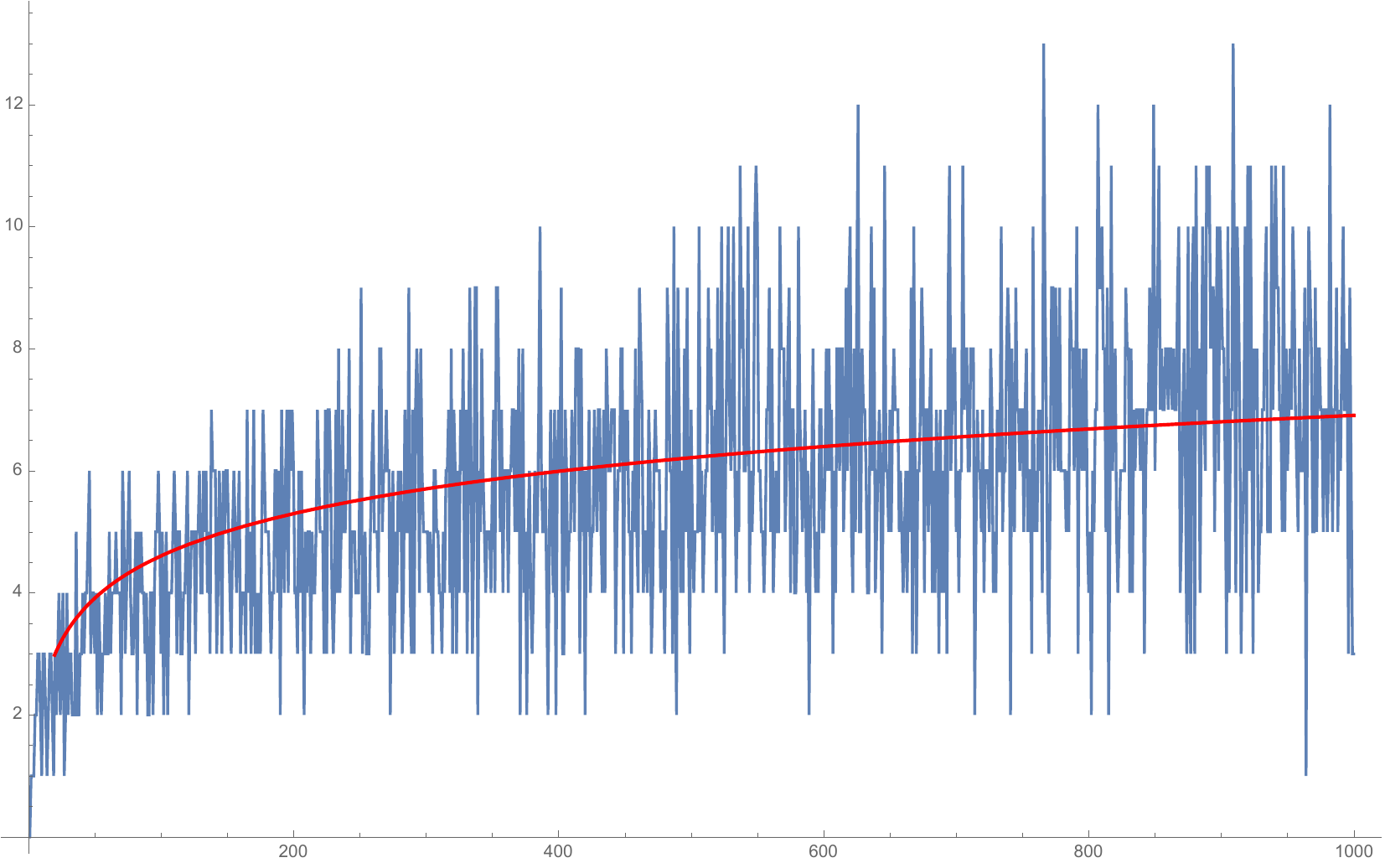}
 \caption{Plot of a simulation of the successive heights $(H_{i} : 0 \leq i \leq 1000)$ against the $\log$ function (in red).}
 \end{center}
 \end{figure}
 
 \noindent \textbf{Proof.} Since \raisebox{.5pt}{\textcircled{\raisebox{-.9pt} {$n$}}} is grafted to a uniform node with label $ <n$ we have the following recursive distributional equation: $H_{0}=0$ and for $n \geq 1$
  \begin{eqnarray}  \label{eq:recursiveheight} H_{n} \overset{(d)}{=} 1 + H_{U_{n-1}},  \end{eqnarray}
 where in the right-hand side $U_{n-1} \in \{0,1,2, \dots , n-1\}$ is independent of the RRT defining $(H_{n} : n \geq 0)$. This type of equality is called a \textbf{recursive distributional equation}. Actually, we saw in the proof of Theorem \ref{thm:PDforcycles} that in a uniform permutation $  \boldsymbol{\sigma}_{n}$, the size $V_{n}$ of the cycle containing $1$ is uniformly distributed over $\{1,2, \dots , n\}$ and conditionally on it the remaining (relabeled) permutation is uniform over $ \mathfrak{S}_{n-V_{n}}$. In particular, $ \mathcal{C}_{n}$ satisfies the same recursive distributional equation as in \eqref{eq:recursiveheight}: $$ \mathcal{C}_{n} \overset{(d)}{=} 1 + \tilde{\mathcal{C}}_{n-V_{n}}, $$
 where on the right-hand side $ (  \tilde{\mathcal{C}}_{i} : i \geq 0)$ are independent variables of law $ \mathcal{C}_{i}$ and also independent of the uniform variable $V_{n} \in \{1,2, \dots , n\}$. With the convention $ \mathcal{C}_{0}=0$, this is sufficient to show that $ \mathcal{C}_{n}$ and $H_{n}$ have the same law since those recursive equations \eqref{eq:recursiveheight} characterize their laws. \qed 
 
 \bigskip

Proposition \ref{thm:numbercycles} directly implies a central limit theorem:
$$ \frac{H_{n}- \log n}{ \sqrt{\log n}} \xrightarrow[n\to\infty]{(d)} \mathcal{N}(0,1),$$ and a weak law of large number $H_{n} / \log n \to 1$ in probability (but not almost surely). 

\subsection{Maximal height}
 
As in the case of vertex degrees, the maximal height $$ \mathrm{Height}(T_n) := \max\{ H_i : 0 \leq i \leq n\}$$ of $T_n$ is much larger than the typical height and is also asymptotically the same as if the heights of different points were independent, that is comparable to $\sup\{  \mathcal{C}_{n}^{(k)} : 1 \leq k \leq n \}$ for independent random variable $ \mathcal{C}_{n}^{(i)}$ of law described in Proposition \ref{prop:degreeroot}.\medskip

\begin{theorem}[Pittel] \label{thm:heightRRT} We have $$ \frac{ \mathrm{Height}( T_n) }{\log n} \xrightarrow[n\to\infty]{ a.s.} \mathrm{e}.$$
\end{theorem}

\noindent \textbf{Proof.} See  Exercise \ref{exo:upperheight} below for the upper bound using the first moment method. The lower bound can in principle be approached by the second moment method  but yield to a very intricate proof. We shall prove this theorem using the continuous time embedding technique in Chapter \ref{chap:AK}. \qed \medskip 

\begin{exo}[Upper bound using Poisson approximation] \label{exo:upperheight} For $p \in (0,1)$ and $\alpha >0$ denote by $ \mathrm{Ber}(p)$ a Bernoulli variable with expectation $p$ and by $ \mathfrak{P}( \alpha)$ a Poisson variable with expectation $\alpha$.
\begin{enumerate}
\item Show that $ \mathrm{Ber}(p) \leq \mathfrak{P}(- \log (1-p))$ for the stochastic order and that 
$$ \mathrm{d_{TV}} ( \mathrm{Ber}(p), \mathfrak{P}(p)) = O (p^2), \quad \mbox{as } p \to 0,$$ where $ \mathrm{d_{TV}}$ is the total variation distance.
\item Deduce that $ H_{n}$ is stochastically dominated by $1$ plus a Poisson variable with expectation $- \sum_{i=2}^{n}\log \left(1-  \frac{1}{i} \right)$.
\item Use \eqref{lem:LDpoisson} to conclude that for all $ \varepsilon>0$ we have  $ \mathbb{P}( \mathrm{Height}(T_{n}) > (\mathrm{e} + \varepsilon) \log n) \to 0$ as $n \to \infty$.
\item Prove that $\mathrm{Height}( T_n) \leq ( \mathrm{e} + \varepsilon) n$ eventually, a.s.
\end{enumerate}
\end{exo}

\paragraph{Bibliographical notes.} The random recursive tree and random uniform permutations over the symmetric group are both very well studied in probability theory. Standard references are Feller  \cite{feller1945fundamental} and  the Saint-Flour lectures of Pitman \cite{Pit06} in particular Section 3.1 or the renowned  \cite{Flajolet:analytic}. Theorem \ref{prop:maxdegreeRRT} is due to Devroye \& Lu \cite{devroye1995strong} and Theorem \ref{thm:heightRRT} to Pittel \cite{pittel1994note}. See \cite{smythe1995survey,goldschmidt2005random,baur2014cutting} for more results about the random recursive tree.  \medskip 

\noindent{\textbf{Hints for Exercises.}}\ \\
\noindent Exercise \ref{exo:polyaassy}: Show that $Z_n\,=\, \frac{R_n(R_n+1)\ldots(R_n+k-1)}{(n+N_0-1)(n+N_0 )\ldots(n+N_0+k-2)}\,$ is a martingale for the canonical filtration and deduce the moments of the limiting proportion of red balls. See Exercise \ref{exo:afresh} for a calculus-free approach.\\
Exercise \ref{exo:upperheight}: For the last question, use the polynomial decay of $ \mathbb{P}( \mathrm{Height}(n) \geq ( \mathrm{e}+ \varepsilon) n)$ obtained in $2)$ along a subsequence $ n = c^{k}$ for $c >1$. Conclude using the fact that $ n \mapsto \mathrm{Height}(T_{n})$ is increasing.

\chapter{Continuous-time branching processes}
\label{chap:AK}
\hfill Randomize to make it simpler!\bigskip

In this chapter, we theorize the Poissonization technique which amounts to transforming a discrete-time process into a continuous-time version which possesses more independence properties. This will be particularly useful for urn processes and random tree growth mechanisms.

 \section{Continuous-time branching trees}
 
Let us first recall the memorylessness property of exponential variables, which will be the crux of the continuous-time embedding technique.  

 \subsection{Properties of exponential laws}
 
In the following, for $ \alpha >0$ we denote by $ \mathcal{E}(\alpha)$ the exponential distribution of expectation $1/ \alpha$, i.e.~given by
$$ \mathcal{E}( \alpha) = \alpha \cdot \mathrm{e}^{- \alpha x} \mathbf{1}_{x >0} \mathrm{d}x,$$ we shall say that $\alpha$ is the \textbf{rate} of the exponential, since by the \textbf{memorylessness} property of the exponential distribution if $X \sim \mathcal{E}(\alpha)$ we have 
 \begin{eqnarray} \label{eq:expomemory} \mathbb{P}( X \in [x , x + \mathrm{d}x] \mid X \geq x) = \alpha \mathrm{d}x,  \end{eqnarray} or equivalently that conditionally on $\{ X \geq t\}$ the variable  $X-t$ has distribution $ \mathcal{E}(\alpha)$.  Recall also that the memorylessness property  is characteristic of the exponential and geometric laws: 
 \begin{exo}[Memorylessness] \label{exo:memory} Let $X$ be a  random variable with values in $ \mathbb{R}_+$ so that for every $a,b  \in \mathrm{Supp}( \mathcal{L}(X))$ we have   $$ \mathbb{P}(X >a+b \mid X >b)= \mathbb{P}(X >a).$$
 Show that $X$ is either an exponential or a multiple of a geometric random variable.
 \end{exo}

\paragraph{Choosing using clocks.} Consider $X_{1} \sim \mathcal{E}( \alpha_{1}), \dots , X_{k}\sim\mathcal{E}( \alpha_{k})$ a family of $k$ independent exponential variables of parameters $ \alpha_{1}, \dots , \alpha_{k}$ and denote by $ \mathscr{M} = \min\{X_{i} : 1 \leq i \leq k\}$ and by $ \mathscr{J}\in \{1, \dots , k \}$ the index at which this minimum is attained.  Then we have:
\begin{proposition}[Choosing with clocks]\label{prop:expomem} The index $\mathscr{J}$ is almost surely well-defined (there is no tie) and we have 
$$  \mathscr{M} \sim \mathcal{E}( \alpha_{1}+  \dots + \alpha_{k}), \quad \mathbb{P}(\mathscr{J} = i) = \frac{\alpha_{i}}{\alpha_{1} + \dots + \alpha_{k}},$$
and conditionally on $\{\mathscr{J}, \mathscr{M}\}$, the remaining variables $  X_{1} - \mathscr{M}, \dots , \widehat{ X_{\mathscr{J}}- \mathscr{M}}, \dots , X_{k}- \mathscr{M}$ are independent and of laws $  \mathcal{E}( \alpha_{1}) , \dots ,\widehat{ \mathcal{E}( \alpha_{\mathscr{J}})}, \dots , \mathcal{E}( \alpha_{k})$.
\end{proposition}
\noindent \textbf{Proof.}  This can can heuristically be explained as follows: by the memorylessness property of the exponential laws \eqref{eq:expomemory}, the variable $  \mathscr{M}$ must follow an exponential law with rate $\alpha_{1} + \dots + \alpha_{k}$ and given that $ \mathscr{M} \in [x, x+ \mathrm{d}x]$, the probability that $ \mathcal{E}( \alpha_{i})$ is the smallest is just proportional to the rate i.e. 
$$ \mathbb{P}(\mathscr{J} = i) = \frac{\alpha_{i}}{\alpha_{1} + \dots + \alpha_{k}}.$$
The remaining statement follows by the memorylessness property. More formally,  since the variables are independent and have a density with respect to the Lebesgue measure, there are a.s. pairwise distincts and so $\mathscr{J}$ is well-defined. Furthermore, for  any positive function $ \phi : \mathbb{R}_+ \times (\mathbb{R}_+)^{k-1}$ we have 
  \begin{eqnarray*} && \mathbb{E}\left[ \phi\Big( \mathscr{M} ; X_{1} - \mathscr{M}, \dots , \widehat{ X_{j})- \mathscr{M}}, \dots , X_{k}- \mathscr{M}\Big) \mathbf{1}_{\mathscr{J} = j}\right] \\ 
  &=& \int_0^\infty \mathrm{d} s_j  \alpha_j \mathrm{e}^{- \alpha_j s_j} \int_{s_j}^\infty  \left(\prod_{i \ne j}\mathrm{d}s_i \alpha_i \mathrm{e}^{-\alpha_i s_i}\right) \phi( s_j ; (s_1-s_j), \dots , \widehat{(s_j-s_j)}, \dots ,( s_k -s_j))\\
  &\underset{ \begin{subarray}{c} m = s_j \\ \tilde{s}_i = s_i-s_j \end{subarray}}{=}& \frac{\alpha_j}{ \sum_i \alpha_i} \int_0^{\infty} \mathrm{d}m \  (\sum_i \alpha_i) \mathrm{e}^{- m ( \sum_i \alpha_{i})} \prod_{i \ne j} \int_0^\infty \mathrm{d} \tilde{s}_i \ \ \phi (m ; \tilde{s}_1, \dots , \widehat{\tilde{s}_j}, \dots , \tilde{s}_k),  \end{eqnarray*}
and this proves the claim.
 \qed \medskip 

A consequence of the above proposition is that if we want to sample from $\{1,2, \dots, k \}$ proportionally to some weights $\alpha_{1}, \dots , \alpha_{k}$; one way, which may seem strange at first glance, is to sample independent exponential clocks $ (X_{i} \sim \mathcal{E}(\alpha_{i}) : 1 \leq i \leq k)$ and consider the index of the first clock that rings.  The advantage of this point of view is that by Proposition \ref{prop:expomem}, the exponential clocks that have not rung can be further used (after subtracting the minimum) to sample according to $(\alpha_{i})$ the remaining items as well! 

We shall use many times the well-know extremal statistics of exponential distribution:
\begin{lemma}[Gumbel distribution] \label{lem:gumbel} Let $ X_1, \dots , X_n$ be i.i.d.~variables of law $ \mathcal{E}(1)$. We denote their maximum by $  \mathrm{M}_n = \max\{ X_i : 1 \leq i \leq n\}$. Then we have the following convergence in distribution towards the Gumbel\footnote{\raisebox{-5mm}{\includegraphics[width=1cm]{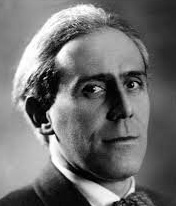}} Emil Julius Gumbel (1891-1966), German} distribution:
 $$ \mathrm{M}_n-\log n \xrightarrow[n\to\infty]{(d)} G \overset{(d)}{=} \mathrm{e}^{- \mathrm{e}^{-x}-x} \mathrm{d}x.$$
\end{lemma}
Remark two useful observations:  First, if $G$ has the Gumbel distribution then $ \mathbb{P}(G \leq x) = \mathrm{e}^{- \mathrm{e}^{-x}}$ so that $ \mathrm{e}^{-G}$ has law $ \mathcal{E}(1)$. Second, by iterating Proposition \ref{prop:expomem} the variable $ \mathrm{M}_{n}$ has the law same as 
$$ \mathrm{M}_{n} \overset{(d)}{=} \mathcal{E}(n) \overset{\coprod}{+} \mathcal{E}(n-1)  \overset{\coprod}{+}\dots \overset{\coprod}{+} \mathcal{E}(1),$$ where the variables are independent. We deduce that the right-hand side of the last display satisfies the same convergence as stated in the lemma.

\noindent \textbf{Proof.} For $x \in \mathbb{R}$, if $X \sim \mathcal{E}(1)$ we have  $$ \mathbb{P}( \mathrm{M}_n \leq \log n +x) = ( \mathbb{P}( X \leq \log n +x))^n= \left(1- \frac{\mathrm{e}^{-x}}{n}\right)^n \xrightarrow[n\to\infty]{} \mathrm{e}^{- \mathrm{e}^{-x}}.$$ \qed

\begin{exo}[Hide and seek]  \label{exo:hideseek}We sample i.i.d.~random variables in a finite set $  \mathbb{X}=\{x_1, \dots , x_n\}$ according to some weights $\alpha_1=1, \alpha_2=1\dots, \alpha_{n-1}=1$ and $\alpha_n=2$  until all elements of $ \mathbb{X}$ have been seen in the sequence. What is the probability that the last element unseen is $x_n$?
\end{exo}

\subsection{Continuous branching trees and their discrete associated Markov chains} \label{sec:AK}
Let $p \in \{1,2, \dots \} \cup \{\infty\}$ and denote by $ \mathfrak{y} = \{1,2, \dots , p\}$ if $p< \infty$  or $ \mathfrak{y} = \mathbb{Z}_{\geq 0}$ be the set of \textbf{discrete types}.  To ease notation, we shall identify the space $(\mathbb{Z}_{\geq 0})^{ \mathfrak{y}}$ with the space of discrete measures   $\sum_{i \in \mathfrak{y}} x_i \delta_i$ with $x_i \in \mathbb{Z}_{\geq 0}$,  for example $(2,0,0,\dots)$ we be written $2\delta_1$. For each type $ i \in \mathfrak{y}$, we are given $\alpha_i >0$ a positive \textbf{weight} and  an \textbf{offspring distribution}  $ (\mu_{i} :  i \in \mathfrak{y})$ over $(\mathbb{Z}_{\geq 0})^{ \mathfrak{y}}$. Finally, let us fix $ \mathbf{x} \in (\mathbb{Z}_{\geq 0})^{ \mathfrak{y}}$ a non-zero \textbf{starting configuration}.

\begin{figure}[!h]
 \begin{center}
 \includegraphics[width=13cm]{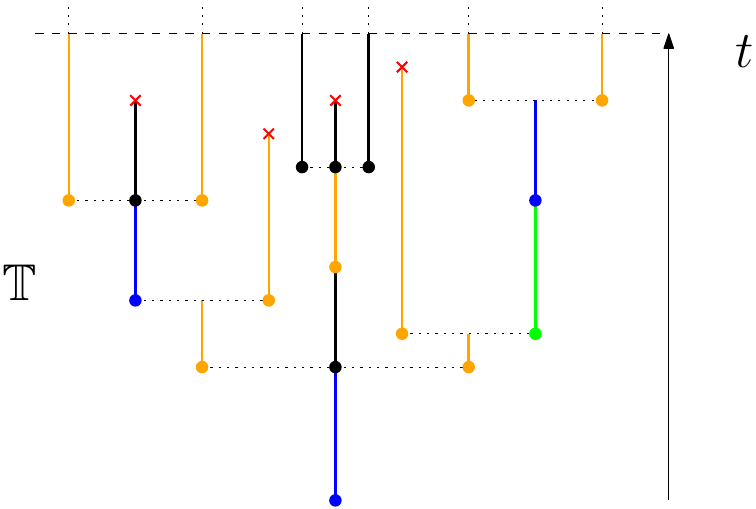}
 \caption{Illustration of the construction of the random tree $ \mathbb{T}$ starting from a single blue individual (the colors represent types of particles): each particle of type $i$ lives for an exponential time of expectation $1/ \alpha_i$, then dies and gives birth to new particles according to the distribution $\mu_i$. All those samplings are made independently of each other. The red crosses represent deaths with no birth. \label{fig:mtgw}}
 \end{center}
 \end{figure}

We now create a random genealogical tree, more precisely a forest of trees, as follows. Under $ \mathbb{P}_{\sum_{i \in \mathfrak{y}} x_i \delta_i}$ the random forest $ \mathbb{F}$ (we shall write $ \mathbb{T}$ if there is a single tree, i.e.~if $ \mathbf{x} = \delta_{i_0}$ for some $ i_0 \in \mathfrak{y}$) is the genealogical forest  of a cloud of particles starting with $x_j$ particles of type $j$, and where subsequently each particule of type $i \in \mathfrak{y}$ behaves independently of the others and lives an exponential time $ \mathcal{E}(\alpha_i)$ of rate  $\alpha_i$ before dying and giving birth to a cloud of particles sampled according to $\mu_{i}$ (independently of the past and of the other particles). The  trees in $ \mathbb{F}$ are  locally finite random rooted (but non-planar) trees with edge lengths as depicted on Figure \ref{fig:mtgw}.

In the case of a single ancestor, it is possible to make a formal definition of $ \mathbb{T}$ as a plane tree with edge lengths, by ordering the children of each particle from left-to-right, so that each particle alive at some time corresponds to a vertex of Ulam's tree. The type and the life time of particles are then additional decorations. We will however not bother to make such construction in general and mostly rely on the intuition of the reader. 
Several limit theorems are available in the literature for the number of particles of each type living at time $t$ in $ \mathbb{T}$, but for the purpose of these lecture notes we shall only deal with the most basic examples, namely Poisson processes and Yule trees, see Section \ref{sec:yules}. But before that, let us connect those random continuous trees to discrete Markov chains using properties of the exponential distributions.

If  $ \mathbb{F}$ is a random forest of law $ \mathbb{P}_{ \mathbf{x}}$ as  above, consider $0 = \tau_0 < \tau_1 < \tau_2 < \dots$ the jump times\footnote{since the exponential distribution has a density and since all particles' life times are independent, it is easy to see that the jump times are a.s. distinct. But we do not exclude the possibility that the jump times accumulate.}, i.e.~the times when a particle dies in $ \mathbb{F}$ and gives birth to a new cloud of particles (possibly empty). Let also introduce $( \mathbb{X}_k : k \geq 0)$  the $ (\mathbb{Z}_{\geq 0})^{ \mathfrak{y}}$-valued process made of the number of particles of each type $i \in \mathfrak{y}$ at time $ \tau_k$.
\begin{lemma}[Athreya--Karlin]  \label{lem:AK} Under $ \mathbb{P}_{ \sum_i x_i \delta_i}$, the process $( \mathbb{X}_k: k \geq 0)$ is a Markov chain starting from $ (x_i : i \in  \mathfrak{y})$ and with transitions described informally as follows: conditionally given 
$ \mathbb{X}_{n} = (x_{n}^{(i)} :  i \in \mathfrak{y})$ we choose a uniform particle of type $i_{0}$ with probability 
$$ \frac{x_{n}^{{(i_{0})}} \cdot \alpha_{i_{0}}}{ \sum_{i \in \mathfrak{y}} x_{n}^{(j)} \alpha_{j}},$$
 then this particle dies and creates new particles $( z^{(i)} : i \in \mathfrak{y})$ with law $\mu_{i_{0}}$. More formally, for any positive function $f : (  \mathbb{Z}_{\geq 0})^p \to \mathbb{R}_+$ we have 
 $$ \mathbb{E}[f( \mathbb{X}_{n+1}) \mid \mathcal{F}_{n}] = \sum_{ i_{0}=1}^{p}  \frac{\alpha_{i_{0}}  \mathbb{X}_{n}^{(i_{0})}}{ \sum_{j } \mathbb{X}_{n}^{(j)} \alpha_{j}} \sum_{ z \in (\mathbb{Z}_{\geq 0})^{p}} \mu_{i_{0}}(z) \cdot f( \mathbb{X}_{n} -\delta_{i_0} + z).$$
 \end{lemma}
 \noindent \textbf{Proof.}  Let us prove by induction on $k \geq 0$ that at time $\tau_k$, conditionally on the past up to time $ \tau_k$, the particles alive at time $\tau_k$ all carry independent exponential clocks of weight $\alpha_i$ for a particle of type $i$.  This is true for $k =0$ and propagates easily by  Proposition \ref{prop:expomem}.  In particular, by Proposition \ref{prop:expomem} again, conditionally on the types of the particles at time $\tau_k$, the next particle to die is chosen proportionally to the rate $\alpha_i$ of its type $i$ and reproduces according to $\mu_i$.  \qed\medskip

 We shall see in Section \ref{sec:examplesAK} several examples of discrete chains which are more efficiently studied via their continuous-time analogs, but before that, let us study the most fundamental cases where particles reproduce at constant rate into a fixed number of new particles.


\section{Yules trees} \label{sec:yules}
In this section, we shall focus on a very special case of continuous-time branching process where there is only one type of particle which reproduce at rate $1$ into exactly $k \geq 1$ particles. When $k=1$ this corresponds to a vanilla constant rate  \textbf{Poisson process} on $ \mathbb{R}_+$ and when $k \geq 2$ we speak of (random) \textbf{Yule} trees.
\subsection{$k=1$ and Poisson process}
\label{sec:poissoncounting}
Fix here $p=1$ (monotype) and $\mu_{1} = \delta_{\delta_1}$, i.e.~when a particle dies, it gives rise to a single particle. In terms of set of particles, nothing is happening. But the temporal death counting process gives the link between exponential variables and Poisson processes. More precisely, consider $ X_{1},  X_{2}, \dots $ a sequence of i.i.d.~exponential variables of rate $1$ and build the counting process  for $t \geq 0$
$$  \mathfrak{P}(t) =  \sup\{i \geq 0 : X_{1} + \dots + X_{i} \leq t\}.$$
This random c\`adl\`ag process turns out to be a Poisson counting process and this connection is the standard way to prove \eqref{eq:lawll}:
\begin{proposition}[Standard Poisson] \label{prop:standardpoisson} For any $0=t_{0}\leq t_{1} \leq t_{2}\leq t_{3} \leq \cdots$ , the variables $  \mathfrak{P}({t_{i+1}})-\mathfrak{P}({t_{i}})$ for $0 \leq i \leq k-1$ are independent and of law 
$$ \mathfrak{P}({t_{i+1}})-\mathfrak{P}({t_{i}}) \sim \mathrm{Poisson}(t_{i+1}-t_{i}).$$
\end{proposition}
\noindent \textbf{Proof.} This is a very classical result whose proof can be found in many textbooks. Let us however sketch the arguments: The independence and stationary of the increments follows by the loss of memory property applied recursively at times $t_{i-1}, t_{i-2}, \dots , t_{1}$. To prove that $\mathfrak{P}({t})$ follows a Poisson distribution one can notice that from Proposition \ref{prop:expomem}  we can write 
$$ \forall t \geq 0,\quad  \mathfrak{P}(t) = \sum_{i=1}^{n} \mathfrak{P}^{(i)}({t/n}),$$ where $\mathfrak{P}^{(i)}(\cdot)$ are i.i.d.~copies of $\mathfrak{P}(\cdot /n)$ i.e.~of the process $ \mathfrak{P}$ constructed with exponentials of mean $n$. For fixed $t>0$, when $n \to \infty$ notice that we have $$\mathbb{P}(\mathfrak{P}(t /n) = 1) \sim  \frac{t}{n} \quad \mbox{ and }\quad \mathbb{P}(\mathfrak{P}(t /n) \geq 2) \leq   \frac{C_{t}}{n^{2}},$$ where $C_{t}>0$ is a positive constant. In particular, the total variation distance $\mathfrak{P}(t /n)$ and $ \mathrm{Bern}(t/n)$ is less than $ 2 C_{t}/n^{2}$ and we deduce that $$ \mathrm{d_{TV}}(\mathfrak{P}({t}), \mathrm{Bin}(n, t/n))  \leq \frac{2 C_{t}}{n^{2}} \cdot n \xrightarrow[n\to\infty]{}0,$$ and since $\mathrm{Bin}(n, t/n) \to \mathfrak{P}(t)$ in distribution we are done. \qed \bigskip

These two visions on the standard Poisson process are already very useful:

\begin{exo}For $n \geq 1$, let $(U_i : 1 \leq i \leq n-1)$ be i.i.d. uniform on $[0,1]$ and denote by $(V_i : 1 \leq i \leq n-1)$ their increasing rearrangement and put $V_0=0$ and $V_n=1$. Let $X_1, \dots , X_{n}$ be i.i.d.~r.v.~of law $ \mathcal{E}(1)$ and denote by $ \mathbb{X} = X_1 + \dots + X_{n}$. Show that 
$$ \left( V_{i}-V_{i-1} : 1 \leq i \leq n\right) \quad \overset{(d)}{=} \quad  \left( \frac{X_i}{ \mathbb{X}} : 1 \leq i \leq n \right).$$
\end{exo}


 \subsection{$k\geq 2$ and Yule trees}
\begin{figure}[!h]
  \begin{center}
  \includegraphics[width=15cm]{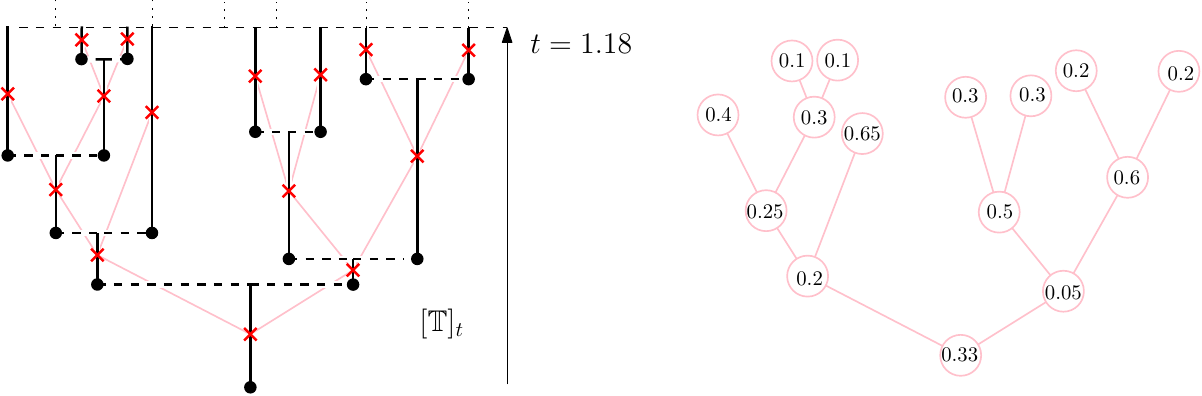}
  \caption{Illustration of the encoding of  (a restriction of) the Yule tree as a (finite) plane $k$-ary tree whose vertices carry positive numbers (in pink on the right). \label{fig:RRTbinary}}
  \end{center}
  \end{figure}

 Another special example of multi-type branching tree is given by setting $p=1$ (monotype), $\alpha_{1}=1$ to fix ideas,  and $\mu_{1} = \delta_{ k\delta_1}$ for some $k \geq 2$, i.e.~each particle creates $k$ new particles when dying. We then speak of the \textbf{Yule\footnote{\raisebox{-5mm}{\includegraphics[width=1cm]{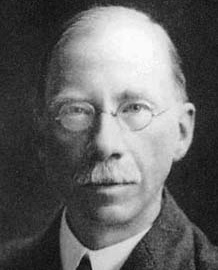}} George Udny Yule (1871--1951) British} tree} of order $k$. In other words, the discrete tree underlying $ \mathbb{T}$ under $ \mathbb{P}_{\delta_1}$ is the full $k$-ary tree whose edge lengths are i.i.d.~distributed according to $ \mathcal{E}(1)$. For later purposes, it will be useful to have a plane ordering of the tree. This can be obtained by starting with the infinite $k$-ary tree  whose vertex set is  $\bigcup_{n \geq 0}\{0,1, \dots , k-1\}^{n}$ and equip each of its vertices with an independent exponential r.v.\ with rate $1$ (the vertex lengths). In this correspondance, the vertices of $\bigcup_{n \geq 0}\{0,1, \dots , k-1\}^{n}$ are associated with the edges of the plane Yule tree $ \mathbb{T}$. For each $t \geq 0$, we denote by $ [ \mathbb{T}]_{t}$ the finite plane tree obtained by cutting $ \mathbb{T}$ at height $t$. By the same procedure as before, it can be seen as a finite plane tree whose vertices have either $0$ or $k$ children and whose vertices are decorated with positive lengths, see Figure \ref{fig:RRTbinary}. In the following, we shall always make such identification without further notice.
 
  We denote by $\# \partial [ \mathbb{T}]_{t}$ the number of leaves of $ [\mathbb{T}]_{t}$ and use 
 $$ \mathcal{Y}_{t}^{(k)} \equiv \# \partial [ \mathbb{T}]_{t},$$ as a short-hand notation.     In this case, the growth of the tree is very well understood since $( \mathcal{Y}^{(k)}_t : t \geq 0)$ is a Markov chain which makes positive  jumps of size $(k-1)$ with rate $ \mathcal{Y}^{(k)}_{t}$. We deduce that $ \mathrm{m}(t) := \mathbb{E}_{\delta_1}[ \mathcal{Y}^{(k)}_{t}]$ satisfies $ \mathrm{m}'(t) = (k-1) \mathrm{m}(t)$  and $ \mathrm{m}(0) =1$ under $ \mathbb{P}_{\delta_1}$ so that 
 $$ \mathbb{E}_{\delta_1}\left[ \mathcal{Y}^{(k)}_{t} \right] = \mathrm{e}^{(k-1) t}, \quad \mbox{ for }t \geq 0.$$ Combined with the Markov property, it follows that 
  \begin{eqnarray}  \label{eq:martingaleyule} \Big(\mathrm{e}^{-(k-1)t}\cdot \mathcal{Y}^{(k)}_{t} : t \geq 0\Big), \quad \mbox{ is a martingale}  \end{eqnarray} for the filtration made of the information up to time $t$ and so converges almost surely (this will be reproved in the following proposition). We can even identify its limit:

 \begin{proposition} \label{prop:yule} Let $( \mathcal{Y}^{(k)}_t : t \geq 0)$ be the counting process in a  Yule tree of order $ k\geq 2$ under $ \mathbb{P}_{\delta_1}$. Then we have 
 $$  \mathrm{e}^{- (k-1)t} \cdot \mathcal{Y}^{(k)}_t \xrightarrow[t\to\infty]{a.s.} \boldsymbol{\Gamma}\left (\frac{1}{k-1}, \frac{1}{k-1}\right),$$
  where $\boldsymbol{\Gamma}(\alpha, \beta)$ is a scaled Gamma random variable, i.e.~with law $ \frac{1}{\Gamma(\alpha)} x^{\alpha-1} \mathrm{e}^{-\beta x} \beta^\alpha  \, \mathrm{d}x \mathbf{1}_{x>0}$ (in particular a standard exponential when $\alpha=\beta=1$).
 \end{proposition}
 \noindent \textbf{Proof.} Let us first prove the proposition in the case $k=2$ for simplicity. Consider the jump times $0 = \tau_0 < \tau_1 < \tau_2 < \dots$ of the process $\mathcal{Y}^{(2)}$ so that we have $ \mathcal{Y}_{\tau_i} ^{(2)} = i+1$ deterministically. By the properties of exponential variables we know that $ (\tau_{i+1}-\tau_i : i \geq 0)$ are independent exponential random variables  with rate $i+1\geq 1$. We write $ \mathrm{h}_n = \sum_{k=1}^n \frac{1}{k}$ for the $n$th harmonic number. Clearly $( \tau_{n}- \mathrm{h}_{n}: n \geq 1)$ is a martingale bounded in $L^{2}$ since 
 $$ \sum_{i=1}^{\infty} (\tau_{i}-\tau_{i-1} - \frac{1}{i})^{2} = \sum_{i=1} \frac{1}{i^2}  \mathrm{Var}( \mathcal{E}(1)) < \infty.$$ Hence $( \tau_{n}- \log n: n \geq 1)$ converges almost surely (and in $L^{2}$) towards some random variable $ \mathcal{X}_{\infty}$. To compute the law of this variable, recall from Lemma \ref{lem:gumbel} and the discussion following it that we have 
 $$\tau_{n} - \log n \xrightarrow[n\to\infty]{(d)} G,$$ where $G$ has the Gumbel distribution. 
 We deduce that 
$$ (n+1) \mathrm{e}^{-\tau_n}  \xrightarrow[n\to\infty]{a.s.} \exp(- \mathcal{X}_\infty)  \overset{(d)}{=}  \mathrm{e}^{-G} \overset{(d)}{=}\mathcal{E}(1),$$ and this proves the statement of the proposition for times $t$ of the form $\tau_n$. Assuming for a moment that $\Delta \tau_n \to 0$ almost surely as $n \to \infty$, a sandwiching argument for times $\tau_n \leq t < \tau_{n+1}$ concludes the proof. To prove that $\Delta \tau_n \to 0$, we use the Borel--Cantelli lemma since for $n \geq 1$
$$ \mathbb{P} \left(\tau_{n+1} - \tau_n \geq \frac{1}{ \sqrt{n}}\right) = \mathbb{P}\left(  n^{-1}\mathcal{E}(1) \geq  \frac{1}{ \sqrt{n}}\right) = \mathrm{e}^{- \sqrt{n}},$$ is summable in $n \geq 1$.\\
 The case $k \geq 3$ is similar: the only trick is to consider the sum of $k-1$ independent Yule trees so that the jump times of the forest are separated by independent variables of law $ \mathcal{E}(k-1), \mathcal{E}(2(k-1)), \mathcal{E}(3(k-1))\dots$  we can reduce to the above problem (and using the fact that $(k-1)$ copies of r.v.~of law $\Gamma( \frac{1}{k-1},\frac{1}{k-1})$ is an exponential of parameter $k-1$). \qed \medskip 

Actually, in the case $k=2$ (and $\alpha_{1}=1$)  the distribution of $ \mathcal{Y}^{(2)}_{t}$ is explicitly given for each $t\geq0$ by a geometric distribution with parameter $  \mathrm{e}^{-t}$, i.e. 
 \begin{eqnarray} \label{eq:yuleexplicit} \mathbb{P}( \mathcal{Y}^{(2)}_{t}=k) =   \mathrm{e}^{-t}(1- \mathrm{e}^{-t})^{k-1}, \quad \mbox{ for }k \in \{1,2,3, \dots\}. \end{eqnarray}  Taking the limit as $t \to \infty$, this recovers the form of the limiting law in the above proposition. Once given, the proof of the above claim is easy  by solving the differential equations satisfied by the probabilities $ p_{k}(t) := \mathbb{P}( \mathcal{Y}^{(2)}_{t}=k)$ for $k \geq 2$
 $$  \frac{ \mathrm{d}}{  \mathrm{d}t} p_{k}(t) = - k p_{k}(t) + (k-1)p_{k-1}(t),$$ with the limiting condition $ p_{1}(t)= \mathrm{e}^{-t}$. See \cite[Chapter III.5 ]{AN72} for analogs when $k \geq 3$.

  \section{Examples} \label{sec:examplesAK}
  We now give a few examples of discrete Markov chains which are easily studied via their continuous time analogs. This includes the classical coupon collector problem, the  pill problem, the O.K. Corral model and the random recursive tree! We shall start with a new look at the Polya urn studied in Section \ref{sec:polyaurn} before moving to the more challenging examples that will require a few results  useful to perform the continuous-time $\to$ discrete time or ``depoissonization'' operation.

  \subsection{Polya Urn, reloaded} \label{sec:revolution}
Let us interpret  the classical Polya urn scheme (Section \ref{sec:polyaurn}) as the counting process of  a continuous time branching process using  Lemma \ref{lem:AK}.  For this we consider the case when $p=2$, i.e.~we have two types of particles (red and blue say) and the offspring mechanisms are deterministic $\mu_{1} = \delta_{2\delta_1}$ and $\mu_{2} = \delta_{2\delta_2}$: each particle reproduces at rate $1$ into two particles of the same color independently of the others. Then the branching forest $ \mathbb{F}$ under $ \mathbb{P}_{\delta_1+\delta_2}$ is made of two trees, one red and one blue, describing the genealogy of the two initial particles. By Proposition \ref{lem:AK}, the discrete Markov chain describing the number of blue and red particles at each jump time is simply given by \eqref{eq:polya}, that is, if we start initially with one particle of each color, we are facing the dynamic of the standard Polya urn! \\
Now, the magic of the continuous time is that, since particles of different colors do not interact, the two trees of $ \mathbb{F}$ under $\mathbb{P}_{\delta_1+\delta_2}$, are independent copies of the standard Yule tree of order $2$ (started with a single particle). If $ \mathcal{B}_{t}$ and $ \mathcal{R}_{t}$ respectively denote the number of blue and red particles alive at time $ t \geq 0$ then from Proposition \ref{prop:yule} we have 
$$ \mathrm{e}^{-t} (\mathcal{B}_{t}, \mathcal{R}_{t}) \xrightarrow[t\to\infty]{}( X, X'),$$ where $ X$ and $ X$ are two independent exponential laws of expectation $1$. In particular, we re-deduce Proposition \ref{prop:polyaunif} on the asymptotic proportion of blue balls:
$$ \frac{ \mathcal{B}_t}{ \mathcal{R}_t + \mathcal{B}_t} \xrightarrow[t\to\infty]{a.s.}  \frac{X}{X+X'}  \overset{(d)}{=} \mathrm{Unif}[0,1].$$

\begin{exo} \label{exo:afresh} Contemplate Exercise \ref{exo:polyaassy} afresh.  \end{exo}

  \subsection{Depoissonization tools}
  We now present two lemmas that we will use repeatedly  below.  The first one is a probabilistic variation on Dini's lemma:
  
  \begin{lemma}[Dini] \label{lem:dini} Let $(D^{{(n)}}_{t} : t \in \mathbb{R})$ be random non-decreasing  c\`adl\`ag processes, i.e such that $ D^{{(n)}}_{s} \leq D^{{(n)}}_{t}$ for every $s \leq t$ and $n \geq 0$. We suppose that $D^{(n)}$  converge point-wise in probability, that is for any $t \in \mathbb{R}$ we have
  $$ D^{(n)}_{t} \xrightarrow[n\to\infty]{( \mathbb{P})} f(t),$$ where $f : \mathbb{R} \to \mathbb{R}$ is a non-decreasing continuous function. Then we also have the stronger convergence  $$ ( D_{t}^{{(n)}} : t \in \mathbb{R}) \xrightarrow[n\to\infty]{( \mathbb{P})} (f(t) : t \in \mathbb{R}),$$ for the topology of uniform convergence over every compact subset of $ \mathbb{R}$.
  \end{lemma}
\noindent \textbf{Proof.} Fix a dense sequence $(t_{i} : i \geq 0)$ in $ \mathbb{R}$. Since $ D^{{(n)}}_{t_{i}} \to f(t_{i})$ in probability as $n \to \infty$ for each $i$, we have $( D_{t_{i}}^{{(n)}} : i \geq 0) \to ( f(t_{i}) : i \geq 0)$ in probability for the topology of point-wise convergence on $ \mathbb{R} ^{ \mathbb{Z}_{\geq 0}}$. By the Skorokhod representation theorem, we can construct a probability space $( \Omega, \mathcal{F}, \mathbb{P})$ and a sequence of processes $( \tilde{D}^{(n)})$ so that $ \tilde{D}^{(n)} {=} D^{{(n)}}$ in law  for  each $n$, and so that we have 
$$ \forall i \geq 0, \quad \tilde{D}^{(n)}_{t_{i}} \xrightarrow[n\to\infty]{} f(t_{i})\quad \mbox{ \textbf{almost surely}}.$$
Since $\tilde{D}^{(n)} \overset{(d)}{=} D^{{(n)}}$, the processes $\tilde{D}^{(n)}$ are non-decreasing, and it follows from (classical) Dini's theorem that we actually have the stronger convergence $(\tilde{D}^{(n)}_{t} : t  \in \mathbb{R}) \to ( f(t) : t \in \mathbb{R})$ for the topology of uniform convergence over every compact subset of $ \mathbb{R}$. We deduce the similar convergence but in probability for $ D^{(n)}$ by equality in law. \qed \medskip

The same result holds true (with the same proof) if we replace convergence in probability by almost sure convergence.  Let us see how we can use such convergences:
\begin{lemma}[Slutsky] \label{lem:slut} Suppose that $(D^{(n)}_{t})_{t \in \mathbb{R}}$ is a sequence of random processes and $T_{n}$ as sequence of random times (which might not be stopping times). Suppose that
$$ D^{{(n)}} \xrightarrow[n\to\infty]{(d)} F \quad \mbox{ and } \quad  T_{n} \xrightarrow[n\to\infty]{(d)} \theta,$$ where $F$ is a random continuous function and $\theta$ is a random variable. The first convergence is in the sense of uniform convergence over every compact of $ \mathbb{R}$. We suppose that either $F \equiv f$ is a fixed continuous function \textbf{or} that $\theta \equiv C$ is a constant (in which case the respective convergence in distribution holds in probability). Then we have 
  \begin{eqnarray*}  D^{{(n)}}_{T_{n}} \xrightarrow[n\to\infty]{ (d)} F(\theta). \end{eqnarray*}
  \end{lemma}
  \noindent \textbf{Proof.} Since one of the limiting variables is deterministic, Slutsky's lemma entails that $( D^{{(n)}}, T_{n})$ converges in distribution towards $((F(t) : t \in \mathbb{R}), \theta)$. We can then use Skorokhod representation again to obtain versions $( \tilde{D}^{{(n)}}, \tilde{T}_{n})$ so that $( \tilde{D}^{{(n)}}, \tilde{T}_{n}) \overset{(d)}{=}( {D}^{{(n)}}, {T}_{n})$ for each $n$ but satisfying 
$$ (\tilde{D}^{{(n)}}_t : t \in \mathbb{R}) \xrightarrow[n\to\infty]{a.s.}(F(t) : t \in \mathbb{R}) \quad \mbox{ and }\quad  \tilde{T}_{n} \xrightarrow[n\to\infty]{a.s.} \theta,$$ where the first arrow holds for the uniform convergence on every compact of $ \mathbb{R}_+$. We deduce the desired convergence in law since 
$$D^{{(n)}}_{T_{n}} \underset{ \mathrm{for\ each \ }n }{\overset{(d)}{ = }} \tilde{D}^{{(n)}}_{\tilde{T}_{n}}  \xrightarrow[n\to\infty]{a.s.} F(\theta).$$\qed 

\subsection{Coupon collector}
The famous coupon collector problem is the following. Fix $n \geq 1$ and let $(X_{k})_{k \geq 1}$ be i.i.d.~uniform variables over $\{1,2, \dots , n\}$. We interpret each $X_{k}$ as a ``coupon'' among a collection of all  $n$ possible ones, and we ask how many coupons we should buy to get the full collection, i.e. 
$$ \mathcal{T}_{n} := \inf\big\{ k \geq 1 : \{X_{1}, \dots , X_{k}\} = \{1,2, \dots , n\}\big\}.$$
Using our continuous time embedding technique we shall prove:
\begin{proposition}[Coupon collector] \label{prop:coupon}We have the following convergence in law
$$ \frac{ \mathcal{T}_{n} -  n \log n}{n} \xrightarrow[n\to\infty]{(d)} G,$$
where $G$ has the Gumbel distribution.
\end{proposition}
\noindent \textbf{Proof.} We pass in continuous time and consider for each $i \in \{1,\dots , n\}$ an independent Poisson processes $  \mathfrak{P}^{{(i)}}$ of unit rate. This is equivalent to considering $p=1$, $\mu_{1} =\delta_{\delta_1}$ and $\alpha_{1}=1$ under $ \mathbb{P}_{n \delta_1}$ in Lemma \ref{lem:AK}. We let $0 < \tau_{1}< \tau_{2} < \dots$ be the jump times of the union of those processes, so that by an application of Lemma \ref{lem:AK} the indices of the corresponding Poisson processes are distributed as $(X_{k})_{ k \geq 1}$. The continuous time analog of $ \mathcal{T}_{n}$ in this setting is thus 
$$ T_{n} + \log n := \max_{1 \leq i \leq n} \min \big\{ {t \geq 0} : \mathfrak{P}^{(i)}(t) \geq 1\big\},$$
which by Proposition \ref{prop:standardpoisson} has the law of the maximum of $n$ independent exponential variables of rate $1$. This is given by Lemma \ref{lem:gumbel} and we have $ T_{n}\to G$ in distribution. Coming back to the discrete setting, the number of coupons bought as time $ T_{n}$ is thus 
$$ \mathcal{T}_{n} \overset{(d)}{=} \sum_{i=1}^{n} \mathfrak{P}^{{(i)}}( T_{n} + \log n).$$
The sum $\sum_{i}\mathfrak{P}^{{(i)}}(\cdot)$ has the same distribution as $ \mathfrak{P}(n \cdot)$, but beware, in this writing $ \mathfrak{P}$ is \textit{not independent} from  $ T_{n}$. To circumvent this problem, notice that for any $t \in \mathbb{R}$ we have the convergence in probability
$$  D^{{(n)}}_{t} \xrightarrow[n\to\infty]{ (\mathbb{P})} t, \quad \mbox{ where }\quad D^{{(n)}}_{t}:=\frac{ \mathfrak{P}( n \log n + n t)- n \log n}{n}.$$
This weak law of large number is easily seen since $ \mathbb{E}[ D^{{(n)}}_{t}]=t$ and $ \mathrm{Var}(D^{{(n)}}_{t}) \leq \mathrm{Cst} n \log n /n^2$ for some $ \mathrm{Cst}>0$. We deduce from Lemma \ref{lem:dini} the stronger version:
$$ \left( D^{{(n)}}_{t}\right)_{t \in \mathbb{R}} \xrightarrow[n\to\infty]{ (\mathbb{P})} (t)_{t \in \mathbb{R}},$$ for the topology of uniform convergence over every compact of $ \mathbb{R}$ and by Lemma \ref{lem:slut} we get  
$$  \frac{ \mathcal{T}_{n} - n \log n}{n} \quad \overset{(d)}{=} D^{{(n)}}_{ T_{n}} =  \frac{ \mathfrak{P}(  n \cdot T_{n})- n \log n}{n} \quad  \xrightarrow[n\to\infty]{(d)} \quad G.$$
\qed 

\subsection{Balls in bins}
The above approach (with the same continuous time process!) can be used to address the balls in bin problem. Let again $(X_k)_{k\geq 1}$ be i.i.d.~r.v.~uniformly distributed over $\{1,2, \dots , n \}$. We interpret this time the $X_k$ as ``balls'' that are thrown uniformly at random in the $n$ ``bins'' numbered $1,2, \dots , n$. The question is: After throwing $n$ balls, what is the maximal load of a bin, i.e. 
$$ \mathrm{ML}_n = \max \{ B_i : 1 \leq i \leq n\}\quad \mbox{where} \quad B_i = \#\{1 \leq k \leq n: X_k =i\}.$$
\begin{proposition}[Balls in bins] We have $$ \mathrm{ML}_n \sim_{ \mathbb{P}} \frac{\log n}{ \log \log n}, \quad \mbox{ as $n \to \infty$}.$$
\end{proposition}
\noindent \textbf{Proof.} We use the same notation as in the proof of Proposition \ref{prop:coupon} and in particular $ \mathfrak{P}^{(i)}$ are independent unit rate Poisson processes carried by each bin, and $\tau_n$ is the time at which $n$ balls have been thrown. We deduce that we have 
$$ \mathrm{ML}_n \overset{(d)}{=} \max_{1\leq i \leq n} \mathfrak{P}^{(i)}( \tau_n).$$
As before, the problem is that $\tau_n$ is not independent from the $ \mathfrak{P}^{(i)}$. However, on the one hand, recalling that sum $\sum_{i}\mathfrak{P}^{{(i)}}(\cdot)$ has the same distribution as $ \mathfrak{P}(n \cdot)$, we clearly have by the law of large numbers that 
$$ \tau_n \xrightarrow[n\to\infty]{ ( \mathbb{P})} 1.$$
On the other hand, for \textit{fixed} $t_0 >0$, the variable $ (\mathfrak{P}^{(i)}(t_0) : 1 \leq i \leq n )$ are independent $ \mathfrak{P}(t_0)$ random variables, so that if we let $ \mathcal{ML}_{n}(t_0) = \max_{1\leq i \leq n} \mathfrak{P}^{(i)}(t_0)$, we have for any $m\geq 1$
$$  \mathbb{P}\left(\mathcal{ML}_{n}(t_0) < m \right) = \big(1- \mathbb{P}( \mathfrak{P}(t_0) \geq m)\big)^n.$$ 
It is easy to see that $\mathbb{P}( \mathfrak{P}(t_0) \geq m)$ is actually equivalent to   $\mathrm{e}^{-t_0}\frac{t_0^m}{m!}$ as $m \to \infty$. So,  for any $  \varepsilon>0$ the above display goes to $0$ for $m \leq (1- \varepsilon) \log n/\log\log n$ and to $1$ for $m \geq (1+ \varepsilon) \log n / \log\log n$ as $n \to \infty$. We deduce that for any $t_0>0$ we have  $$ \left( D^{{(n)}}_{t_{0}}\right) \xrightarrow[n\to\infty]{( \mathbb{P})} 1 \quad \mbox{ where }\quad D^{{(n)}}_{t_{0}} :=\frac{\mathcal{ML}_{n}(t_0)}{\log n / \log \log n}.$$
This convergence  is reinforced using monotonicity and Lemma \ref{lem:dini} into 
$$ \left(D^{{(n)}}_{t}\right)_{ t \in [0.9,1.1]} \xrightarrow[n\to\infty]{(\mathbb{P})} (1)_{t \in [0.9, 1.1]},$$ for the topology of uniform convergence over $[0.9, 1.1]$. Since $ \tau_n \to 1$ in probability, we can then apply Lemma \ref{lem:slut} to deduce as desired that $$ \frac{\mathrm{ML}_n}{\log n / \log \log n} \overset{(d)}{=} D^{{(n)}}_{ \tau_n} \xrightarrow[n\to\infty]{( \mathbb{P})}1.$$  \qed  

   \subsection{Pill problem} \label{sec:pills}
  From Wikipedia: \begin{quote}
  The pill jar puzzle is a probability puzzle, which asks the  value of the number of half-pills remaining when the last whole pill is popped from a jar initially containing $n$ whole pills and the way to proceed is by removing a pill from the bottle at random. If the pill removed is a whole pill, it is broken into two half pills. One half pill is consumed and the other one is returned to the jar. If the pill removed is a half pill, then it is simply consumed and nothing is returned to the jar.\end{quote}
  
This problem (attributed to Knuth and McCarthy) can be approached using the Athreya--Karlin embedding. Indeed, suppose we have two types of particles: those of type $1$ corresponding to half-pills and those of type $2$ corresponding to whole pills. We set the rates $\alpha_{1}= \alpha_{2}=1$ and suppose that when a particle of type $2$ dies, it gives rise to a single particle of type $1$, whereas particle of type $1$ have no descendance. Formally $\mu_{2} = \delta_{\delta_1}$ and $\mu_{1} = \delta_{\varnothing}$. If we start initially from $ \mathbb{P}_{n \delta_2}$ i.e.~a forest $ \mathbb{F}$ with $n$ particles of type $2$ (whole pills) then by Lemma \ref{lem:AK} the evolution of the underlying discrete time Markov chain corresponds to the evolution of the content of the jar in the pill puzzle above. If $L_{n}$ is the number of half-pills remaining when all whole pills have been consumed we can then easily prove:

\begin{proposition}[Pill problem] \label{prop:pills}Under $ \mathbb{P}_{n \delta_2}$ the random variable $ \frac{L_{n}}{\log n}$ converges in law towards an exponential variable of mean $1$.
\end{proposition}

\noindent \textbf{Proof.} Under $ \mathbb{P}_{n \delta_2}$ the evolution of the $n$ genealogies starting from the $n$ particles of type $2$ are independent and are described by a sequence $( X_{2}^{(i)}, X_{1}^{(i)} : 1 \leq i \leq n)$ of i.i.d.~r.v.~of law $ \mathcal{E}(1)$ giving the life time of the particles of type $2$ and of their only child of type $1$. If for every $t >0$ we introduce the number of particles still alive at time $t + \log n$
$$  D_{t}^{{(n)}} =  \frac{1}{\log n} \sum_{i=1}^{n} \mathbf{1}_{ X_{2}^{(i)} + X_{1}^{(i)} > t + \log n}, $$ then by Athreya--Karlin embedding we have 
$$ \frac{L_{n}}{\log n} \overset{(d)}{=} D_{T_{n}}^{{(n)}}\quad \mbox{ where } \quad T_{n} = \max_{1 \leq i \leq n} X_{2}^{(i)}-\log n.$$
By Lemma \ref{lem:gumbel} we have  the convergence to a Gumbel distribution $T_{n} \xrightarrow{(d)}  G$ as $n \to \infty$. On the other hand, since $ \mathbb{P}(X_{2}^{(1)} + X_{1}^{(2)}>t) = \mathrm{e}^{-t}(t+1)$, an easy law of large number  (proved using first and second moment for example)  shows that for deterministic times $\log n +t$ for $t \in \mathbb{R}$ we have 
$$ D_{t}^{{(n)}} \xrightarrow[n\to\infty]{ (\mathbb{P})}  \mathrm{e}^{-x}.$$ This convergence is as usual reinforced using Lemma \ref{lem:dini} and monotonicity into a functional one. We can then couple the previous three displays to deduce using Lemma \ref{lem:slut} that 
$$ \frac{L_{n}}{\log n}  \overset{(d)}{=} D^{(n)}_{T_{n}} \xrightarrow[n\to\infty]{(d)}  \mathrm{e}^{-G} \overset{(d)}{=} \mathcal{E}(1).$$ \qed  

\subsection{O.K. Corral}

\begin{figure}[!h]
 \begin{center}
 \includegraphics[width=14cm]{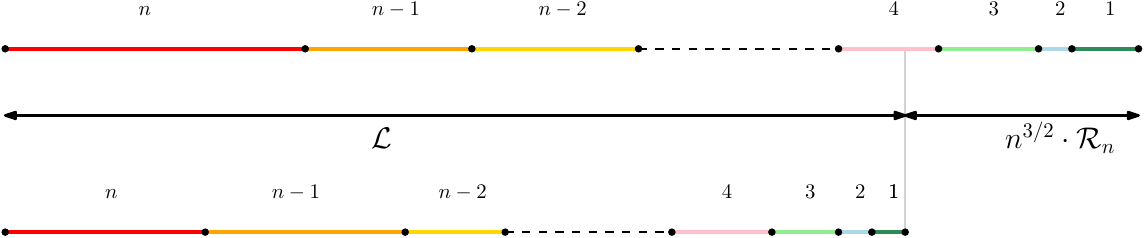}
 \caption{Illustration of the proof: the two independent lineages of particles of types $n \to n-1 \to \dots \to 2 \to 1$. The type of the particle still standing at the death of the other lineage (here $4$) is studied through $ \mathcal{R}_n$. \label{fig:okcorral}}
 \end{center}
 \end{figure}

Imagine two groups of $n$ people facing each other. At each time step, one individual is chosen uniformly and shouts a person of the other group. The question is: ``How many people are still standing when one of the group dies out''. This riddle is usually named the O.K. Corral \footnote{The gunfight at the O.K. Corral took place on October 26, 1881. Four lawmen were facing five outlaws.  During that brief battle (less than a minute), three men were killed, three were wounded, two ran away, and one was unharmed.} problem. Formally, let $( O_1(k),O_2(k) : k \geq 0)$ a Markov chain on $ \{0,1,\dots,n\}^2$ starting from $O_1(0)=O_2(0) =n$ and with transition probabilities
$$ \mathbb{P}\left( \displaystyle \begin{subarray}{l} O_1(k+1)\\ O_2(k+1) \end{subarray} = \begin{subarray}{l} O_1(k)-1\\ O_2(k) \end{subarray} \left| \begin{subarray}{l} O_1(k)\\ O_2(k) \end{subarray} \right)\right. = 1- \mathbb{P}\left( \begin{subarray}{c} O_1(k+1)\\ O_2(k+1) \end{subarray} = \begin{subarray}{l} O_1(k)\\ O_2(k)-1 \end{subarray} \left| \begin{subarray}{l} O_1(k)\\ O_2(k) \end{subarray} \right)\right. = \frac{ O_2(k)}{O_1(k)+O_2(k)}.$$
We then let $\theta_n = \inf\{ k \geq 0 : O_1(k)=0 \mbox{ or }O_2(k)=0\}.$
\begin{theorem} We have the following convergence in distribution
$$ n^{-3/4}\left(O_1(\theta_n) + O_2(\theta_n)\right) \xrightarrow[n\to\infty]{(d)}  \left( \frac{8}{3}\right)^{1/4} \sqrt{| \mathcal{N}|},$$ where $  \mathcal{N}$ is a standard normal variable.
\end{theorem}
\noindent \textbf{Proof.} We shall embed the discrete Markov chain in continuous time using the Athreya--Karlin lemma. Specifically suppose that we start from two particles of type $n$. Each particle of type $i \in \{1,2,\dots  ,n\}$ behave independently of each other and lives for an exponential variable of parameter $ \frac{1}{i}$ (or mean $i$) and then gives rise to a particle of type $i-1$. If $i=1$ then the lineage dies out when the particle of type $1$ dies out. Formally, this is obtained by taking an infinite number of types $p=\infty$, with rates $\alpha_{i} = \frac{1}{i}$  and offspring distribution $\mu_{i} = \delta_{\delta_{i-1}}$ for $i \geq 1$ and $\mu_{1}= \delta_{\varnothing}$, see Figure \ref{fig:okcorral}. Then under $ \mathbb{P}_{2 \delta_n}$, we have two independent lineages of particles of type $n \to n-1 \to \dots \to 2 \to 1$. We denote by $ \mathcal{L}_1, \mathcal{L}_2$ the lengths of the lineages and put $ \mathcal{L} = \mathcal{L}_1 \wedge \mathcal{L}_2$. By Lemma \ref{lem:AK}, the discrete evolution of the types of particles at the jump times $0=\tau_0<\tau_1<\dots$  has the same law as $(O_1(k), O_2(k) : 0 \leq k \leq \theta_n)$. The quantity we are looking for is the type of the remaining particle at time $ \mathcal{L}$ and we shall observe this through its remaining life time:
$$  \mathcal{L}_{1} + \mathcal{L}_{2} - \mathcal{L} = | \mathcal{L}_{1}- \mathcal{L}_{2}|.$$
We have 
$$ \mathbb{E}[\mathcal{L}_{1}- \mathcal{L}_{2}] =0, \quad \mathrm{Var}(\mathcal{L}_{1}- \mathcal{L}_{2}) = 2 \sum_{i=1}^{n} n^{2} \mathrm{Var}( \mathcal{E}(1)) \sim \frac{2 n^{3}}{3},$$ and we leave to the reader verify (using e.g.~Lindenberg CLT, or using characteristic functions) that we have $ n^{-3/2}(\mathcal{L}_{1}- \mathcal{L}_{2}) \to (2/3)^{-1/2} \mathcal{N}$ in law so that
$$  \mathcal{R}_{n} = n^{-3/2}|\mathcal{L}_{1}- \mathcal{L}_{2}| \xrightarrow[n\to\infty]{(d)}  \sqrt{\frac{2}{3}} | \mathcal{N}|,$$ where $ \mathcal{N}$ is a standard normal. 

We therefore know that the remaining life time of the lineage of the particle still standing at time $ \mathcal{L}$ is of order $ \sqrt{\frac{2}{3}} |\mathcal{N}| n^{3/2}$, to connect this variable with the type of the particle in question, we use the following: Let $ {D}^{(n)}_{t}$ be the type of the particle still alive in the first lineage  at time $ \mathcal{L}_{1} - n^{3/2}t$, renormalized by $n^{3/4}$. We will show that 
$$ (D^{{(n)}}_{t})_{t \geq 0} \xrightarrow[n\to\infty]{ (\mathbb{P})} ( \sqrt{2 t})_{t \geq 0}.$$ Once this is done, since the same result holds for the second lineage where the process is denoted $(\tilde{D}^{{(n)}}_{t})_{t \geq 0}$, the result is again a consequence of Lemma \ref{lem:AK} since we have $$n^{-3/4}\left(O_1(\theta_n) + O_2(\theta_n)\right)  \overset{(d)}{=} D^{(n)}_{ \mathcal{R}_{n}} \mathbf{1}_{ \mathcal{L}_{2}< \mathcal{L}_{1}} + \tilde{D}^{(n)}_{ \mathcal{R}_{n}} \mathbf{1}_{ \mathcal{L}_{1}< \mathcal{L}_{2}}   \xrightarrow[n\to\infty]{(d)}  \sqrt{2 \sqrt{2/3} |\mathcal{N}|}.$$
To prove the penultimate display, we shall rather focus on the inverse function of $D^{(n)}_{\cdot}$ and consider for $ x \geq 0$ the remaining time $ n^{3/2}\cdot \mathcal{H}^{{(n)}}_x$ in the lineage starting from a particle of type $\lfloor x n^{{3/4}} \rfloor$. It is thus sufficient to show that $ (\mathcal{H}^{{(n)}}_{ x })_{x\geq 0} \to  (\frac{x^{2}}{2})_{x \geq 0}$, or by monotonicity and Lemma \ref{lem:dini} that for each $x\geq 0$ we have 
  \begin{eqnarray} \mathcal{H}^{{(n)}}_x \xrightarrow[n\to\infty]{( \mathbb{P})} \frac{x^{2}}{2}.   \label{eq:hnxlgn}\end{eqnarray}
Since $\mathcal{H}^{{(n)}}_x \overset{(d)}{=} n^{-3/2}\cdot  \sum_{i=1}^{\lfloor x n^{3/4} \rfloor} X_{i} $ where the variables $X_{i}$ are independent and of law $ \mathcal{E}( 1/i)$, the expectation and variance of $ \mathcal{H}^{(n)}_{x}$ are easily estimated:
$$ \mathbb{E}[\mathcal{H}^{{(n)}}_x ] = n^{-3/2}\sum_{i=1}^{\lfloor x n^{3/4} \rfloor} i  \xrightarrow[n\to\infty]{} \frac{x^{2}}{2}, \quad \mbox{ and }\quad   \mathrm{Var}(\mathcal{H}^{{(n)}}_x) = n^{-3}\sum_{i=1}^{\lfloor x n^{3/4} \rfloor} i^{2} \xrightarrow[n\to\infty]{}0.$$
Our goal \eqref{eq:hnxlgn} then follows by Markov's inequality. \qed

\section{Back to the Random Recursive Tree}

Our last example is the random recursive tree process (Chapter \ref{chap:RRT}) which we will construct from a standard Yule tree of order $2$. This will enables us to give quick proofs of (stronger) results about the geometry of the RRT. As we will see in the next chapter, the Athreya-Karlin embedding will give independence properties that make life much simpler when proving the deep Theorems \ref{prop:maxdegreeRRT} and \ref{thm:heightRRT}.

\subsection{Construction of the RRT from a Yule process}
\label{sec:rrtfromyule}

\begin{figure}[!h]
 \begin{center}
 \includegraphics[width=14cm]{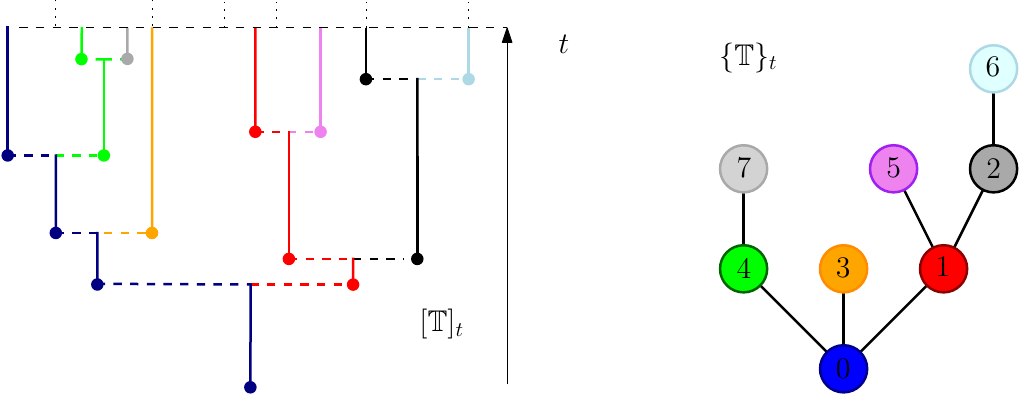}
 \caption{Constructing the random recursive tree (Right) from a standard Yule process (Left): each particle gives rise to a particle of a new type at an exponential rate and this is interpreted as an attachment in the RRT. \label{fig:yuleRRT}}
 \end{center}
 \end{figure}

Let us consider the  plane version of the Yule tree $ \mathbb{T}$ of order $2$ started from a single particle and recall the notation  $( [\mathbb{T}]_{t} : t \geq 0)$ for the tree cut at height $t$.  In the plane version of $ [ \mathbb{T}]_{t}$ we contract all the edges going to the left: we obtain a plane genealogical tree whose vertices are labeled by $0,1, 2, \dots$ by order of appearance in the Yule  tree, see Figure \ref{fig:yuleRRT}. We denote by $ \{  \mathbb{T}\}_t$ the increasing tree obtained after forgetting the plane ordering. The following is easily proved using the same  techniques as in the proof of Lemma \ref{lem:AK}:
\begin{proposition}[From Yule to RRT]  \label{prop:RRTyule} If $0= \tau_0 < \tau_{1}< \dots < \tau_{n}< \dots$ are the first times at which $\# \partial [ \mathbb{T}]_{\tau_n} = n+1$ then conditionally on $( \tau_{n} : n \geq 0)$ the process $( \{ \mathbb{T}\}_ {\tau_{n}} : n \geq 0)$ is a random recursive tree. \end{proposition}

 \noindent \textbf{Proof.} Let us prove by induction on $n \geq 0$ that at time $\tau_n$, conditionally on the past up to time $ \tau_n$, the Yule tree has $n+1$ alive particles carrying independent exponential clocks, the first one that rings inducing a splittings into two particles. This is true for $n =0$ and propagates easily by the memorylessness property of the exponential laws \eqref{eq:expomemory}. In particular, by Proposition \ref{prop:expomem}, conditionally on the past up to time $\tau_n$, the next particle to split is a uniform particle of $\partial [ \mathbb{T}]_{\tau_n}$. Translating the dynamics in terms of $\{\mathbb{T}\}_{\tau_n}$ directly shows that this chain evolves as a random recursive tree. \qed  \medskip 
 
 \subsection{Degree statistics}
 Let us use Proposition \ref{prop:RRTyule} to give streamlined proofs of basic results on degree distribution in the RRT. Recall in particular from Proposition	\ref{prop:yule}  that we have   \begin{eqnarray} \label{eq:taynb} \frac{\tau_{n}}{ \log n} \xrightarrow[n\to\infty]{a.s.}  1,  \end{eqnarray} and  more precisely $ \tau_n - \log n \to G$ where $ G$ has the Gumbel distribution. By the above construction and Proposition \ref{prop:standardpoisson}, for all $t \geq 0$, the degree of the root vertex $\raisebox{.5pt}{\textcircled{\raisebox{-.9pt} {$0$}}}$ in $ \{\mathbb{T}\}_{t} $ is given by $ \mathfrak{P}(t)$ where  $( \mathfrak{P}(t) : t \geq 0)$ is a unit-rate Poisson counting process. This enables us to deduce a stronger version of  \eqref{eq:smalldegreelog} given in the last chapter:
 \begin{proposition} \label{prop:degreeroot}We have the following convergences
$$ \frac{\mathrm{deg}_{\{ \mathbb{T}\}_{\tau_{n}}}( \raisebox{.5pt}{\textcircled{\raisebox{-.9pt} {$0$}}})}{\log n} \xrightarrow[n\to\infty]{a.s.} 1 \quad \mbox{ and }\quad \frac{\mathrm{deg}_{\{ \mathbb{T}\}_{\tau_{n}}}( \raisebox{.5pt}{\textcircled{\raisebox{-.9pt} {$0$}}})- \log n}{ \sqrt{\log n}} \xrightarrow[n\to\infty]{(d)} \mathcal{N}(0,1).$$
\end{proposition}
\noindent \textbf{Proof.} Since the degree of the root in $ \raisebox{.5pt}{\textcircled{\raisebox{-.9pt} {$0$}}} $ in $ \{ \mathbb{T}\}_t$ is given by the Poisson counting process $ \mathfrak{P}(t)$ along the left-most branch, using \eqref{eq:lawll} we deduce that 
$$ \frac{\mathrm{deg}_{ \{ \mathbb{T} \}_{t}}( \raisebox{.5pt}{\textcircled{\raisebox{-.9pt} {$0$}}})}{t} \xrightarrow[t\to\infty]{a.s.} 1 \quad \mbox{ and }\quad \left(\frac{\mathrm{deg}_{ \{\mathbb{T}\}_{tu}}( \raisebox{.5pt}{\textcircled{\raisebox{-.9pt} {$0$}}})- tu}{ \sqrt{t}}\right)_{u \geq 0} \xrightarrow[t\to\infty]{(d)}  (B_{u})_{u \geq 0},$$ where $B$ is a standard linear Brownian motion and where the convergence in the right-hand side holds with respect to the topology of uniform convergence for every compact subset of $ \mathbb{R}_{+}$. From Proposition 	\ref{prop:yule} it follows that $\tau_{n}  \sim \log n$ a.s. as $ n\to \infty$ and the desired statement follows by combining those observations and using  Lemma \ref{lem:slut}.\qed \medskip

\subsection{A new look at Goncharov \& Kolchin's result} \label{sec:goncharovback}
Let us now use the link between uniform  permutations and the RRT, and the construction of the latter from a standard Yule tree, to give a fresh look at Goncharov \& Kolchin's result (Theorem \ref{thm:poissoncycle}) on the Poisson statistics of small cycle counts. More precisely, we shall give a direct proof of Lemma \ref{lem:shepplloyd} due to Loyd \& Shepp without relying on Cauchy formula: \medskip 

\noindent \textbf{Proof of Lemma \ref{lem:shepplloyd}, second version.} Consider the increasing tree $\{ \mathbb{T}\}_t$ and let us denote by $\sigma$ the random permutation associated with it thanks to Section \ref{sec:CRP-RRT}. In particular, conditionally on its size, the permutation $\sigma$ is uniformly distributed. Recall also that the cycle lengths of $\sigma$ correspond to the sizes of the subtrees above $\noeud{0}$ in $\{ \mathbb{T}\}_t$, the later corresponding via the construction of Figure \ref{fig:yuleRRT} to the size (number of individuals living at time $t$) of the subtrees branching of from the left-most branch in $[ \mathbb{T}]_t$. By \eqref{eq:yuleexplicit}, the process of points on the left-most branch, identified with $[0,t]$, at which branches a subtree reaching $k \in \{1,2, \dots \}$ individuals at time $t$ is Poisson with intensity 
$$ \mathrm{e}^{-(t-s)}(1- \mathrm{e}^{-(t-s)})^{k-1} \mathbf{1}_{s \in [0,t]} \mathrm{d}s,$$ and furthermore, by Poisson thinning, those processes are independent for different values of $k$. We deduce that the number of cycles of length $k \in  \{1,2, \dots \}$ in $ \sigma$ are given by independent Poisson variables with mean 
$$ \int_0^t \mathrm{d}s\ \mathrm{e}^{-(t-s)}(1- \mathrm{e}^{-(t-s)})^{k-1} =  \frac{(1- \mathrm{e}^{-t})}{k}.$$
This is exactly the statement of Lemma \ref{lem:shepplloyd} with $x = (1- \mathrm{e}^{-t})$. \qed 
\bigskip 

\subsection{Concentration of local statistics}
The continuous time embedding and its independence properties can also be used to efficiently prove concentration of local statistics in the RRT. Let us focus on the degree to illustrate the method: 
For $k \geq 0$ and $t \geq 0$ introduce the variable 
$$ \mathrm{D}_k([  \mathbb{T}]_t) := \# \Big\{ u \in \{ \mathbb{T}\}_t \backslash \noeud{0} : \mathrm{deg}^{+}_{\{ \mathbb{T}\}_t}(u) =k\Big\},$$ which counts the number of vertices (except the root) in the contraction of $[ \mathbb{T}]_t$ whose out-degree is $k$. Then we have 
\begin{proposition}[Concentration of local statistics] \label{prop:concentrationlocal} We have 
$$ \frac{\mathrm{D}_k([  \mathbb{T}]_t)}{ \# \partial [ \mathbb{T}]_t} \xrightarrow[t\to\infty]{ ( \mathbb{P})}  \lim_{t \to \infty} \mathrm{e}^{-t} \mathbb{E}[ \mathrm{D}_k([  \mathbb{T}]_t)],$$ where the limit exists.
\end{proposition}
It will follow from the forthcoming Theorem \ref{prop:manyto1} that the limit above is equal to $2^{-k-1}$, thus proving Proposition  \ref{prop:empirical}, see Section \ref{sec:maxdegreeyule}. A little more effort in the proof enables to prove an almost sure convergence. \medskip 

\noindent \textbf{Proof.} 
The proof crucial relies on the Markov property of the Yule tree: Recall that conditionally on $[ \mathbb{T}]_{t}$ the tree $  [\mathbb{T}]_{t+s}$ is obtained by grafting $\# \partial [ \mathbb{T}]_t$ i.i.d.~copies of $[ \mathbb{T}]_s$ on the leaves of $[ \mathbb{T}]_t$. This enables us to write for any $s,t \geq 0$ the stochastic inequalities
  \begin{eqnarray}\sum_{i=1}^{\# \partial [ \mathbb{T}]_t} \mathrm{D}_k( [ \mathbb{T}^{(i)}]_{s})-\# \partial [ \mathbb{T}]_t \leq  \mathrm{D}_k( [ \mathbb{T}]_{t+s}) \leq  \sum_{i=1}^{\# \partial [ \mathbb{T}]_t} \mathrm{D}_k( [ \mathbb{T}^{(i)}]_{s}) + \# \partial [ \mathbb{T}]_t,   \label{eq:stosandwhich}\end{eqnarray}
where $ \mathbb{T}^{(i)}$ are i.i.d.~standard Yule trees of order $2$ independent of $ \mathbb{T}$. Taking expectation and dividing by $ \mathrm{e}^{-(t+s)}$ we deduce with the shorthand notation $ {d}_{k}(t) = \mathbb{E}[D_{k}( [ \mathbb{T}]_{t})] \mathrm{e}^{	-t}$ that
$$ d_{k}(s) - \mathrm{e}^{-s}\leq d_{k}(t+s) \leq  d_{k}(s) + \mathrm{e}^{-s}.$$ Taking $t>> s>>1$, this shows that ${d}_{k}(t)$ converges as $t \to \infty$ and we denote its limit by $ d_{k}(\infty)$.  Since $  \mathrm{e}^{-t}\# \partial [ \mathbb{T}]_{t} \to X$ almost surely where $X \sim \mathcal{E}(1)$, for any $ \varepsilon>0$, the weak law of large numbers applied twice in \eqref{eq:stosandwhich} shows that with a probability tending to $1$ as $ t \to \infty$ we have
  \begin{eqnarray} X ( d_{k}(s) - \mathrm{e}^{-s})(1 - \varepsilon) \leq \frac{ \mathrm{D}_k( [ \mathbb{T}]_{t+s})}{ \mathrm{e}^{t+s}} \leq (1+ \varepsilon) ( d_{k}(s) + \mathrm{e}^{-s}) X, \end{eqnarray}
and taking again $ t >> s >>1$ large, this implies the convergence in probability claimed in the lemma. \qed \bigskip

\noindent \textbf{Bibliographical notes.} Passing discrete processes into continuous time to get more independence properties is usually called ``randomization'', ``Poissonization'' or ``continuous time embedding'' \cite{athreya1968embedding}. Background on Yule process can be found in \cite{AN72}. Actually, Proposition \ref{prop:yule} is stated there but with a wrong limit law. This has been corrected in \cite[Lemma 3]{bertoin2004dual} with a  proof different from the one presented here. The continuous time-embedding of the O.K. Corral model is taken from \cite{levine2007internal}. The connection between Yule tree and the random growing trees has already been exploited many times in the literature, see e.g.~\cite[Section 3]{janson2021convergence} and the reference therein. The pill problem (Proposition \ref{prop:pills}) has been solved in \cite{kuba2012limiting} using analytic combinatoric. Our solution based on continuous time seems to be new. Proposition \ref{prop:concentrationlocal} (in a more general local version) implies that the random recursive tree converges in the Benjamini--Schramm sense (quenched), see \cite{Ald91c} or \cite[Example 6.1]{holmgren2015limit} for details. \medskip 

\noindent{\textbf{Hints for Exercises.}}\ \\
\noindent Exercise \ref{exo:memory}: The cumulative function $g(s) =  \mathbb{P}(X > s)$ satisfies $g(s/n)^n = g(s)$ for any $s/n$ in the support of the law of $X$. If $ \mathrm{Supp}( \mathcal{L}(X)) = \mathbb{R}_+$, and since $g$ is decreasing, this forces $g(s) = \mathrm{e}^{-\alpha s}$ for some $\alpha > 0$.\\
\noindent Exercise \ref{exo:hideseek}: If $ X_1, \dots , X_n$ are independent exponential variables of parameter $1$ the probability is given by 
$$ \mathbb{P}\left(  \frac{1}{2} \cdot X_n > \max X_1, \dots , X_{n-1}\right) = \int_0^\infty  \mathrm{d}t \, \mathrm{e}^{-t}(1- \mathrm{e}^{-t/2})^{n-1} = \frac{2}{n(n+1)}.$$

\chapter{Spine decomposition and applications}
\label{chap:spine}
\hfill Grow a spine!\bigskip

We describe in this chapter the spine decomposition of Yule trees which will be a key ingredient in our forthcoming applications to the random recursive and Barab\'asi--Albert trees. In particular, it will enable us to prove Theorems \ref{prop:maxdegreeRRT} and \ref{thm:heightRRT} on the max degree and max height in a random recursive tree of size $n$.

\section{Spine decomposition of Yule trees}

 We fix $k \geq 2$ and consider under $ \mathbb{P} \equiv \mathbb{P}_{\delta_1}$ the plane Yule tree $ \mathbb{T}  $ of order $k$ started from a single particle with rates equal to $1$ (see Section \ref{sec:yules}). Recall that for any $t\geq 0$ we denote by $ [\mathbb{T}]_{t}$ the tree $ \mathbb{T}$ cut at level $t$ and write $ \partial [\mathbb{T}]_{t}$ for the boundary of $ [\mathbb{T}]_{t}$ made of all particles alive at time $t$. If $u \in \partial [\mathbb{T}]_t$ is a particle living at time $t$ on the Yule tree, we denote by $[\mathbb{T}]_t^u$ the tree obtained from $[\mathbb{T}]_t$ by distinguishing the branch going from the root to the particle $u$ living at height $t$. We also use the notation $ \# \partial [ \mathbb{T}]_{t}$ for the number of particles alive at time $t$ in $ \mathbb{T}$ (this was abbreviated by $ \mathcal{Y}^{(k)}_{t}$ in the previous chapter).

 \subsection{Martingale transform}
 
This section, rather abstract, can be skipped at first reading. It presents the spine construction in a broader context, that of \textbf{martingale transformation}. We do not aim at the same level of rigor as in the rest of these pages and just hope to pique the reader's interest. Those willing to proceed with the applications should take Theorem \ref{prop:manyto1} (the many-to-one formula) as granted. \medskip

 In general, a positive martingale $ (M_{n} : n \geq 0)$ over a filtered $( \mathcal{F}_{n} : n \geq 0)$ probability space enables us to change the underlying measure $ \mathbb{P}$ by biasing with the martingale $(M)$, see Exercise \ref{exo:martingalebiais} for a toy model. This is the essence of the famous ``Girsanov transformation'' in continuous stochastic calculus, and let us see the effect of this transformation when applied to Yule trees with the martingale identified in the previous chapter.

 Recall from \eqref{eq:martingaleyule} that the process $ M_{t} := \mathrm{e}^{-(k-1)t} \# \partial [\mathbb{T}]_{t}$ is a martingale starting from $1$ for the filtration $ \mathcal{F}_{t} = \sigma ( [\mathbb{T}]_{s} : 0 \leq s \leq t)$. When in possession of such a positive martingale, one can perform a change of measure by biasing the underlying random variables by this martingale. 
 Specifically, this is obtained by considering the probability $ \mathbb{Q}_{t}$ whose Radon--Nikodym  derivative with respect to the underlying probability $ \mathbb{P}$ is 
$$ \left. \frac{ \mathrm{d} \mathbb{Q}_{t}}{ \mathrm{d} \mathbb{P}}\right|_{ \mathcal{F}_{t}} = M_{t}.$$
Actually, since $ M_{t}$ is a martingale, this change of measure is coherent in the sense that for $0 \leq s \leq t$ we have $ \left.\mathbb{Q}_{t}\right|_{ \mathcal{F}_{s}} = \mathbb{Q}_{s}$. This can be checked by a one-line calculation using the martingale property: for any positive measurable function $F$ we have 
$$ \mathbb{E}[M_{t} F( [[\mathbb{T}]_{t}]_{s})] = \mathbb{E}[ \mathbb{E}[M_{t} F( [\mathbb{T}]_{s}) \mid \mathcal{F}_{s}]] =  \mathbb{E}[M_{s} F( [\mathbb{T}]_{s})].$$
By coherence of the restrictions (and leaving the details of the topology, restriction ... to the courageous reader) one can thus define a probability measure $ \mathbb{Q}$ under which the random infinite tree $ \mathbb{T}$ has the property that 
$$ [ \mathbb{T}]_{t}\quad  \mbox{ under } \mathbb{Q} \qquad \overset{(d)}{=} \qquad  [ \mathbb{T}]_{t} \quad \mbox{ under } M_{t} \cdot  \left. \mathrm{d} \mathbb{P}\right|_{ \mathcal{F}_{t}} = \mathbb{Q}_{t}.$$
Now, if $ \mathbb{T}_{t}^{\bullet}$ is obtained under $ \mathbb{Q}_{t}$ by distinguishing a particle of $\partial [  \mathbb{T}]_{t}$ uniformly at random (this actually distinguishes a branch in $[ \mathbb{T}]_t$), the same calculation as above enables us to see that the tree with distinguished branch obtained by restricting up to height $s$ has the same law as $ \mathbb{T}_{s}^{\bullet}$. By coherence of the restriction (and again leaving the details   to the courageous reader) one can thus define a probability measure $ \mathbb{Q}$ and a random infinite tree $ \mathbb{T}^{\bullet}$ with an \textbf{infinite line of descent} so that for each $t$ the finite tree $ [ \mathbb{T}^{\bullet}]_{t}$ obtained by restricting to height $t$ and keeping the distinguished branch, has the distribution of $ \mathbb{T}_{t}^{\bullet}$ under $ \mathbb{Q}_{t}$.

\begin{exo}[An example of martingale transform] \label{exo:martingalebiais} Let $(S_n : n \geq 0)$ be a simple symmetric random walk \textit{started from $1$}. We denote by $ \mathcal{F}_n$ its canonical filtration such that if $\tau_0 = \inf \{ k \geq 0 : S_k=0\}$ then the process $M_n = S_{n \wedge \tau_0}$ is a non-negative martingale. As above define the law $ \mathbb{Q}$ so that 
$$ \left. \frac{ \mathrm{d} \mathbb{Q}}{ \mathrm{d} \mathbb{P}}\right|_{ \mathcal{F}_{n}} = M_{n}.$$
Show that under $ \mathbb{Q}$ the process $(S_n : n \geq 0)$ is a Markov chain with probability transitions 
$$ \mathbb{Q}(S_{n+1}= S_n +1\mid S_n = i ) = \frac{i+1}{2i}, \quad  \mathbb{Q}(S_{n+1}= S_n -1\mid S_n = i ) = \frac{i-1}{2i}, \quad i \geq 1.$$
\end{exo} 	

\subsection{Spine decomposition}
The law of $ \mathbb{T}^{\bullet}$ under $ \mathbb{Q}$ is actually quite simple to describe. Consider a continuous time branching tree as in Section \ref{sec:AK} with \textbf{two types} of particles: standard particles of type $1$ which reproduce at rate $1$ and mutant particles of type $2$ which reproduce at rate $k\geq 2$. When a standard particle dies, it gives rise to $k$ standard particles, but when a mutant particle dies it gives rise to $k-1$  standard particles (type $1$) and a single mutant particle (type $2$). Actually, since we shall consider them as plane trees, we need to prescribe an ordering in the case of reproduction of a mutant, by placing the mutant descendant  uniformly among its children. We then consider the random plane tree  $ \mathbb{T}$ under the measure $ \mathbb{P}_{\delta_2}$ started with only one mutant: it is clear that there is a single line of descent composed of mutant particles and this defines a random tree with a distinguished ray $  \mathbb{T}^{ \bullet}$.

\begin{proposition}[Description of the law of $ \mathbb{T}^\bullet$]  \label{prop:spinedecomp}The law of $ \mathbb{T}^{\bullet}$ under $\mathbb{Q}$ is that of $ \mathbb{T}^{\bullet}$ under $ \mathbb{P}_{\delta_2}$. 
\end{proposition}
\begin{figure}[!h]
 \begin{center}
 \includegraphics[width=13cm]{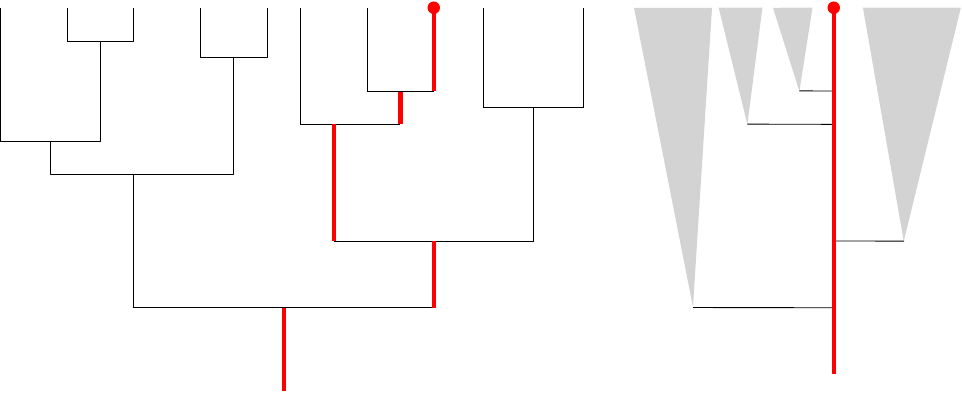}
 \caption{The law of the pointed tree $[\mathbb{T}^{\bullet}]_{t}$ under $ \mathbb{Q}$ is the same as that of the Yule tree started with a mutant particle. In particular, when $k=2$ the ancestral line (Right on the figure) from the distinguished point to the root in $[\mathbb{T}^{\bullet}]_{t}$ under $ \mathbb{Q}$ is obtained by superimposing two independent Poisson processes of intensity $1$ for each side of the spine.}
 \end{center}
 \end{figure}
Before giving the proof, let us provide the reader with an equivalent formulation, the so-called ``\textbf{Many-to-one formula}'', which can be read without reference to the measure $ \mathbb{Q}$. It will be very practical for applications as it  enables us to perform first-moment calculation over all branches:
\begin{theorem}[Many-to-one formula] \label{prop:manyto1}\noindent For any positive and measurable function $F$ we have 
$$ \mathbb{E}_{\delta_1}\left[ \sum_{u \in \partial  [\mathbb{T}]_{t}} F(  [\mathbb{T}]_{t}^{u})\right] = \mathrm{e}^{(k-1)t} \cdot  \mathbb{E}_{\delta_2}\left[  F([\mathbb{T}^{\bullet}]_{t}) \right].$$
\end{theorem}
\noindent \textbf{Proof.} By the definition of the objects we have with $ \mathbb{E}\equiv \mathbb{E}_{\delta_1}$
  \begin{eqnarray*}\mathbb{E}\left[ \sum_{u \in \partial  [\mathbb{T}]_{t}} F(  [\mathbb{T}]_{t}^{u})\right] &=& \frac{\mathbb{E}[\# \partial [ \mathbb{T}]_t]}{ \mathbb{E}[\# \partial [ \mathbb{T}]_t]} \cdot \mathbb{E}\left[ \frac{\# \partial [ \mathbb{T}]_t}{\# \partial [ \mathbb{T}]_t} \sum_{u \in \partial  [\mathbb{T}]_{t}} F(  [\mathbb{T}]_{t}^{u})\right] \\
  &=& \mathbb{E}[\# \partial [ \mathbb{T}]_t] \cdot \mathbb{E}_{ \mathbb{Q}}[F(\mathbb{T}_t^\bullet)] \\
  &=& \mathrm{e}^{(k-1)t} \mathbb{E}_{\delta_2}\left[  F([\mathbb{T}^{\bullet}]_{t}) \right].  \end{eqnarray*}\qed  \medskip 
   
\noindent \textbf{Proof of the Proposition \ref{prop:spinedecomp}.} We consider the set $ \mathcal{S}_{t}$ (resp. $ \mathcal{S}^\bullet_{t}$)   of all plane trees  $\tau$ (resp $\tau^\bullet$) where each vertex has $0$ or $k$ children and endowed with vertex lengths $ (\ell(u) : u \in \mathrm{Vertices}(\tau))$ so that the $\ell$-height (the sum of the vertex lengths from a vertex to the root) of all its leaves is exactly $t$ (resp. with a distinguished leaf $\bullet$). Recall that $[ \mathbb{T}]_t$ (resp. $[ \mathbb{T}]_t^u$ or $[ \mathbb{T}^\bullet]_t$) can be seen as an element of $ \mathcal{S}_t$ (resp. $ \mathcal{S}_t^\bullet$), see Figure \ref{fig:RRTbinary}. There is a natural measure on $ \mathcal{S}_{t}$ (resp. $ \mathcal{S}_{t}^\bullet$) obtained as the sum for each finite plane tree $\tau$ as above of the product of the Lebesgue measure  for each $\ell(u)\geq 0$ for all non leaves $u \in \tau$, subject to the condition that the sum of all $\ell(u)$ for all $u$ from the root to a leaf stays below $t$ (the label of a leaf $v$ is then obtained as $t - \sum \ell(u)$ where the sum runs of all ancestors of $v$). By construction of the (plane) Yule tree, the law of  $[\mathbb{T}]_{t}$  under $ \mathbb{P}_{\delta_1}$ is absolutely continuous with respect to the above measure on $ \mathcal{S}_t$ with density given by 
 $$   \prod_{u \in \tau \backslash \mathrm{Leaves}(\tau)} \mathrm{e}^{-\ell(u)} \prod_{v \in \mathrm{Leaves}(\tau)} \mathbb{P}( \mathcal{E} \geq \ell(v)) = \mathrm{\exp}\left(- \sum_{  u \in \tau} \ell(u)\right),$$ so that the law of $ [\mathbb{T}^\bullet]_t$ under $ \mathbb{Q}$ has density with respect to the above measure on $ \mathcal{S}_t^\bullet$ given by
 $$   \mathrm{\exp}\left(- \sum_{  u \in \tau^\bullet} \ell(u)\right) \times  \mathrm{e}^{-(k-1)t}  \times \frac{ \# \mathrm{Leaves}(\tau^\bullet)}{\# \mathrm{Leaves}(\tau^\bullet)} = \mathrm{\exp}\left(- \sum_{  u \in \tau^\bullet} \ell(u)\right) \mathrm{e}^{-(k-1)t},$$
On the other hand, the law of $[ \mathbb{T}^\bullet]_t$ under the two-type measure $ \mathbb{P}_{\delta_2}$ is also absolutely continuous with respect to the above measure: taking separately the behavior of the mutant particles along $ \mathrm{Spine}(\tau^\bullet)$, the  path going from the root to the distinguished leaf, this density is seen to be
  \begin{eqnarray*} &=& \exp\left(- \sum_{u \in \tau^\bullet  \backslash \mathrm{Spine}(\tau^\bullet) } \ell(u)\right) \times \left(\prod_{u \in \mathrm{Spine} \backslash \bullet}  k \mathrm{e}^{- k \ell(u)} \cdot \frac{1}{k}\right) \times \underbrace{  \mathbb{P}( \mathcal{E}(k) \geq \ell(\bullet))}_{ \mathrm{e}^{-k \ell( \bullet)}}\\
  &=&\mathrm{\exp}\left(- \sum_{  u \in \tau^\bullet} \ell(u)\right) \mathrm{e}^{-(k-1)t}.   \end{eqnarray*} Since the last two displays agree we have proved the proposition. \qed \medskip

 \section{Application to extreme geometric properties of the RRT}
Recall the construction of the random recursive tree $(T_n : n \geq 0)$ from the plane Yule tree $ \mathbb{T}$ described in Proposition \ref{prop:RRTyule}: in this section we shall suppose that $T_n = \{ \mathbb{T}\}_{\tau_n}$ where $(\tau_i : i \geq 0)$ are the jump times of the particle counting process  and where $\{ \mathbb{T}\}_t$ is the increasing labeled tree obtained from $[ \mathbb{T}]_t$ by ``contracting'' the edges going to the left and numbering the vertices by order of appearance. We use the spinal decomposition to give quick proofs of the two results that were left unproven in Chapter \ref{chap:RRT}.
 \subsection{Maximal Height in RRT} \label{sec:heightRRT}
We recall Theorem \ref{thm:heightRRT}  here: For the random recursive tree $(T_{n} : n \geq 0)$ we have $$ \frac{\mathrm{Height}( T_{n})}{\log n} \xrightarrow[n\to\infty]{ a.s.} \mathrm{e}.$$
 
\noindent \textbf{Proof of Theorem \ref{thm:heightRRT}.} From Proposition \ref{prop:RRTyule} we can write $(T_n : n \geq 0) = (\{ \mathbb{T}\}_{\tau_n}: n \geq 0)$ where $\tau_n$ is the first time when there are $n+1$ particles alive in the Yule tree. Recall from Proposition \ref{prop:yule} and Eq.\ \eqref{eq:taynb} that $\tau_{n} \sim \log n$ almost surely as $n \to \infty$. Hence, by Lemma \ref{lem:slut}, the above theorem is a consequence of the previous two remarks provided that we prove 
$$ \left(\frac{\mathrm{Height}(\{ \mathbb{T}\}_{x \log n})}{\log n} : x \geq 0\right)   \xrightarrow[t\to\infty]{a.s.} (x \cdot \mathrm{e} : x \geq 0), $$ for the uniform convergence over every compact subset of $ \mathbb{R}_+$.  Since the height of $\{ \mathbb{T} \}_t$ is increasing with $t$, by Lemma \ref{lem:dini} it  suffices to prove that    \begin{eqnarray} \label{eq:goalRRTheight} \frac{\mathrm{Height}(\{ \mathbb{T}\}_t)}{t}  \xrightarrow[t\to\infty]{ a.s.} \mathrm{e}.  \end{eqnarray} 
Now, recall from the construction of Section \ref{sec:rrtfromyule} that each particle $u \in \partial [ \mathbb{T}]_{t} $ is associated with a vertex in $\{ \mathbb{T}\}_{t}$, which we still denote $u$ by abuse of notation, whose distance to the root $\raisebox{.5pt}{\textcircled{\raisebox{-.9pt} {$0$}}}$ of $\{ \mathbb{T}\}_t$ satisfies 
 \begin{eqnarray} \label{eq:distancespine} \mathrm{Dist}_{\{ \mathbb{T}\}_{t}}(\raisebox{.5pt}{\textcircled{\raisebox{-.9pt} {$0$}}}, u) =  \begin{array}{l} \mbox{ number of ancestral lineages pointing to the left} \\ \mbox{ in the spine between } u \mbox{ and the root in } [\mathbb{T}]_t. \end{array}  \end{eqnarray}
 
 Let us start with the easy upper bound for \eqref{eq:goalRRTheight}.\\
\noindent \texttt{Upper bound}. By the many to one formula  (Theorem \ref{prop:manyto1}) we have 
  \begin{eqnarray*} \mathbb{P}\left( \mathrm{Height}(\{ \mathbb{T}\}_t) \geq x\}\right) &\leq & \mathbb{E}\left[ \# \{ u \in \partial [\mathbb{T}]_{t} : \mathrm{Dist}_{\{ \mathbb{T}\}_{t}}(\raisebox{.5pt}{\textcircled{\raisebox{-.9pt} {$0$}}}, u) \geq x\}\right]\\
   &\underset{ \mathrm{Thm.} \ref{prop:manyto1}}=& \mathrm{e}^{t} \cdot \mathbb{P}_{\delta_2}( \mathrm{Dist}_{\{ \mathbb{T}^{\bullet}\}_{t}}(\raisebox{.5pt}{\textcircled{\raisebox{-.9pt} {$0$}}}, \bullet) \geq x) \\ 
&\underset{\eqref{eq:distancespine}}{=}&   \mathrm{e}^{t} \cdot  \mathbb{P}( \mathfrak{P}(t) \geq x),  \end{eqnarray*} where $( \mathfrak{P}(t) : t \geq 0)$ is a standard Poisson counting process. When $x  = (  \mathrm{e}+ \varepsilon)t$ for $ \varepsilon>0$ small, Lemma \ref{lem:LDpoisson} entails that the above probability decays to $0$ exponentially fast in $t$. By Markov's inequality and the Borel-Cantelli Lemma we deduce that $\mathrm{Height}(\{\mathbb{T}\}_{n}) \leq  ( \mathrm{e}+ \varepsilon) n$ eventually for $n \in \mathbb{Z}_{ >0}$ large enough   $ \mathbb{P}$-a.s. Since $ t \mapsto  \mathrm{Height}(\{\mathbb{T} \}_{t})$ is increasing, the same holds true when the integer $n$ is replaced by $t >0$.\\
\noindent \texttt{Lower bound}. By the previous calculation, we know that the expected number of branches $u \in \partial [\mathbb{T}]_{t}$ corresponding to a vertex $u$ at height $ \geq (  \mathrm{e} - \varepsilon)t$ in $\{ \mathbb{T}\}_t$ tends to $\infty$ exponentially fast with $t$. As usual, this does not imply right away that the number of such branches is non zero with high probability. However, this fact can be used together with the branching property of $  \mathbb{T}$: Fix $t_{0}>0$ large enough so that 
  \begin{eqnarray} \label{eq:meangeq2}  \mathbb{E}\left[ \# \big\{ u \in \partial [\mathbb{T}]_{t_0} : \mathrm{Dist}_{\{ \mathbb{T}\}_{t_0}}(\raisebox{.5pt}{\textcircled{\raisebox{-.9pt} {$0$}}}, u) \geq (  \mathrm{e} - \varepsilon)t_0\big\}\right] =\mathrm{e}^{t_0} \cdot  \mathbb{P}( \mathfrak{P}(t_0) \geq ( \mathrm{e}- \varepsilon)t_0) \geq 2.  \end{eqnarray} We now consider the branching process obtained by restricting the Yule tree to times $k \cdot t_{0}$ for $ k \in \{1,2, \dots\}$ and considering those particles $ u \in \partial [ \mathbb{T}]_{k t_{0}}$ for which there are at least $( \mathrm{e}- \varepsilon)t_0$ ancestral lineages pointing to the left between time $k t_{0}$ and time $ (k-1)t_{0}$ in $ \mathbb{T}$. By the Markov property of the Yule tree, those ``particles'' form a Bienaymé--Galton--Watson tree in discrete time $k \geq 0$ whose mean offspring is larger than $2$ by \eqref{eq:meangeq2}, so it survives with positive probability. Hence, there exists a random generation $0 \leq M < \infty$ from which the branching process survives on. For $k \geq M$, a particle $u \in \partial [\mathbb{T}]_{k t_{0}}$ in this branching process   has the property that 
$$ \mathrm{Dist}_{\{ \mathbb{T}\}_{kt_0}}(\raisebox{.5pt}{\textcircled{\raisebox{-.9pt} {$0$}}}, u) \geq ( \mathrm{e}- \varepsilon) (k-M) t_{0}\quad  ,$$ which easily entails the lower bound $\mathrm{Height}(\{ \mathbb{T}\}_t) \geq ( \mathrm{e}-2 \varepsilon)t$ for $t$ large enough a.s. \qed 

 \subsection{Maximal degree in RRT} \label{sec:maxdegreeyule}
 We now prove Theorem \ref{prop:maxdegreeRRT} on the maximum degree in $T_n$ which we also recall for the reader's convenience:  Let $  \mathrm{MaxDegree}(T_n)$ be the largest vertex (out)-degree in the random recursive tree $T_{n}$. Then as $n \to \infty$ we have 
$$ \frac{ \mathrm{MaxDegree}(T_n)}{ \log_{2} n} \xrightarrow[n\to\infty]{a.s.} 1.$$

\noindent \textbf{Proof of Theorem \ref{prop:maxdegreeRRT}.} As in the previous section, since $ t \mapsto \mathrm{MaxDegree}(\{ \mathbb{T}\}_t)$ is increasing in $t$ and by virtue of \eqref{eq:taynb} it is sufficient to prove that 
  \begin{eqnarray} \label{eq:goalddegreeyule} \frac{\mathrm{MaxDegree}( \{ \mathbb{T}\}_t)}{t} \xrightarrow[t\to\infty]{a.s.} \frac{1}{\log(2)}.  \end{eqnarray}
As for the height,  if $u \in \partial  [\mathbb{T}]_{t}$, we can read on $[ \mathbb{T}]_t$ the degree of $u$ inside $\{ \mathbb{T}\}_t$: it is easy by looking at Figure \ref{fig:yuleRRT}  to convince oneself that we have 
 \begin{eqnarray} \label{eq:degreespine} \mathrm{deg}^+_{\{ \mathbb{T}\}_{t}}( u) =  \begin{array}{l} \mbox{ number of ancestral lineages pointing to the right} \\ \mbox{ in the spine between } u \mbox{ and the root in } \mathbb{T} \\ \mbox{ before the first ancestral lineage pointing to the left} . \end{array}  \end{eqnarray}
 We now proceed separately with the upper and lower bounds for \eqref{eq:goalddegreeyule}. We set 
 $$ \beta = \frac{1}{\log 2},$$ to ease notation.

\noindent \texttt{Upper bound}.  In the two-type tree $ [ \mathbb{T}^\bullet]_t$ under $ \mathbb{P}_{\delta_2}$, the branching events to the left and right of the mutant branch are independent and appear as Poisson processes with intensity $1$. The number of lineages branching to right before encountering a lineage branching to the left is then stochastically bounded from above by a geometric random variable with parameter $1/2$. By the many to one formula we thus have for $x \geq 1$ 
\begin{eqnarray*} \mathbb{P}\left( \exists u  \in [ \mathbb{T}]_{t}: \mathrm{deg}^+_{\{ \mathbb{T} \}_t}(u) \geq x \right) 
&\leq& \mathbb{E}\left[ \# \{ u \in \partial [\mathbb{T}]_{t} : \mathrm{deg}^+_{\{ \mathbb{T}  \}_t}(u) \geq x\}\right] \\ 
& \underset{ \mathrm{Thm.} \ref{prop:manyto1}}= & \mathrm{e}^{t} \mathbb{P}( \mathrm{deg}^+_{\{ \mathbb{T}^\bullet  \}_t}(\bullet) \geq x ) \\ 
& \underset{ \eqref{eq:degreespine}}{\leq} & \mathrm{e}^{t} \mathbb{P}( \mathrm{Geo}(1/2) \geq x) = \mathrm{e}^{t} \,2^{-x}.  
\end{eqnarray*}
If $x = (\beta +\varepsilon)t$ the above display goes to $0$ exponentially fast in $t$. We conclude using the Borel-Cantelli lemma and monotonicity as in the previous proof that $\mathrm{MaxDegree}(\{ \mathbb{T}  \}_t) \leq ( \frac{1}{\log 2} + \varepsilon) t$ for all $t$ large enough a.s.\\
\texttt{Lower bound}.  Let us consider all particles alive at time $t' =  t (1-  \frac{\beta}{2}) \approx t \times  0,27\dots$ inside $ [\mathbb{T}]_t$. 
Using the independence property of the Yule tree, and by considering only the monochromatic branches going from time $t'$ to time $t$ in $ [\mathbb{T}]_t$ (always turning left) we deduce that 
$$ \mathrm{MaxDegree}( \{ \mathbb{T}  \}_t) \geq \max_{1 \leq i \leq \# \partial [ \mathbb{T}]_{t'}} X_i,$$ where conditionally on $ \# \partial [ \mathbb{T}]_{t'}$ the variables $X_i$ are independent and of law $ \mathfrak{P}(t-t') = \mathfrak{P}(t \beta/2)$.  By Proposition \ref{prop:yule} we have $\# \partial [ \mathbb{T}]_{t'} \sim_{a.s.} \mathcal{E}(1) \mathrm{e}^{ (  1-\beta/2)t}$. In the notation of Lemma \ref{lem:LDpoisson},  an easy computation shows that with $\beta = \frac{1}{\log 2}$ we have 
$$ (1- \frac{\beta}{2}) - \frac{\beta}{2} I \left( \frac{\beta}{\beta/2}\right) =  0,$$ so that for $\varepsilon>0$ there exists $\delta >0$ with $(1- \frac{\beta}{2}-\delta)  - \frac{\beta}{2} I \left( \frac{\beta - \varepsilon}{\beta/2}\right)>0$. In particular  
  \begin{eqnarray*} &&\mathbb{P}\Big(\mathrm{MaxDegree}(\{ \mathbb{T}  \}_t) \leq t( \beta- \varepsilon) \ \big| \  \# \partial [ \mathbb{T}]_{t'} \geq \mathrm{e}^{ (1- \frac{\beta}{2}-\delta) t}\Big) \\
  & \leq &   \Big(1- \mathbb{P}( \mathfrak{P}( t\beta/2) > (\beta- \varepsilon) t)\Big)^{\mathrm{e}^{ (1- \frac{\beta}{2}-\delta) t}} \\
   & \underset{ \mathrm{Lem.} \ref{lem:LDpoisson}}{\leq} & \exp \left( - \mathrm{e}^{ (1- \frac{\beta}{2}-\delta) t}  \mathrm{e}^{ -\frac{\beta}{2} I \left( \frac{\beta - \varepsilon}{\beta/2}\right) t}\right) \leq \exp(  - \mathrm{e}^{ \mathrm{cst}\, t}),\end{eqnarray*} for some $ \mathrm{cst}>0$. Since the right-hand side is summable for $t \in \mathbb{Z}_{>0}$ and since eventually $\# \partial [ \mathbb{T}]_{t'} \geq \mathrm{e}^{ (1- \frac{\beta}{2}-\delta) t}$ with probability one, we deduce from the Borel--Cantelli lemma that 
$\mathrm{MaxDegree}( \{ \mathbb{T}  \}_t) \geq  t(\beta - \varepsilon)$ eventually along integer values of $t$. By monotonicity the same holds for all $t$ large enough and this concludes the proof.\qed

\begin{remark}  \label{rek:ouestgros} The proof of Theorem \ref{prop:maxdegreeRRT} actually shows that the maximal degree in the random recursive tree $T_n$ is attained by a vertex  $\raisebox{.5pt}{\textcircled{\raisebox{-.9pt} {$i$}}}$ with $i \approx n^{(1-\beta/2) + o(1)} = n^{0.27\dots}$. This may be seem counterintuitive since the vertex $\raisebox{.5pt}{\textcircled{\raisebox{-.9pt} {$0$}}}$ clearly has the largest degree for the stochastic order. 
\end{remark}

The many-to-one formula and Equation \eqref{eq:degreespine}  directly show that the limit appearing
in Proposition \ref{prop:concentrationlocal} is equal to $2^{-k-1}$ as announced after the proposition.
\bigskip 

\noindent \textbf{Bibliographical notes.} Spinal decomposition (and the associated many-to-one formula) is a very important tool in the theory of branching processes. Although it had precursors e.g. \cite{chauvin1988kpp}, this method has  been popularized by Lyons, Pemantle and Peres \cite{LPP95b}. See also \cite{shi2015branching} for its numerous applications in branching random walk or \cite{AD15} for discrete Bienaym\'e--Galton--Watson trees. In general, martingale change of measures are frequently met in probability theory ($h$-transforms, Girsanov formula, exponential tiltings...). See \cite{addario2018high} and \cite{addario2013poisson} for recent results about maximal degree and height of random recursive trees. \medskip

\noindent{\textbf{Hints for Exercises.}}\ \\

\noindent Exercise \ref{exo:martingalebiais}: Biaising by the martingale is equivalent to performing a $h$-transformation with the function $h : i \mapsto i$ which is harmonic for the walk killed at $\tau_0$. See \cite[Appendix A.3]{CurStFlour} for more.\\

\chapter{Barab\'asi-Albert preferential attachment tree}
\hfill Rich get richer.\bigskip

In this chapter we modify the RRT construction using a preferential attachment rule:
\begin{definition}[BA] The \textbf{Barab\'asi--Albert (BA) preferential attachment tree} is the Markov chain with values in the set of unoriented  labeled trees such that $ \mathsf{T}_1 =$  \raisebox{.5pt}{\textcircled{\raisebox{-.9pt} {$0$}}}--\raisebox{.5pt}{\textcircled{\raisebox{-.9pt} {$1$}}} and so that for $n \geq 2$, conditionally on $ \mathsf{T}_{n-1}$, the labeled tree $ \mathsf{T}_{n}$ is obtained by attaching the new vertex  \raisebox{.5pt}{\textcircled{\raisebox{-.6pt} {$n$}}}  onto the vertex \raisebox{.5pt}{\textcircled{\raisebox{-.9pt} {$k$}}} of $ \mathsf{T}_{n-1}$ with probability 
$$ \mathbb{P}\left( \raisebox{.5pt}{\textcircled{\raisebox{-.9pt} {$n$}}} \to \raisebox{.5pt}{\textcircled{\raisebox{-.9pt} {$k$}}} \big| \mathsf{T}_{n-1}\right) = \frac{ \mathrm{deg}_{ \mathsf{T}_{n-1}}(\raisebox{.5pt}{\textcircled{\raisebox{-.9pt} {$k$}}})}{2(n-1)}.$$ 
\end{definition}

 \begin{figure}[!h]
 \begin{center}
 \includegraphics[height=3cm]{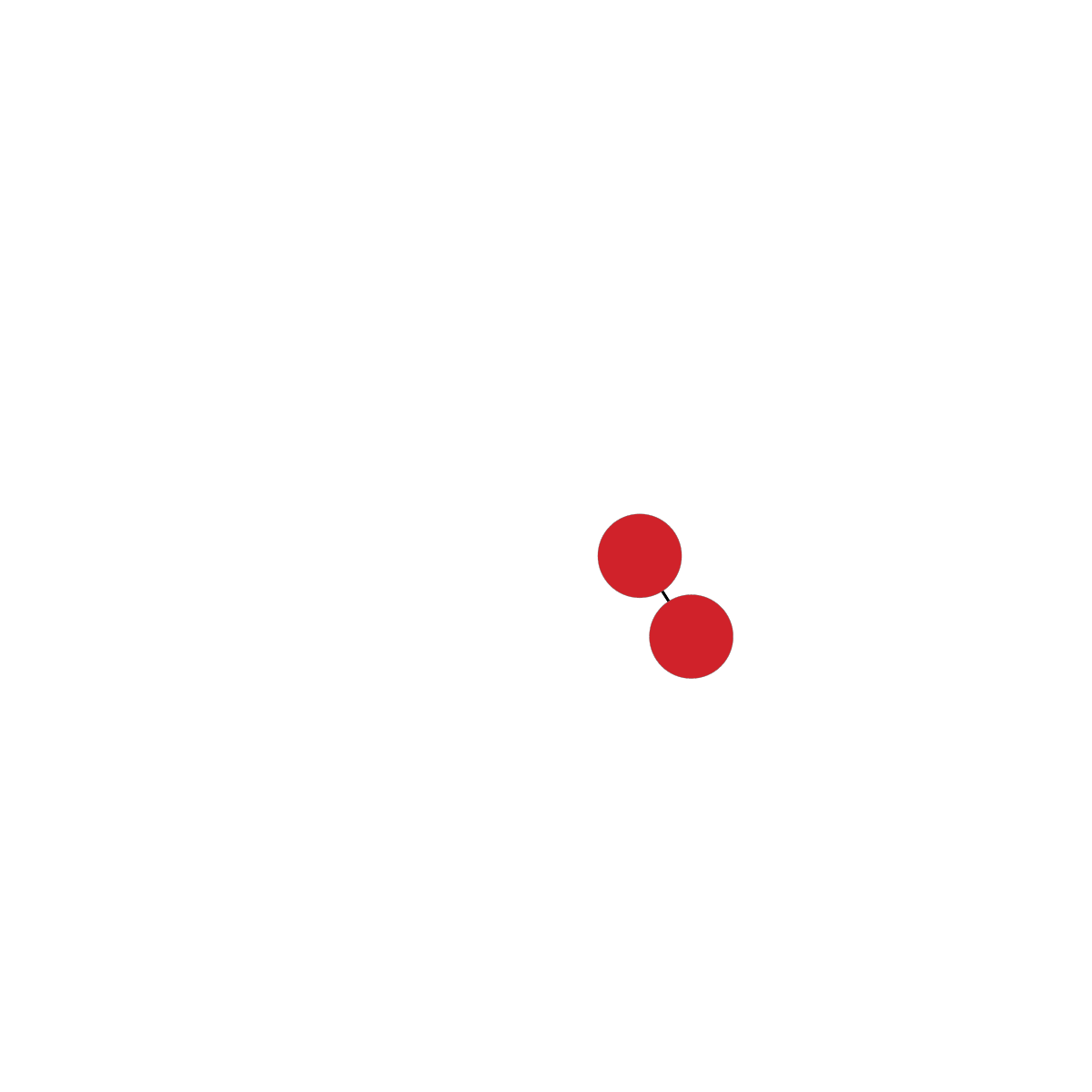}
     \includegraphics[height=2.95cm]{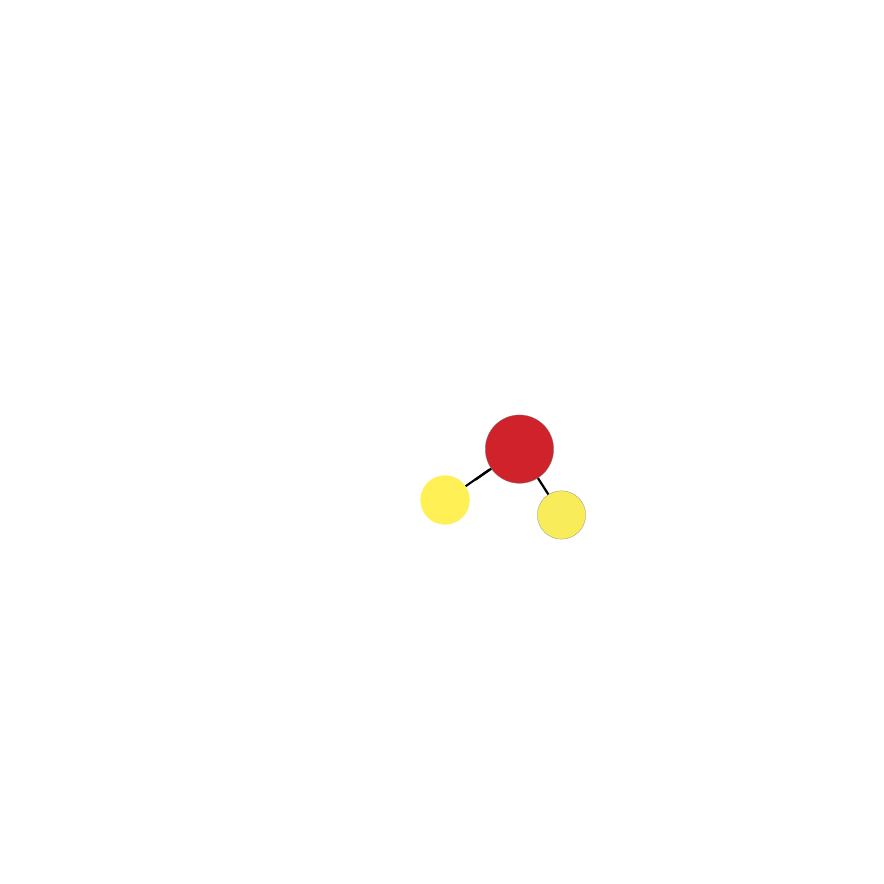}
  \includegraphics[height=3cm]{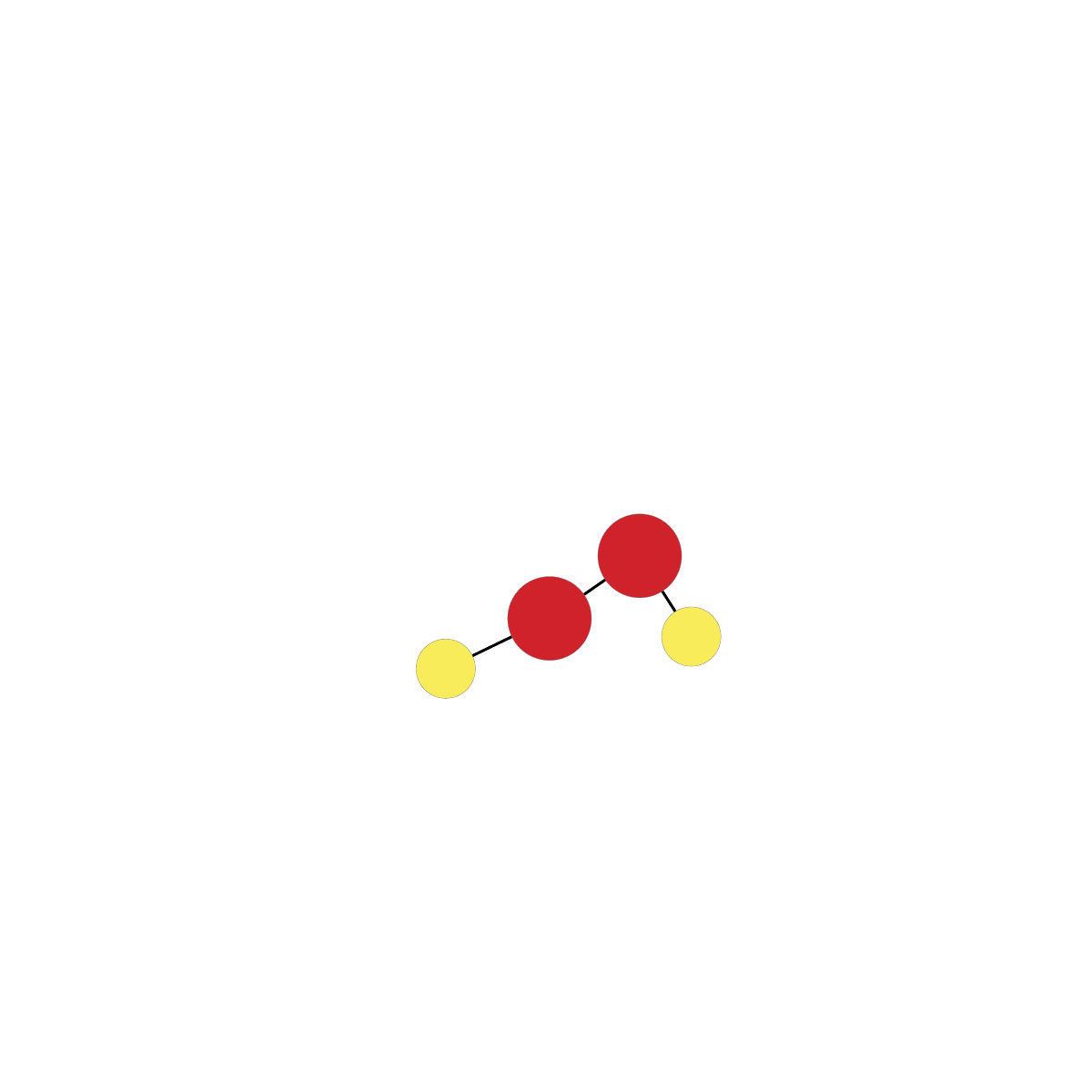}
   \includegraphics[height=3cm]{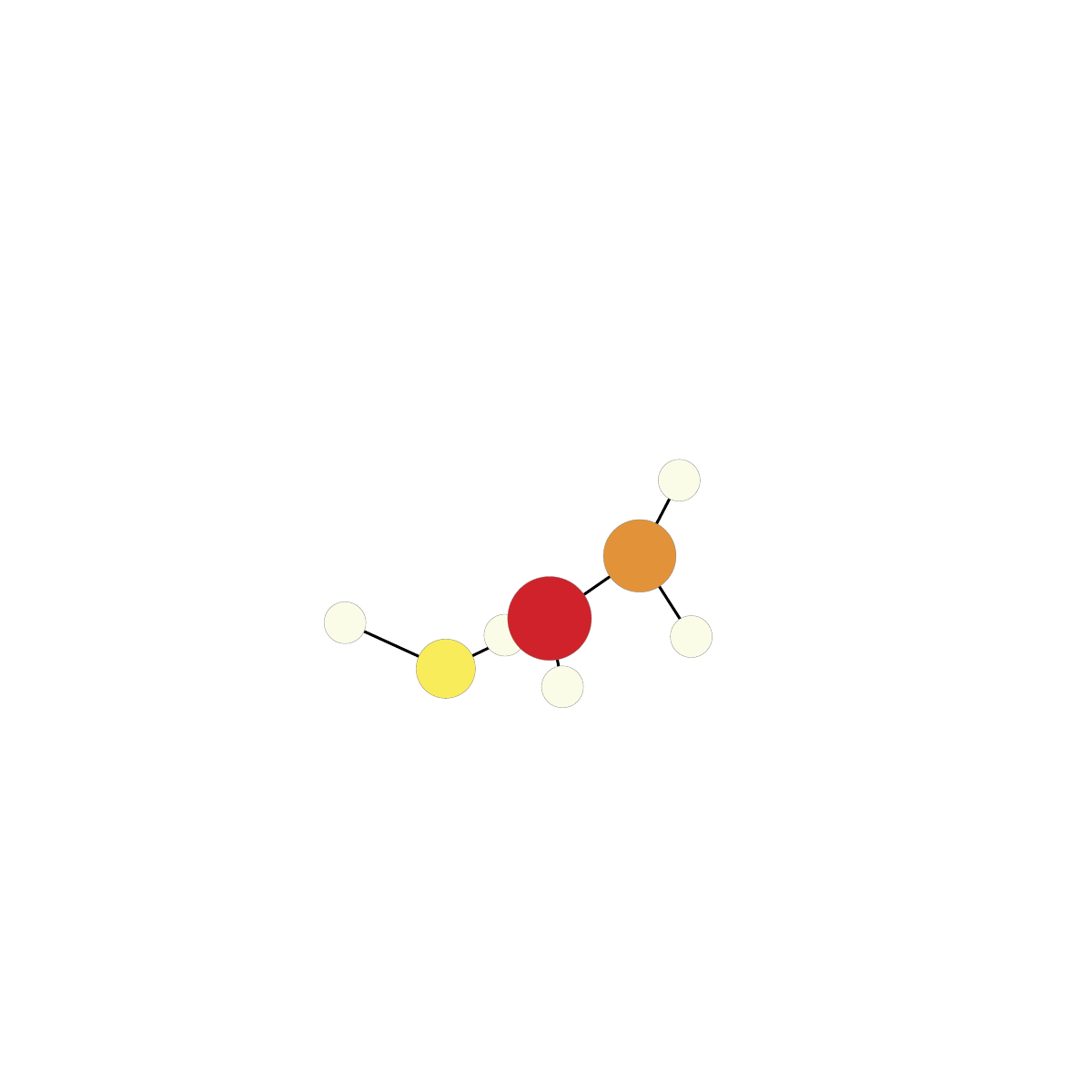}
    \includegraphics[height=3cm]{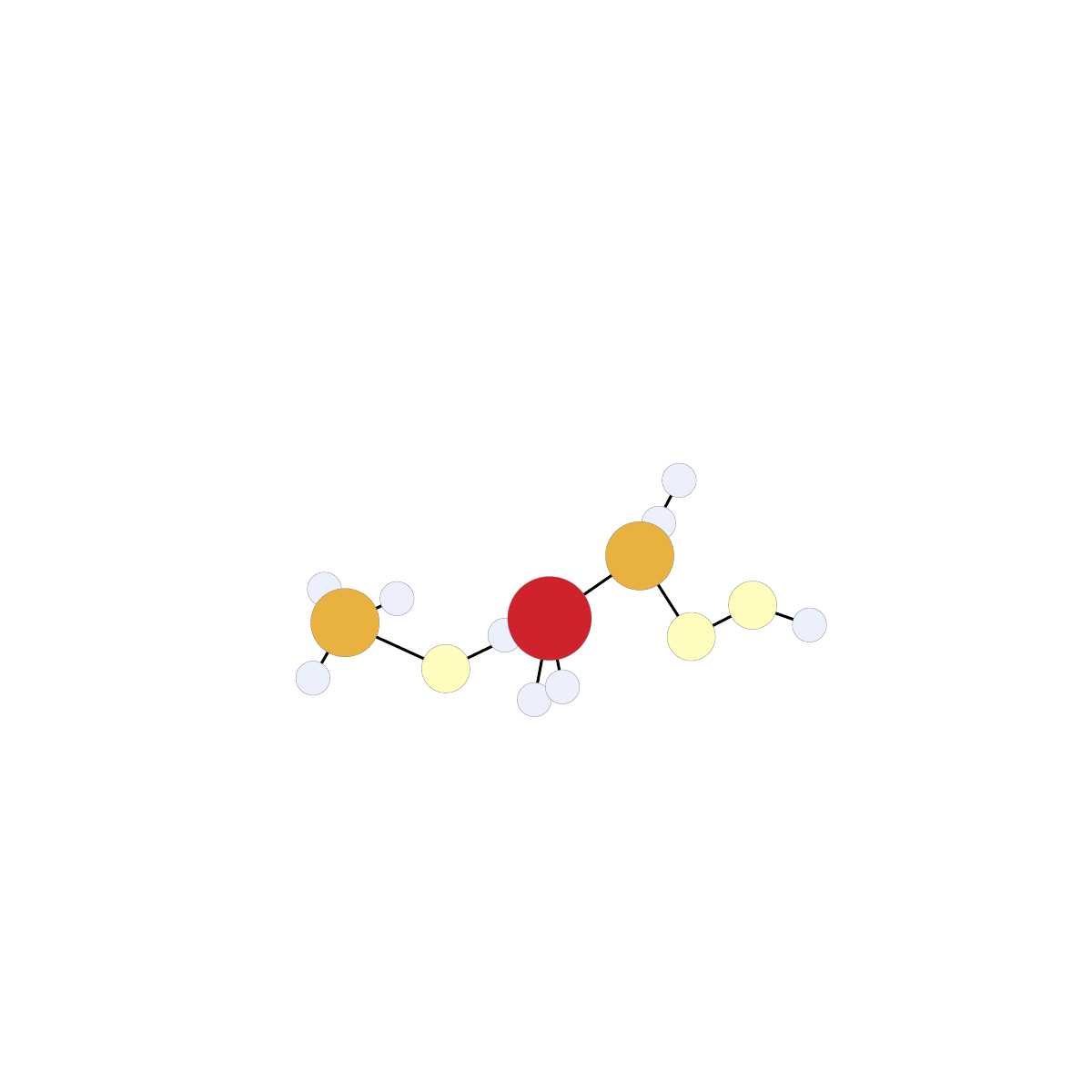}
     \includegraphics[height=3cm]{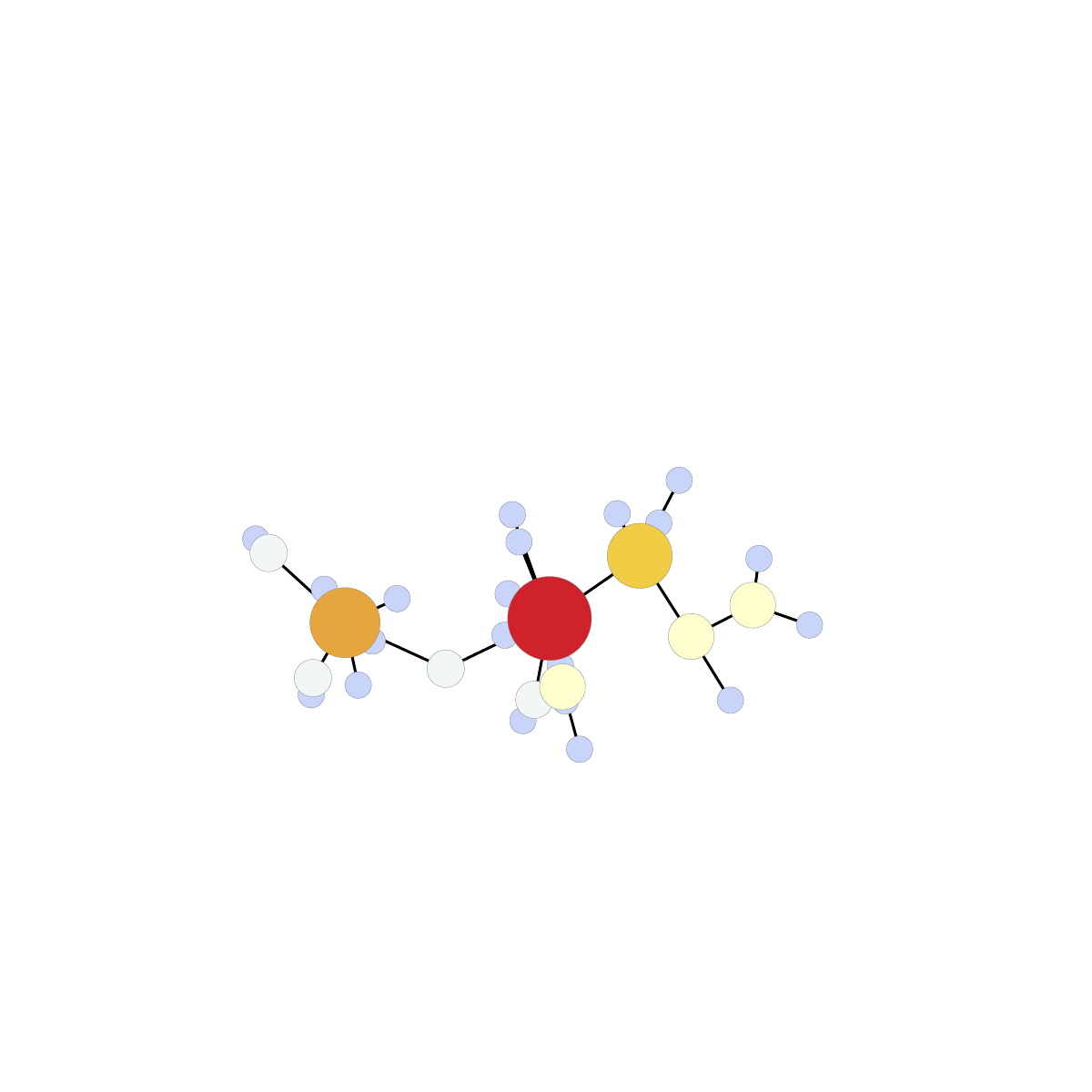}
      \includegraphics[height=3cm]{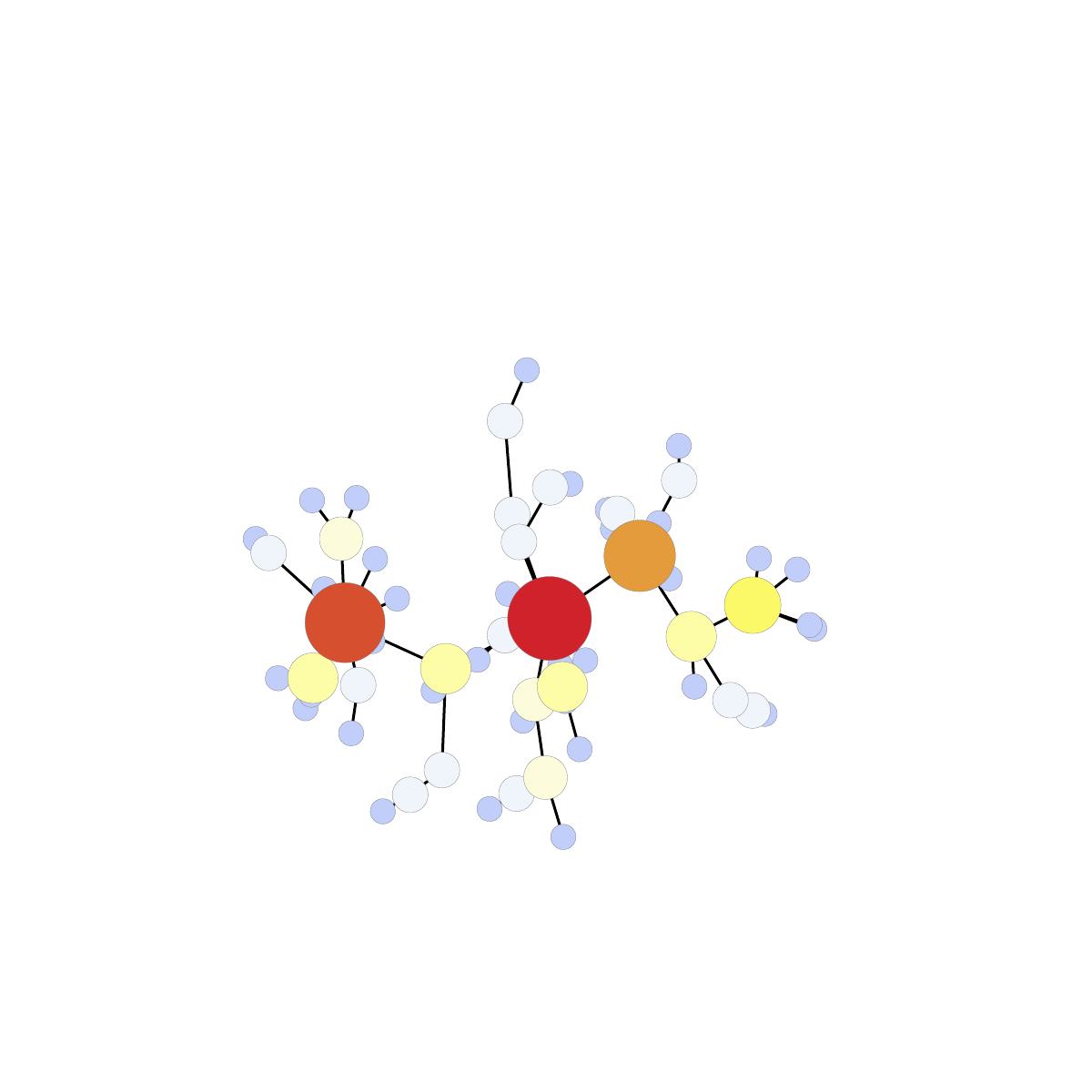}
       \includegraphics[height=3cm]{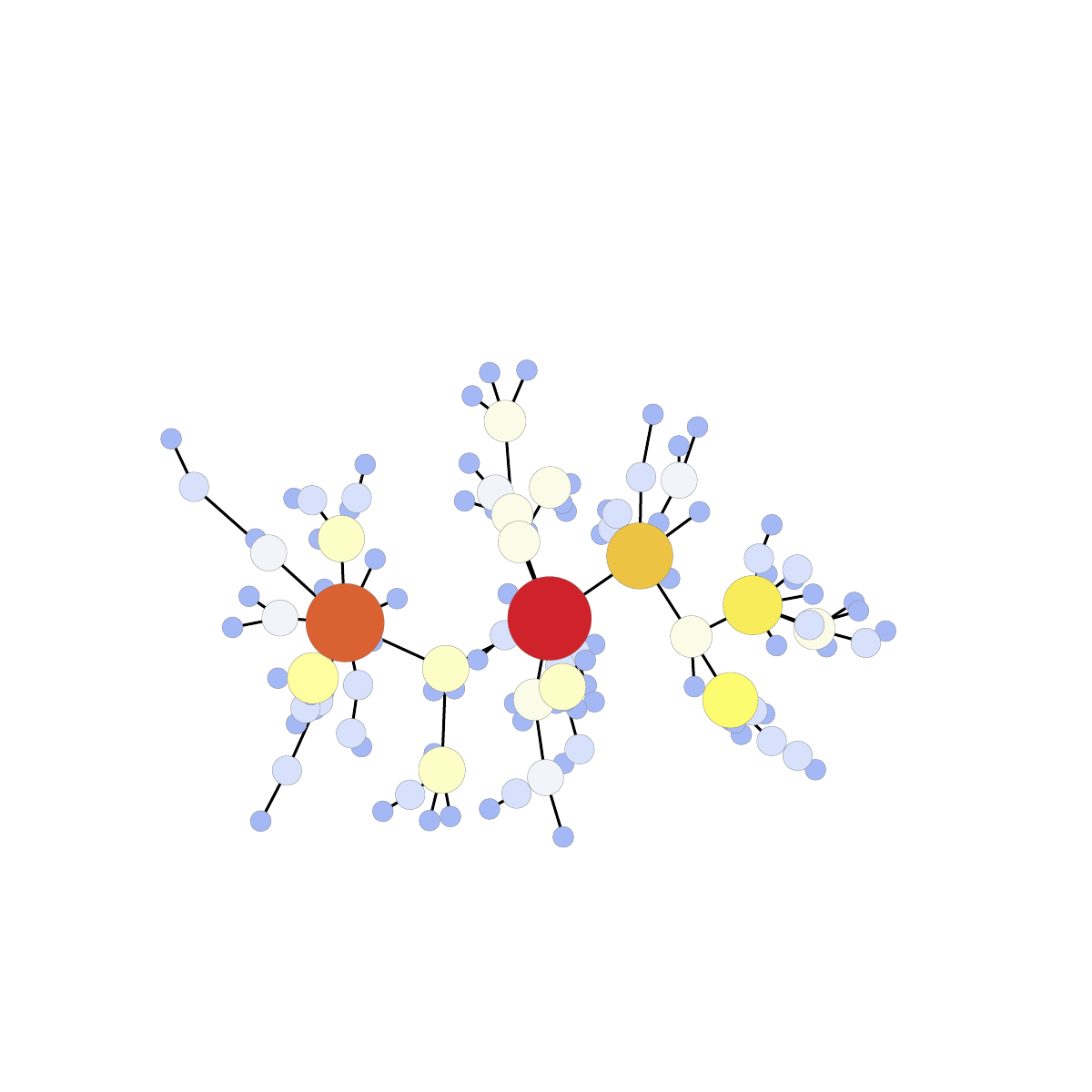}
        \includegraphics[height=3cm]{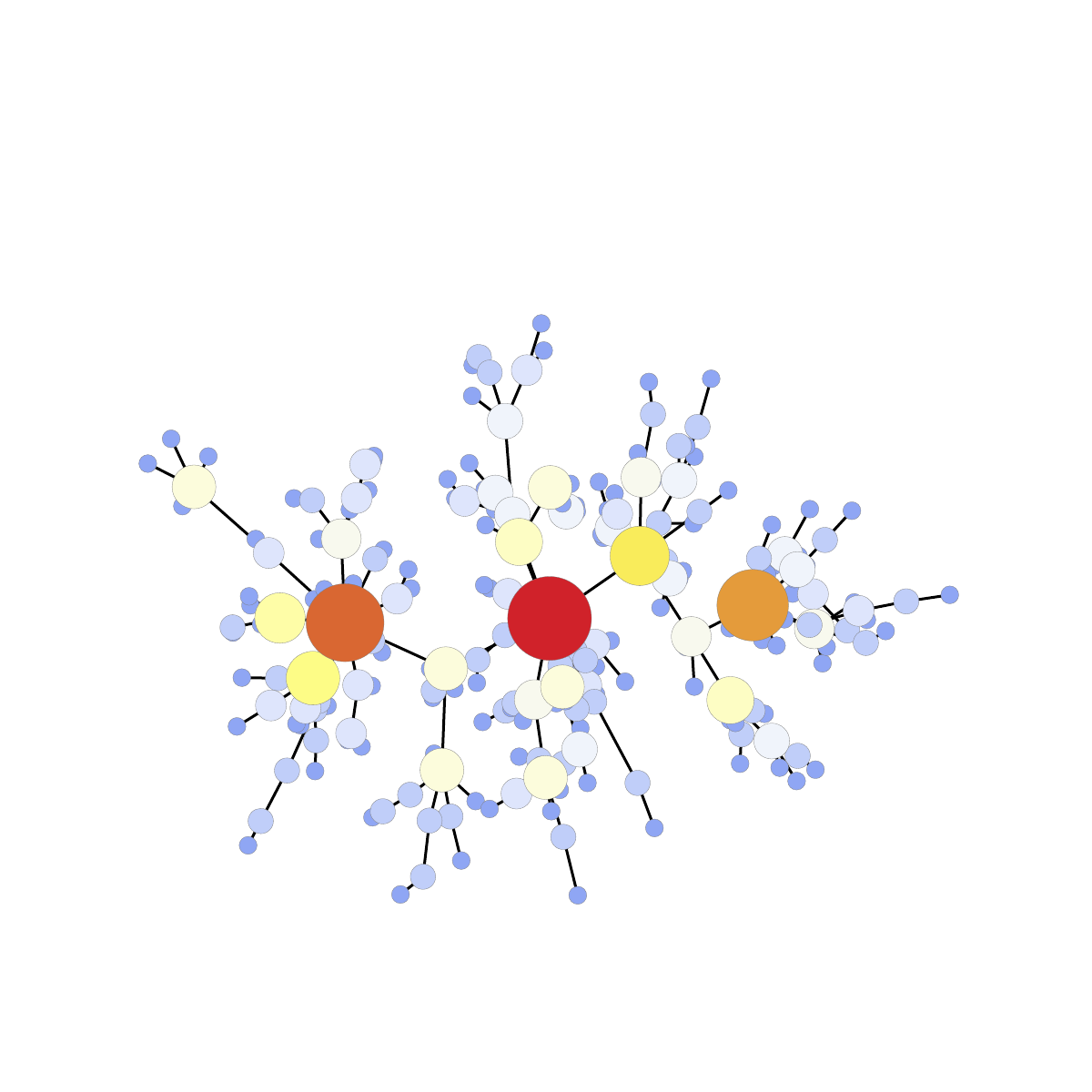}
                \includegraphics[height=3cm]{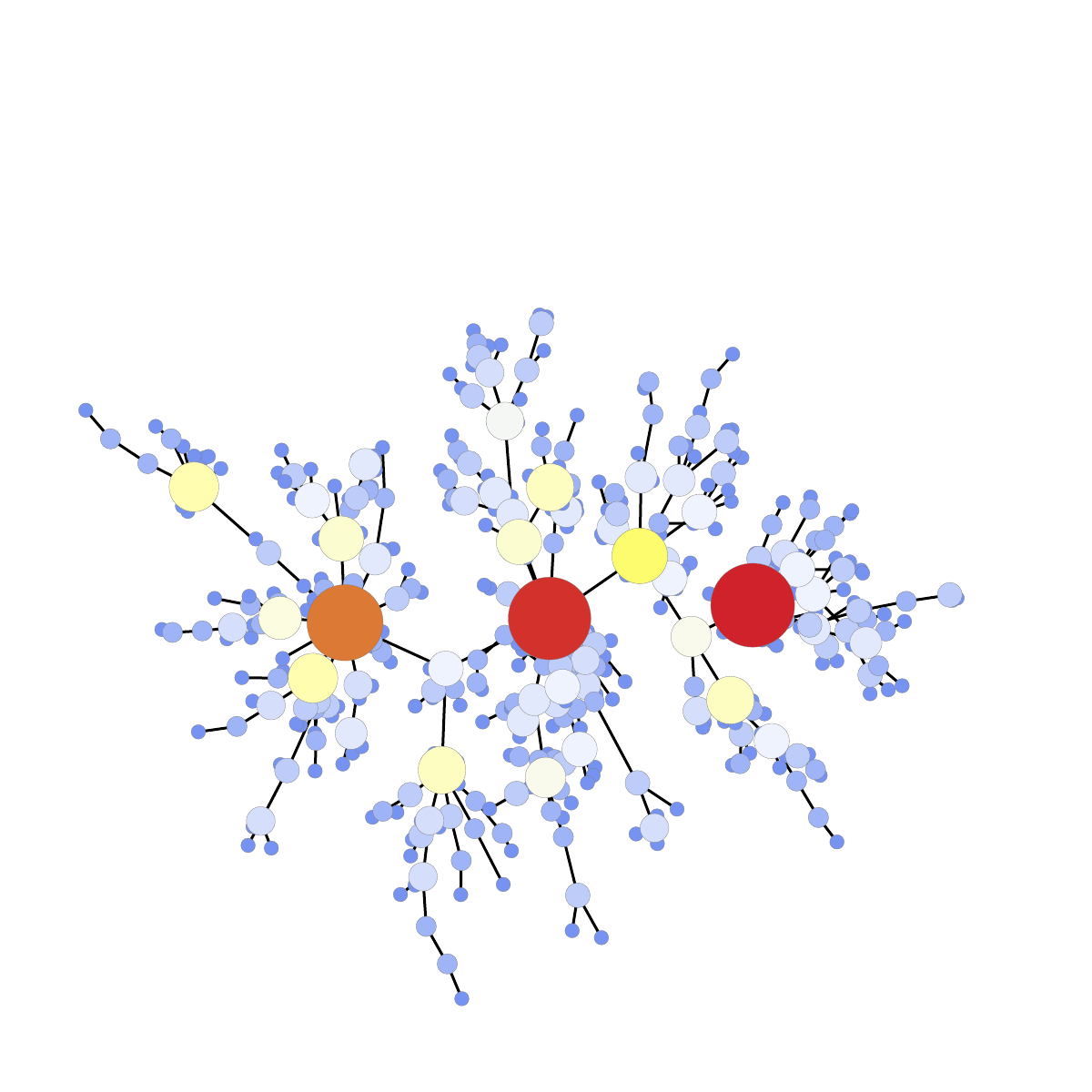}
                        \includegraphics[height=3cm]{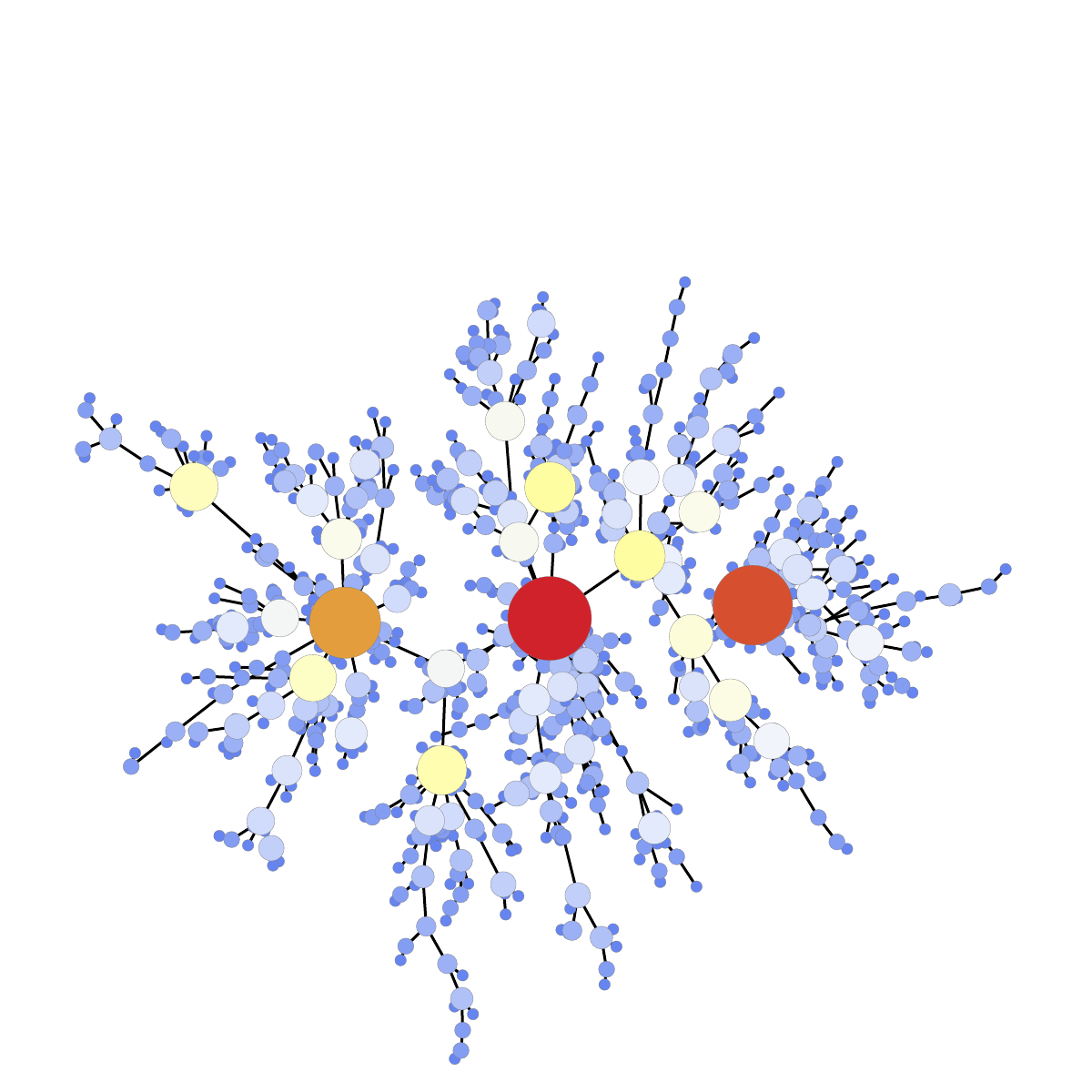}
                                \includegraphics[height=3cm]{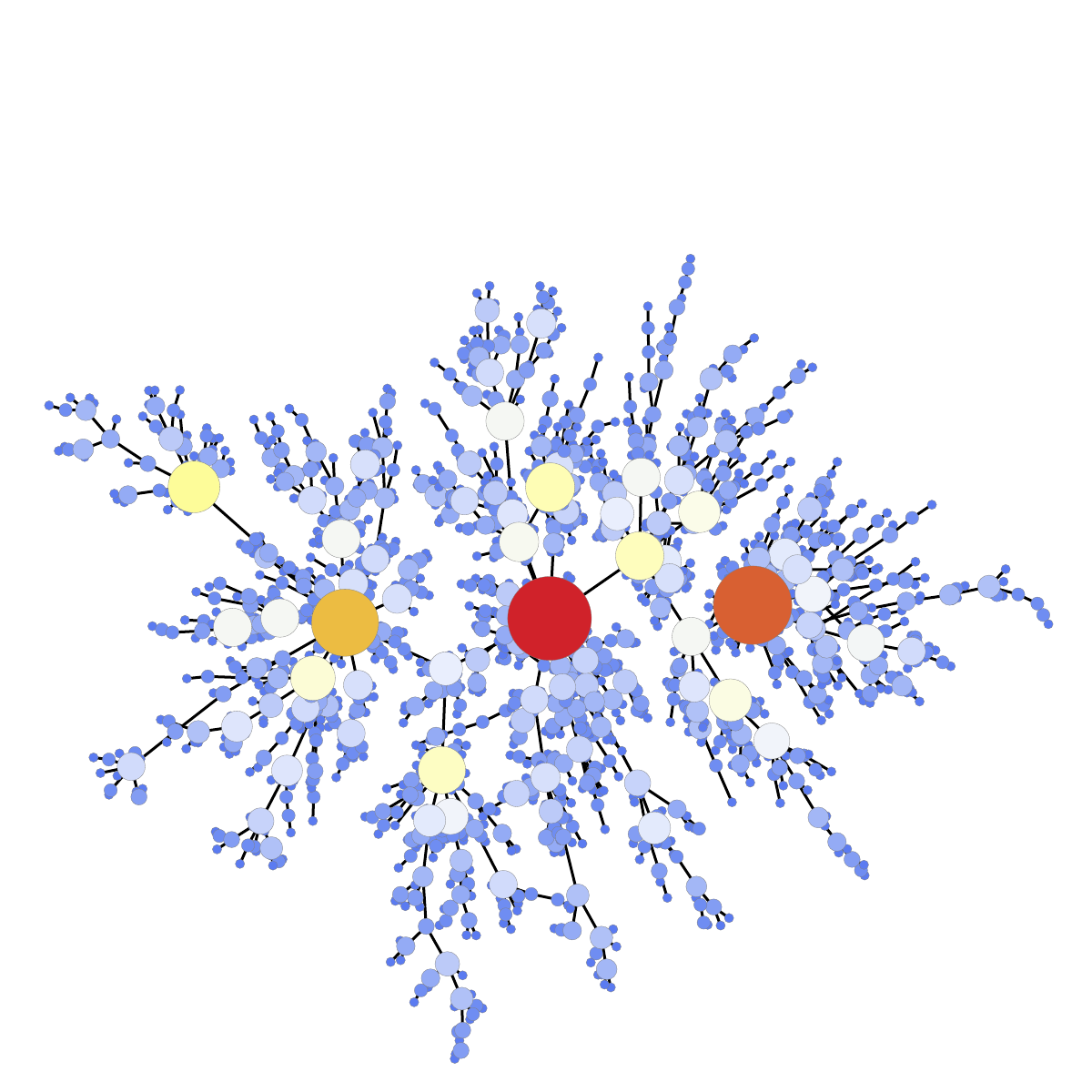}
                                        \includegraphics[height=3cm]{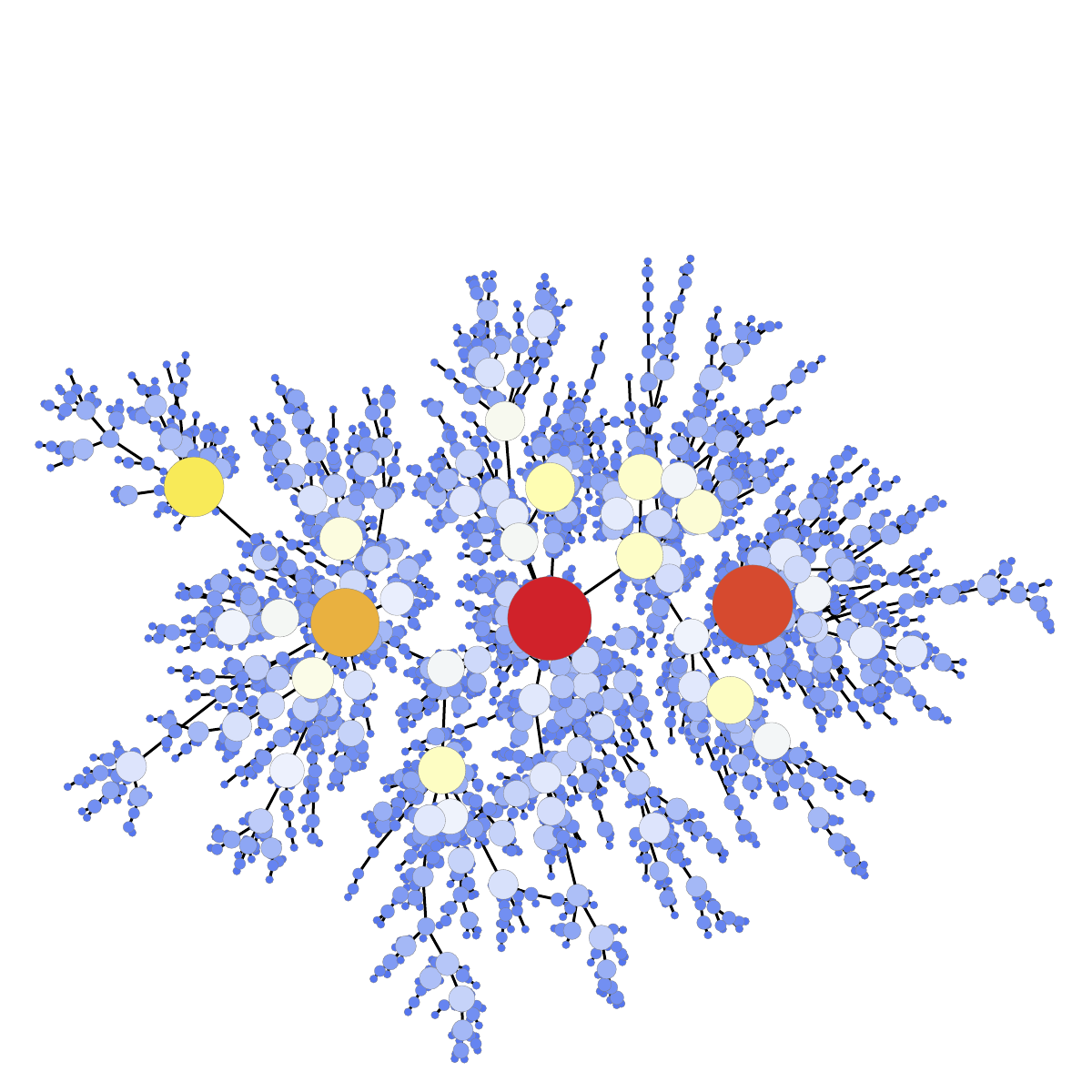}
                                                                                \includegraphics[height=3cm]{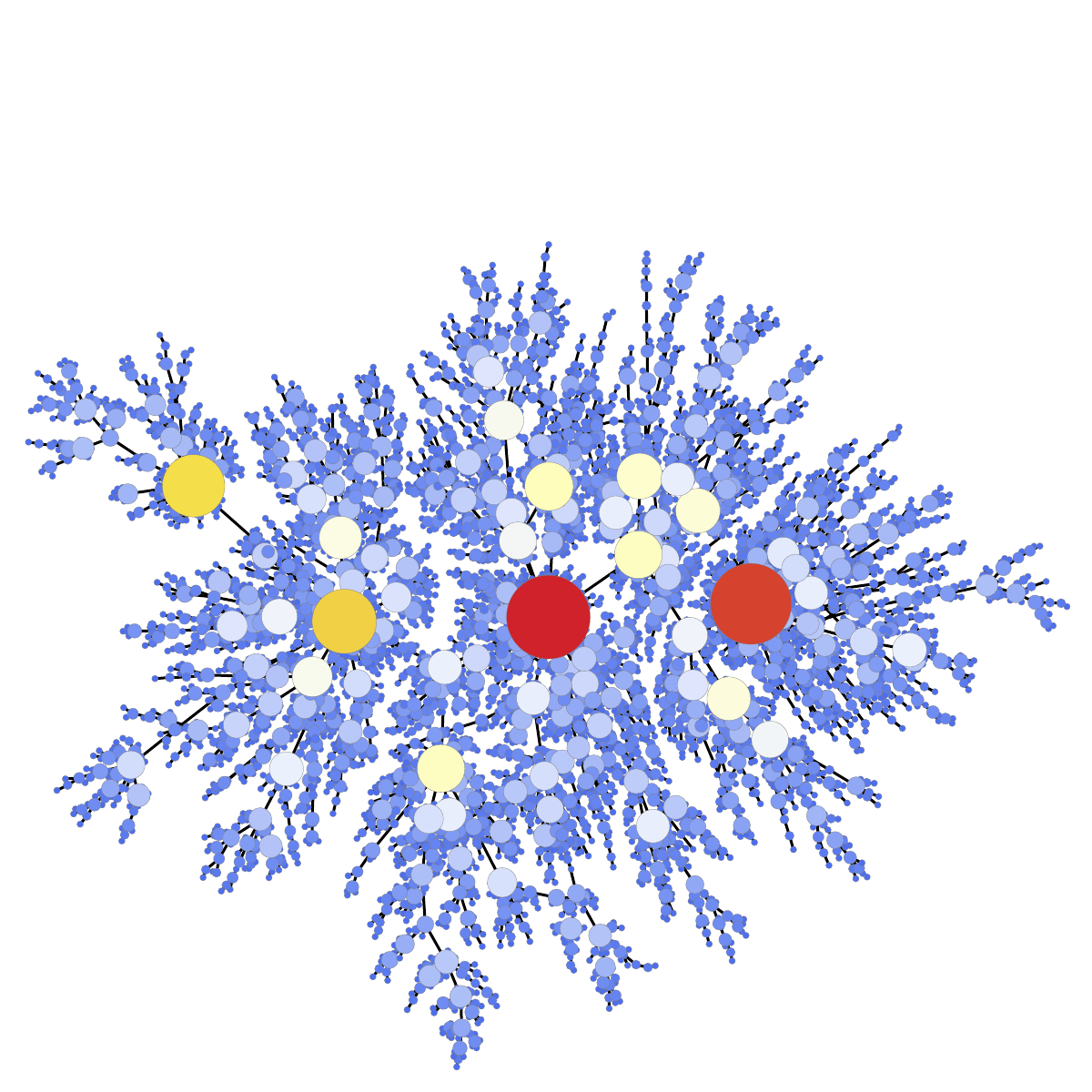}
                                                                                                                        \includegraphics[height=3cm]{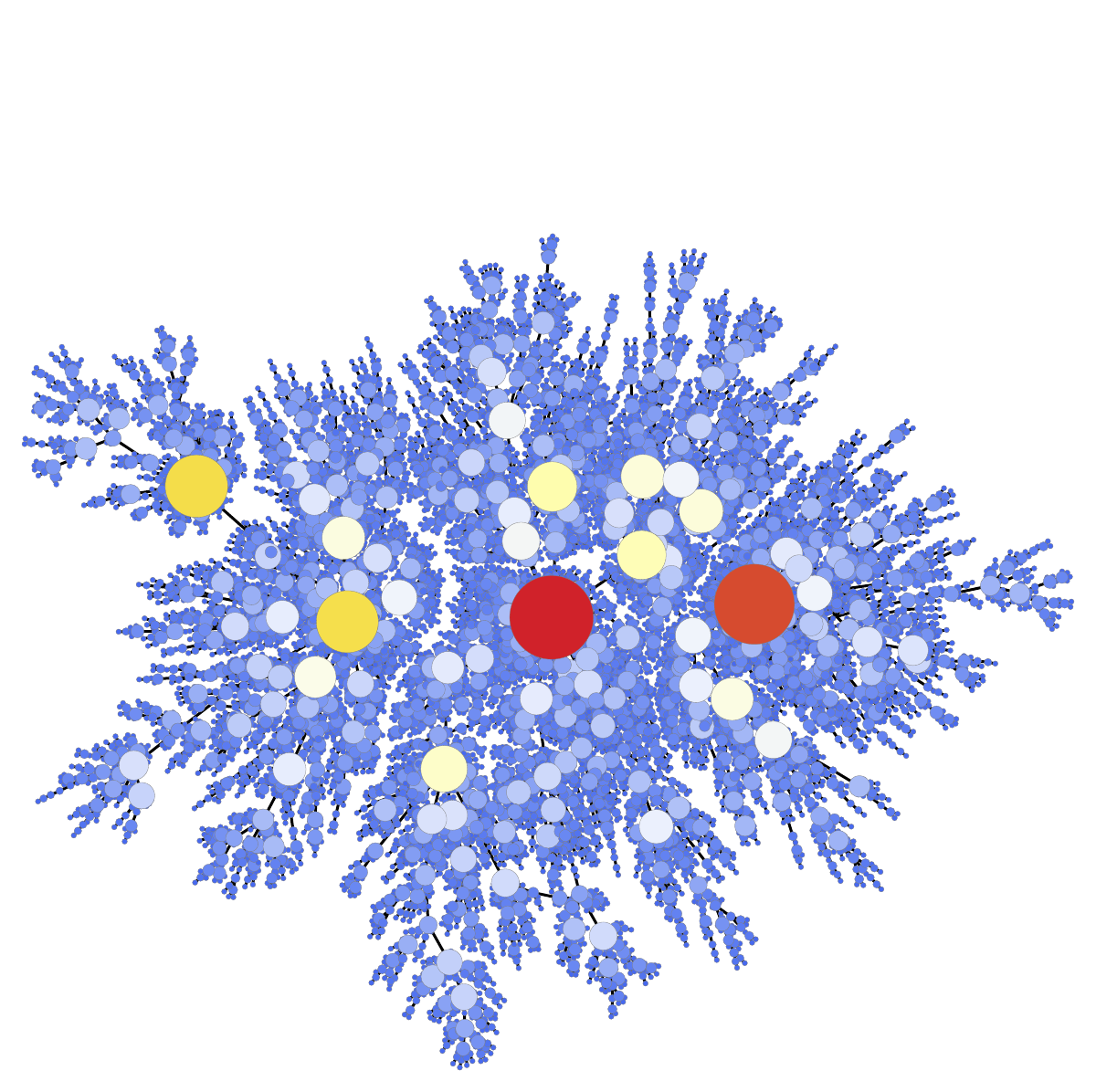}

 \caption{A sampling of the process $ \mathsf{T}_n$ for $n=1,2,3,4,8,16,32, \dots , 2^{14}$. The colors and sizes of the vertices indicate their degrees.\label{fig:BA}}
 \end{center}
 \end{figure}
 
 Since $ \mathsf{T}_{n}$ has $n$ edges, the sum of its vertex degrees is equal to $2n$, so that the normalization in the above definition indeed produces probability transitions. Compared to the random recursive tree, the preferential attachment model has a reinforcement of large degrees ``the rich get richer'' paradigm. This mechanism has been popularized by  Barab\'asi \& Albert \footnote{\raisebox{-3mm}{\includegraphics[height=1.5cm]{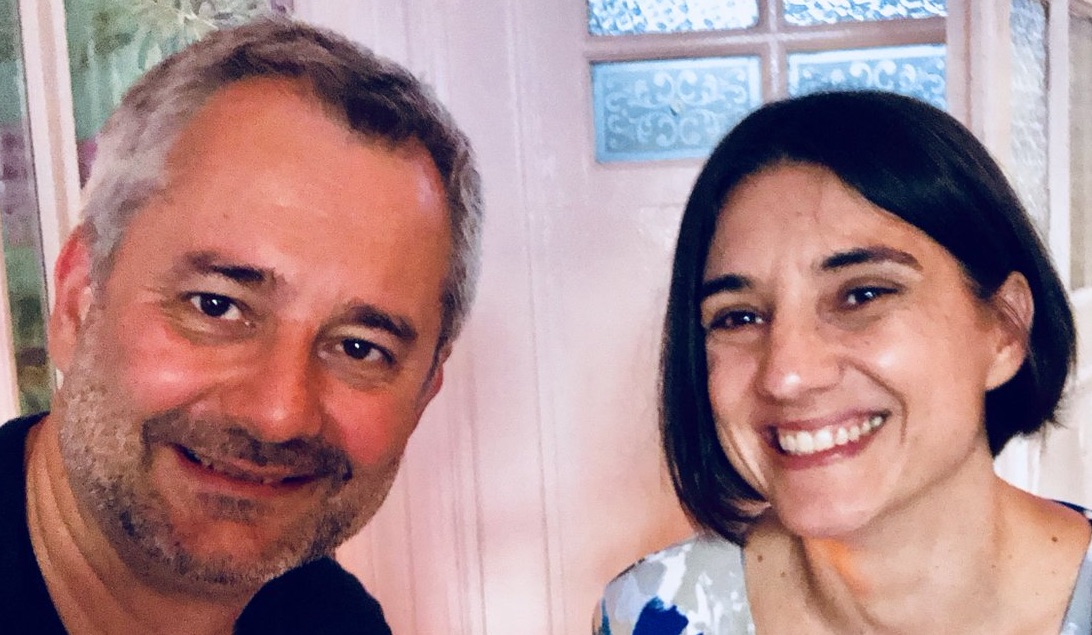}} Albert-László Barabási (1967--), and Réka Albert (1972--), Romanian} as a tractable model for real-world networks. It is possible to analyze this random tree growth using combinatorics as we did in Chapter \ref{chap:RRT} but we shall rather use the convenient tools developed in the previous two chapters.

 \section{Equivalent constructions}
As in Section \ref{sec:rrtfromyule} we shall see that the Barab\'asi--Albert tree process $( \mathsf{T}_n : n \geq 1)$ can  be constructed from a Yule process. But before that, let us reinterpret it as a random \textbf{plane} recursive trees.

\subsection{Plane recursive tree} 
Let us consider a plane variant of the random recursive tree construction in which we consider a Markov chain $( {T}^{ \mathrm{plan}}_{n} : n \geq 1)$ of labeled \textit{plane} trees where $T^{\mathrm{plan}}_{1} = \noeud{0}-\noeud{1}$ and where for $n \geq 2$, conditionally on $T^{\mathrm{plan}}_{n-1}$ the tree $T^{\mathrm{plan}}_{n}$ is obtained by grafting $ -\noeud{n}$ in one of the $2(n-1)$ corners (an angular sector made by two consecutive edges around a vertex) of $T^{\mathrm{plan}}_{n-1}$ uniformly at random.

\begin{figure}[!h]
 \begin{center}
 \includegraphics[width=14cm]{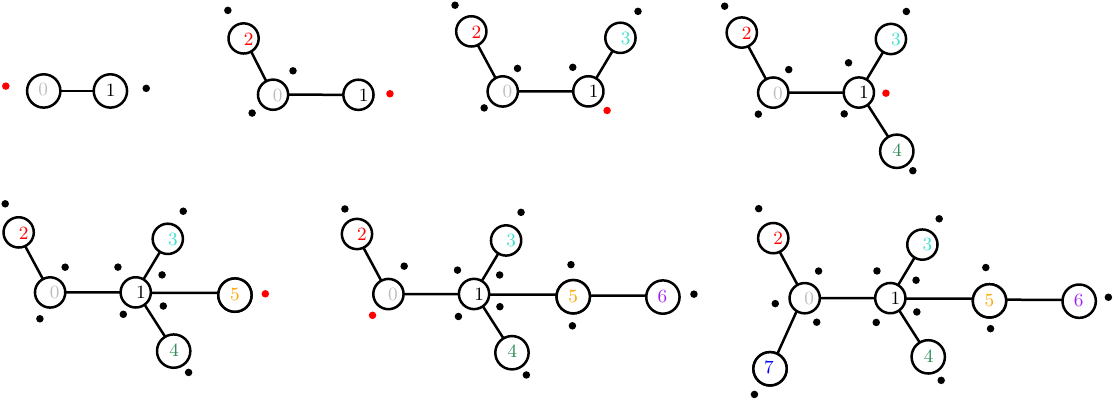}
 \caption{Illustration of the construction of $( {T}^{ \mathrm{plan}}_{n} : n \geq 1)$. The corners are represented by dots and the one selected for the grafting at the next step is in red.}
 \end{center}
 \end{figure}
 
 The tree $ {T}^{ \mathrm{plan}}_{n} $ is thus a plane tree (the root edge being the oriented edge $\noeud{0} \to \noeud{1}$) whose $n+1$ vertices are labeled by $\{0,1,2, \dots , n \}$ and such that the labels are increasing along branches starting from $\noeud{0}$. There are exactly $2^{n-1} (n-1)!$ such discrete tree structures and $ {T}^{ \mathrm{plan}}_{n}$ is, for each $n$, uniformly distributed over them. The following should then be clear:
 \begin{proposition}[Random plane recursive tree] The sequence of unlabeled non-plane trees obtained from $({T}^{ \mathrm{plan}}_{n} : n \geq 1)$ by forgetting the plane ordering is distributed as $ ( \mathsf{T}_{n} : n \geq 1)$.\end{proposition}
 
 It is also possible to obtain (a small variant of the) Barab\'asi--Albert tree process $( \mathsf{T}_n : n \geq 1)$ by modifying the uniform attachment rule:
 \begin{exo}[From  RRT to BA]  \label{exo:pawel} Consider the following attachment mechanism for labeled increasing trees starting with $ \mathfrak{T}_1 = \noeud{0}-\noeud{1}$: for $n \geq 2$ pick a uniform node $ \noeud{i}$ of $ \mathfrak{T}_{n-1}$ and attach $\noeud{n}$ with probability $1/2$ to $\noeud{i}$ or with probability $1/2$ to the first ancestor of $\noeud{i}$ (when going back towards $\noeud{0}$). If $\noeud{i} = \noeud{0}$, just attach $\noeud{n}$ to $\noeud{0}$. Show that the chain $ (\mathfrak{T}_n : n \geq 1)$ is very close to $( \mathsf{T}_n : n \geq 1)$.
 \end{exo}
 
\subsection{Construction via Yule tree of order $3$}\label{sec:yuletoBA}

Consider two independent plane Yule trees $ \mathbb{F} =  \mathbb{T}_0 \cup \mathbb{T}_1$ of order $3$ with rates equal to $1$, that is, in Section \ref{sec:AK} take $p=1$, $\alpha_1 =1$ and $\mu_1 = \delta_{3 \delta_1}$ and work under $ \mathbb{P}_{2 \delta_1}$. To ease notation in the rest of this section, under $ \mathbb{P}$ the forest $ \mathbb{F}$ has law $ \mathbb{P}_{2 \delta_1}$ whereas $ \mathbb{T}, \mathbb{T}_0, \mathbb{T}_1$ have law $ \mathbb{P}_{\delta_1}$. As in the previous chapter, we shall suppose that those trees are obtained by labeling the vertices of the full ternary tree $ \bigcup_{n \geq 0} \{0,1,2\}^n$ with i.i.d.~random exponential variables with mean $1$. For $t \geq 0$, we shall perform a contraction operation on $[ \mathbb{F}]_t = [\mathbb{T}_0]_t \cup  [\mathbb{T}_1]_t$ similar to that introduced in Section \ref{sec:rrtfromyule}: at each branch point of $[ \mathbb{F}]_t$,  we shall separate the right-most particle created from its two brothers. This creates a partitioning of $[ \mathbb{F}]_t$  into smaller ``Yule trees of order $2$''. Contracting each of these smaller subtrees into a single node and labeling them in their time-order of apparition\footnote{with the convention that the subtree  associated to the root of the first tree of $  \mathbb{F}$ corresponds to $\raisebox{.5pt}{\textcircled{\raisebox{-.9pt} {$0$}}}$} yields to an increasing (non-plane) labeled tree which we denote by $\llbrack  \mathbb{F}\rrbrack_t$.


\begin{figure}[!h]
 \begin{center}
 \includegraphics[width=16cm]{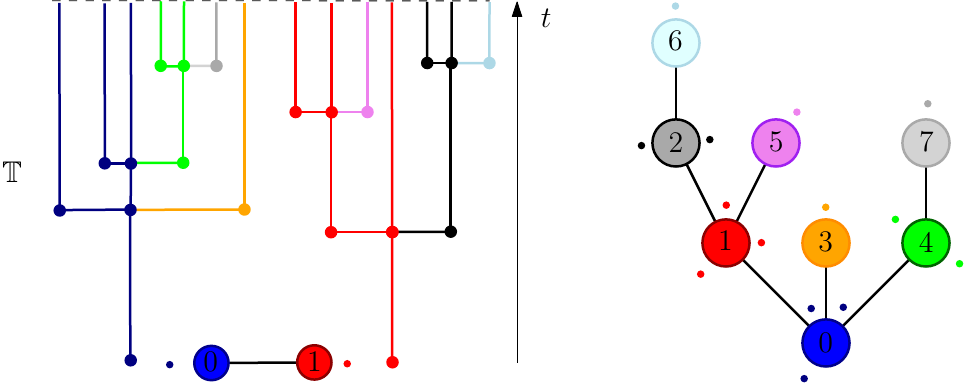}
 \caption{Constructing the increasing labeled tree $\llbrack \mathbb{F} \rrbrack_t$ by contracting all ``sub Yule trees of order 2" obtained by forgetting the right-most particle at each branch point.}
 \end{center}
 \end{figure}
We then have the analog of Proposition \ref{prop:RRTyule} which is proved using the same techniques:

\begin{proposition}[From Yule to BA]  \label{prop:BAyule} If $0= \tau_1 < \tau_{2}< \dots < \tau_{n}< \dots$ are the first times at which $\# \partial [ \mathbb{F}]_{\tau_n} = 2n$ then conditionally on $( \tau_{n} : n \geq 1)$ the process $( \llbrack \mathbb{F} \rrbrack_ {\tau_{n}} : n \geq 1)$ is a Barab\'asi--Albert preferential attachment tree. \end{proposition}

As in the preceding chapter, we will use the above construction together with our knowledge on Yule process  to deduce interesting geometric properties of the Barab\'asi-Albert tree, in particular on its maximal degree and its height.

\section{Degrees}
We denote by $\mathrm{deg}_{\llbrack \mathbb{F} \rrbrack_{t}}( \noeud{i})$ the degree of the $i$th vertex in the contraction of $[ \mathbb{F}]_t$ so that by Proposition \ref{prop:BAyule} we have the equality in terms of processes $$\left(\mathrm{deg}_{\llbrack \mathbb{F} \rrbrack_{\tau_n}}( \noeud{i}),  \mbox{ for }  i \leq n\right)_{n \geq 1}  = \left(\mathrm{deg}_{ \mathsf{T}_n}( \noeud{i}),  \mbox{ for }  i \leq n\right)_{n \geq 1}.$$

\subsection{Almost sure convergence}
Let us focus first on the degree of the root vertex $\noeud{0}$ inside $ \llbrack \mathbb{F} \rrbrack_{t}$. On the one hand, the variable $ \mathrm{deg}_{\llbrack \mathbb{F} \rrbrack_{t}}( \noeud{0})$ is equal to the number $ \mathcal{Y}^{{(2)}}_{t}$ of particles alive at time $t$ in the ``sub Yule process''  of order $2$ obtained by keeping only the first two children at each branching point.  On the other hand, the total number of particles $\mathcal{Y}^{{(3)}}_{t} = \#
 \partial [ \mathbb{F}]_{t}$ alive at time $t$ in the forest  is the sum of two independent Yule processes of order $3$. We deduce from Proposition \ref{prop:yule} the following almost sure convergences 
 $$ \mathrm{e}^{-t} \mathcal{Y}^{{(2)}}_{t} \xrightarrow[n\to\infty]{a.s.} \mathcal{E} \quad  \mbox{ and } \quad \mathrm{e}^{-2t} \mathcal{Y}^{{(3)}}_{t} \xrightarrow[n\to\infty]{a.s.}  \frac{1}{2} \cdot \mathcal{E}',$$ where $ \mathcal{E}$ and $ \mathcal{E}'$ are two exponential variables of mean $1$ which are \textbf{not independent}. In particular, it follows from the last display together with  Proposition \ref{prop:BAyule}  and Proposition	\ref{prop:yule} that $  n^{-1/2} \cdot \mathrm{deg}_{\llbrack \mathbb{F} \rrbrack_{\tau_n}}( \noeud{0})$ converges almost surely towards $ \mathcal{E}/ \sqrt{ \mathcal{E}'/4}$ and more generally that: 
 \begin{proposition}[Almost sure convergence of degrees]  \label{prop:asdegree}There exists a vector of almost surely positive and finite random variables $( \mathcal{X}_{i} : i \geq 0)$ so that for each $i \geq 0$ we have the following almost sure convergences 
 $$ \frac{\mathrm{deg}_{ \mathsf{T}_{n}}( \noeud{i})}{n^{1/2}} \xrightarrow[n\to\infty]{a.s.} \mathcal{X}_{i}.$$
 Moreover the $( \mathcal{X}_i : i \geq 0)$ are almost surely distinct.
 \end{proposition}
\noindent \textbf{Proof.} Recall that  $\tau_i$ is the first time at which the particle $\noeud{i}$ appears in $ \llbrack \mathbb{F} \rrbrack_t$. By the Markov property of the Yule process, for $t \geq \tau_i$ the degree $ \mathrm{deg}_{\llbrack \mathbb{F} \rrbrack_t}( \noeud{i})$ can be expressed as a counting process in a Yule tree of order $2$, whereas the total number of corners is given by the sum of two independent Yule process of order $3$ (the number of vertices is half of it). Using  Proposition \ref{prop:yule} three times, we deduce the almost sure convergence towards positive r.v. $ \mathcal{X}_i$. Let us now explain why  $ \mathcal{X}_0 \ne \mathcal{X}_1$ with probability one, leaving the general case $ \mathcal{X}_i \ne \mathcal{X}_j$ to the reader. For $ a \in \{0,1\}$,  denote by $ \mathcal{D}_a$ (resp. $ \mathcal{M}_a$) the limit of the renormalized size of the Yule tree of order $2$ (resp. of order $3$) obtained by keeping the first two children in each branching  (resp. keeping all children) in the tree $ \mathbb{T}_a$. By the above discussion, we have 
$$ \mathcal{X}_0 = \frac{ \mathcal{D}_0}{ \sqrt{ (\mathcal{M}_0 + \mathcal{M}_1)/2}} \quad \mbox{ and }\quad \mathcal{X}_1 = \frac{ \mathcal{D}_1}{ \sqrt{ (\mathcal{M}_0 + \mathcal{M}_1)/2}}.$$ Remark now that $ ( \mathcal{D}_0, \mathcal{M}_0)$ and $( \mathcal{D}_1, \mathcal{M}_1)$ are independent and $ \mathcal{D}_a$ have no atoms (they are exponentially distributed). Hence the probability that $ \mathcal{D}_0 = \mathcal{D}_1$ is $0$ implying that $ \mathcal{X}_0 \ne \mathcal{X}_1$ a.s. \qed \medskip 

 \begin{exo}[A martingale approach]  \label{exo:martingale} Here is a way to prove the  almost sure  convergence of renormalized degrees without the continuous-time embedding. Let $D_{n} = \mathrm{deg}_{ \mathsf{T}_{n}}( \noeud{0})$ for $n \geq 1$. Show that we have 
 $$ \mathbb{E}[D_{n+1} \mid \sigma( \mathsf{T}_{k} : 1 \leq k \leq n)] = D_{n} \cdot \left( 1 + \frac{1}{2n}\right).$$
 Conclude that $D_{n} \cdot \prod_{k=1}^{n-1} (1+ \frac{1}{2k})^{-1}$ is positive martingale which converges almost surely and recover the first part of the previous proposition.
 \end{exo}
 
   \begin{figure}[!h]
  \begin{center}
  \includegraphics[width=13cm]{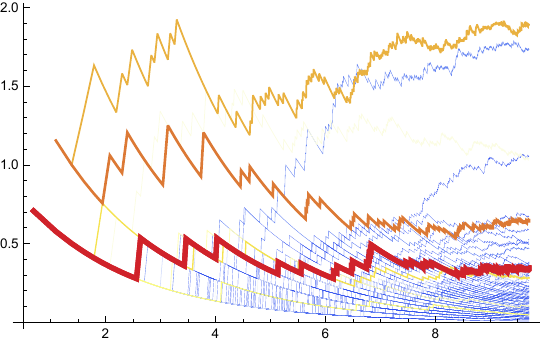}
  \caption{\textsc{Degree's race in the evolution of Figure \ref{fig:BA}.} Curves of the renormalized vertex degrees in the scale  $ \sqrt{n}$ (in $y$ axis)  in a logarithmic scale for $n$  (in the $x$ axis). The renormalized degree of  \noeud{0} is in red, that of  $\noeud{1}$ in orange and that of  $\noeud{2}$ (this vertex asymptotically has the largest degree) in dark yellow. \label{fig:degrees}}
  \end{center}
  \end{figure}
  
  \subsection{Maximal degree}
  As for the case of the random recursive tree, one can wonder about the maximal degree in the Barab\'asi--Albert tree process. In the RRT, the largest degree after $n$ steps turned out not to be among the first nodes of the network but among the nodes arrived at time $ \approx n^{0,27\dots}$, see Remark \ref{rek:ouestgros}. Here, the fast decay of the degrees enables us to show that the largest degree actually belongs to the first few nodes of the network. More precisely we have: \clearpage
  
  \begin{theorem}[Mori] \noindent  With the notation of Proposition \ref{prop:asdegree}, the random vector $ ( \mathcal{X}_i : i \geq 0)$ almost surely satisfies $ \mathcal{X}_i \to 0$ as $i \to \infty$ and the pointwise a.s.\ convergence can be reinforced into an almost sure convergence for the $\ell^\infty$ metric: \label{thm:mori}
   $$ \left( \frac{\mathrm{deg}_{ \mathsf{T}_{n}}( \noeud{i})}{n^{1/2}} : i \geq 1 \right) \xrightarrow[n\to\infty]{a.s. \mathrm{ \ for \ } \ell^\infty} (\mathcal{X}_{i} : i \geq 0).$$
   \end{theorem}
 Combining the previous result with the fact (proved in Proposition  \ref{prop:asdegree}) that the $  \mathcal{X}_i$ are positive and almost surely distinct, we deduce that the relative position $ \mathrm{Pos}( \noeud{i}, n)$  of the degree of node $ \noeud{i}$ among $\{ \noeud{0}, \dots , \noeud{n}\}$ converges almost surely as $n \to \infty$ towards $ \mathfrak{S}_i$ where $ \mathfrak{S} : \{0,1,2, \dots \} \to \{1,2,3, \dots \}$ is a bijection. This implies in particular the convergence of the index of the largest vertex's degree in $ \mathsf{T}_n$. 
   \medskip 
   
The main technical input for the proof of Theorem \ref{thm:mori} is  a maximal inequality based on Proposition \ref{prop:yule}:
   \begin{lemma} \label{lem:techdeg} Let $ (\mathcal{Y}^{(2)}_t : t \geq 0)$ be the counting process of a standard Yule tree of order $2$, rate $1$, and starting from $1$ particle. For all $x \geq 2$ we have
$$ \mathbb{P}\left( \sup_{t \geq 0}  \mathrm{e}^{-t} \mathcal{Y}^{(2)}_t  \geq x\right) \leq  2 \mathrm{exp}(- x/2).$$
\end{lemma}
\noindent \textbf{Proof.} Fix $x \geq 2$ and denote by $\theta = \inf\{ t \geq 0 : \mathrm{e}^{-t} \mathcal{Y}^{(2)}_t  \geq x\}$. On the event where the stopping time $\theta$ is finite, the strong Markov property entails that conditionally on $[ \mathbb{T}]_\theta$, the $\mathcal{Y}^{(2)}_\theta = \# \partial[ \mathbb{T}]_\theta=:N$ particles alive at time $\theta$ will have independent offsprings distributed according to a standard Yule tree of order $2$. Recalling from  Proposition \ref{prop:yule} that $\mathrm{e}^{-t} \mathcal{Y}^{(2)}_t \to \mathcal{E}$ a.s., on the event $\{ \theta < \infty\}$ we can write
$$  \mathcal{E} = \lim_{t \to \infty} \mathrm{e}^{-t} \mathcal{Y}^{(2)}_t = \mathrm{e}^{-\theta} \sum_{i=1}^N \lim_{t \to \infty }  \mathrm{e}^{-t} \mathcal{Y}^{(2),i}_t = \mathrm{e}^{-\theta} \sum_{i=1}^N \mathcal{E}_i,$$  where on the right-hand side, the variables $ (\mathcal{E}_i : i \geq 1)$ are i.i.d.~exponential variables of rate $1$ independent of $N$. Using the easy fact that  $\inf_{k \geq 1} \mathbb{P}( \sum_{i=1}^k \mathcal{E}_i \geq k/2) \geq \frac{1}{2}$ we have 
  \begin{eqnarray*}  \mathrm{e}^{-x/2} = \mathbb{P}( \mathcal{E} \geq x/2) &\geq& \mathbb{E}\left[ \mathbf{1}_{\theta< \infty} \mathbb{P}\left( \mathrm{e}^{-\theta} \sum_{i=1}^N \mathcal{E}_i \geq  \frac{x}{2} \right) \right]\\
  & \underset{N \mathrm{e}^{-\theta} \geq x}{\geq} & \mathbb{E}\left[ \mathbf{1}_{\theta< \infty} \mathbb{P}\left( \sum_{i=1}^N \mathcal{E}_i \geq  \frac{N}{2}\right) \right]\\
  & \geq &  \inf_{k \geq 1} \mathbb{P}\left( \sum_{i=1}^k \mathcal{E}_i \geq k/2\right) \cdot \mathbb{P}(\theta < \infty) \geq \frac{1}{2} \mathbb{P}( \theta < \infty).  \end{eqnarray*}\qed \medskip

  \noindent \textbf{Proof of Theorem \ref{thm:mori}.} Given the work done in the proof of Proposition \ref{prop:asdegree}, the convergence for the $\ell^\infty$ metric follows if we can show that 
  $$ \lim_{m \to \infty }\sup_{i \geq m} \sup_{n \geq 1} \frac{ \mathrm{deg}_{ \mathsf{T}_n}(\noeud{i})}{ \sqrt{n}} =0, \quad a.s.$$ or via the continuous time representation that 
  \begin{eqnarray} \lim_{T \to \infty }\sup_{{\footnotesize \noeud{i}} \mathrm{ \ created \  after\  }T} \sup_{t \geq T} \frac{ \mathrm{deg}_{ \llbrack  \mathbb{F} \rrbrack_t}(\noeud{i})}{  \mathrm{e}^t} =0, \quad a.s.   \label{eq:goalt0}\end{eqnarray}
For $t > A$, a new splitting appears in $ \mathbb{F}$ with intensity $ \# \partial [ \mathbb{F}]_t \ \mathrm{d}t$, this creates a new vertex in  $ \llbrack \mathbb{F} \rrbrack_t$ and the probability that such a vertex gets a degree larger than $ \varepsilon \mathrm{e}^u$ at some later time $u = s+t \geq t$ is upper bounded by
$$ \mathbb{P}(\sup_{s \geq 0}  \mathrm{e}^{-s} \mathcal{Y}^{(2)}_s \geq \varepsilon \mathrm{e}^{t}) \underset{ \mathrm{Lem.\ } \ref{lem:techdeg} }{\leq} 2 \exp(- \varepsilon \mathrm{e}^{t} /2).$$
We deduce that 
 \begin{eqnarray*} &&   \mathbb{E}\left[ \sum_{{\footnotesize \noeud{i}} \mathrm{ \ created \  after\  }T} \mathbf{1}\left\{\sup_{t \geq A} \frac{ \mathrm{deg}_{ \llbrack  \mathbb{F} \rrbrack_t}(\noeud{i})}{  \mathrm{e}^t} \geq \varepsilon \right\}\right] \\ &\leq& \int_{A}^\infty \mathrm{d}t\  \mathbb{E}[ \# \partial [ \mathbb{F}]_t ] \cdot  2 \exp(- \varepsilon \mathrm{e}^{t} /2) \\ &=&  \int_{A}^\infty \mathrm{d}t\ 2 \mathrm{e}^{2t} \cdot 2 \exp(- \varepsilon \mathrm{e}^{t} /2).  \end{eqnarray*}
For $ \varepsilon>0$ fixed,  the above integral can be made arbitrarily small provided that $A>0$ is chosen large enough. This implies \eqref{eq:goalt0}. \qed

\subsection{Empirical degree distribution}
As in Section \ref{sec:empiricalRRT} we can also study the \textbf{empirical degree distribution} in $ \mathsf{T}_{n}$: We let $ \nu_{n} $ be the (random) empirical distribution of the out-degrees defined by 
$$ \nu_{n} =  \frac{1}{n+1}\sum_{i=0}^{n} \delta_{ \mathrm{deg}_{ \mathsf{T}_{n}}^+( {\footnotesize \noeud{i}})}.$$
As for Proposition \ref{prop:empirical},  the empirical degree distribution converges towards a deterministic distribution which now has an interesting polynomial tail behavior:
   \begin{theorem}[Convergence of the empirical degree distribution] \label{prop:empiricalBA}The empirical distribution of the out-degrees in $ \mathsf{T}_{n}$ converges in probability towards an explicit deterministic law: for each $k \geq 0$ we have 
   $$ \nu_{n}(\{k\})  \xrightarrow[n\to\infty]{( \mathbb{P})}  \frac{4}{(k+1)(k+2)(k+3)}.$$
   \end{theorem}
\noindent  \textbf{Proof.} We obviously use the construction of $ \mathsf{T}_n = \llbrack \mathbb{F} \rrbrack_{\tau_n}$ valid for all $n \geq 1$ simultaneously. The same proof as for Proposition \ref{prop:concentrationlocal} shows that 

$$ \frac{\mathsf{D}_k([  \mathbb{F}]_t)}{ \# \partial [ \mathbb{F}]_t} \xrightarrow[t\to\infty]{ ( \mathbb{P})}  \lim_{t \to \infty} \mathrm{e}^{-2t} \mathbb{E}[ \mathsf{D}_k([  \mathbb{F}]_t)],$$ where  

$$ \mathsf{D}_k([  \mathbb{F}]_t) := \# \Big\{ u \in \llbrack \mathbb{F} \rrbrack_t \backslash \noeud{0} : \mathrm{deg}^{+}_{ \llbrack \mathbb{F}\rrbrack_t}(u) =k\Big\},$$ and where the limit exists. We compute the expectation of $ { \mathsf{D}}_k( [ \mathbb{T}]_t)$, the number of vertices different from $ \noeud{0}$ and of degree $k\geq 1$ in $\llbrack \mathbb{T} \rrbrack_t$, in a single contracted Yule tree of order $3$. As in Section \ref{sec:goncharovback}, recall that a new vertex is created at time $s$ with intensity $ \# \partial [ \mathbb{T}]_s$ and by \eqref{eq:yuleexplicit}, this vertex has degree $k$ at time $t$ with probability $ \mathrm{e}^{-(t-s)}(1- \mathrm{e}^{-(t-s)})^{k-1}$. Recalling that $ \mathbb{E}[\# \partial [ \mathbb{T}]_t]=  \mathrm{e}^{2t}$ we have 
 \begin{eqnarray*}  \mathrm{e}^{-2t} \cdot \mathbb{E}\left[{ \mathsf{D}}_k( [ \mathbb{T}]_t) \right] &=& \mathrm{e}^{-2t} \int_0^t \mathrm{d}s \, \mathbb{E}[\# \partial [ \mathbb{T}]_s] \cdot  \mathrm{e}^{-(t-s)}(1- \mathrm{e}^{-(t-s)})^{k-1}\\
 &\underset{u=t-s}{=}&  \int_0^t \mathrm{d}u \, \mathrm{e}^{-3u} (1- \mathrm{e}^{-u})^{k-1} \\
 & \xrightarrow[t\to\infty]{} &  \int_0^\infty \mathrm{d}u \, \mathrm{e}^{-3u} (1- \mathrm{e}^{-u})^{k-1}\\
 & \underset{x= \mathrm{e}^{-u}}{ = }&  \int_0^1 \mathrm{d}x\, x^2(1-x)^{k-1}= \frac{2}{k(k+1)(k+2)} \end{eqnarray*}
In the case of two trees, we also have for any $k \geq 1$ 
$$ \frac{ { \mathrm{D}}_k( [ \mathbb{F}]_t)}{ \# \partial [ \mathbb{F}]_t} \xrightarrow[t\to\infty]{( \mathbb{P})}\frac{2}{k(k+1)(k+2)},$$ which proves the result since the number of vertices in $ \llbrack \mathbb{F} \rrbrack_t$ is half of 	$\# \partial [ \mathbb{F}]_t$. \qed

\begin{remark}[Scale-free property] The fact that the empirical degree distribution $\nu_n$ converges towards a limiting law $\nu_\infty$ with a polynomial tail behavior $\nu_\infty(\{k\}) \approx k^{-\alpha}$ with $\alpha > 0$ is usually refer to as the \textbf{scale-free} property. In the case of the Barabi--Albert trees the tail with exponent $\alpha=3$ is coherent with the fact that the largest degree in $ \mathsf{T}_n$ is of order $ \sqrt{n} = n^{1/(\alpha-1)}$ which is the order of magnitude of the maximum of $n$ i.i.d.~samplings according to $\nu_\infty$.  \end{remark}

\section{Height}
We finish by studying the maximal height in $ \mathsf{T}_n$. The preferential attachment mechanism do yield to smaller trees compared to the uniform attachment case, but they stay of logarithmic order:
\begin{theorem}[Pittel] We have  \label{thm:heightBA}
$$ \frac{ \mathrm{Height}( \mathsf{T}_n)}{\log n} \xrightarrow[n\to\infty]{a.s.} c\approx 1.79\dots,$$ where $c=(2 \gamma)^{-1}$ for $\gamma$ the solution to $ \gamma \mathrm{e}^{1+\gamma}= 1$.
\end{theorem}

\noindent \textbf{Sketch of proof.} The proof follows the same strategy as that of Theorem \ref{thm:heightRRT} presented in Section \ref{thm:heightRRT}. Similar to \eqref{eq:distancespine}, a particle $u \in \partial [ \mathbb{F}]_{t} $ is associated with a vertex in $\llbrack \mathbb{F}\rrbrack_{t}$ whose distance to the root $\raisebox{.5pt}{\textcircled{\raisebox{-.9pt} {$0$}}}$ in $\llbrack \mathbb{F} \rrbrack_t$ is equal to the number of branch points along the spine for which the lineage to $u$ is the right-most. When $u = \bullet$ is the distinguished particle of $[ \mathbb{T}^\bullet]_t$ under $ \mathbb{P}_{\delta_2}$, branchings happens at rate $3$ and a third of them is of the above form. By the many to one formula (Theorem \ref{prop:manyto1}) we then have 
$$ \mathbb{E}\left[\sum_{u \in  \partial [ \mathbb{T}]_t}  \mathbf{1}\{\mathrm{dist}_{\llbrack \mathbb{T} \rrbrack_t}(u, \noeud{0}) \geq \alpha t \}\right] = \mathrm{e}^{2t}  \mathbb{P}( \mathfrak{P}(t) \geq \alpha t).$$
By Lemma \ref{lem:LDpoisson}, when $\alpha = \gamma^1 + \varepsilon$ the previous display converges to $0$ exponentially fast with $t \to \infty$ (notice that $a=\gamma^{-1}$ is solution to $a \log a - (a-1)=2$). Since there are roughly $ n \approx \mathrm{e}^{2t}$ vertices at time $t$ in $ \llbrack \mathbb{F} \rrbrack_t$ we deduce using the same arguments as in Section \ref{sec:heightRRT} that the height of $ \mathsf{T}_n$ is eventually less than $  \frac{\gamma^{-1}}{2} n$ as $n \to \infty$ a.s. The lower bound follows mutatis mutandis the same lines as in Section \ref{sec:heightRRT} and we leave it as an exercise for the reader. \qed

%

\paragraph{Bibliographical notes.} Although generally attributed to Albert \& Barab\'asi \cite{AB99} which is one of the most cited papers in mathematics with more than 45 000 citations up to 2023, the model of linear preferential attachment tree has been studied before (at least) by Szymanski \cite{szymanski1987nonuniform} and Mahmoud \cite{mahmoud1992distances}. This is a very good example of Stigler's law of eponymy.  Exercise \ref{exo:pawel} was suggested by Pavel Krapivsky. The almost sure convergence of the largest degrees (Theorem \ref{thm:mori}) is due to Mori \cite{mori2005maximum}. Theorem \ref{thm:heightBA} is first proved in \cite{pittel1994note} using the continuous time embedding technique.\medskip 

\noindent{\textbf{Hints for Exercises.}}\ \\

\noindent Exercise \ref{exo:pawel}: A new vertex attaches to $\noeud{i}$ with probability proportional $ \mathrm{deg}^+_{ \mathfrak{T}_n}(\noeud{i})+1$ which is equal to $ \mathrm{deg}_{ \mathfrak{T}_n}(\noeud{i})$ except for the root $\noeud{0}$ which has a small bias. 
\chapter*{Appendix}
\hfill So the last shall be first.\\

\hfill (Matthew 20:16)
\bigskip

\section*{Large deviations for Poisson random variables}

Let us state a simple lemma on Poisson random variables which we used many times in these lecture notes. Recall that $ \mathfrak{P}(a)$ is a Poisson variable of mean $a >0$.
  \begin{lemma}[Large deviations and maximum of i.i.d.~Poisson random variables]  \label{lem:LDpoisson} For $a >0$ denote by $I(a) := a\log a-(a-1)$ then for all $t \geq 0$ we have 
 $$\mathbb{P}( \mathfrak{P}(t) > at) \leq \mathrm{e}^{-t I(a)} \quad \mbox{ if } a>1 \quad \mbox{ and } \quad \mathbb{P}( \mathfrak{P}(t) < at) \leq \mathrm{e}^{-t I(a)} \quad \mbox{ if } 0<a<1.$$
Fix $c \geq 0$ and let $ X_1,\dots , X_{ \lfloor n^c \rfloor}$ be i.i.d.~random variables with Poisson law of expectation $\log n$. Then we have 
$$ \frac{\log{\max_{1 \leq i \leq n^c} X_i }}{\log n} \xrightarrow[n\to\infty]{(  \mathbb{P})} x_c, \quad \mbox{ with $x_c \geq 1$ solution to } I(x_c) = c,$$ and furthermore $ \mathbb{P}(\max_{1 \leq i \leq n^c} X_i  \leq n^{x_c - \varepsilon})$ tends to $0$ stretched-exponentially fast.\end{lemma}
\noindent \textbf{Proof.} Suppose $a >1$ and let us apply a standard exponential Markov's inequality to write for $\lambda >0$
$$ \mathbb{P}( \mathfrak{P}(t) > at) \leq \frac{ \mathbb{E}[ \mathrm{e}^{\lambda  \mathfrak{P}(t) }]}{ \mathrm{e}^{\lambda a t}} =   \mathrm{exp}(t(( \mathrm{e}^\lambda -1)-\lambda a)) \underset{\lambda = \log a}{\leq} \exp(-t I(a)).$$ The case $a<1$ is dealt with similarly using negative $\lambda$. For the second point notice that for $x \geq 1$ we have using the first point 
$$ \mathbb{P}( \max_{1 \leq i \leq n^c} X_i \leq n^{x}) = (1- \mathbb{P}( \mathfrak{P}(\log n) > x \log n))^{n^c} = \exp\left(- n^{c -I(x) + o(1)}\right),$$ and so the above probability tends to $0$ stretched-exponentially fast if $I(x) < c$ and to $1$ if $I(x)>c$. \qed \medskip 

First part of Lemma \ref{lem:LDpoisson} is known under the name of ``Bennett’s inequality'' (see Terence Tao's blog for a nice sharpening of it).

  \bibliographystyle{siam}
\bibliography{bibli}

            \end{document}